\documentclass[a4paper,11pt,reqno]{amsart}

\usepackage{amsmath}
\usepackage{amssymb}
\usepackage{amsfonts}
\usepackage{epsfig}
\usepackage{epigraph}
\usepackage[hang]{caption}
\usepackage{mathtools}

\usepackage{color}
\usepackage{comment}
\usepackage{mathrsfs}

\usepackage{MnSymbol}

\graphicspath{{eps_1A/}{eps_1B/}{eps_1B_mirror/}}

\def\colorcom#1#2#3{
 \def#1##1 {
  \begin{color}{#3}
      \bf[\,#2: ##1\,]
  \end{color}}
} \colorcom\RS{RS}{red}

\def\colorcom#1#2#3{
 \def#1##1 {
  \begin{color}{#3}
      \bf[\,#2: ##1\,]
  \end{color}}
} \colorcom\DL{DL}{blue}

\makeatletter

\@addtoreset{equation}{chapter}

\newtheorem{theorem}{Theorem}[section]
\@addtoreset{equation}{section}

\def\botG{\stackrel{}{\Gamma}}

\newtheorem{corollary}[theorem]{Corollary}
\newtheorem{remark}[theorem]{Remark}
\newtheorem{proposition}[theorem]{Proposition}
\newtheorem{lemma}[theorem]{Lemma}
\newtheorem{example}[theorem]{Example}

\def\PerfProof{{\it Proof.\ }}

\DeclareFontFamily{U}{matha}{\hyphenchar\font45}
\DeclareFontShape{U}{matha}{m}{n}{
      <5> <6> <7> <8> <9> <10> gen * matha
      <10.95> matha10 <12> <14.4> <17.28> <20.74> <24.88> matha12
      }{}
\DeclareSymbolFont{matha}{U}{matha}{m}{n}

\DeclareMathSymbol{\notperp}{3}{matha}{"4D}

 \DeclarePairedDelimiter\abs{\lvert}{\rvert}%
 \DeclarePairedDelimiter\norm{\lVert}{\rVert}%

 \makeatletter
 \let\oldabs\abs
 \def\abs{\@ifstar{\oldabs}{\oldabs*}}
 \let\oldnorm\norm
 \def\norm{\@ifstar{\oldnorm}{\oldnorm*}}

\makeatother

\makeindex

\begin{document}

\title[Root systems and diagram calculus. I]
{\qquad\qquad Root systems and diagram calculus. \newline
  I. Regular extensions of Carter diagrams and the uniqueness of
       conjugacy classes}
         \author{Rafael Stekolshchik}

\date{}

\begin{abstract}
 In $1972$, R.~Carter introduced admissible diagrams to
 classify conjugacy classes in a finite Weyl group $W$.
 We say that an admissible diagram $\Gamma$ is a {\it Carter diagram}
 if any edge $\{\alpha, \beta\}$ with
 inner product $(\alpha, \beta) > 0$ (resp. $(\alpha, \beta) < 0$)
 is drawn as dotted (resp. solid) edge.
 We construct an explicit transformation of any Carter diagram containing long cycles
 (with the number of vertices $l > 4$) into another Carter diagram containing only $4$-cycles.
 Thus, all Carter diagrams containing long cycles can be eliminated from the classification list.

 There exist diagrams determining two conjugacy classes in $W$.
 It is shown that any connected Carter diagram $\Gamma$ containing a $4$-vertex pattern
 $D_4$ or $D_4(a_1)$ determines a single conjugacy class.
 The main approach is studying different extensions of Carter diagrams.
 Let $\widetilde{\Gamma}$ be the Carter diagram obtained from a certain Carter
 diagram $\Gamma$ by adding a single vertex $\alpha$
 connected to $\Gamma$ at $n$ points, $n \leq 3$.
 Let a {\it socket} be the set of vertices of $\Gamma$ connected to $\alpha$.
 If the number of sockets available for extensions is equal to $2$, then there is a pair of extensions
 $\Gamma < \widetilde{\Gamma}_L$ and $\Gamma < \widetilde{\Gamma}_R$,
 called {\it mirror extensions} and the pair elements $w_L$ and $w_R$
 associated with $\widetilde{\Gamma}_L$ and $\widetilde{\Gamma}_R$.
 We show that $w_R = T^{-1}w_L{T}$ for some $T \in W$, where the map $T$ is explicitly constructed for all mirror extensions.

 In Carter's description of the conjugacy classes in a Weyl group
 a key result (Carter's theorem) states that every element in a Weyl group is a product of two involutions.
 One of the goals of this paper and its sequels is to prepare the notions and framework
 in which we give the proof of this fact without appealing to the classification of conjugacy classes.

 \end{abstract}

\maketitle

\tableofcontents

\setlength{\epigraphwidth}{95mm}

\clearpage

\clearpage

\epigraph{The use of trees as diagrams for groups was anticipated in
1904, when C.~Rodential \cite{Rod04} was commenting on a set of
models of cubic surfaces. He was analyzing the various rational
double points that can occur on such a surface. In 1931, I used
these diagrams in my enumeration of kaleidoscopes, where the dots
represent mirrors. E.B.Dynkin re-invented the diagrams in 1946 for
the classification of simple Lie algebras.}{H.~S.~M.~Coxeter, The
evolution of Coxeter-Dynkin diagrams, \cite[p.224]{Cox91}, 1991}

\section{\sc\bf Preview}

\subsection{Surprising cycles and dotted edges}
   \label{sec_surprising}
   \index{root system}
   \index{Coxeter element}
   \index{Coxeter plane}
   \index{Weyl group}
 Let $W$ be a Weyl group and $\varPhi$  the root system associated
 with $W$. Let us connect the non-orthogonal simple roots in
 $\varPhi$ with each other. We get a graph called a \emph{Dynkin
 diagram}. One may want to connect all (not only simple)
 non-orthogonal roots with each other. How does the graph thus
 obtained look like?

 The graphs thus obtained are the beautiful color computer-generated
 pictures given on John Stembridge's home page (based on Peter
 McMullen's drawings). These pictures are projections of the root
 system of $\varPhi$ into the Coxeter plane\footnotemark[1], see
 \cite{Stm07}. Though beautiful, these graphs are not easy to grasp: For example,
 in the picture of the root system $E_8$, there are $6720$ edges, see \cite{Ma10}.

 \footnotetext[1]{The {\it Coxeter plane} $P$ is the
 span of the real and imaginary parts of an eigenvector for the
 Coxeter element {\bf C} with eigenvalue
 $\cos(\frac{2\pi}{h})+i\cdot\sin(\frac{2\pi}{h})$, where $h$ is the
 Coxeter number associated with the root system $\varPhi$.}

 To see some details in Stembridge's pictures, one can confine
 oneself to only connected subsets of linearly independent roots.
 Then the graphs simplify drastically: There are $\leq n + 4$ edges
 in the connected diagram associated with every subset of
 linearly independent roots lying in the root system $\varPhi$ of
 rank $n$, see Fig. \ref{diagram_tree}. Essentially, such diagrams
 were presented by Carter in 1972, in \cite{Ca70}, \cite{Ca72}. These
 graphs are said to be {\it admissible diagrams} and are designed to
 characterize elements of the Weyl group, see definition in \S\ref{sec_adm_diagr}.

 Each element $w \in W$ can be expressed in the form
 \begin{equation}
   \label{any_roots_0}
    w  = s_{\alpha_1} s_{\alpha_2} \dots s_{\alpha_k}, \text{ where } \alpha_i \in \varPhi,
 \end{equation}
 and $s_{\alpha_i} \in W$ are reflections
 corresponding to not necessarily simple roots $\alpha_i \in \varPhi$.

Carter proved that $k$ in the decomposition \eqref{any_roots_0} is
the smallest if and only if the subset of roots $\{\alpha_1,
\alpha_2, \dots, \alpha_k\}$ is linearly independent; such a
decomposition is said to be {\it reduced}. The admissible diagram
corresponding to the given element $w$ is not unique, since the
reduced decomposition of the element $w$ is not unique.

 When I first got acquainted with admissible diagrams I was surprised
 by the fact that these diagrams contain cycles, though the extended Dynkin diagram $\widetilde{A}_l$,
 cannot be a part of any admissible diagram (Lemma \ref{lem_must_dotted}).
 It turned out that the cycles in
 admissible diagrams essentially differ from the cycle $\widetilde{A}_l$.
 Namely, in such a cycle, there is necessarily two pairs of roots: A pair with a
 positive inner product together with a pair with a negative
 inner product. This does not happen for $\widetilde {A}_l$.

 This observation motivated me to distinguish such pairs of roots:
 Let us draw the {\it dotted} (resp. {\it solid}) edge $\{\alpha,
 \beta\}$ if $(\alpha, \beta) > 0$ (resp. $(\alpha, \beta) < 0$), see Fig.  \ref{two_types_diagram}.
 Let the diagrams with properties of admissible diagrams and containing
 dotted edges be called {\it Carter diagrams}. Up to dotted edges,
 the classification of Carter diagrams coincides with the
 classification of admissible diagrams. Recall that $(\alpha, \beta) > 0$
 (resp. $(\alpha, \beta) < 0$) means that the angle between
 roots $\alpha$ and $\beta$ is acute (resp. obtuse). For the Dynkin
 diagrams, all angles between simple roots are obtuse, thus all edges
 are solid.

\subsection{The Carter theorem and the theorem on a single conjugacy class}
  \label{sec_two_th}
 Admissible and Carter diagrams $\Gamma$ involve also the following
 restriction: If a subdiagram of $\Gamma$ is a cycle, then it contains
 an even number of nodes. Denote the set of elements $w \in W$, each of
 which corresponds to an admissible diagram, by $W_0$. The
 existence of an admissible diagram for the element $w$ means that
 $w$ can be decomposed into the product of two involutions as
 follows:
 \begin{equation}
   \label{two_invol_0}
      w = w_1{w}_2, \quad \text{ where } \quad
      w_1 = s_{\alpha_1} s_{\alpha_2} \dots s_{\alpha_k}, \quad
      w_2 = s_{\beta_1} s_{\beta_2} \dots s_{\beta_h},
 \end{equation}
 the roots $\{\alpha_i\mid i = 1,\dots,k\}$ being mutually
 orthogonal, and the roots $\{\beta_j\mid j = 1,\dots,h\}$ being also
 mutually orthogonal.

 Thus, $W_0$ is the subset of elements $w \in W$ that can be
 decomposed as \eqref{two_invol_0}. It turns out that $W_0 = W$. This
 fact is one of the main results of \cite[Theorem C]{Ca72}. I call
 this result {\bf the Carter theorem}.  The proof of the Carter
 theorem is based on the classification of conjugacy classes.
 (By \cite[p. G-21]{Ca70}, N.~Burgoyne carried out a check of classification for
 $E_7$ and $E_8$ with a computer aid).

 \index{Coxeter groups}
 \index{Weyl group}
 I would like to quote Carter  from \cite[p. 4]{Ca00}:

{\it \lq\lq One remarkable feature of the theory of Coxeter groups
and Iwahori-Hecke algebras is the number of key properties which at
present have no uniform proof and can only be proved in a
case-by-case manner. We mention three examples of such properties.
In the first place every element of a Weyl group is a product of two
involutions. This is a key property in Carter's description of the
conjugacy classes in the Weyl group \cite{Ca72}. Secondly, every
element of a finite Coxeter group can be transformed into an element
of minimal length in its conjugacy class by a sequence of
conjugations by simple reflections such that the length does not
increase at any stage. This property is basic to the Geck-Pfeiffer
approach to the conjugacy classes of Coxeter groups. Thirdly, there
are basic properties of Lusztig's a-function which at present have
only case-by-case proofs. The a-function is an important invariant
of irreducible characters of Coxeter groups. One cannot be satisfied
with the theory of Coxeter groups until such case-by-case proofs are
replaced by uniform proofs of a conceptual nature.\rq\rq}

 One of the goals of this paper, and its sequels \cite{St10} and
 \cite{St11}, is a \emph{conceptual} proof of the Carter theorem, and
 a number of issues associated with this theorem.

 I would like not to use the classification of conjugacy classes
 either for proving the Carter theorem or for other reasoning. My
 proof of the Carter theorem presented in \cite{St10}, \cite{St11}
 does not rely on the classification of conjugacy classes; instead it is
 based on the classification of the graphs I introduced in \cite{St10}
 and called {\it linkage diagrams}. This proof is also a
 case-by-case checking, but a much simpler one, and is obtained
 without aid of computers, see \cite{St10}.

 {\bf In the current paper, we construct the framework for studies of
 basic patterns and facts concerning the Carter diagrams.}

As I said above, each element $w \in W_0$,  together with its
conjugacy class $\{w\}$, is associated with either a certain Carter diagram
or an admissible diagram. {\bf One of the central questions} considered
in this paper is: {\it Do all elements associated with a given Carter
diagram constitute a single conjugacy class?} The answer is that,
generally speaking, this is not so.

\begin{figure}[h]
\centering
\includegraphics[scale=1.5]{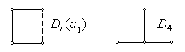}
\caption{Characterizing patterns $D_4(a_1)$, $D_4$}
\label{tetris}
\end{figure}

We show that {\bf any Carter diagram containing the square
$D_4(a_1)$ or the Dynkin diagram $D_4$ determines a single conjugacy
class (Theorem \ref{th_uniq_diagr})}, see Fig. \ref{tetris}. We call
this result {\bf the uniqueness theorem for conjugacy classes}.

Of course, this theorem can be derived from the classification of
conjugacy classes. In accordance with the above, we do not use this
classification. We use another approach based on the consideration
of regular extensions of Carter diagrams and the induction in the
number of vertices.

Consider the connection diagram $\widetilde{\Gamma}$ obtained from a
certain Carter diagram $\Gamma$ by adding only one vertex $\alpha$,
where $\alpha$ is connected to $\Gamma$ at $v$ points, where $v=1,2$
or $3$. If $\widetilde{\Gamma}$ is also a Carter diagram, this
extension $\Gamma$ to $\widetilde{\Gamma}$ is said to be a {\it
regular extension} and is denoted by
\begin{equation*}
 \Gamma < \widetilde{\Gamma} \qquad \text{ or } \qquad  \Gamma \stackrel{\alpha}{<} \widetilde{\Gamma}.
\end{equation*}

The uniqueness theorem  is proved by induction: If $\Gamma$ contains
$D_4(a_1)$ or $D_4$, and $\Gamma$ determines a single conjugacy
class, then $\widetilde{\Gamma}$ also determines a single conjugacy
class (Proposition \ref{prop_induct_step}).

\subsection{Theorem on eliminating long cycles}
There are different decompositions \eqref{any_roots_0} of $w$: They
can be obtained from each other by some transformations.
Transforming one Carter diagram $\Gamma_1$ to another Carter diagram
$\Gamma_2$ we can get a certain intermediate diagram $\Gamma'$ which
is not necessarily a Carter diagram. Such an intermediate diagram
will be called a {\it connection diagram}.  This term is motivated
by the fact that this diagram describes only the connectivity between roots,
nothing more.

{\bf The study of certain properties of connection and Carter
diagrams is one of the goals of this paper}.

Consider an example of basic properties of connection and Carter
diagrams. Let $\{\alpha_1, \alpha_2, \alpha_3\}$ be $3$ linearly
independent and mutually orthogonal roots. There do not exist two
non-connected roots $\beta$ and $\gamma$ connected to every
$\alpha_i$ in such a way that $\{\alpha_1, \alpha_2, \alpha_3,
\beta, \gamma\}$ is a linearly independent quintuple. First of
all, any cycle of linearly independent roots contains an odd number
of dotted edges. Let $n_1$, $n_2$, $n_3$ be the odd numbers of
dotted edges in every cycle $\{\alpha_i, \beta, \alpha_j, \gamma
\}$, where $1 \le  i < j \le 3$. Therefore, $n_1 + n_2 + n_3$ is
odd, contradicting the fact that every dotted edge appears twice, so
$n_1 + n_2 + n_3$ is even (Corollary \ref{cor_numb_ep}), see Fig.
\ref{3-cells}.
 This and similar properties allow us to simplify the classification of
 Carter diagrams. The main result obtained in this direction is the
 following one: {\bf Any Carter diagram containing $n$-cycles ($n >
 4$) is equivalent to another Carter diagram containing only
 $4$-cycles} (Theorem \ref{th_get_rid}).
 This is the theorem on eliminating Carter diagrams with long cycles,
 i.e., $l$-cycles with $l > 4$. To exclude long cycles {\bf we
 construct explicit transformations mapping every Carter diagram with
 long cycles into a certain Carter diagram containing only
 $4$-cycles}, see Table \ref{tab_eq_pairs}.

\subsection{Mirror extensions}  In the current paper,
we use what we call regular extensions, see \S\ref{sec_two_th}.

{\bf Extensions of Carter diagrams constitute the essential part of
the study. Among different types of regular extensions there is one
difficult case which we dubbed {\it mirror extensions}.}

Let $\Gamma \stackrel{\alpha}{<}  \widetilde{\Gamma}$ be a regular
extension. The Carter diagram $\widetilde{\Gamma}$ extends $\Gamma$
by adding the vertex $\alpha$; the set of vertices connected to
$\alpha$ is called a {\it socket}. The number of sockets available
to get $\widetilde{\Gamma}$ is said to be a {\it sockets number}. If
the sockets number is $n = 2$, the corresponding regular extensions are called
{\it mirror extensions}.
\begin{figure}[h]
\centering
\includegraphics[scale=0.7]{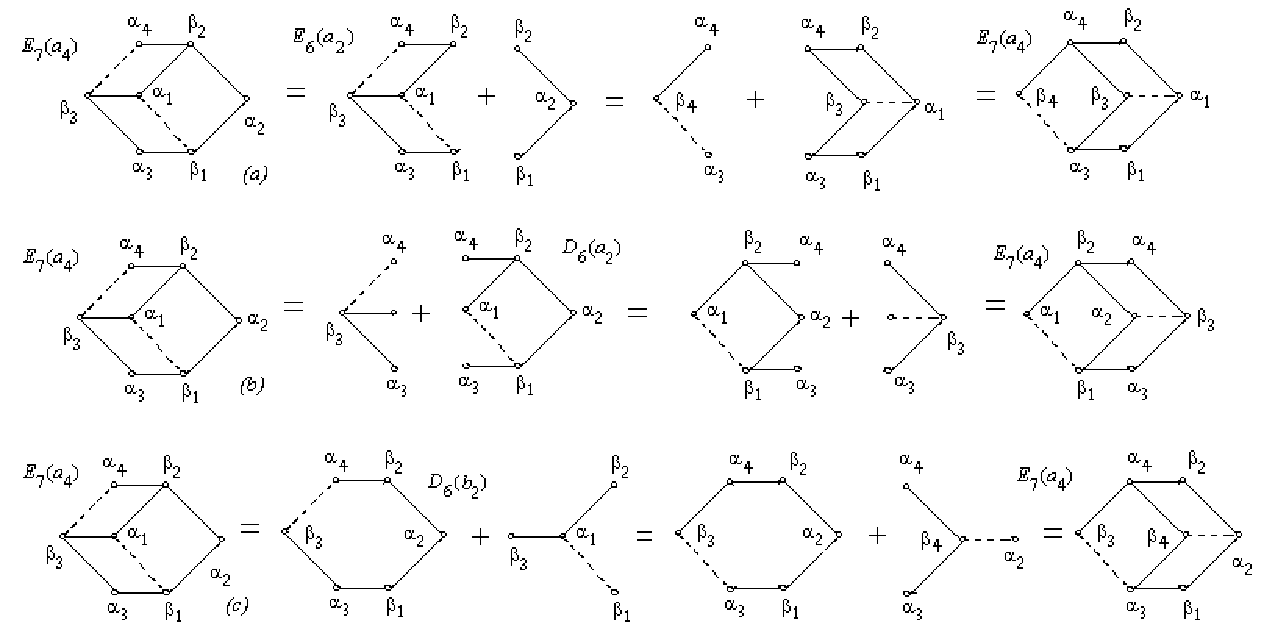}
 \caption{Three pairs of mirror extensions $\Gamma < \widetilde{\Gamma}$, where
 $\widetilde{\Gamma} = E_7(a_4)$, see Table \ref{tab_2option_E7a4}}
\label{multi_sockets}
\end{figure}
For example, look at the three pairs of mirror extensions $\Gamma <
\widetilde{\Gamma}$ depicted in Fig. \ref{multi_sockets}. Speaking
of mirror extensions, we mean two extensions: $\Gamma <
\widetilde\Gamma_L$ and $\Gamma < \widetilde\Gamma_R$ that are said
to be \emph{left} and \emph{right} extensions; they are presented in
Tables \ref{tab_mirror_1} -- \ref{tab_mirror_2}. The main relevant
fact concerning the mirror extensions is the equivalence of
$\widetilde\Gamma_L$ and $\widetilde\Gamma_R$, i.e., conjugacy of
elements of the Weyl group corresponding to $\widetilde\Gamma_L$ and
$\widetilde\Gamma_R$: {\bf Any $\widetilde\Gamma_L$-associated
element $w_L$ and any $\widetilde\Gamma_R$-associated element $w_R$
are conjugate}  (Theorem \ref{th_mirror_ext}):
\begin{equation*}
  w_R = T^{-1}{w}_L{T} \text{ for some } T \in W.
\end{equation*}

\subsection{Extensions and the uniqueness theorem}

 The regular extensions are distinguished by the sockets number $n =
 1,2,3$. The regular extensions with $n=1$ (resp. $n=2$, resp. $n=3$)
 are said to be {\it single-track} (resp. {\it mirror}, resp. {\it
 threefold}) extensions. Tables \ref{tab_single_track} -- \ref{tab_threefold} contain sets of single-track, mirror and
 threefold extensions, see Fig. \ref{multi_sockets}.

One more parameter characterizing extensions is the {\it pinholes
number} meaning the number of vertices in the socket. The socket
with the pinholes number equal to $v$ is said to be a {\it
$v$-socket}; a regular extension with a $v$-socket is said to be a
$\mathsf{P}v$-extension. The first pair of mirror extensions $\Gamma
< \widetilde{\Gamma}$ depicted in Fig. \ref{multi_sockets}
constitutes a pair of $2$-socket extensions, the second and third
pairs are $3$-socket extensions, see Table \ref{tab_2option_E7a4}.
Any Carter diagram can be obtained as the extension of a smaller
Carter diagram by means of one of $\mathsf{P}v$-extension with
$v=1,2,3$. {\bf The proof of the uniqueness theorem is derived from
the consideration of different cases of the sockets number $n$ and
the pinholes number $v$.} Actually, there are fewer cases, not all
cases are realized.
 \index{class of Carter diagrams $\mathsf{C4}$}
 \index{class of Carter diagrams $\mathsf{DE4}$}
 \index{$\mathsf{C4}$ (class of Carter diagrams)}
 \index{$\mathsf{DE4}$ (class of Carter diagrams)}
 \index{$\mathsf{C4} \coprod \mathsf{DE4}$}
\begin{table}[h]
 \Small
  \centering
  \renewcommand{\arraystretch}{1.4}
  \begin{tabular} {|c|c|c|c|}
  \hline
    $E_6(a_2) < E_7(a_4)$
      & $E_6(a_2) ~\stackrel{\{\alpha_3, \alpha_4\}}\rightpitchfork~ E_7(a_4)$
      & $E_6(a_2) ~\stackrel{\{\beta_1, \beta_2\}}\rightpitchfork~ E_7(a_4)$
                                                & $2$-sockets  \\
  \hline
    $D_6(a_2) < E_7(a_4)$
       & $D_6(a_2) ~\stackrel{\{\alpha_1, \alpha_3, \alpha_4\}}\rightpitchfork~ E_7(a_4)$
       & $D_6(a_2) ~\stackrel{\{\alpha_2, \alpha_3, \alpha_4\}}\rightpitchfork~ E_7(a_4)$
                                                & $3$-sockets \\
  \hline
    $D_6(b_2) < E_7(a_4)$
       & $D_6(b_2) ~\stackrel{\{\beta_1, \beta_2, \beta_3\}}\rightpitchfork~ E_7(a_4)$
       & $D_6(b_2) ~\stackrel{\{\alpha_2, \alpha_3, \alpha_4\}}\rightpitchfork~ E_7(a_4)$
                                                & $3$-sockets  \\
  \hline
\end{tabular}
  \vspace{1mm}
  \caption{Mirror extensions $\Gamma < \widetilde{\Gamma}$, where $\widetilde{\Gamma} = E_7(a_4)$,
  see Fig. \ref{multi_sockets}}
  \label{tab_2option_E7a4}
\end{table}
From the point of view of the sockets number any regular extension $\Gamma <
\widetilde{\Gamma}$ is a certain single-track, mirror or threefold
extension, while from the point of view of pinholes numbers any regular
extension is a certain $\mathsf{P}1$-, $\mathsf{P}2$-,
$\mathsf{P}3$-extension (Lemma \ref{lem_123_socket}).
The process of adding a new root to a certain root subset $S$
associated with a given Carter diagram $\Gamma$ is depicted by
different extensions. {\bf This process is one of the
representations of what I'd like to call the Diagram Calculus.}

\begin{figure}[h]
\centering
\includegraphics[scale=0.62]{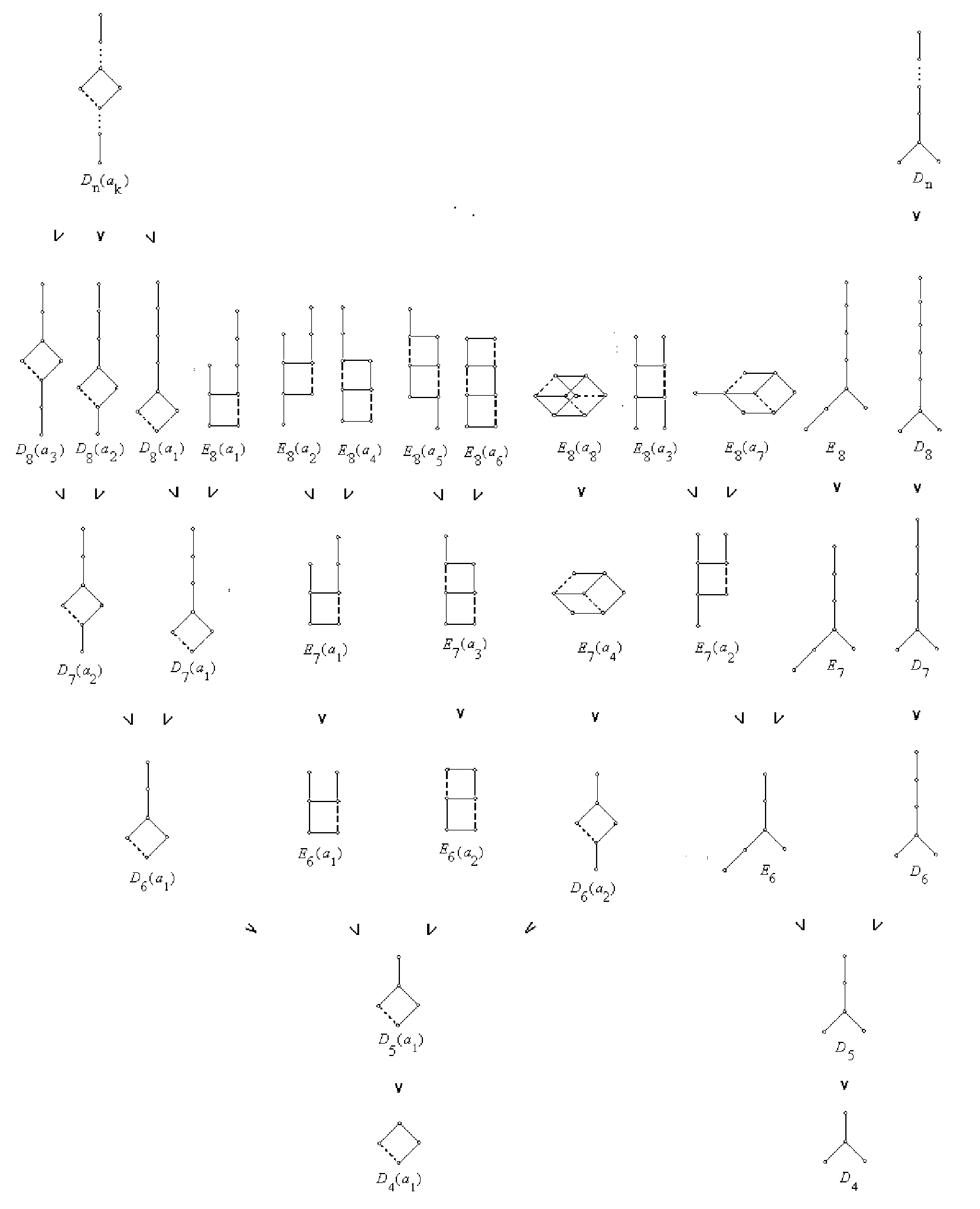}
\caption{The tree of Carter diagrams from $\mathsf{C4} \coprod \mathsf{DE4}$}
\label{diagram_tree}
\end{figure}

In the proof of the uniqueness theorem we see that some symmetric
Carter diagrams have several regular extensions caused by the
symmetry of order $2$ (resp. order $3$) of a given Carter
diagram, see Fig. \ref{multi_sockets}; when this happens, mirror
extensions (resp. threefold extensions) arise.
  Let $\mathsf{C4}$ be the class of connected simply-laced Carter
  diagrams each of which {\it contains a $4$-cycle $D_4(a_1)$}.
  Let $\mathsf{DE4}$ be the class of connected simply-laced Carter diagrams each of which is
  a {\it Dynkin diagram and contains $D_4$.}
  The regular extensions determine a partial order on the class of
  Carter diagrams $\mathsf{C4} \coprod \mathsf{DE4}$.
  {\bf The partially ordered tree of Carter diagrams is depicted in Fig. \ref{diagram_tree}.}

\subsection{The partial Cartan matrix and linear dependence}
  \index{linearly dependent root}
Considering the regular extensions of Carter diagrams we frequently
need a criterion that tells if a given root performing the extension
is linearly dependent on a certain root subset $S$ or not.

 In this paper, I introduce a matrix called a \emph{partial Cartan matrix} $B_{\Gamma}$.
 {\bf There is a simple criterion that
 tells if $\gamma$ connected with the root $\tau_i \in S$
 is linearly dependent on $S$ or not.
 This condition is formulated in terms of diagonal elements  of the
 matrix $B_{\Gamma}^{-1}$}, see \S\ref{sec_partial_B}. In \cite{St10}, we use
 partial Cartan matrices to derive a criterion that tells if a
 given root $\gamma$ is linearly independent of $S$ or not.

 For any diagram $\Gamma$, we consider a certain subset $S \subseteq
 \varPhi$ of linearly independent roots such that there is a
 one-to-one correspondence between roots of $S$ and vertices of
 $\Gamma$. The subset $S$ is said to be {\it $\Gamma$-associated}.
 Let $S = \{\tau_1, \dots, \tau_l\}$. Consider the matrix
 $B_{\Gamma}$ defined analogously to the Cartan matrix associated
 with a given Dynkin diagram:
 \begin{equation*}
  \scriptsize
   B_{\Gamma} :=
      \left (
        \begin{array}{cccccc}
         (\tau_1, \tau_1) & \dots & (\tau_1, \tau_l) \\
          & \dots  & \\
         (\tau_l, \tau_1) & \dots & (\tau_l, \tau_l) \\
        \end{array}
      \right ).
 \end{equation*}
The matrix $B_{\Gamma}$ will be called a \emph{partial Cartan
matrix}. The off-diagonal elements of the partial Cartan matrix
$B_{\Gamma}$ might be positive integers. For example, $B_{\Gamma}$
for $D_5$ and $D_5(a_1)$ are as follows:
 \index{Cartan matrix! - of hyperbolic type}
 \index{$b^{\vee}_{\eta, \eta}$ (diagonal element of $B^{-1}_{\Gamma}$)}
\begin{equation*}
 \scriptsize
 \begin{split}
  & \begin{array}{c}  \qquad \includegraphics[scale=0.62]{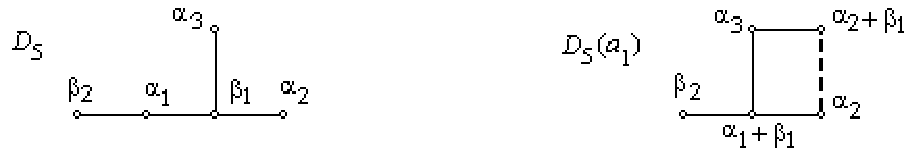}  \end{array} \\
  & \begin{array}{c}
   \left [
   \begin{array}{cccccc}
     2 & 0 & 0  &   -1 & -1  \\
     0 & 2 & 0  &   -1 & 0 \\
     0 & 0 & 2  &   -1 & 0  \\
    -1 & -1 & -1 &   2 & 0  \\
    -1 & 0 & 0 &    0 & 2 \\
  \end{array}
   \right ]
   \begin{array}{c}
      \alpha_1 \\
      \alpha_2 \\
      \alpha_3 \\
      \beta_1 \\
      \beta_2 \\
   \end{array}
  \end{array} 
  \qquad \qquad \qquad
    \begin{array}{c}
   \left [
   \begin{array}{cccccc}
     2 & 0 & 0  &   -1 & 1  \\
     0 & 2 & 0  &   -1 & -1 \\
     0 & 0 & 2  &   -1 & 0  \\
    -1 & -1 & -1 &   2 & 0  \\
     1 & -1 & 0 &    0 & 2 \\
  \end{array}
   \right ]
   \begin{array}{c}
      \alpha_2 \\
      \alpha_3 \\
      \beta_2 \\
      \alpha_1 + \beta_1 \\
      \alpha_2 + \beta_1 \\
   \end{array}
  \end{array} \\ 
   & \qquad \qquad \qquad B_{\Gamma} \text{ for } D_5
   \qquad \qquad \qquad \qquad \qquad \qquad \qquad \qquad B_{\Gamma} \text{ for } D_5(a_1)
  \end{split}
\end{equation*}
 In the case $D_5(a_1)$, the matrix $B_{\Gamma}$ contains $1$ in the
 slot $\{\alpha_2, \alpha_2 + \beta_1\}$ corresponding to the dotted
 edge.\footnotemark[1]
 Let $V$ be spanned by the subset of simple roots $\Delta^+ \subset \varPhi$.
 The subspace $L \subseteq V$ spanned by the root subset is said
 to be \emph{$S$-associated} or \emph{$\Gamma$-associated}. The root
 system $\varPhi$ for $D_5$, and the $D_5(a_1)$-associated root
 subset $S$ are as follows:
 \begin{equation*}
   \Delta^+ = \{\alpha_1, \alpha_2, \alpha_3, \beta_1, \beta_2\}, \qquad
   S = \{\alpha_2, \alpha_3, \beta_2, \alpha_1 + \beta_1, \alpha_2 + \beta_1\}.
 \end{equation*}
Roots $\alpha_2 + \beta_1$ and $\alpha_1 + \beta_1$ are not simple
in $\varPhi = \varPhi(D_5)$. In this case, we have $L = V$.

Let $\gamma$ be a root \underline{linearly dependent on $S$}, let $\gamma$ be
connected with only one root $\tau_i \in S$. We have
\footnotetext[1]{D.~Leites pointed out that there are a number
 of other cases, where some off-diagonal elements of the Cartan matrix are
 positive integers. In particular, this is so in the case of Lorentzian algebras, see
 \cite{GN02}, \cite{CCLL10}. However, note that in these cases the Cartan matrices are of {\it hyperbolic type},
 whereas the partial Cartan matrices for Carter diagrams are {\it positive definite},
 see Proposition \ref{restr_forms_coincide}.}
\begin{equation}
   \label{eq_sum_tau_j}
    \gamma = \sum\limits_{j=1}^l{t_j\tau_j}.
 \end{equation}
Set $\gamma^{\vee}:=\{(\gamma, \tau_i)\mid i = 1,\dots,l\}$. By
\eqref{eq_sum_tau_j}, we have
 \begin{equation*}
   \Small
       \gamma^{\vee} := \left (
    \begin{array}{c}
      (\gamma, \tau_1) \\
      \dots \\
      (\gamma, \tau_i) \\
      \dots \\
      (\gamma, \tau_l) \\
    \end{array}
  \right ) =
    \left (
    \begin{array}{c}
      \sum t_j(\tau_j, \tau_1) \\
      \dots \\
      \sum t_j(\tau_j, \tau_i) \\
      \dots \\
      \sum t_j(\tau_j, \tau_l) \\
    \end{array}
    \right ) =
    B_{\Gamma}
    \left (
    \begin{array}{c}
      t_1 \\
      \dots \\
      t_i \\
      \dots \\
      t_l \\
    \end{array} \right ) =
    B_{\Gamma}\gamma =
    \left (
    \begin{array}{c}
      0 \\
      \dots \\
      \pm{1} \\
      \dots \\
      0 \\
    \end{array}
    \right )
    \begin{array}{l}
       \\
       \\
     \hspace{-2mm}i \\
       \\
       \\
    \end{array},
    \text{ then }
   \gamma = B_{\Gamma}^{-1}\gamma^{\vee}.
 \end{equation*}
{\bf The root $\gamma$ is linearly dependent on $S$
 and connected only with $\tau_i$ if and only if}
\begin{equation*}
  b^{\vee}_{\tau_i,\tau_i} = 2,
\end{equation*}
where $b^{\vee}_{\tau_i,\tau_i}$ is the $i$th diagonal element of
$B_{\Gamma}^{-1}$. For $D_5$, we have $b^{\vee}_{\alpha_1,\alpha_1}
= 2$; for $D_5(a_1)$, only $b^{\vee}_{\tau_1,\tau_1} = 2$, where
$\tau_1 = \alpha_1 + \beta_1$; the corresponding slots are boxed:
\begin{equation*}
 \scriptsize
  \begin{split}
  & \frac{1}{4} \left [
   \begin{array}{cccccc}
     \fbox{8} & 4 & 4 & 8 & 4 \\
     4 & 5 & 3 & 6 & 2 \\
     4 & 3 & 5 & 6 & 2 \\
     8 & 6 & 6 & 12 & 4 \\
     4 & 2 & 2 & 4 & 4 \\
   \end{array} \right ]
   \begin{array}{c}
      \alpha_1 \\
      \alpha_2 \\
      \alpha_3 \\
      \beta_1 \\
      \beta_2 \\
   \end{array}
  \qquad \qquad \qquad
   \frac{1}{4} \left [
   \begin{array}{ccccc}
     5 & 1 &   2 & 4 & -2  \\
     1 & 5 &   2 & 4 & 2 \\
     2 & 2 &   4 & 4 & 0  \\
     4 & 4 &   4 & \fbox{8} & 0  \\
    -2 & 2 &   0 & 0 & 4 \\
  \end{array}
   \right ]
   \begin{array}{c}
      \alpha_2 \\
      \alpha_3 \\
      \beta_2 \\
      \alpha_1 + \beta_1 \\
      \alpha_2 + \beta_1 \\
   \end{array} \\
   & \qquad \qquad B_{\Gamma}^{-1} \text{ for } D_5
   \qquad \qquad \qquad \qquad \qquad \qquad \qquad B_{\Gamma}^{-1} \text{ for } D_5(a_1)
  \end{split}
\end{equation*}

~\\

\index{$\mathsf{C4} \coprod \mathsf{DE4}$}

\clearpage

\section{\sc\bf Introduction}

\subsection{Diagrams}
  \label{sec_intro}

\subsubsection{Admissible and Carter diagrams}
  \label{sec_adm_diagr}
  \index{root system}
 Let $\varPhi$ be the root system associated with
 a Weyl group $W$; let $s_{\alpha_i}$ be the reflection in $W$
 corresponding to \underline{not necessarily simple root} $\alpha_i \in \varPhi$.
 Each element $w \in W$ can be expressed in the form
 \begin{equation}
   \label{any_roots}
    w  = s_{\alpha_1} s_{\alpha_2} \dots s_{\alpha_k}, \text{ where } \alpha_i \in \varPhi,
 \end{equation}
 We denote by $l_C(w)$ the smallest value $k$ in any expression like
 \eqref{any_roots}, see \cite[p. 3]{Ca72}.
 We always have $l_C(w) \leq l(w)$. Recall that $l(w)$ is the smallest value
 $k$ in any expression like \eqref{any_roots} such that all roots $\alpha_i$ are \underline{simple}.
 The decomposition \eqref{any_roots} is called {\it reduced} if
 $l_C(s_{\alpha_1} s_{\alpha_2} \dots s_{\alpha_k}) = k$.

 \begin{lemma}{\cite[Lemma 3]{Ca72}}
  \label{lem_lin_indep}
   Let $\alpha_1, \alpha_2, \dots, \alpha_k \in \varPhi$.
   Then $s_{\alpha_1} s_{\alpha_2} \dots s_{\alpha_k}$ is reduced
   if and only if $\alpha_1, \alpha_2, \dots, \alpha_k$ are linearly
   independent. \qed
 \end{lemma}

\index{admissible diagram}
\index{diagram! - admissible diagram}
 A diagram $\Gamma$ is said to be {\it admissible}, see \cite[p. 7]{Ca72}, if
\begin{equation}
 \label{eq_def_adm}
 \begin{split}
 & (a)  \text{ The nodes of } \Gamma \text{ correspond to a set of linearly
  independent roots in } \varPhi. \\
 & (b) \text{
  If a subdiagram of } \Gamma \text{ is a cycle, then it contains an even number of nodes. } 
 \end{split}
\end{equation}

 \index{Carter diagram}
 \index{diagram! - Carter diagram}
 \index{$\Gamma$-associated set of roots}
 \index{$\Gamma$-associated element}
 \index{$S$-associated element}
 \index{$S$ (set of linearly independent roots)}
 \index{$S_{\alpha}$ (subset $\{\alpha_i \mid i = 1,\dots,k\}$)}
 \index{$S_{\beta}$ (subset $\{\beta_j \mid j = 1,\dots,h\}$)}
 \index{$\alpha$-set}
 \index{$\beta$-set}
 \index{Weyl group}
  Any admissible diagram $\Gamma$ is said to be a {\it Carter diagram}
  if any edge connecting a pair of roots $\{\alpha, \beta\}$ with
  inner product $(\alpha, \beta) > 0$ (resp. $(\alpha, \beta) < 0$)
  is drawn as dotted (resp. solid) edge. Let
\begin{equation}
   \label{eq_alpha_bet}
     S = \{\alpha_1, \alpha_2, \dots, \alpha_k, \beta_1, \beta_2, \dots, \beta_h \}
 \end{equation}
 be any set of linearly independent, not necessarily simple, roots associated with $\Gamma$,
 where roots of the set $S_{\alpha} := \{\alpha_i \mid i = 1,\dots,k\}$ are
 mutually orthogonal, roots of the set $S_{\beta} := \{\beta_j \mid j = 1,\dots,h\}$ are also mutually orthogonal.
 According to (\ref{eq_def_adm}$(a)$), there exists the set \eqref{eq_alpha_bet} of linearly independent roots.
 Thanks to (\ref{eq_def_adm}$(b)$), such a partitioning into the sum of two mutually orthogonal sets $S_{\alpha}$
 and $S_{\beta}$ is possible.
 The set $S$ is said to be a {\it $\Gamma$-associated} set of roots. Let
 \begin{equation}
   \label{two_invol}
      w = w_1{w}_2, \quad \text{ where } \quad
      w_1 = s_{\alpha_1} s_{\alpha_2} \dots s_{\alpha_k}, \quad
      w_2 = s_{\beta_1} s_{\beta_2} \dots s_{\beta_h}.
 \end{equation}
 Since $S$ is linearly independent, the decomposition
 \eqref{two_invol} is reduced, see Lemma \ref{lem_lin_indep}, and
 $k + h = l_C(w)$.
 The element $w$ is said to be {\it $\Gamma$-associated},
 and also {\it $S$-associated}. The decomposition \eqref{two_invol}
 is said to be a {\it bicolored decomposition}.  The set of roots $S_{\alpha}$
 (resp. $S_{\beta}$) is said to be the {\it $\alpha$-set}
 (resp. {\it $\beta$-set}) of roots  corresponding to the bicolored
 decomposition \eqref{two_invol}.

 \begin{remark}[On the Carter theorem]
 {\rm
   {\it \lq\lq For any  $w \in W$, there is a Carter diagram $\Gamma$ such that $w$ is
   the $\Gamma$-associated element.\rq\rq} The existence of such a Carter diagram means that any $w \in W$
   can be decomposed into the product of two involutions.  This statement is equivalent to the well-known
   theorem proved by Carter in \cite[Theorem C]{Ca72}. Carter's proof
   is based on the description of all conjugacy classes for any
   Weyl group. For the Weyl groups $G_2$, $F_4$, $E_6$, $E_7$, $E_8$, the conjugacy classes are
   presented in \cite[Tables 7 -- 11]{Ca72}.
   For another proof of the Carter theorem based on the
   classification of so-called {\it linkage diagrams} extending Carter
   diagrams, see \cite{St10}, \cite{St11}.
   The classification of linkage diagrams presented in \cite{St10}
   does not involve the use of a computer.
 }
 \end{remark}

 \index{primitive element}
 \index{semi-Coxeter element}
\begin{remark}[On the semi-Coxeter conjugacy class]
  \label{rem_primitive}
{\rm
  A conjugacy class of $W$ that can be determined by a connected
Carter diagram with number of nodes equal to the rank of $W$ is
called a {\it semi-Coxeter conjugacy class}, \cite{CE72}. The
conjugacy class, whose Carter diagram is the Dynkin diagram of $W$,
is called the {\it Coxeter conjugacy class}.
 Any representative of a semi-Coxeter conjugacy class
 is called a {\it primitive element}, or a {\it semi-Coxeter element} \cite{KP85}, \cite{B89}, \cite{St11}.
 }
\end{remark}

\subsubsection{Connection diagrams}
  \label{sec_conn_diagr}

 \index{connection diagram}
 \index{$(\Gamma, o)$ (connection diagram)}
 \index{diagram! - connection diagram}
 \index{semi-Coxeter element associated with  the connection diagram $(\Gamma, o)$}
 \index{$(\Gamma, o)$-semi-Coxeter element}
 \index{Dynkin diagram associated with a Weyl group (customary use)}
 \index{Dynkin diagram representing a conjugacy class}
 \index{Carter diagram looked like a Dynkin diagram}
 \index{CCl (conjugacy class)}
 \index{solid edge}
 \index{dotted edge}
 \index{acute angle between roots}
 \index{obtuse angle between roots}

  Let $\Gamma$ be the diagram characterizing
  connections between roots of a certain set $S$ of linearly
  independent and not necessarily simple roots,
  $o$
  be the order of reflections in the decomposition \eqref{any_roots}.
  The pair $(\Gamma, o)$ is said to be a {\it connection diagram}.
  We omit indicating order $o$ in the description of the connection diagram if
  the order of reflections in the decomposition \eqref{any_roots} is clear.
  The connection diagram determines the element $w$ (and its inverse $w^{-1}$)
  obtained as the product of all reflections associated with the diagram,
  while the order $o$ (resp. $o^{-1}$) describes the order of reflections
  in the decomposition of $w$ (resp. $w^{-1}$).
  Similarly to Remark \ref{rem_primitive}, we call $w$  the {\it semi-Coxeter element
  associated with  the connection diagram} $(\Gamma, o)$, or {\it $(\Gamma, o)$-semi-Coxeter element}.

  Connection diagrams describe connected sets with any cycles in the diagram, not necessarily even.
  Converting a Carter diagram $\Gamma_1$ into another
  Carter diagram $\Gamma_2$ we sometimes get connection  diagrams (but not Carter diagrams),
  and the \lq\lq{evenness}\rq\rq of cycles is violated during this conversion, see \S\ref{sec_get_rid}.

  The Dynkin diagrams in this paper appear in two ways: (1)
  associated with Weyl groups (customary use); (2) representing
  conjugacy classes (CCl), i.e, a Carter diagram which looks
  like (and actually is) a Dynkin diagram. In a few cases Dynkin diagrams represent
  two (and even three!) conjugacy classes.

  For the Carter diagrams and connection diagrams,
  we distinguish acute and obtuse angles between roots. Recall that a {\it solid edge} indicates an obtuse angle
  between the roots exactly as for Dynkin diagrams.
  A {\it dotted edge} indicates an acute angle between the roots considered, see \S\ref{sec_surprising}
  and Fig. \ref{two_types_diagram}.

\subsubsection{The $4$-cycles in Carter diagrams and connection diagrams}
\begin{figure}[h]
\centering
\includegraphics[scale=0.77]{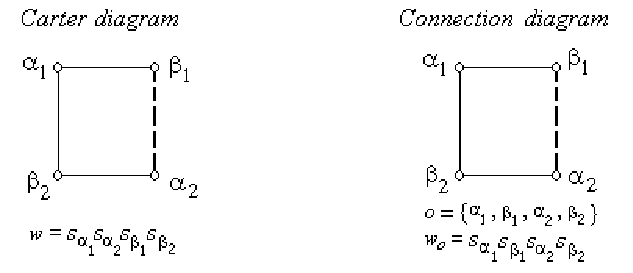}
\caption{The Carter diagram $D_4(a_1)$ and connection
diagram $(\mathcal{G}_4, o)$}
\label{two_types_diagram}
\end{figure}

The Carter diagram for a $4$-cycle in Fig. \ref{two_types_diagram}
determines a bicolored decomposition:
 \index{diagram! - characterizing pattern $D_4(a_1)$}
 \index{trivial order}
 \begin{equation*}
   \label{index_c}
    w = s_{\alpha_1}s_{\alpha_2}s_{\beta_1}s_{\beta_2}.
    \vspace{2mm}
 \end{equation*}
 Here, $w$ is the $D_4(a_1)$-associated element, where $D_4(a_1)$ denotes a $4$-cycle, see \cite{Ca72}.
 The diagrams in Fig. \ref{two_types_diagram} differ in the order.
 In the case of Carter diagrams, the order is {\it trivial} (related with a given bicolored decomposition)
 and we do not indicate it.
 The connection diagram in Fig. \ref{two_types_diagram}
 has order $o = \{ \alpha_1, \beta_1, \alpha_2, \beta_2 \}$:
 \begin{equation}
   \label{index_omega}
    w_o = s_{\alpha_1}s_{\beta_1}s_{\alpha_2}s_{\beta_2}.
 \end{equation}
 In \eqref{index_omega}, $w_o$ is the $(\mathcal{G}_4, o)$-associated
 element, where $\mathcal{G}_4$ is a $4$-cycle.
 We will omit the index $o$ of the element $w_o$
 if the order $o$ is clear from the context.
\begin{remark}
  {\rm
   Hereafter, we suppose that every cycle contains only one dotted edge.
   Otherwise, we apply reflections  $\alpha \longmapsto -\alpha$.
   These operations do not change the element $w$ since
   $s_\alpha = s_{-\alpha}$. In this case,
   every dotted edge with an endpoint vertex $\alpha$ is changed to the solid one,
   the cycle with all edges solid cannot occur, see Lemma \ref{lem_must_dotted}.
   Note also that the dotted edge can be moved
   to any other edge of the cycle by means of reflections.
  }
\end{remark}

 \index{semi-Coxeter element}

The semi-Coxeter elements generated by reflections
 $s_{\alpha_1}, s_{\alpha_2}, s_{\beta_1}, s_{\beta_2}$
 constitute exactly two conjugacy classes, $w$ and $w_o$
 being their representatives. In the basis
 $\{ \alpha_1, \alpha_2, \beta_1, \beta_2 \}$, we have:
\begin{equation}
   \label{2_Coxeter_elem}
   w =   \left ( \begin{array}{cccc}
          1  & 0  &  -1 & -1 \\
          0  & 1  &   1 & -1 \\
          1  & -1 &  -1 &  0 \\
          1  & 1  &  0 &  -1
          \end{array} \right ), \qquad \qquad \qquad
   w_o =   \left ( \begin{array}{cccc}
          0  & 1  &  0 & 0 \\
          1  & 0  &   -1 & -1 \\
          0  & 0  &  0 &  1 \\
          1  & 1  &  0 &  -1
          \end{array} \right )
 \end{equation}
 and their characteristic polynomials are:
 \begin{equation}
   \label{char_2_elem}
    \chi(w) = x^4 + 2x^2 + 1, \qquad \qquad \qquad
    \chi(w_o) = x^4 + x^3 + x + 1.
 \end{equation}

\subsubsection{Transformation of $4$-cycles}

\index{diagram! - characterizing pattern $D_4$}
\index{$\stackrel{u}{\simeq}$ (conjugation $w \longrightarrow u^{-1}w{u}$)}

\begin{figure}[h]
\centering
\includegraphics[scale=0.7]{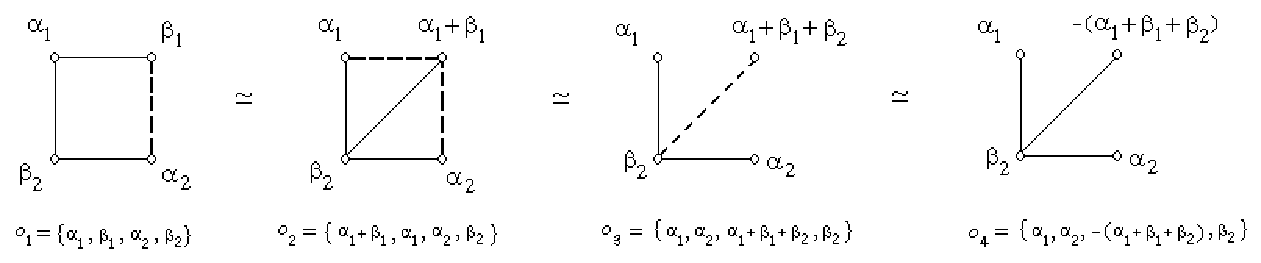}
\caption{Eliminating of the cycle. The equivalence of ($\mathcal{G}_4$, $o_1$) and ($D_4$, $o_4$)}
\label{equiv_diagrams_2}
\end{figure}

Denote by $\stackrel{u}{\simeq}$ the conjugation $w \longrightarrow u^{-1}w{u}$.
Let us transform the element $w_o$ from \eqref{index_omega}:

\begin{equation}
  \label{transf_4cycle}
  \begin{split}
     w_o ~=~
        & s_{\alpha_1}s_{\beta_1}s_{\alpha_2}s_{\beta_2} =
        s_{\alpha_1 + \beta_1}s_{\alpha_1}s_{\alpha_2}s_{\beta_2}
        \stackrel{s_{\alpha_1 + \beta_1}}{\simeq}
        s_{\alpha_1}s_{\alpha_2}s_{\beta_2}s_{\alpha_1 + \beta_1} =
        \\
        & s_{\alpha_1}s_{\alpha_2}s_{\alpha_1 + \beta_1 + \beta_2}s_{\beta_2} =
          s_{\alpha_1}s_{\alpha_2}s_{-({\alpha_1 + \beta_1 + \beta_2})}s_{\beta_2}.
   \end{split}
\end{equation}
 We have:
\begin{equation}
 \Small
   \begin{split}
   & (\alpha_1 + \beta_1 + \beta_2, \alpha_1) =
     (\alpha_1, \alpha_1)  + (\beta_1, \alpha_1)  + (\beta_2, \alpha_1) =
     1 - \frac{1}{2} - \frac{1}{2} = 0,   \\
   & (\alpha_1 + \beta_1 + \beta_2, \alpha_2) =
     (\beta_1, \alpha_2) + (\beta_2, \alpha_2) = - \frac{1}{2} +  \frac{1}{2} = 0, \\
   & (\alpha_1 + \beta_1 + \beta_2, \beta_2) = 1 - \frac{1}{2} = \frac{1}{2}.
   \end{split}
\end{equation}
  Hence, the roots $\{ \alpha_1, \alpha_2,  -(\alpha_1 + \beta_1 + \beta_2)\}$
  are mutually orthogonal, so
  in eq. \eqref{transf_4cycle}, we obtained a bicolored decomposition.
  Thus, the connection diagram
  ($\mathcal{G}_4$, $o_1$) is reduced to the Carter diagram ($D_4$, $o_4$),
  which is also the Dynkin diagram $D_4$, see Fig.~\ref{equiv_diagrams_2}.
  That is why, in eq.~\eqref{char_2_elem}
  the characteristic polynomial $\chi(w_{o_1}) = x^4 + x^3 + x + 1 = (x^3 + 1)(x + 1)$
  is equal to the characteristic polynomial of the $D_4$-associated
  element, see \cite[Table 3]{Ca72}, or \cite[Table 1]{St08}.

\subsection{Equivalence}
 \label{sec_equiv_diagrams}

\subsubsection{Three transformations} 

 \index{equivalence transformations}
 \index{three equivalence transformations}
 \index{unipotent conjugacy classes}
 Talking about a certain diagram $\Gamma$ we actually have in
 mind a set of roots with orthogonality relations as it is prescribed by the diagram $\Gamma$. 
 We try to find some common properties of sets of roots (from the root systems associated with
 the simple Lie algebras) and diagrams associated with these sets. 
 These diagrams are not necessarily Dynkin diagrams since sets of roots we study
 are not necessarily sets of simple roots and are not
 root subsystems. 
 We use the term \lq{\lq}Dynkin diagram\rq\rq to describe connected sets of
 linearly independent simple roots in the root system. 
 Similarly, \lq{\lq}Carter diagrams{\rq\rq} describe connected sets of linearly independent roots,
 not necessarily simple, and such that any cycle is even. 

 Same as Dynkin diagrams  describe simple Lie algebras, Carter diagrams describe
 conjugacy classes in Weyl groups\footnotemark[1].
 \footnotetext[1]{There is an interesting relationship between conjugacy classes in the Weyl group
 and unipotent conjugacy classes in the reductive algebraic group studied by G.~Lusztig in \cite{Lu08},
 \cite{Lu11}, \cite{Lu12}, \cite{Lu12a}.}
  First of all, in this paper we will see that
 any Carter diagram with cycles of any length can be transformed into an {\it equivalent Carter diagram}
 with cycles of length $4$. The equivalence of connection diagrams (and, in particular, of Carter diagrams)
 is discussed in \S\ref{sec_eq_conn_diagr}.
 Below we consider a rather natural set of three transformations operating on connection diagrams:
 Similarities, conjugations and $s$-permutations.

\begin{figure}[h]
\centering
\includegraphics[scale=0.7]{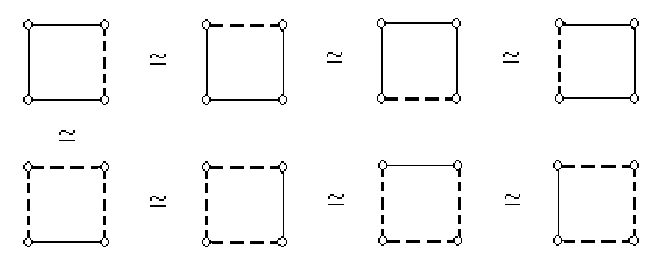}
\caption{Eight similar $4$-cycles equivalent to $D_4(a_1)$}
\label{8-equiv_4cycles}
\end{figure}

 \index{similarity of Carter diagrams}
 \underline{Similarity.}
 This is replacing a root with the opposite one:
\begin{equation}
  \label{eq_1_equiv}
   \alpha \longmapsto -\alpha.
\end{equation}
 Two connection diagrams obtained from each other by a sequence
 of reflections \eqref{eq_1_equiv}, are said to be {\it similar} connection diagrams, see Fig. \ref{8-equiv_4cycles}.
 An equivalence transformation of connection diagrams obtained by a sequence of reflections \eqref{eq_1_equiv}
 is said to be a {\it similarity transformation} or {\it similarity}.
\index{similar Carter diagrams}

 \index{projective roots}
 By applying similarity \eqref{eq_1_equiv} any solid edge with an endpoint vertex being $\alpha$
 can be changed to a dotted one and vice versa; this does not change, however, the corresponding reflection\footnotemark[1]:
 \footnotetext[1]{Recently, for the case of Dynkin diagrams, Dynkin and Minchenko in  \cite{DM10}
 considered so-called {\it projective roots}, i.e., pairs $\{\alpha, -\alpha\}$.}
\begin{equation*}
   s_{\alpha} = s_{-\alpha}.
\end{equation*}

  \begin{remark}[On trees]
   \label{rem_tree}
  {\rm For the set
  $\{ \alpha_1, \dots, \alpha_i, \alpha_{i+1}, \dots \alpha_n \}$
  forming a tree, we may assume that, up to the similarity, all non-zero inner products
  $(\alpha_i, \alpha_j)$ are negative. Indeed, if $(\alpha_i, \alpha_j) >  0$,
  we apply similarity transformation $\alpha_j \longmapsto -\alpha_j$, consider all
  inner products $(\alpha_k, \alpha_j) > 0$ and repeat similarity transformations
  $\alpha_k  \longmapsto -\alpha_k$ if necessary. This process converges
  since the diagram is a tree.
  }
 \end{remark}

 \index{conjugation of Carter diagrams}
  \underline{Conjugation.} Let $(\Gamma, o)$ be a connection diagram,
  $S = \{\alpha_1 \dots, \alpha_n \}$ a $\Gamma$-associated set.
  A conjugation sends all roots of a given set $S$ to another set by means of
  the same element $T$ from the Weyl group:
\begin{equation}
  \label{eq_1_equiv_2}
   \alpha_1 \longmapsto T\alpha_1, \quad \dots, \quad  \alpha_n \longmapsto T\alpha_n.
\end{equation}
Then
\begin{equation*}
   s_{\alpha_i} \longmapsto s_{T\alpha_i} = T^{-1}s_{\alpha_i}T \text{ for } i = 1,\dots,n, \quad \text{ and } \quad
     \prod\limits_{i} s_{\alpha_i} \longmapsto \prod\limits_{i} s_{T\alpha_i}. \\
\end{equation*}
 If $o = \{\alpha_{i_1} \dots, \alpha_{i_n} \}$ is an order of roots, then
 the conjugation $T$ sends $o$ into $T{o} = \{T\alpha_{i_1} \dots, T\alpha_{i_n} \}$.
 Let $\Gamma$ be a Carter diagram.
 Since $T$ preserves relations between roots, $T$ preserves $\Gamma$
 and the $\Gamma$-associated conjugacy class.

 \index{$s$-permutation of Carter diagrams}
\underline{$s$-Permutation.}
 The \lq{\lq}evenness\rq\rq of cycles is not violated by similarities
 \eqref{eq_1_equiv} and conjugations \eqref{eq_1_equiv_2}.
 It can be violated by the transformations of the third type, we call them {\it $s$-permutations}:
 \begin{equation}
  \label{eq_1_equiv_4}
  \begin{split}
   s_{\alpha}s_{\beta} =
    \begin{cases}
      s_{\beta}s_{\alpha + \beta} = s_{\alpha + \beta}s_{\alpha}  & \text{ for } (\alpha, \beta)  < 0,
     \\
      s_{\beta}s_{\alpha - \beta} = s_{\alpha - \beta}s_{\alpha} & \text{ for } (\alpha, \beta)  > 0,
     \\
      s_{\beta}s_{\alpha} & \text{ for } (\alpha, \beta) = 0.
    \end{cases}
  \end{split}
\end{equation}
Relations \eqref{eq_1_equiv_4} take place only for a simply-laced
connection between vertices $\alpha$ and $\beta$. In the general
case, the  $s$-permutations satisfy the following relation:
\begin{equation*}
   s_{\alpha}s_{\beta} = s_{\beta}s_{s_{\beta}(\alpha)} = s_{s_{\alpha}(\beta)}s_{\alpha}.
\end{equation*}
 Clearly, the $s$-permutation \eqref{eq_1_equiv_4} is non-trivial only
 if $\alpha$ and $\beta$ are connected.
 A non-trivial $s$-permutation \eqref{eq_1_equiv_4} yelds a new set
 of roots in which $\alpha$ (or $\beta$) is changed to $\alpha + \beta$ or
 $\alpha - \beta$ according to whether the edge $\{\alpha, \beta\}$
 is solid or dotted. For the new set, we also draw the diagram
 which is not necessarily a Carter diagram anymore but is a certain connection
 diagram.

\index{similarity of Carter diagrams}
  The set of transformations \eqref{eq_1_equiv}, \eqref{eq_1_equiv_2} and \eqref{eq_1_equiv_4}
  operates on a connection diagram $\Gamma$ and the root subset $S$ associated with the diagram $\Gamma$.
  Similarities \eqref{eq_1_equiv} change a given connection diagram to a similar one;
  conjugations \eqref{eq_1_equiv_2} preserve connection diagrams;
  $s$-permutations \eqref{eq_1_equiv_4} essentially change connection diagrams.
  However, both similarities and $s$-permutations preserve the element $w$ associated with the given diagram.
  Transformations \eqref{eq_1_equiv}, \eqref{eq_1_equiv_2} and \eqref{eq_1_equiv_4}
  \underline{preserve the conjugacy class} containing $w$ and also \underline{preserve the linear independence}
  of the roots constituting the subset $S$.

 \index{equivalent connection diagrams}
 \subsubsection{The equivalence of connection diagrams}
   \label{sec_eq_conn_diagr}
   Similarities, conjugations and $s$-permutations are said to be
   {\it equivalence transformations}. The equivalence transformations preserve associated conjugacy classes.
   Connection diagrams $(\Gamma_1, o_1)$
   and $(\Gamma_2, o_2)$ are said to be {\it equivalent} if for any
   $(\Gamma_1, o_1)$-associated element $w_1$, there exists
   a $(\Gamma_2, o_2)$-associated element $w'_2$ such that
   $w'_2$ can be obtained from $w_1$ by means of equivalence transformations,
   and for any $(\Gamma_2, o_2)$-associated element $w_2$, there exists
   a $(\Gamma_2, o_1)$-associated element $w'_1$ such that
   $w_2$ can be obtained from $w'_1$ by means of equivalence transformations.
   In this case, we will write
 \begin{equation*}
    (\Gamma_1, o_1) \simeq (\Gamma_2, o_2).
 \end{equation*}
 Such a definition of the equivalence of connection diagrams does not require
 the uniqueness of the conjugacy class associated
 with $\Gamma_1$ (resp. $\Gamma_2$). However, if one of diagrams $\Gamma_1$ and $\Gamma_2$
 determines a single conjugacy class, the same holds for another diagram.
 Indeed, let $\{w_1\}$ be a single $\Gamma_1$-associated
 conjugacy class and $w_2$, $w_2'$ be arbitrary $\Gamma_2$-associated elements, i.e.,
 $w_2 \simeq w_1$, and  $w'_2 \simeq w_1$. Then by transitivity,  we have $w_2 \simeq w_2'$.
 For example, it will be shown in \S\ref{th_get_rid} that
 \begin{equation}
   \label{eq_A_B_1}
   E_8(b_3) \simeq E_8(a_3), \quad E_7(b_2) \simeq E_7(a_2), \quad D_6(b_2) \simeq D_6(a_2),
   \quad E_8(b_5) \simeq E_8(a_5).
 \end{equation}
 By to Theorem \ref{th_uniq_diagr}, each of the diagrams $E_8(a_3)$, $E_7(a_2)$, $D_6(a_2)$ and $E_8(a_5)$
 determine a single conjugacy class. Therefore, the same holds for diagrams
 $E_8(b_3)$, $E_7(b_2)$, $D_6(b_2)$ and $E_8(b_5)$.

   Some of admissible and Carter diagram  may be equivalent to a connection diagram and vice versa.
   In \S\ref{sec_get_rid}, we use this fact in the process of excluding diagrams with cycles of length $>4$
   from Carter's list \cite[p. 10]{Ca72}, see Theorem \ref{th_get_rid}.
   We exclude a number of diagrams from possible candidates for the role of admissible or Carter diagram, since
   they have a subdiagram equivalent to an extended Dynkin diagram, a case which cannot be
   (Proposition \ref{prop_non_ext_Dynkin}, Lemma \ref{lem_2_impos_diagr}).

 \subsubsection{Two $\Gamma$-associated conjugacy classes}
   \label{sec_two_diff_calsses}
\begin{figure}[h]
\includegraphics[scale=1.8]{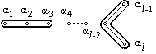}
 \caption{Elements $s_{\alpha_1}s_{\alpha_3}s_{\alpha_2}$ and
     $s_{\alpha_{l-1}}s_{\alpha_l}s_{\alpha_{l-2}}$ are not conjugate}
\label{Dl_w1_w2}
\end{figure}

 There exist $\Gamma$-associated elements
 $w_1$ and $w_2$ such that $w_1 \not\simeq w_2$. For example, the Carter diagram $A_3$
 determines two different conjugacy classes in $D_l$, see Fig.
 \ref{Dl_w1_w2}; for details, see \S\ref{sec_A3}.

\index{conjugate sets}
\index{non-conjugate $\Gamma$-associated sets}

 \subsubsection{Two non-conjugate $\Gamma$-associated sets}
   Let $S_1 = \{\varphi_1,\dots,\varphi_n\}$ and $S_2 = \{\delta_1,\dots,\delta_n\}$
   be two $\Gamma$-associated sets of roots. The sets $S_1$ and $S_2$ are
   said to be {\it conjugate} if there exists an element $T \in W$
   such that $T: \varphi_i \longmapsto \delta_i$ for $i = 1,\dots,n$. In this case, we write
 \begin{equation*}
    S_1 \simeq S_2 \text{ and } T{S}_1 = S_2.
 \end{equation*}
   Let $w_1$ (resp. $w_2$) be any $S_1$-associated (resp. $S_2$-associated) element.
   If $S_1 \simeq S_2$, then $w_1 \simeq w_2$.
~\\

\begin{figure}[h]
 \centering
\includegraphics[scale=0.9]{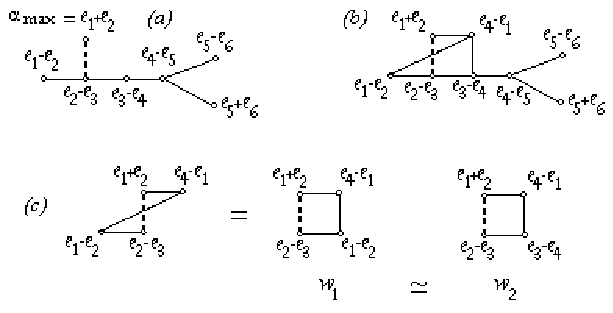}
\caption{Equivalence of the $\mathcal{C}_1$-associated element $w_1$ and
         the $\mathcal{C}_2$-associated element $w_2$}
\label{example_D6_2sq}
\end{figure}

   There exist, however, conjugate elements $w_1$ and $w_2$ such that $S_1 \not\simeq S_2$.
   Consider two $4$-cycles in $D_6$:
 \begin{equation*}
   \begin{split}
     & \mathcal{C}_1 = \{ e_1 + e_2,  e_4 - e_1, e_1 - e_2, e_2 - e_3  \}, \\
     & \mathcal{C}_2 = \{ e_1 + e_2,  e_4 - e_1, e_3 - e_4, e_2 - e_3  \}.
   \end{split}
 \end{equation*}
 These sets are non-conjugate: $\mathcal{C}_1 \not\simeq \mathcal{C}_2$, see Fig. \ref{example_D6_2sq}
 and \S\ref{sec_example_4cycles},
 but the $\mathcal{C}_1$-associated element $ w_1 = s_{e_1 + e_2}s_{e_1 - e_2}s_{e_4 - e_1}s_{e_2 - e_3}$
 and the $\mathcal{C}_2$-associated element $w_2 = s_{e_1 + e_2}s_{e_3 - e_4}s_{e_4 - e_1}s_{e_2 - e_3}$ are
 conjugate.

\subsection{Regular extensions}
   \label{sec_regular}
\subsubsection{Single-track, mirror and threefold extensions}
  \label{sec_1_2_option}
   For any diagram $\Gamma$, let $\alpha$ be a certain new vertex
   connected to $\Gamma$ at several vertices.
   We say that  $\widetilde\Gamma$ is the {\it extension of $\Gamma$ by the vertex $\alpha$}.
   If $\alpha$ is considered as a root, not a vertex, the phrase
   \lq\lq{\it extension of $\Gamma$ by the root $\alpha$}\rq\rq~means that we
   attach $\alpha$ to a $\Gamma$-associated subset $S$ and we get the $\widetilde\Gamma$-associated subset
   $\widetilde{S} = S \cup\alpha$.
   Recall that the extension $\widetilde\Gamma$ of $\Gamma$ is said to be a
   {\it regular} one if the initial diagram $\Gamma$ and the extended
   diagram $\widetilde\Gamma$ are both connected Carter diagrams and $\alpha$ is
   connected to $\Gamma$ at not more than three vertices.
   For the regular extension, we write
 \index{regular extensions}
 \index{regular extension $\Gamma < \widetilde{\Gamma}$ (or $\Gamma \stackrel{\alpha}{<} \widetilde{\Gamma}$)}
 \index{extension! - regular extension}
 \index{socket of the extension $\Gamma <\widetilde\Gamma$}
 \index{socket $\Gamma \rightpitchfork \widetilde\Gamma$ or $\Gamma ~\stackrel{\{v_1,\dots,v_n\}}{\rightpitchfork}~  \widetilde\Gamma$}
 \index{socket}
 \index{single-track extension}
 \index{extension! - single-track extension}
 \index{mirror map $T$}
 \index{$\{\Gamma\}$ (set of vertices of $\Gamma$)}
  \begin{equation}
    \label{eq_ext_diagr}
     \Gamma < \widetilde\Gamma \quad \text{ or } \quad \Gamma \stackrel{\alpha}{<} \widetilde\Gamma.
  \end{equation}
   If $\widetilde\Gamma$ extends $\Gamma$ by adding the vertex
   $\alpha$, the set of vertices connected to $\alpha$ is called the
   {\it socket of the extension $\Gamma <\widetilde\Gamma$}.
   If no ambiguity arises, this set will be simply called a {\it socket}.
   We denote the socket by
  \begin{equation}
    \label{eq_socket}
      \Gamma \rightpitchfork \widetilde\Gamma \quad \text{ or }
      \quad \Gamma ~\stackrel{\{v_1,\dots,v_n\}}{\rightpitchfork}~  \widetilde\Gamma,
  \end{equation}
  where $\{v_1,\dots,v_n\}$ is the set of socket's vertices.
  For the extension \eqref{eq_ext_diagr}, there can be several options for sockets.
  In this case, we call the
  extension \eqref{eq_ext_diagr} the {\it multi-option extension}.
  If only one socket is available for the extension \eqref{eq_ext_diagr}
  of the given $\Gamma$, we call this extension the {\it single-track extension},
  see Table \ref{tab_single_track}. Let
  \begin{equation}
    \label{eq_G1_G2_ext}
       \Gamma ~\stackrel{\mathscr{P}_1}\rightpitchfork~ \widetilde\Gamma_1,
        \qquad \Gamma ~\stackrel{\mathscr{P}_2}\rightpitchfork~ \widetilde\Gamma_2
  \end{equation}
  be two extensions with different sockets $\mathscr{P}_1$ and $\mathscr{P}_2$.
  For any diagram $\Gamma$, denote by $\{\Gamma\}$ the set of vertices of $\Gamma$.
  Let $M$ act on the set of vertices $\{\widetilde\Gamma_1\} \bigcup \{\widetilde\Gamma_2\}$
  as the mirror map, translating one socket into another, i.e,
\index{mirror extensions}
\index{extension! - mirror extensions}
\index{threefold extensions}
\index{extension! - threefold extension}
  \begin{equation}
    \begin{array}{lll}
      & M:\widetilde\Gamma_1 \longrightarrow \widetilde\Gamma_2,
        & M:\widetilde\Gamma_2 \longrightarrow \widetilde\Gamma_1, \\
      & M:\mathscr{P}_1 \longrightarrow \mathscr{P}_2,
        & M:\mathscr{P}_2 \longrightarrow \mathscr{P}_1, \quad M^2 = I.\\
     \end{array}
  \end{equation}
 Extensions \eqref{eq_G1_G2_ext} are said to be a pair of {\it mirror extensions},
 see Tables \ref{tab_mirror_1},  \ref{tab_mirror_2}.
 Let
  \begin{equation}
    \label{eq_G1_G2_G3_ext}
    \begin{split}
      & \Gamma ~\stackrel{\mathscr{P}_1}\rightpitchfork~ \widetilde\Gamma_1,
        \qquad \Gamma ~\stackrel{\mathscr{P}_2}\rightpitchfork~ \widetilde\Gamma_2,
        \qquad \Gamma ~\stackrel{\mathscr{P}_3}\rightpitchfork~ \widetilde\Gamma_3
    \end{split}
  \end{equation}
  be three extensions with different sockets $\mathscr{P}_1$, $\mathscr{P}_2$, $\mathscr{P}_3$.
  If any pair of extensions from \eqref{eq_G1_G2_G3_ext} forms a pair of mirror extensions,
  the triple \eqref{eq_G1_G2_G3_ext} is said to be {\it threefold extensions},
  see Table \ref{tab_threefold}.

\begin{remark}
  \label{rem_abuse_not}
{\rm
 By abuse of notation, we sometimes write
\begin{equation*}
    mirror~extensions~\Gamma < \widetilde\Gamma \quad (\text{resp. } threefold~extensions~\Gamma < \widetilde\Gamma)
\end{equation*}
instead of
\begin{equation*}
    mirror~extensions~\Gamma < \widetilde\Gamma_1 ~and~ \Gamma < \widetilde\Gamma_2 \quad
    (\text{resp. } threefold~extensions~\Gamma < \widetilde\Gamma_1,
     \Gamma < \widetilde\Gamma_2 ~and~ \Gamma < \widetilde\Gamma_3)
\end{equation*}
 Here, $\widetilde\Gamma_1$, $\widetilde\Gamma_2$ (resp. $\widetilde\Gamma_1$, $\widetilde\Gamma_2$, $\widetilde\Gamma_3$,)
 are isomorphic diagrams which differ by roots extending $\Gamma$. In this case, the notation $\Gamma < \widetilde\Gamma$
 omits the roots corresponding to the vertices.
 \qed
}
\end{remark}

  \begin{remark}
  \label{rem_all_ext}
  {\rm
   Tables \ref{tab_single_track}, \ref{tab_mirror_1},   \ref{tab_mirror_2}, \ref{tab_threefold}
   do not contain {\it all} possible single-track, mirror and threefold extensions, but only a certain set of extensions
   sufficient to get every Carter diagram  from $\mathsf{C4} \coprod \mathsf{DE4}$ with number of vertices
   $l \geq 5$ from a certain smaller Carter diagram, see Lemma \ref{lem_123_socket}.
   } \qed
 \end{remark}

  Single-track, mirror and threefold extensions describe extensions
  from the point of view of the number of sockets available. This number is said to be the {\it sockets number}.
  In addition, it is important for us to describe extensions from the point of view of the number of vertices in the socket,
  this is said to be the {\it pinholes number}. The socket with the pinholes number equal to $v$ is said to be
  a {\it $v$-socket}.
  For example, we can get the diagram $E_7(a_2)$
  by the single-track extension $E_6 < E_7(a_2)$ with the $2$-socket or by mirror extensions
  $D_6(a_2) < E_7(a_2)$ with $1$-sockets, or by mirror extensions
  $E_6(a_2) < E_7(a_2)$ with $1$-sockets, see Fig. \ref{ext_E7a2}, Table \ref{tab_2option_E7a2}.

 \index{$v$-socket}
 \index{sockets number}
 \index{$\mathsf{P}1$-, $\mathsf{P}2$-, and $\mathsf{P}3$-extensions}
 \index{$\mathsf{P}v$-extensions $(n=1,2,3)$}
 \index{extension! - with a $v$-socket}
 \index{extension! - $\mathsf{P}1$-, $\mathsf{P}2$-, and $\mathsf{P}3$-extensions}
 \index{extension! - $\mathsf{P}v$-extensions $(v=1,2,3)$}
 \index{sockets number}
 \index{pinholes number}
 \index{$1$-socket}

 \subsubsection{$\mathsf{P}1$-, $\mathsf{P}2$-, and $\mathsf{P}3$-extensions}

  The regular extension with the $v$-socket is said to be {\it $\mathsf{P}v$-extension}.
  Any Carter diagram can be obtained as an extension of a smaller Carter diagram
  by means of one of $\mathsf{P}v$-extensions, where $v=1,2,3$.

  \begin{lemma}
    \label{lem_123_socket}
    For any Carter diagram $\widetilde{\Gamma}$ from $\mathsf{C4} \coprod \mathsf{DE4}$ with number of vertices
    $l \geq 5$, there exists a smaller Carter diagram $\Gamma$ such that,
    from the point of view of the \underline{sockets number},
    any regular extension $\Gamma < \widetilde{\Gamma}$ is
    a certain single-track, mirror or threefold extension while from
    the point of view of the \underline{pinholes number} any regular extension is
    a certain $\mathsf{P}1$-, $\mathsf{P}2$-, $\mathsf{P}3$-extension.
  \end{lemma}
   \PerfProof Lemma follows from comparing Table \ref{tab_class_Carter} with
   Tables \ref{tab_single_track}, \ref{tab_mirror_1},   \ref{tab_mirror_2}, \ref{tab_threefold},
   see Remark \ref{rem_all_ext}.
   \qed
~\\
 \index{$\mathsf{P}1$-extension of the $D$-joint type}
 \index{$\mathsf{P}1$-extension of the $A$-joint type}
 \index{extension! - $\mathsf{P}1$-extension of the $D$-joint type}
 \index{extension! - $\mathsf{P}1$-extension of the $A$-joint type}

  We will also distinguish two types of $\mathsf{P}1$-extensions to be used only
  in one specific case during the proof of Proposition \ref{prop_induct_step}.
  The $\mathsf{P}1$-extension with an $1$-socket belonging to $D_4(a_1)$
  will be called {\it $\mathsf{P}1$-extension of the $D$-joint type}.
  The $\mathsf{P}1$-extension with an $1$-socket belonging only to one edge
  will be called {\it $\mathsf{P}1$-extension of the $A$-joint type},
  see Tables \ref{tab_single_track}, \ref{tab_mirror_1}.

\subsubsection{Examples of single-track and mirror extensions}
\begin{figure}[h]
\centering
\includegraphics[scale=0.6]{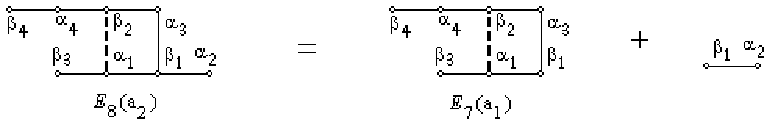}
 \caption{The single-track extension $E_7(a_1) \stackrel{\alpha_2}{<} E_8(a_2)$}
\label{ext_E7a1_to_E8a2}
\end{figure}

  In Fig. \ref{ext_E7a1_to_E8a2}, for $\widetilde\Gamma = E_8(a_2)$,
  we have the single-track extension $E_7(a_1) \stackrel{\alpha_2}{<}  E_8(a_2)$, the
  $1$-socket is as follows:
  \begin{equation*}
      E_7(a_1) ~\stackrel{\{\beta_1\}}\rightpitchfork~ E_8(a_2).
  \end{equation*}

  In Fig. \ref{multi_sockets}, for $\widetilde\Gamma = E_7(a_4)$, we have
  three pairs of mirror extensions. There is no any other single-track extension for getting $E_7(a_4)$.

 \begin{figure}[h]
\centering
\includegraphics[scale=0.7]{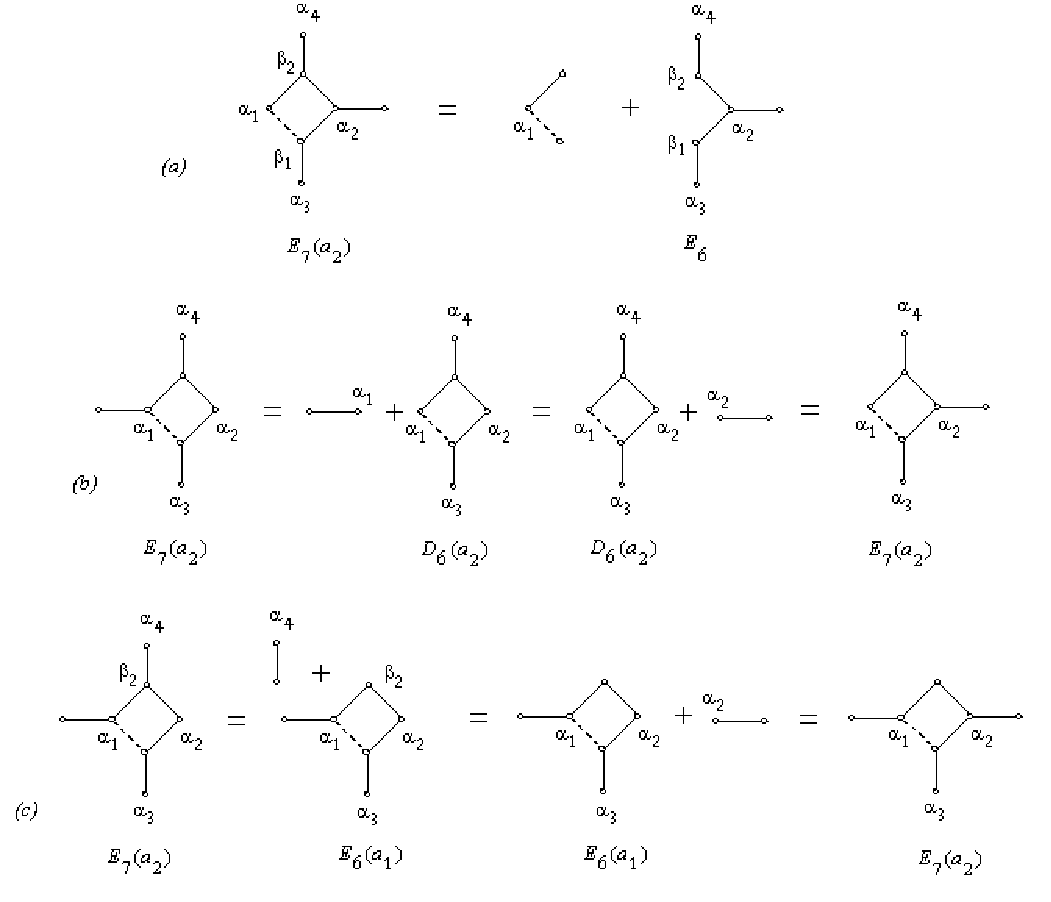}
 \caption{The single-track extension and mirror
          ~extensions $\Gamma < E_7(a_2)$, see Table \ref{tab_2option_E7a2}}
\label{ext_E7a2}
\end{figure}

\begin{table}[h]
 \Small
  \centering
  \renewcommand{\arraystretch}{1.4}
  \begin{tabular} {|c|c|c|c|}
  \hline
    $E_6 < E_7(a_2)$
        & \multicolumn{2}{|c|}{$E_6 ~\stackrel{\{\beta_1, \beta_2\}}\rightpitchfork~ E_7(a_2)$}
        & $2$-socket \\
  \hline
    $D_6(a_2) < E_7(a_2)$
        & $D_6(a_2) ~\stackrel{\{\alpha_1\}}\rightpitchfork~ E_7(a_2)$
        & $D_6(a_2) ~\stackrel{\{\alpha_2\}}\rightpitchfork~ E_7(a_2)$
                                                & $1$-sockets \\
  \hline
    $E_6(a_1) < E_7(a_2)$
        & $E_6(a_1) ~\stackrel{\{\alpha_2\}}\rightpitchfork~ E_7(a_2)$
        & $E_6(a_1) ~\stackrel{\{\beta_2\}}\rightpitchfork~ E_7(a_2)$
                                                & $1$-sockets \\
  \hline
\end{tabular}
  \vspace{2mm}
  \caption{The single-track extension and mirror extensions, Fig. \ref{ext_E7a2}}
  \label{tab_2option_E7a2}
\end{table}

\subsubsection{Single-track Condition, Mirror Condition and Threefold Condition}
  We consider the following three conditions which are fundamental for the proof of
  the uniqueness of the conjugacy classes, see Theorem  \ref{th_uniq_diagr}.
~\\

\begin{figure}[h]
\centering
\includegraphics[scale=0.8]{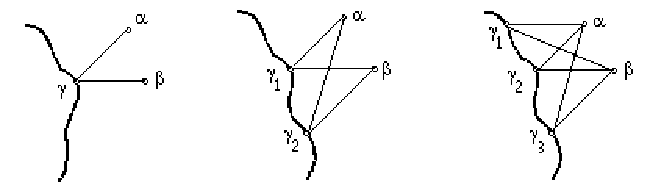}
 \caption{For $n=1,2,3$: The same $v$-socket for $\Gamma \stackrel{\alpha}{<} \widetilde\Gamma$ and $\Gamma
\stackrel{\beta}{<} \widetilde\Gamma$}
\label{the_same_socket}
\end{figure}

 \underline{Single-track Condition.}
  Let $\Gamma$ be a Carter diagram. Let $\Gamma \stackrel{\alpha}{<} \widetilde\Gamma$
  and $\Gamma \stackrel{\beta}{<} \widetilde\Gamma$ be
  two regular extensions with \underline{the same socket} $\mathscr{P}$, where $\alpha$ and $\beta$ are two
  different roots connected with socket $\mathscr{P}$.
  Let $w$ be a $\Gamma$-associated element.
  If $ws_\alpha \simeq ws_\beta$ for any such roots $\alpha$ and $\beta$,
  then we say that \emph{the Single-track Condition holds}, see Fig. \ref{the_same_socket}.
~\\

\index{Single-track Condition}
\index{$\widetilde\Gamma_L$ (left extension)}
\index{$\widetilde\Gamma_R$ (right extension)}
\index{extension! - left extension $\widetilde\Gamma_L$}
\index{extension! - right extension $\widetilde\Gamma_R$}
\index{$\widetilde\Gamma_L$-associated element}
\index{$\widetilde\Gamma_R$-associated element}
\index{mirror map $T$}
\index{Mirror Condition}

  \underline{Mirror Condition.}
      Let $\Gamma$ be a Carter diagram with the mirror map $M$, see \S\ref{sec_1_2_option}.
      Let $\mathscr{P}_1$ and $\mathscr{P}_2$ be two sockets for mirror extensions $\Gamma < \widetilde\Gamma$.
      We denote extensions of $\Gamma$ corresponding to $\mathscr{P}_1$ and $\mathscr{P}_2$
      by $\widetilde\Gamma_L$ (left extension) and $\widetilde\Gamma_R$ (right extension),
      see Fig. \ref{multi_sockets}, Fig. \ref{ext_E7a2}, Tables \ref{tab_mirror_1} -- \ref{tab_mirror_2}.
      Let $\widetilde{\Gamma}_L$ (resp. $\widetilde{\Gamma}_R$)
      be obtained from $\Gamma$ by adding a root $\alpha$ (resp. $\beta$)
      connected to $\mathscr{P}_1$ (resp. $\mathscr{P}_2$).
      Let $w$ be a certain $\Gamma$-associated element.
      If for any root $\alpha$ connected to $\mathscr{P}_1$,
      there exists a root $\beta$ connected to $\mathscr{P}_2$
      such that $w{s}_\alpha$ (resp. $w{s}_\beta$) is
      the $\widetilde\Gamma_L$-associated (resp. $\widetilde\Gamma_R$-associated) element,
      and $w{s}_\alpha \simeq w{s}_\beta$, then we say that \emph{the Mirror Condition holds}.
~\\

 \index{Threefold Condition}
 \index{triality}
 \index{Weyl group}

  \underline{Threefold Condition.}
      Let $\mathscr{P}_1$, $\mathscr{P}_2$ and $\mathscr{P}_3$ be three sockets for threefold extensions
      $\Gamma < \widetilde\Gamma$, see \S\ref{sec_1_2_option}.
      We denote the three extensions of $\widetilde\Gamma$ corresponding to $\mathscr{P}_1$, $\mathscr{P}_2$ and $\mathscr{P}_3$,
      by $\widetilde\Gamma_1$ $\widetilde\Gamma_2$ and $\widetilde\Gamma_3$, respectively,
      see Tables \ref{tab_threefold}.
      Let $\widetilde{\Gamma}_i$ (resp. $\widetilde{\Gamma}_j$), where $i,j \in \{1,2,3\}$ and $i \ne j$,
      be obtained from $\Gamma$ by adding roots $\alpha_i$ 
      connected to the socket $\mathscr{P}_i$. 
      Let $w$ be a $\Gamma$-associated element.
      If for any root $\alpha_i$ connected to $\mathscr{P}_i$,
      there exists a root $\alpha'_j$ connected to $\mathscr{P}_j$
      such that $w{s}_{\alpha_i}$ (resp. $w{s}_{\alpha'_j}$) is
      $\widetilde\Gamma_i$-associated (resp. $\widetilde\Gamma_j$-associated) element,
      and $w{s}_{\alpha_i} \simeq w{s}_{\alpha'_j}$,
      then we say that \emph{the Threefold Condition holds}\footnotemark[1].
      \footnotetext[1]{The special features of the Dynkin diagram $D_4$
       and objects associated with $D_4$ arise because the corresponding Weyl group
       has an outer automorphism of order three, see \cite{A96}.
       The simple Lie group Spin$(8)$ has the most symmetrical Dynkin diagram $D_4$.
       Outer automorphisms of Spin$(8)$ were discovered in 1925 by \'{E}lie Cartan, who called symmetries of $D_4$
       {\it triality}, see \cite{C25}.}

   \begin{proposition} 
    \label{prop_extensions}
  Let $\Gamma$ be a Carter diagram determining a single conjugacy class.
  Let $w_1$, $w_2$ be $\Gamma$-associated, and $w_1 \simeq w_2$.

  {\rm (i)} Let $\alpha$, $\beta$ be roots extending $\Gamma$ to $\widetilde{\Gamma}$ in \underline{the same socket}.
     (Note that $\Gamma < \widetilde\Gamma$ is not necessarily a single-track extension.)
     If the Single-track Condition holds, then
  \begin{equation}
     \label{eq_w1_w2_alp}
       w_1{s}_\alpha ~\simeq~ w_2{s}_\beta.
  \end{equation}
  In particular, it can be that $\alpha = \beta$.

  {\rm (ii)} Let $\Gamma < \widetilde\Gamma$ be a \underline{single-track extension}.
     If the Single-track Condition holds, then $\widetilde{\Gamma}$ determines a single conjugacy class.

  {\rm (iii)} Let $\Gamma < \widetilde\Gamma_1$ and $\Gamma < \widetilde\Gamma_2$
    be \underline{mirror extensions} corresponding to
    sockets $\mathscr{P}_1$ and $\mathscr{P}_2$, respectively. Let $\alpha$ (resp. $\beta$)
    be a root connected to the socket $\mathscr{P}_1$ (resp. $\mathscr{P}_2$)
    in the diagram $\Gamma$. If the Single-track Condition and the Mirror Condition hold
    for $\Gamma$, then
  \begin{equation}
     \label{eq_w1_w2_alp_2}
       w_1  s_\alpha ~\simeq~ w_2  s_\beta.
  \end{equation}
    In other words, any two $\widetilde{\Gamma}$-associated elements $w_1{s}_\alpha$ and $w_2{s}_\beta$
    are conjugate, i.e., $\widetilde{\Gamma}$ determines a single conjugacy class.

  {\rm (iv)} Let $\Gamma < \widetilde\Gamma_i$  be \underline{threefold extensions}
  corresponding to  sockets $\mathscr{P}_i$, where $i=1,2,3$.
  Let $\alpha_i$ be roots connected to sockets $\mathscr{P}_i$
  in the diagram $\Gamma$. If the Single-track Condition and the Threefold Condition hold for $\Gamma$,
  then
  \begin{equation}
     \label{eq_w1_w2_alp_3}
       w_1{s}_{\alpha_i} ~\simeq~ w_2{s}_{\alpha_j}.
  \end{equation}
    In other words, any two $\widetilde{\Gamma}$-associated elements $w_1{s}_{\alpha_i}$ and $w_2{s}_{\alpha_j}$
    are conjugate, i.e., $\widetilde{\Gamma}$ determines a single conjugacy class.
  \end{proposition}

   \PerfProof (i) Since  $P{w}_1{P}^{-1} = w_2$, then
  \begin{equation*}
         w_2{s}_\alpha =  P{w_1}{P}^{-1}s_\alpha \simeq w_1{P}^{-1}s_\alpha{P} = w_1{s}_{P\alpha}.
  \end{equation*}
         Since the Single-track Condition holds for $w_1$, then $w_1  s_{P\alpha} \simeq w_1  s_\beta$
         and $w_2  s_\alpha \simeq w_1  s_\beta$, i.e.,
         \eqref{eq_w1_w2_alp} holds.

       (ii) Let $\Gamma < \widetilde\Gamma$ be a single-track
       extension with the socket $\mathscr{P}$, let $\widetilde{w}$, $\widetilde{v}$ be two
       $\widetilde\Gamma$-associated elements.
       Then $\widetilde{w}$ (resp. $\widetilde{v}$) is the product of a certain
       $\Gamma$-associated element $w_1$ (resp. $w_2$) and some
       reflection $s_\alpha$  (resp. $s_\beta$), where $\alpha$ (resp. $\beta$)
       is connected to the same socket $\mathscr{P}$.
       By (i) we have
  \begin{equation*}
       \widetilde{w} = w_1  s_\alpha \simeq w_2  s_\beta = \widetilde{v}.
   \end{equation*}

       (iii) By the Mirror Condition, for any $\beta$ connected to $\mathscr{P}_1$,
   there exists a root $\gamma$ connected to $\mathscr{P}_2$
   such that $w{s}_\gamma \simeq w{s}_\beta$. If $\alpha$ is another root connected to the socket $\mathscr{P}_2$
   in $\Gamma$, then by the Single-track Condition we have
  \begin{equation}
     \label{eq_w1_w2_alp3}
      w{s}_\alpha ~\simeq~ w{s}_\gamma ~\simeq~ w{s}_\beta.
  \end{equation}
  Then
  \begin{equation*}
      w_1{s}_\alpha \quad\stackrel{by \eqref{eq_w1_w2_alp3}}\simeq\quad w_1{s}_\beta
                    \quad\stackrel{by \eqref{eq_w1_w2_alp}}\simeq\quad  w_2{s}_\beta.
  \end{equation*}


\begin{table}[h]
\scriptsize
  \centering
  \renewcommand{\arraystretch}{1.3}
  \begin{tabular} {|c|c|c|}
  \hline
  \hline
  \multicolumn{3}{|c|}{$\mathsf{P}1$-extensions $\Gamma \stackrel{\alpha_2}{<} \widetilde\Gamma$ \qquad ($D$-joint type)} \\
  \hline
    $E_7(a_1) < E_8(a_2)$ &
     $\begin{array}{c} \\ \qquad \qquad \includegraphics[scale=0.7]{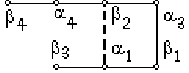}  \qquad \qquad  \end{array}$ &
     $\begin{array}{c} \\ \qquad \qquad \includegraphics[scale=0.7]{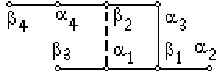} \qquad     \end{array}$ \\
    $E_7(a_2) < E_8(a_3)$ &
     $\begin{array}{c} \quad \includegraphics[scale=0.7]{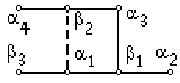}   \end{array}$ &
     $\begin{array}{c} \quad \includegraphics[scale=0.7]{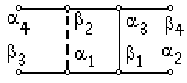}  \end{array}$ \\
    $E_7(a_3) < E_8(a_5)$ &
     $\begin{array}{c} \quad \includegraphics[scale=0.7]{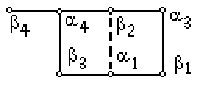}   \end{array}$ &
     $\begin{array}{c} \includegraphics[scale=0.7]{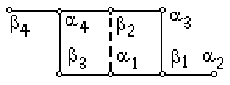}   \end{array}$ \\
    $D_l(a_1) < D_{l+1}(a_2)$ &
    $\begin{array}{c} \qquad \includegraphics[scale=0.7]{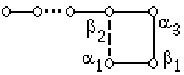}  \end{array}$ &
    $\begin{array}{c} \quad \includegraphics[scale=0.7]{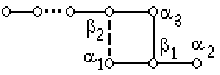} \end{array}$ \\
  \hline
  \multicolumn{3}{|c|}{$\mathsf{P}1$-extensions $\Gamma \stackrel{\varphi}{<} \widetilde\Gamma$ \qquad ($A$-joint type)} \\
  \hline
  $D_l(a_k) < D_{l+1}(a_k)$ &
    $\begin{array}{c} \\ \includegraphics[scale=0.5]{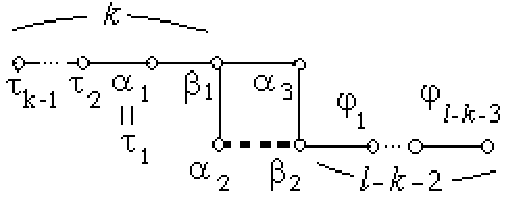}    \end{array}$ &
    $\begin{array}{c}  \\ \includegraphics[scale=0.5]{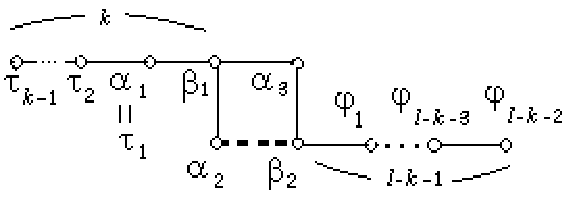} \\
    \quad \{ \delta = \varphi_{l-k-3}, ~\varphi = \varphi_{l-k-2} \} \end{array}$ \\
  $D_l < D_{l+1}$ &
    $\begin{array}{c}  \includegraphics[scale=1.1]{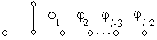}   \end{array}$ &
    $\begin{array}{c}  \includegraphics[scale=1.1]{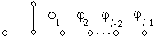}   \\
    \qquad \qquad \{ \delta = \varphi_{l-2}, ~\varphi = \varphi_{l-1} \} \end{array}$ \\
  $E_7 < E_8$ &
    $\begin{array}{c}  \includegraphics[scale=0.5]{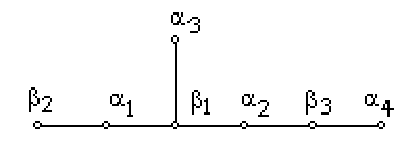}   \end{array}$ &
    $\begin{array}{c}  \includegraphics[scale=0.5]{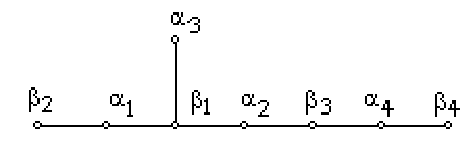} \\
       \qquad \qquad \{ \delta = \alpha_4, ~\varphi = \beta_4 \} \end{array}$ \\
 \hline
 \hline
  \multicolumn{3}{|c|}{$\mathsf{P}2$-extensions $\Gamma \stackrel{\alpha_2}{<} \widetilde\Gamma$} \\
  \hline
    $E_6 < E_7(a_2)$ &
      $\begin{array}{c} \\ \includegraphics[scale=0.7]{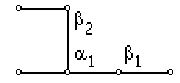}   \end{array}$ &
      $\begin{array}{c} \\ \qquad  \includegraphics[scale=0.7]{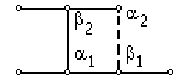}  \end{array}$ \\
    $E_7(a_3) < E_8(a_6)$ &
     $\begin{array}{c} \quad \includegraphics[scale=0.7]{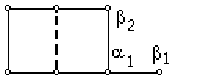}   \end{array}$ &
     $\begin{array}{c} \quad \includegraphics[scale=0.7]{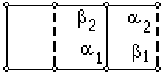}   \end{array}$ \\
   $E_7(a_1) < E_8(a_4)$ &
     $\begin{array}{c} \quad \qquad \includegraphics[scale=0.7]{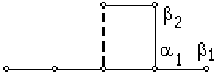} \qquad \qquad   \end{array}$ &
     $\begin{array}{c} \qquad \qquad \includegraphics[scale=0.7]{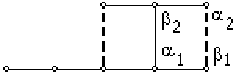} \qquad  \end{array}$ \\
 \hline
 \hline
  \multicolumn{3}{|c|}{$\mathsf{P}3$-extensions $\Gamma \stackrel{\alpha}{<} \widetilde\Gamma$} \\
  \hline
  $E_7(a_4) < E_8(a_8)$ &
    $\begin{array}{c}  \\ \includegraphics[scale=1.2]{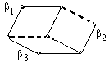}   \end{array}$ &
    $\begin{array}{c}  \\ \includegraphics[scale=1.2]{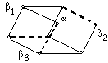}   \end{array}$ \\
  $E_7(a_2) < E_8(a_7)$ &
    $\begin{array}{c} \qquad \qquad  \includegraphics[scale=0.7]{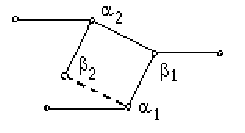} \qquad  \end{array}$ &
    $\begin{array}{c} \qquad \qquad   \includegraphics[scale=0.7]{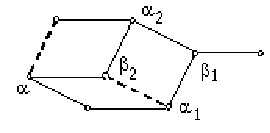} \quad \end{array}$ \\
 \hline
\end{tabular}
  \vspace{2mm}
   \caption{\hspace{3mm}
     Single-track extensions $\Gamma < \widetilde\Gamma$}
  \label{tab_single_track}
\end{table}

 \index{single-track $\mathsf{P}1$-extension ($D$-joint type)! - $E_7(a_1) < E_8(a_2)$}
 \index{single-track $\mathsf{P}1$-extension ($D$-joint type)! - $E_7(a_2) < E_8(a_3)$}
 \index{single-track $\mathsf{P}1$-extension ($D$-joint type)! - $E_7(a_3) < E_8(a_5)$}
 \index{single-track $\mathsf{P}1$-extension ($D$-joint type)! - $D_l(a_1) < D_{l+1}(a_2)$}
 \index{single-track $\mathsf{P}1$-extension ($A$-joint type)! - $D_l(a_k) < D_{l+1}(a_k)$}
 \index{single-track $\mathsf{P}1$-extension ($A$-joint type)! - $D_l < D_{l+1}$}
 \index{single-track $\mathsf{P}1$-extension ($A$-joint type)! - $E_7 < E_8$}
 \index{single-track $\mathsf{P}2$-extension! - $E_6 < E_7(a_2)$}
 \index{single-track $\mathsf{P}2$-extension! - $E_7(a_3) < E_8(a_6)$}
 \index{single-track $\mathsf{P}2$-extension! - $E_7(a_1) < E_8(a_4)$}
 \index{single-track $\mathsf{P}3$-extension! - $E_7(a_4) < E_8(a_8)$}
 \index{single-track $\mathsf{P}3$-extension! - $E_7(a_2) < E_8(a_7)$}


\index{mirror $\mathsf{P}1$-extension ($D$-joint type)! - $D_5(a_1) < E_6(a_1)$}
\index{mirror $\mathsf{P}1$-extension ($D$-joint type)! - $D_6(a_1) < E_7(a_1)$}
\index{mirror $\mathsf{P}1$-extension ($D$-joint type)! - $D_7(a_1) < E_8(a_1)$}
\index{mirror $\mathsf{P}1$-extension ($D$-joint type)! - $D_4(a_1) < D_5(a_1) (D_5(a_1) \subset E_l)$}
\index{mirror $\mathsf{P}1$-extension ($D$-joint type)! - $D_4(a_1) < D_5(a_1) (D_5(a_1) \subset D_l, E_l)$}

\index{mirror $\mathsf{P}1$-extension ($A$-joint type)! - $D_5 < E_6$}
\index{mirror $\mathsf{P}1$-extension ($A$-joint type)! - $E_6 < E_7$}
\index{mirror $\mathsf{P}1$-extension ($A$-joint type)! - $D_{2k+2}(a_k) < D_{2k+3}(a_k)$}

\begin{table}[h]
 \scriptsize
  \centering
  \renewcommand{\arraystretch}{1.1}
  \begin{tabular} {|c|c|c|}
   \hline
   \hline
  \multicolumn{3}{|c|}{Mirror $\mathsf{P}1$-extensions, $\Gamma \stackrel{\varphi}{<} \widetilde{\Gamma}_L$,
  $\Gamma \stackrel{\gamma}{<} \widetilde{\Gamma}_R$ \qquad ($D$-joint type)} \\
  \hline
     $D_5(a_1) < E_6(a_1)$
       & $\begin{array}{c} \\ \includegraphics[scale=1.3]{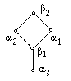}   \end{array}$
       & $\begin{array}{c} \\ \includegraphics[scale=1.3]{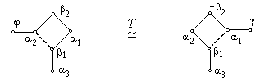}  \end{array}$  \\
       & & \\
     $D_6(a_1) < E_7(a_1)$ & $\begin{array}{c} \includegraphics[scale=1.3]{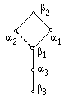}   \end{array}$
       & $\begin{array}{c} \includegraphics[scale=1.3]{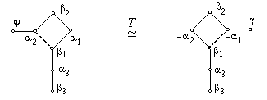}  \end{array}$ \\
       & & \\
     $D_7(a_1) < E_8(a_1)$ & $\begin{array}{c} \includegraphics[scale=1.3]{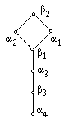}  \end{array}$
       & $\begin{array}{c} \includegraphics[scale=1.3]{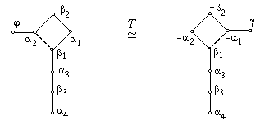}  \end{array}$ \\
       & & \\
      $\begin{array}{c}
       D_4(a_1) < D_5(a_1), \\
       (D_5(a_1) \subset E_l)
      \end{array}$
       & $\begin{array}{c} \includegraphics[scale=1.4]{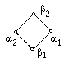}  \end{array}$
       & $\begin{array}{c} \includegraphics[scale=0.6]{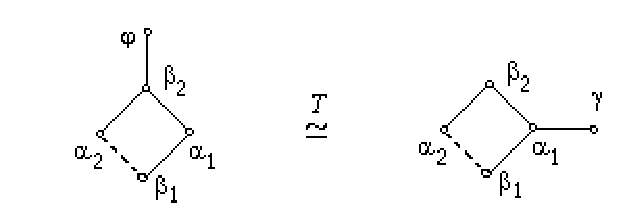} \end{array}$ \\
       & & \\
     $\begin{array}{c}
       D_4(a_1) < D_5(a_1), \\
       (D_5(a_1) \subset D_l, E_l)
      \end{array}$
       & $\begin{array}{c} \includegraphics[scale=1.4]{D4a1_2.eps}  \end{array}$
       & $\begin{array}{c} \quad  \includegraphics[scale=0.8]{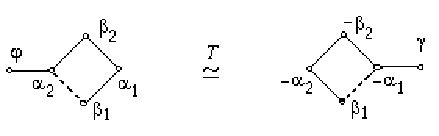}   \end{array}$ \\
       & & \\
  \hline
  \multicolumn{3}{|c|}{Mirror $\mathsf{P}1$-extensions, $\Gamma \stackrel{\varphi}{<} \widetilde{\Gamma}_L$,
  $\Gamma \stackrel{\gamma}{<} \widetilde{\Gamma}_R$ \qquad ($A$-joint type)} \\
  \hline
      & & \\
     $D_5 < E_6$
       & $\begin{array}{c} \includegraphics[scale=1.4]{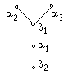}  \end{array}$
       & $\begin{array}{c} \includegraphics[scale=0.6]{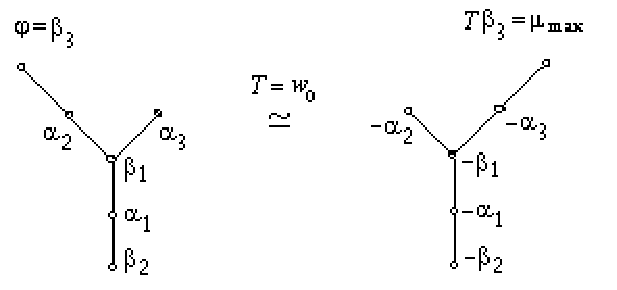} \end{array}$  \\
      & & \\
     $E_6 < E_7$
       & $\begin{array}{c} \includegraphics[scale=1.4]{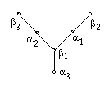}  \end{array}$
       & $\begin{array}{c} \includegraphics[scale=0.6]{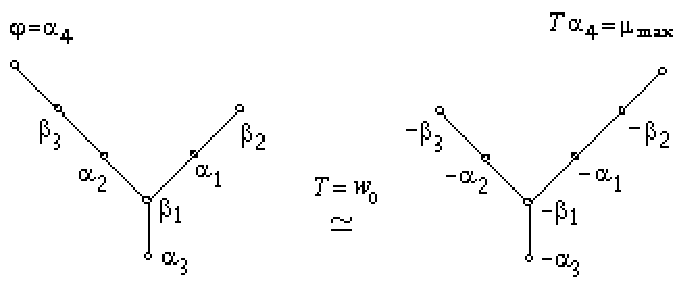}  \end{array}$ \\
      & & \\
     $\begin{array}{c} D_{2k+2}(a_k) < \\
       D_{2k+3}(a_k)
       \end{array}$
       & $\begin{array}{c} \includegraphics[scale=1.2]{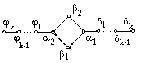}   \end{array}$
       & $\begin{array}{c} \includegraphics[scale=1.2]{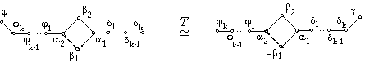}   \end{array}$  \\
  \hline
\end{tabular}
   \caption{Mirror extensions, $\Gamma \stackrel{\varphi}{<} \widetilde{\Gamma}_L$,
  $\Gamma \stackrel{\gamma}{<} \widetilde{\Gamma}_R$, where $\gamma = T\varphi$
   (or $\mu_{max} = T\varphi$), see Proposition \ref{prop_mirror_ext}}
  \label{tab_mirror_1}
\end{table}

\index{mirror $\mathsf{P}2$-extension! - $D_5(a_1) < E_6(a_2)$}
\index{mirror $\mathsf{P}2$-extension! - $D_6(a_1) < E_7(a_3)$}
\index{mirror $\mathsf{P}2$-extension! - $D_7(a_1) < E_8(a_4)$}

\index{mirror $\mathsf{P}3$-extension! - $D_6(a_2) < E_7(a_4)$}

\begin{table}[h]
 \scriptsize
  \centering
  \renewcommand{\arraystretch}{1.2}
  \begin{tabular} {|c|c|c|}
  \hline
  \hline
  \multicolumn{3}{|c|}{Mirror $\mathsf{P}2$-extensions, $\Gamma \stackrel{\varphi}{<} \widetilde{\Gamma}_L$,
  $\Gamma \stackrel{\gamma}{<} \widetilde{\Gamma}_R$} \\
  \hline
     $D_5(a_1) < E_6(a_2)$
       & $\begin{array}{c} \\ \includegraphics[scale=1.3]{D5a1_toExt.eps}   \end{array}$
       & $\begin{array}{c} \\ \includegraphics[scale=1.3]{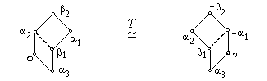}   \end{array}$ \\
     & & \\
     $D_6(a_1) < E_7(a_3)$
       & $\begin{array}{c} \includegraphics[scale=1.3]{D6a1_toExt.eps}   \end{array}$
       & $\begin{array}{c} \includegraphics[scale=1.3]{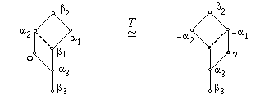}  \end{array}$ \\
     & & \\
     $D_7(a_1) < E_8(a_4)$
       & $\begin{array}{c} \includegraphics[scale=1.3]{D7a1_toExt.eps}  \end{array}$
       & $\begin{array}{c} \includegraphics[scale=1.3]{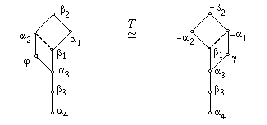} \end{array}$ \\
  \hline
  \hline
  \multicolumn{3}{|c|}{Mirror $\mathsf{P}3$-extensions, $\Gamma \stackrel{\varphi}{<} \widetilde{\Gamma}_L$,
  $\Gamma \stackrel{\gamma}{<} \widetilde{\Gamma}_R$} \\
  \hline
      & & \\
     $D_6(a_2) < E_7(a_4)$  & $\begin{array}{c} \includegraphics[scale=1.3]{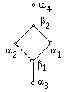}  \end{array}$
        & $\begin{array}{c} \includegraphics[scale=0.62]{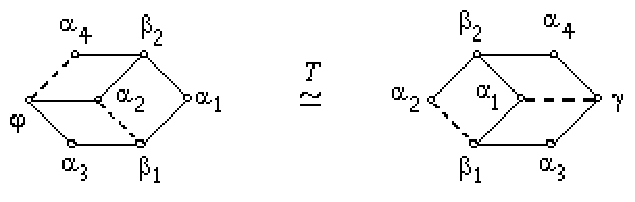}   \end{array}$  \\
  \hline
\end{tabular}
  \caption{Mirror extensions (cont.), $\Gamma \stackrel{\varphi}{<} \widetilde{\Gamma}_L$,
  $\Gamma \stackrel{\gamma}{<} \widetilde{\Gamma}_R$, where $\gamma = T\varphi$, see Proposition \ref{prop_mirror_ext} }
  \vspace{3cm}
  \label{tab_mirror_2}
\end{table}

\begin{table}[h]
 \scriptsize
  \centering
  \renewcommand{\arraystretch}{1.2}
  \begin{tabular} {|c|c|c|}
 \hline
 \hline
  \multicolumn{3}{|c|}{Threefold extensions $D_4 < D_5$} \\
 \hline
      & & \\
     $\begin{array}{c} D_4 < D_5 \\ (D_5 \subset E_l) \end{array}$
        & $\begin{array}{c} \includegraphics[scale=0.7]{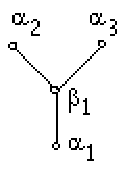}  \end{array}$
        & $\begin{array}{c} \includegraphics[scale=0.62]{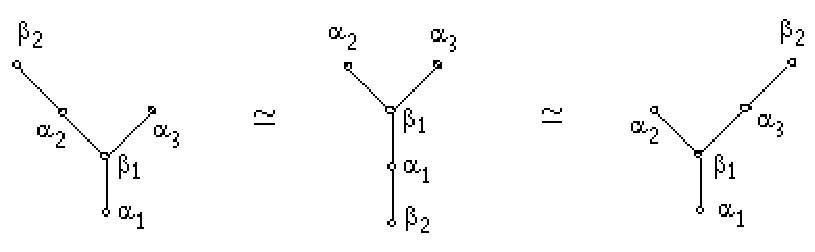} \end{array}$  \\
  \hline
  \end{tabular}
   \caption{Threefold extensions}
  \label{tab_threefold}
\end{table}

\index{class of Carter diagrams $\mathsf{C4}$! - $E_6(a_1)$, $E_6(a_2)$}
\index{class of Carter diagrams $\mathsf{C4}$! - $E_7(a_1)$, $E_7(a_2)$, $E_7(a_3)$, $E_7(a_4)$}
\index{class of Carter diagrams $\mathsf{C4}$! - $E_8(a_1)$, $E_8(a_2)$, $E_8(a_3)$, $E_8(a_4)$}
\index{class of Carter diagrams $\mathsf{C4}$! - $E_8(a_5)$, $E_8(a_6)$, $E_8(a_7)$, $E_8(a_8)$}
\index{class of Carter diagrams $\mathsf{C4}$! - $D_4(a_1)$, $D_5(a_1)$, $D_6(a_1)$, $D_6(a_2)$}
\index{class of Carter diagrams $\mathsf{C4}$! - $D_7(a_1)$, $D_7(a_2)$, $D_8(a_1)$, $D_8(a_2)$, $D_8(a_3)$}
\index{class of Carter diagrams $\mathsf{C4}$! - $D_n(a_k)$, ($n > 8$)}

\index{class of Carter diagrams $\mathsf{DE4}$! - $D_4$, $D_5$, $D_6$, $D_7$, $D_8$, $E_6$, $E_7$, $E_8$}
\index{class of Carter diagrams $\mathsf{DE4}$! - $D_n$, ($n > 8$)}

\clearpage
~\\

  (iv) By the Threefold Condition, for any $\alpha_i$ connected to $\mathscr{P}_i$,
   there exists a root $\alpha'_j$ connected to $\mathscr{P}_j$, where $i \neq j$,
   such that $w{s}_{\alpha_i} ~\simeq~ w{s}_{\alpha'_j}$. If $\beta_j$ another root connected to the socket $\mathscr{P}_j$
   in $\Gamma$, then by the Single-track Condition
  \begin{equation}
     \label{eq_w1_w2_alp3_2}
      w{s}_{\alpha_i} ~\simeq~ w{s}_{\alpha'_j} ~\simeq~ w{s}_{\beta_j}.
  \end{equation}
  Then
  \begin{equation*}
      w_1{s}_{\alpha_i} \quad\stackrel{by \eqref{eq_w1_w2_alp3_2}}\simeq\quad w_1{s}_{\beta_j}
                    \quad\stackrel{by \eqref{eq_w1_w2_alp}}\simeq\quad  w_2{s}_{\beta_j}.
  \end{equation*}
   \qed

 \begin{remark}{\rm
  By Proposition \ref{prop_extensions},
  in order to prove Theorem \ref{th_uniq_diagr}, it suffices for any Carter
  diagram $\widetilde\Gamma$, to find a smaller Carter diagram $\Gamma$
  such that $\widetilde\Gamma$ can be obtained by extending $\Gamma$ in one, two or three sockets,
  and extensions  $\Gamma < \widetilde\Gamma$ (see Remark \ref{rem_abuse_not}) have the following property:
~\\

    (a) if $\Gamma < \widetilde\Gamma$ is a single-track extension, then \underline{the Single-track Condition holds}.
~\\

    (b) if $\Gamma < \widetilde\Gamma_1$ and $\Gamma < \widetilde\Gamma_2$
    are mirror extensions ($\widetilde\Gamma_1 \simeq \widetilde\Gamma_2 \simeq \widetilde\Gamma$,
    see Remark \ref{rem_abuse_not}),
    then \underline{the Single-track} \underline{Condition and the Mirror Condition hold}.
~\\

    (c) if $\Gamma < \widetilde\Gamma_1$,  $\Gamma < \widetilde\Gamma_2$ and $\Gamma < \widetilde\Gamma_3$
    ($\widetilde\Gamma_i \simeq \widetilde\Gamma_j \simeq \widetilde\Gamma$, where $i,j=1,2,3$, see Remark \ref{rem_abuse_not}),
    are threefold extensions, then
    \underline{the Single-track Condition and the Threefold Condition hold}.
  }
 \end{remark}

\subsection{The main results}
  \label{sec_main_res}
 \index{bridge}

\subsubsection{Bridges}
\begin{figure}[h]
 \centering
 \includegraphics[scale=0.8]{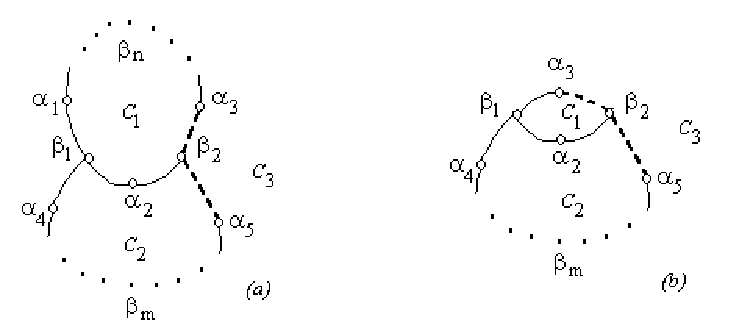}
 \caption{Intersecting cycles}
 \label{inters_cycles}
 \end{figure}

  Consider Carter diagrams containing intersecting cycles, i.e., cycles having a common path, see Fig. \ref{inters_cycles}$(a)$.
  There are three cycles in this figure, see eq. \eqref{eq_3_cycles}.
  To speak about intersecting cycles we choose the \underline{two shortest ones}.
  In the case of Fig. \ref{inters_cycles}$(a)$,
  we throw away from consideration the cycle $\mathcal{C}_3$, where
  \begin{equation}
   \label{eq_3_cycles}
   \begin{split}
    & \mathcal{C}_1 = \{\alpha_1, \beta_1, \alpha_2, \beta_2, \alpha_3, \beta_n\}, \\
    & \mathcal{C}_2 = \{\beta_1, \alpha_4,  \beta_m,   \alpha_5, \beta_2, \alpha_2\}, \\
    & \mathcal{C}_3 = \{\alpha_1, \beta_1, \alpha_4, \beta_m,   \alpha_5, \beta_2, \alpha_3, \beta_n\}.
   \end{split}
  \end{equation}
  Then $\mathcal{C}_1$ and $\mathcal{C}_2$
  have the common path $\{\beta_1, \alpha_2, \beta_2\}$. We denote this path by $\mathcal{C}_1 \cap \mathcal{C}_2$.
  It remains to consider the case, where $2$ cycles have the same length, see Fig. \ref{inters_cycles}$(b)$
  and eq. \eqref{eq_3_cycles_2}.
  \begin{equation}
   \label{eq_3_cycles_2}
   \begin{split}
    & \mathcal{C}_1 = \{\alpha_3, \beta_1, \alpha_2, \beta_2\}, \\
    & \mathcal{C}_2 = \{\beta_1, \alpha_4,  \beta_m,   \alpha_5, \beta_2, \alpha_2\}, \\
    & \mathcal{C}_3 = \{\alpha_3, \beta_1, \alpha_4, \beta_m, \alpha_5, \beta_2\}.
   \end{split}
  \end{equation}
  In Fig. \ref{inters_cycles}$(b)$, lengths of $\mathcal{C}_2$ and $\mathcal{C}_3$ coincide.
  Then the choice of $C_2$ or $C_3$ does not matter. The common path will be called a {\it bridge}.
  For the pair $\{\mathcal{C}_1, \mathcal{C}_2\}$ (resp. $\{\mathcal{C}_1, \mathcal{C}_3\}$), the
  bridge is as follows:
  \begin{equation*}
      \mathcal{C}_1 \cap \mathcal{C}_2 = \{\beta_1, \alpha_2, \beta_2\}, \qquad
      (\text{resp. } \mathcal{C}_1 \cap \mathcal{C}_3 = \{\beta_1, \alpha_3, \beta_2\}).
  \end{equation*}

\subsubsection{Refining the classification of Carter diagrams}
  In {\bf \S\ref{sec_classif}}, we revise Carter's classification of his diagrams.
  We add new arguments to get Carter's classification, namely we use
  the method of equivalent diagrams. The following proposition makes the way leading to the classification shorter:
~\\
~\\
  {\bf Proposition} (On intersecting cycles and bridges (Proposition \ref{prop_two_cycles})).
      {\rm (i)} {\em Let $\Gamma$ be any Carter diagram, or connection diagram,
         containing two cycles with bridge $\mathcal{P}$.
         Then $\mathcal{P}$ consists of exactly one edge.}

\begin{figure} [h]
 \centering
 \includegraphics[scale=0.55]{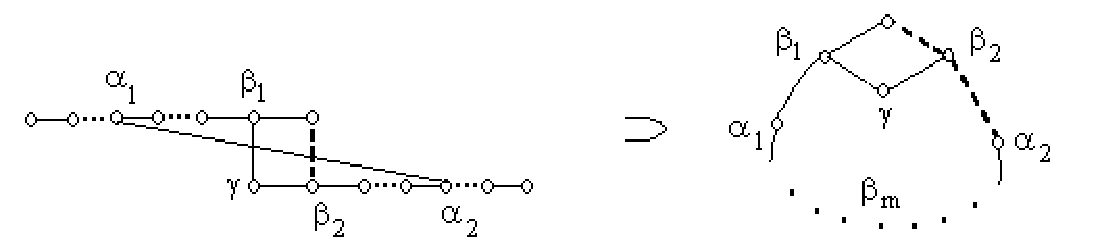}
 \caption{Two paths stemming from the opposite vertices of a $4$-cycle}
\label{2paths}
 \end{figure}

      {\rm (ii)} {\em Let $\Gamma$ be any Carter diagram.
         Let $P_1$, $P_2 \subset \Gamma$ be two paths stemming from the opposite vertices of
         a $4$-cycle in $\Gamma$; let $\alpha_1$ (resp. $\alpha_2$) be the vertex lying in $P_1$ (resp. $P_2$).
         The diagram obtained from $\Gamma$ by adding the edge $\{ \alpha_1, \alpha_2 \}$
         is not a Carter diagram, see Fig. $\ref{2paths}$.}

       {\rm (iii)} {\em Let $\Gamma$ be any Carter diagram containing two intersecting cycles.
         Then one of the cycles consists of $4$ vertices, and the other one can contain only $4$ or $6$
         edges.}

  We repeatedly use the following corollary describing basic
  restrictions for possible configurations of linearly independent
  quintuples of roots.
\index{quintuple of roots}
~\\
~\\
   {\bf Corollary} (On restrictions on quintuples of roots (Corollary \ref{cor_numb_ep})).
    {\rm (i)} {\em Let an $\alpha$-set contain $3$ roots $\{\alpha_1, \alpha_2, \alpha_3 \}$.
    There does not exist two non-connected roots $\beta$ and $\gamma$
    connected to every $\alpha_i$ so that the vectors of the quintuple $\{\alpha_1, \alpha_2, \alpha_3, \beta, \gamma \}$
    are linearly independent, see Fig.} \ref{3-cells}$(a)$.

    {\rm (ii)} {\em Let $\{\alpha_1, \beta_1, \alpha_2, \beta_2 \}$ be
    a square in a connection diagram. There does not exist
    a root $\gamma$ connected to all vertices of the square
    so that the vectors of the quintuple $\{ \alpha_1, \beta_1, \alpha_2, \beta_2, \gamma \}$
    are linearly independent, see Fig.} \ref{3-cells}$(b)$.
~\\

  The classification of the simply-laced Carter diagrams with cycles
  is given in \S\ref{sec_class_s-l}, see Table \ref{tab_class_Carter}.
~\\

\begin{table}[h]
  \centering
  \renewcommand{\arraystretch}{0.4}
  \begin{tabular} {|c|c|}
  \hline  
       & \\
       4-5 &
       $\begin{array}{c}
       \includegraphics[scale=0.5]{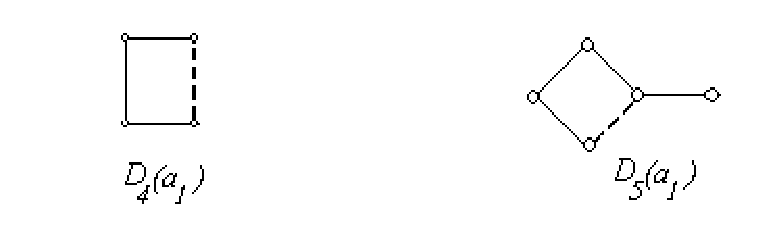}
       \end{array}$ \\
    \hline
       & \\
       6 &
       $\begin{array}{c}
       \includegraphics[scale=0.6]{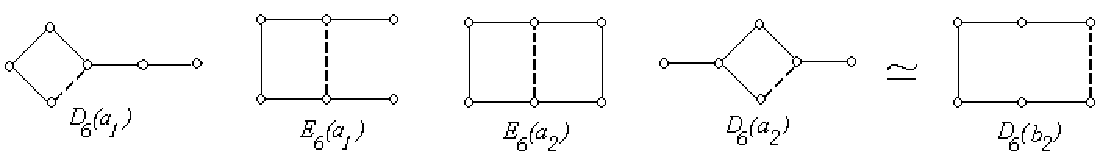}
       \end{array}$ \\
   \hline
       & \\
       7 &
       $\begin{array}{c}
       \includegraphics[scale=0.6]{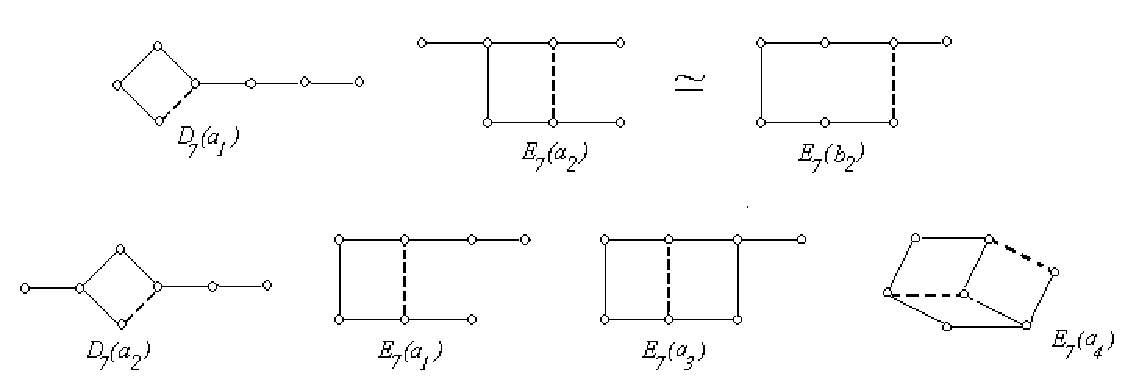}
       \end{array}$ \\
   \hline
       & \\
       8 &
       $\begin{array}{c}
       \includegraphics[scale=0.6]{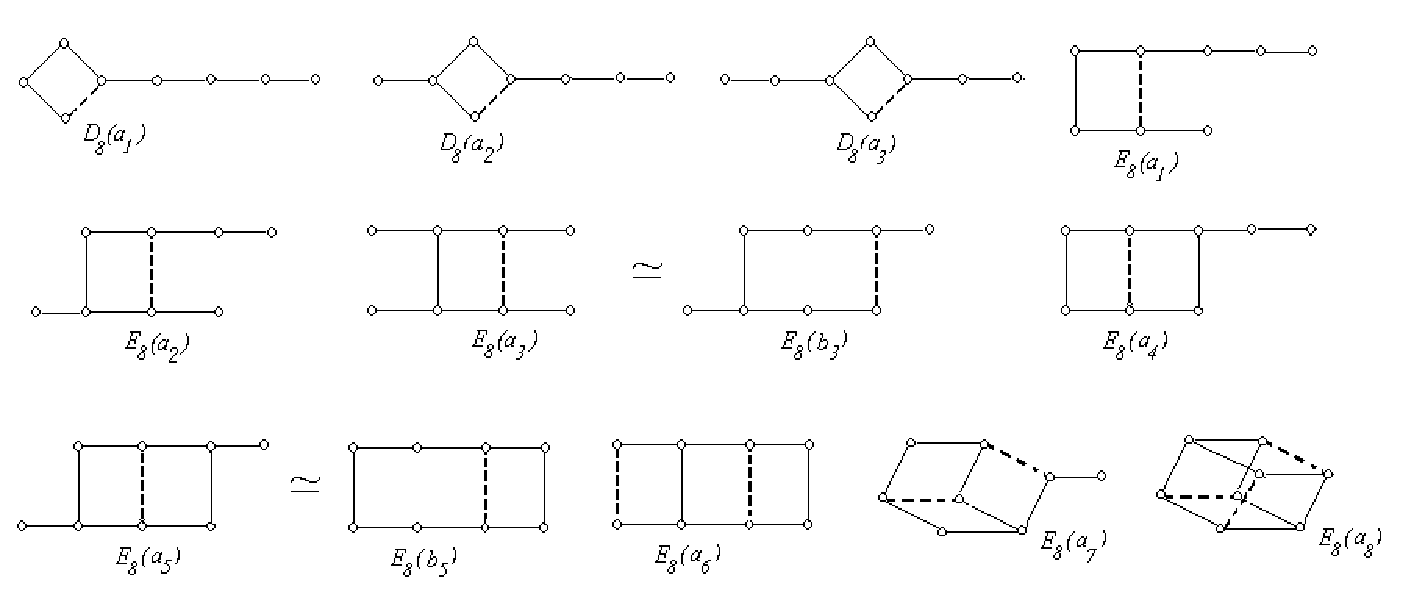}
       \end{array}$ \\
   \hline
       & \\
       $l > 8$ &
       $\begin{array}{c}
       \includegraphics[scale=0.6]{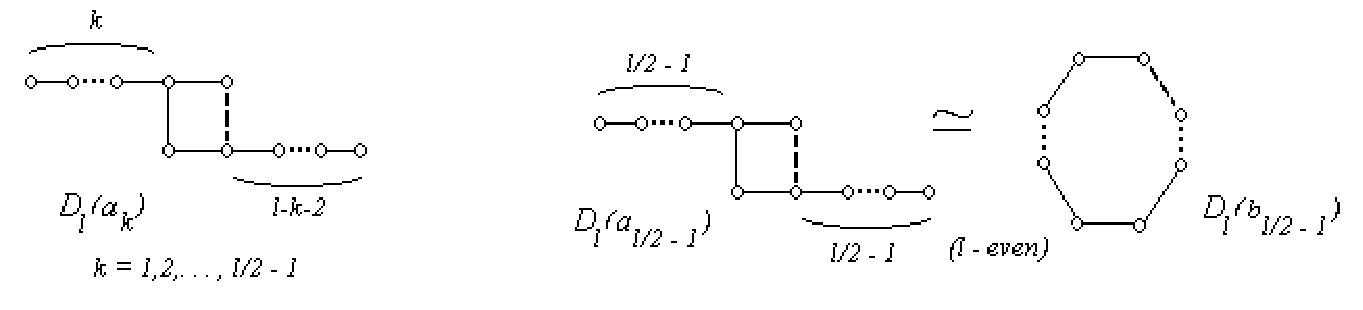}
       \end{array}$ \\
  \hline
\end{tabular}
  \vspace{2mm}
  \caption{The simply-laced Carter diagrams with cycles}
  \label{tab_class_Carter}
\end{table}
\vspace{2mm}
\index{explicit transformation}

  In {\bf \S\ref{sec_get_rid}}, we exclude all diagrams  with cycles of length $>4$.
  (This fact is taken into account in the classification of \S\ref{sec_class_s-l}).
~\\
~\\
  {\bf Theorem} (On exclusion of long cycles (Theorem  \ref{th_get_rid})).
   {\em Any Carter diagram containing $l$-cycles, where $l > 4$, is
   equivalent to another Carter diagram containing only $4$-cycles.}

 For the proof, we construct an explicit transformation
 of the element $w$ determined by the Carter diagram $\Gamma$ containing long cycles, i.e., $l$-cycles
 with $l > 4$, into an element $w'$ from the same conjugacy class described by another
 Carter diagram $\Gamma'$ containing only $4$-cycles.

 \begin{remark}
  \label{rem_excl_long_cycles}
 {\rm
   The possibility of excluding long cycles (Theorem \ref{th_get_rid}) and the uniqueness of the
   conjugacy class (Theorem \ref{th_uniq_diagr}) can be derived from \cite{Ca72} by means of
   classification of all conjugacy classes obtained by case-by-case computations.
   In the current paper, the proof of Theorem \ref{th_get_rid} and Theorem \ref{th_uniq_diagr}
   is also divided into a number of cases, however:

     (i) we do not use the classification of conjugacy classes,

     (ii) reasoning is short enough for each case, see \S\ref{sec_get_rid} and \S\ref{sec_uniq_C4DE4}.
  }
 \end{remark}

  \index{partial Cartan matrix $B_{\Gamma}$}
  \index{Cartan matrix! - partial}
  \index{quadratic form $\mathscr{B}_{\Gamma}$}
  \index{$B^{-1}_{\Gamma}$, the inverse of the partial Cartan matrix}

\subsubsection{Cases of linearly dependent roots}
  \index{linearly dependent root}
  \index{Coxeter element}
  \index{orbit of $\mu_{max}$}
  \index{$\mu_{max}$ (maximal root)}
  \index{$b^{\vee}_{\eta, \eta}$ (diagonal element of $B^{-1}_{\Gamma}$)}
  \index{Weyl group}
  In {\bf \S\ref{sect_basic_pat}}, we introduce {\it partial Cartan
  matrix} $B_{\Gamma}$ corresponding to the Carter diagram $\Gamma$,
  the corresponding quadratic form is denoted by
  $\mathscr{B}_{\Gamma}$. Unlike Cartan matrices,
  a partial Cartan matrix $B_{\Gamma}$ may have positive
  non-diagonal elements. These elements of $B_{\Gamma}$ are associated with dotted edges of the Carter diagram.
  Let $S$ be a certain $\Gamma$-associated subset of roots
  in $W$, see \S\ref{sec_adm_diagr}, and $L \subseteq V$ the subspace spanned by the subset $S$.
  For any $v \in L$, we have
  \begin{equation*}
       \mathscr{B}_{\Gamma}(v) = \mathscr{B}(v),
 \end{equation*}
  where $\mathscr{B}$ is the quadratic form associated with the Weyl group $W$,
  see Proposition \ref{restr_forms_coincide}.

  There are two frequently occurring cases where the vector $\gamma$
  is linearly dependent on roots $S$:

  $(a)$ The root $\gamma$  is connected with only one $\tau_i \in S$.
  Then, the following condition is necessary:
  \begin{equation*}
    b^{\vee}_{i,i} = 2,
  \end{equation*}
  where $b^{\vee}_{i,i}$ is the $i$th diagonal element of $B^{-1}_{\Gamma}$, see Remark \ref{rem_max_root}(ii).

  $(b)$ This case involves the use of the element $w_0$, the longest element in $W$, see \S\ref{sec_max_root}.
  Let $\Gamma$ be the Dynkin diagram $D_l$ or $E_l$,
  let $h$ be the Coxeter number, i.e., the order of the Coxeter element $w$
  for cases $\Gamma = D_l$ or $E_l$.
  Let $w \in W$ be any $\Gamma$-associated element, and
  \begin{equation*}
    S = \{\varphi_1,\dots,\varphi_k, \delta_1,\dots,\delta_m\}
  \end{equation*}
  the $\Gamma$-associated subset of linearly independent roots corresponding to the bicolored decomposition of $w$.
  Here, roots of the set $S_{\varphi} = \{\varphi_i \mid i = 1,\dots,k\}$
  (resp. $S_{\delta} = \{\delta_i \mid i = 1,\dots,m\}$)
  are mutually orthogonal, see \S\ref{sec_adm_diagr}.
~\\

\begin{table} 
 \Small
 \centering
 \renewcommand{\arraystretch}{1.1}
  \begin{tabular} {|c|c|c|c|}
  \hline
    \multicolumn{4}{|c|}{The orbit of $\mu_{max} \in \Phi(D_5)$ under
      $w = s_{\alpha_1}s_{\alpha_2}s_{\alpha_3}s_{\beta_1} \in  W(D_4)$} \\
    \multicolumn{4}{|c|}{The Coxeter number $h(D_4) = 6$} \\
  \hline
      $\mu_{\max}$   &
      $w_{\beta}\mu_{max}$ &
      $w_{\alpha}w_{\beta}\mu_{max} = w\mu_{max}$ &
      $w_{\beta}w\mu_{\max}$    \\
  \hline
      $\begin{array}{ccccc}
        1 & 2 & 2 & 1 \\
          & 1 &   &
      \end{array}$
&     $\begin{array}{ccccc}
        1 & 2 & 2 & 1 \\
          & 1 &   &
      \end{array}$
&     $\begin{array}{ccccc}
        1 & 2 & 1 & 1 \\
          & 1 &   &
      \end{array}$
&     $\begin{array}{ccccc}
        1 & 1 & 1 & 1 \\
          & 1 &   &
      \end{array}$
      \\
  \hline
      $w_{\alpha}w_{\beta}w\mu_{max} =w^2\mu_{max}$ &
      $w_{\beta}w^2\mu_{max}$  &
      $w_{\alpha}w_{\beta}w^2\mu_{max} = w^3\mu_{\max}$ &  \\
  \hline
     $\begin{array}{ccccc}
        0 & 1 & 1 & 1 \\
          & 0 &   &
      \end{array}$
&     $\begin{array}{ccccc}
        0 & 0 & 1 & 1 \\
          & 0 &   &
      \end{array}$
&      $\begin{array}{ccccc}
        0 & 0 & 0 & 1 \\
          & 0 &   &
      \end{array}$
&      \\
  \hline
 \end{tabular}

 \begin{tabular} {|c|c|c|c|c|}
    \hline
    \multicolumn{5}{|c|}{The orbit of $\mu_{max} \in \Phi(E_6)$ under
       $w = s_{\alpha_1}s_{\alpha_2}s_{\alpha_3}s_{\beta_1}s_{\beta_2} \in W(D_5)$} \\
    \multicolumn{5}{|c|}{The Coxeter number $h(D_5) = 8$} \\
  \hline
      $\mu_{\max}$   &
      $w_{\beta}\mu_{max}$ &
      $w\mu_{max}$ &
      $w_{\beta}w\mu_{\max}$ &
      $w^2\mu_{max}$ \\
  \hline
      $\begin{array}{ccccc}
        1 & 2 & 3 & 2 & 1 \\
        &     & 2 &     &
      \end{array}$
&     $\begin{array}{ccccc}
        1 & 2 & 3 & 2 & 1 \\
        &     & 2 &     &
      \end{array}$
&     $\begin{array}{ccccc}
        1 & 2 & 3 & 2 & 1 \\
        &     & 1 &     &
      \end{array}$
&     $\begin{array}{ccccc}
        1 & 2 & 2 & 2 & 1 \\
        &     & 1 &     &
      \end{array}$
&      $\begin{array}{ccccc}
        1 & 1 & 2 & 1 & 1 \\
        &     & 1 &     &
      \end{array}$ \\
  \hline
      $w_{\beta}w^2\mu_{max}$ &
      $w^3\mu_{\max}$   &
      $w_{\beta}w^3\mu_{max}$  &
      $w^4\mu_{\max}$ & \\
  \hline
      $\begin{array}{ccccc}
        0 & 1 & 1 & 1 & 1 \\
        &     & 1 &     &
      \end{array}$
&      $\begin{array}{ccccc}
        0 & 0 & 1 & 1 & 1 \\
        &     & 0 &     &
      \end{array}$
&     $\begin{array}{ccccc}
        0 & 0 & 0 &  1 & 1 \\
        &     & 0 &     &
      \end{array}$
&      $\begin{array}{ccccc}
        0 & 0 & 0 & 0 & 1 \\
        &     & 0 &     &
      \end{array}$
      & \\
  \hline
 \end{tabular}

\renewcommand{\tabcolsep}{0.3mm}
  \begin{tabular} {|c|c|c|c|c|}
    \hline
    \multicolumn{5}{|c|}{The orbit of $\mu_{max} \in \Phi(E_7)$ under
        $w = s_{\alpha_1}s_{\alpha_2}s_{\alpha_3}s_{\beta_1}s_{\beta_2}s_{\beta_3} \in W(E_6)$} \\
    \multicolumn{5}{|c|}{The Coxeter number $h(E_6) = 12$} \\
  \hline
      $\mu_{\max}$   &
      $w_{\beta}\mu_{max}$ &
      $w\mu_{max}$ &
      $w_{\beta}w\mu_{max}$ &
      $w^2\mu_{max}$ \\
  \hline
      $\begin{array}{cccccc}
        2 & 3 & 4 & 3 & 2 & 1\\
          &   & 2 &   &   &
      \end{array}$
&     $\begin{array}{cccccc}
        1 & 3 & 4 & 3 & 2 & 1 \\
          &   & 2 &   &   &
      \end{array}$
&     $\begin{array}{cccccc}
        1 & 2 & 4 & 3 & 2 & 1 \\
        &     & 2 &     &
      \end{array}$
&      $\begin{array}{cccccc}
        1 & 2 & 3 & 3 & 2 & 1\\
          &   & 2 &   &   &
      \end{array}$
&      $\begin{array}{cccccc}
        1 & 2 & 3 & 2 & 2 & 1 \\
          &   & 1 &   &   &
      \end{array}$
      \\
  \hline
     $w_{\beta}w^2\mu_{max}$ &
     $w^3\mu_{max}$ &
     $w_{\beta}w^3\mu_{max}$ &
     $w^4\mu_{max}$ &
     $w_{\beta}w^4\mu_{max}$ \\
  \hline
      $\begin{array}{cccccc}
        1 & 2 & 2 & 2 & 1 & 1 \\
        &     & 1 &     &
      \end{array}$
&      $\begin{array}{cccccc}
        1 & 1 & 2 & 1 & 1 & 1\\
          &   & 1 &   &   &
      \end{array}$
&     $\begin{array}{cccccc}
        0 & 1 & 1 & 1 & 1 & 1 \\
          &   & 1 &   &   &
      \end{array}$
&      $\begin{array}{cccccc}
        0 & 0 & 1 & 1 & 1 & 1 \\
        &     & 0 &     &
      \end{array}$
&      $\begin{array}{cccccc}
        0 & 0 & 0 & 1 & 1 & 1\\
          &   & 0 &   &   &
      \end{array}$
      \\
  \hline
     $w^5\mu_{max}$ &
     $w_{\beta}w^5\mu_{max}$ &
     $w^6\mu_{max}$ & & \\
  \hline
      $\begin{array}{cccccc}
        0 & 0 & 0 & 0 & 1 & 1 \\
          &   & 0 &   &   &
      \end{array}$
&     $\begin{array}{cccccc}
        0 & 0 & 0 & 0 & 0 & 1 \\
        &     & 0 &     &
      \end{array}$
&      $\begin{array}{cccccc}
        0 & 0 & 0 & 0 & 0 & 1\\
          &   & 0 &   &   &
      \end{array}$
    &  &   \\
 \hline
\end{tabular}

 \vspace{2mm}
 \caption{Orbits of the maximal roots, see Fig. \ref{D5_E6_E7_ext}}
 \label{tab_orbit_max_roots}
\end{table}

 \index{root system}

 {\bf Proposition} (On the maximal root $\mu_{max}$ and the longest element $w_0 \in W$ (Proposition \ref{prop_conj_max_root})).
  {\em
  Let $\Gamma$ be extended to another Dynkin diagram $\widetilde{\Gamma}$
  by adding a root $\delta_{m+1}$ connected to $S$ only at $\varphi_k$ and linearly independent of $S$;
  let $\mu_{max}$ be the maximal root in the root system $\Phi(\widetilde\Gamma)$.}

  {\rm (i)} {\em For $w_0 = w^{\frac{h}{2}} =
  (w_{\delta}w_{\varphi})^{\frac{h}{2}}$, we have}
  \begin{equation*}
     w_0 \mu_{max}  ~=~  \delta_{m+1}, \quad
     \mu_{max}  ~=~ w_0\delta_{m+1}.
  \end{equation*}

  {\rm (ii)} {\em The following conjugacy relation holds:}
  \begin{equation*}
     w{s}_{\mu_{max}}  ~\simeq~ w{s}_{\delta_{m+1}}.
  \end{equation*}

  See Table \ref{tab_orbit_max_roots}, Fig. \ref{D5_E6_E7_ext}, and Fig. \ref{Dl_extending}.
~\\

 \index{dipole}
 \index{triangle}
 \index{square}
 \index{diamond}
 \index{Coxeter number}

\begin{figure}[h]
\centering
\includegraphics[scale=0.9]{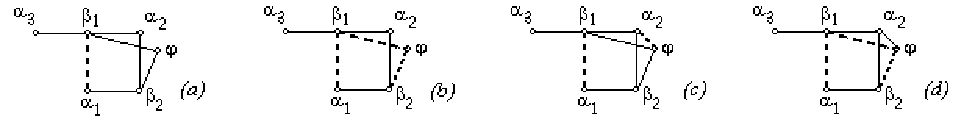}
\caption{Two $D_5(a_1)$-associated elements differing in one vertex}
\label{D5a1_with_fi}
\end{figure}

 A number of basic patterns such as {\it dipoles}, {\it triangles}, {\it squares}, {\it diamonds}
 and other is presented in \S\ref{sect_basic_pat} where some facts pertaining to these patterns are discussed.
 We just present one more frequently used lemma
 related to the pattern consisting of two $D_5(a_1)$-associated subsets differing only in one vertex.
~\\
~\\
 {\bf Lemma} (On necessarily connected roots (Lemma \ref{lem_glue_1})).
   {\em Let $\Gamma = D_5(a_1)$; let
   $S_1 = \{\alpha_1, \alpha_2, \alpha_3, \beta_1, \beta_2\}$ and
   $S_2 = \{\alpha_1, \varphi, \alpha_3, \beta_1, \beta_2 \}$ be two
   $\Gamma$-associated subsets, the vectors of each of which being linearly
   independent, let $\{\alpha_1, \alpha_2, \beta_1, \beta_2\}$
   and $\{\alpha_1, \varphi, \beta_1, \beta_2\}$ be
   $D_4(a_1)$-associated subsets, let the root $\alpha_3$ be connected only with
   $\beta_1$, see Fig. $\ref{D5a1_with_fi}$.}

   {\rm (i)} {\em Configurations of Fig. $\ref{D5a1_with_fi}$$(a)$,$(b)$ are impossible:
   roots $\varphi$ and $\alpha_2$ are necessarily connected, see Fig. $\ref{D5a1_with_fi}$$(c)$,$(d)$}.

   {\rm (ii)} {\em Let $w_1$ (resp. $w_2$) be $S_1$-associated (resp. $S_2$-associated). Then $w_1 \simeq w_2$}.

\subsubsection{The minimal elements of the partial order on Carter diagrams}

\index{partial ordering of Carter diagrams}
  In \S\ref{sec_regular}, the regular extensions of Carter diagrams are introduced and a number of examples is given.
  In {\bf \S\ref{sec_unique}}, we consider regular extensions of Carter diagrams in detail.
  The regular extensions of Carter diagrams determine a {\it partial ordering} $<$
  on the class of Carter diagrams $\mathsf{C4} \coprod \mathsf{DE4}$:
  If $\widetilde\Gamma$ is a regular extension of $\Gamma$,
  we write $\Gamma < \widetilde\Gamma$, see Fig. \ref{diagram_tree}.
  Two $4$-vertex diagrams $D_4(a_1)$ and $D_4$ are \underline{minimal elements} in
  the partially ordered tree of Carter diagrams depicted in Fig. \ref{diagram_tree}.
  In \S\ref{sect_base_case}, we prove the uniqueness
  theorem (Theorem \ref{th_uniq_diagr}) for base diagrams $D_4(a_1)$ and $D_4$.
  The uniqueness theorem does not hold for
  the Carter diagram $A_3$ being the subdiagram of $D_4(a_1)$ and $D_4$, see \S\ref{sec_two_diff_calsses}.
  The pair of opposite vertices of any diagonal in $D_4(a_1)$
  or the pair of endpoints of $D_4$ is said to be the {\it dipole}, see \S\ref{sec_mut_orth}.
  For $4$-vertex diagrams $D_4(a_1)$ and $D_4$, not every two dipoles are conjugate; however, for any two
  diagrams $\Gamma_1$ and $\Gamma_2$ which are both of type $D_4(a_1)$ (resp. $D_4$), there exist dipoles
  $d_1 \in \Gamma_1$ and $d_2 \in \Gamma_2$ such that $d_1$ and $d_2$ are conjugate dipoles.
  We ascend from conjugate dipoles to conjugate triples,
  and from conjugate triples to conjugate quadruples, see  Lemmas \ref{lem_united_2}, \ref{lem_united_3}, \ref{lem_united_4}.
~\\
~\\
  {\bf Corollary} (The base of induction (Corollary \ref{cor_base_case_ind})).
  {\em For $\Gamma = D_4(a_1)$ (resp. $\Gamma = D_4$), let}
   \begin{equation*}
       \begin{cases}
           C_1 = \{\varphi_1, \varphi_2, \delta_1, \delta_2\}, \\
           C_2 = \{\beta_1, \beta_2, \alpha_1, \alpha_2\},
       \end{cases}
       \quad \text{ resp. } \quad
       \begin{cases}
           C_1 = \{\delta_1, \delta_2, \delta_3, \varphi_1\}, \\
           C_2 = \{\alpha_1, \alpha_2, \alpha_3, \beta_1\}
       \end{cases}
   \end{equation*}
  {\em be two $D_4(a_1)$-associated subsets (resp. $D_4$-associated subsets) and}
   \begin{equation*}
       \begin{cases}
          w_1 = s_{\varphi_1}s_{\varphi_2}s_{\delta_1}s_{\delta_2}, \\
          w_2 = s_{\beta_1}s_{\beta_2}s_{\alpha_1}s_{\alpha_2},
       \end{cases}
       \quad \text{ resp. } \quad
       \begin{cases}
          w_1 = s_{\delta_1}s_{\delta_2}s_{\delta_3}s_{\varphi_1}, \\
          w_2 = s_{\alpha_1}s_{\alpha_2}s_{\alpha_3}s_{\beta_1}
       \end{cases}
   \end{equation*}
  {\em be $C_i$-associated elements, where $i = 1,2$.
  Then there exists an element $T \in W$ such that $T^{-1}w_1{T}  = w_2$.
  In other words, any two $\Gamma$-associated elements $w_1$ and $w_2$ are conjugate.}

\begin{figure}[h]
 \centering
 \includegraphics[scale=0.6]{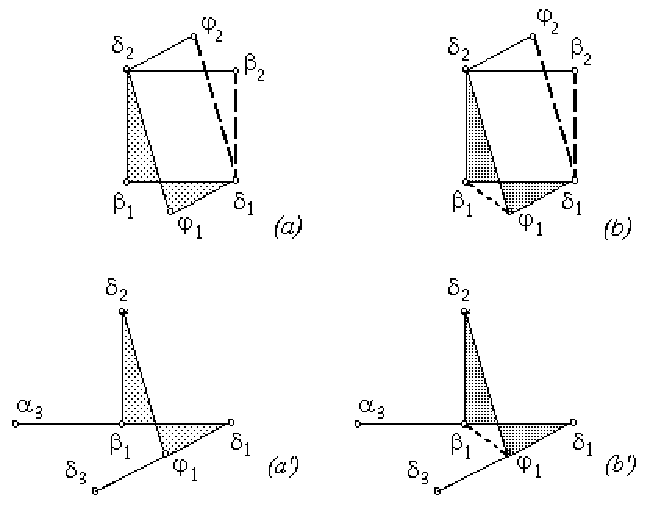}
 \caption{Compare the shaded tetragon $\{\delta_1, \beta_1, \delta_2, \varphi_1\}$ in $D_4(a_1)$ and $D_4$}
\label{same_diagon_2b}
\end{figure}

\begin{remark}{\rm
Observe that although diagrams $D_4(a_1)$ and $D_4$ do not look alike,
occasionally they do have common features.
Comparing diagrams $(a)$ and $(a')$, or
$(b)$ and $(b')$ in Fig. \ref{same_diagon_2b}, see also Fig. \ref{same_diagon_2},
we can determine the key shape which is responsible for the pass from conjugate dipoles to conjugate triples
in both cases: $D_4(a_1)$ and $D_4$.
This shape is the shaded tetragon $\{\delta_1, \beta_1, \delta_2, \varphi_1\}$.
For diagrams $(a)$ and $(a')$, the equality
$$
  \varphi_1 + \beta_1  + \delta_1 + \delta_2 = 0
$$
holds, see Lemma \ref{lem_united_3}. This property allows one to pass from conjugate dipoles to conjugate triples.
For diagrams $(b)$, $(b')$, we can accomplish this pass since $\beta_1 - \varphi_1$ is the root.
For details, see \S\ref{sect_dipoles_triples}.
}
\end{remark}

\subsubsection{The three Principal Cases for the three extension types}
   \index{linearly dependent root}
Two $\widetilde\Gamma$-associated elements are said to be {\it homogeneous $\widetilde\Gamma$-associated elements} if they
are both $\widetilde\Gamma_R$-associated or both $\widetilde\Gamma_L$-associated.
In {\bf \S\ref{sec_unique}}, we show that homogeneous $\widetilde\Gamma$-associated elements
are conjugate:
~\\
~\\
 {\bf Proposition} (On conjugacy of homogeneous $\Gamma$-associated elements (Proposition \ref{prop_induct_step})).
    {\em Let $\Gamma$ be a Carter diagram such that all $\Gamma$-associated elements are
    conjugate.}

    {\rm (i)} {\em  For any \underline{single-track extension} $\Gamma < \widetilde\Gamma$,
    all $\widetilde\Gamma$-associated elements are also conjugate.}

    {\rm (ii)} {\em For any \underline{mirror extensions} $\Gamma < \widetilde\Gamma_L$ and $\Gamma < \widetilde\Gamma_R$,
    all homogeneous $\widetilde\Gamma_R$-associated (resp. $\widetilde\Gamma_L$-associated)
    elements are also conjugate.}

    {\rm (iii)} {\em For any \underline{threefold extensions} $\Gamma < \widetilde\Gamma_1$,
    $\Gamma < \widetilde\Gamma_2$ and $\Gamma < \widetilde\Gamma_3$
     all homogeneous $\widetilde\Gamma_1$-associated (resp. $\widetilde\Gamma_2$-associated,
    resp. $\widetilde\Gamma_3$-associated) elements are also conjugate.}
~\\

    The proof of this proposition is based on considering the following {\it three Principal Cases}:
~\\

    {\it Principal Case $1$: $\alpha$ is connected to $\gamma$.

         Principal Case $2$: $\alpha \perp \gamma$, where $\gamma$ is linearly independent of $\widetilde\Gamma$.

         Principal Case $3$: $\alpha \perp \gamma$, where $\gamma$ is linearly dependent on $\widetilde\Gamma$. }
~\\

    Here, $\alpha$ (resp. $\gamma$) is the root extending $\Gamma$ to $\widetilde\Gamma$.
    For the exact wording of the Principal Cases $2$ and $3$, see eqs. \eqref{eq_prin_case_2} and \eqref{eq_prin_case_3}.
    For mirror extensions, the extending diagram
    $\widetilde\Gamma$ is $\widetilde\Gamma_L$ (resp. $\widetilde\Gamma_R$) for both roots ($\alpha$ and $\gamma$).
    For threefold extensions,  the extending diagram
    $\widetilde\Gamma$ is $\widetilde\Gamma_1$ (resp. $\widetilde\Gamma_2$, resp. $\widetilde\Gamma_3$)
    for both roots, see \S\ref{sec_1_2_option}.

    The proof is the three Principal Cases is given for the three extension types
    (relative the sockets number):
    single-track extension, mirror extensions and threefold extensions, each of which is
    considered for the following $\mathsf{P}v$-extensions:
    $\mathsf{P}1$-extensions, ($D$-joint type and $A$-joint type), $\mathsf{P}2$-extensions and $\mathsf{P}3$-extensions.
    Actually, there are fewer cases, some cases are similar.
    The proof of Principal Cases $1$ and $2$ is relatively simple, see \S\ref{sec_case_1} and \ref{sec_case_2}.
    Principal Case $3$ is the longest and it is proved separately for each type of $\mathsf{P}v$-extensions,
    see \S\ref{sec_case_lin_dep}.

  \index{uniqueness of conjugacy classes}
  \index{conjugacy class! - uniqueness}
  \index{class of Carter diagrams $\mathsf{C4}$}
  \index{class of Carter diagrams $\mathsf{DE4}$}
  \index{connected simply-laced Carter diagrams}
   \index{$\mathsf{C4}$ (class of Carter diagrams)}
   \index{$\mathsf{DE4}$ (class of Carter diagrams)}

 Conjugate elements of $W$ are associated with the same Carter diagram
 $\Gamma$, or with one equivalent to $\Gamma$. In other words, a given conjugacy
 class is associated with the class of diagrams equivalent to $\Gamma$.
 The diagram  $\Gamma$ does not determine a single conjugacy class in $W$, see \cite[Lemma 27]{Ca72}.

 The following theorem  gives a sufficient condition of the uniqueness of the conjugacy class
 determined by the Carter diagram $\Gamma$:
~\\
~\\
 {\bf Theorem} (On the conjugacy class of the diagram (Theorem \ref{th_uniq_diagr})).
 {\em  Let $\Gamma$ be a connected Carter diagram from $\mathsf{C4} \coprod \mathsf{DE4}$.
   Then $\Gamma$ determines a single conjugacy class.}


\subsubsection{Mirror extensions and mirror maps}  In {\bf \S\ref{sec_mirror_threefold}},
we consider mirror and threefold extensions. The following theorem is the main result of the section:
~\\
~\\
{\bf Theorem} (On conjugacy of $\widetilde\Gamma_L$- and $\widetilde\Gamma_R$-associated elements
 (Theorem \ref{th_mirror_ext})).
 {\em Let $\Gamma$ be a Carter diagram, $\Gamma < \widetilde\Gamma_L$ and
   $\Gamma < \widetilde\Gamma_R$ left and right extensions from Tables $\ref{tab_mirror_1}$ -- $\ref{tab_mirror_2}$.
   Then any $\widetilde\Gamma_L$-associated element $w_L$ and any $\widetilde\Gamma_R$-associated element $w_R$
   are conjugate, i.e.,
   \begin{equation*}
      w_R = T^{-1}{w}_L{T} \text{ for some } T \in W.
   \end{equation*}
 }
 ~\\
  For any mirror extensions $\Gamma < \widetilde\Gamma_L$ and $\Gamma < \widetilde\Gamma_R$,
  we construct the map $T$ for any pair elements $w_L$ and $w_R$.
  In most cases, the map $T$ is the composition of the longest element $w_0 \in W(\Gamma_0)$,
  where $\Gamma_0$ is a certain Dynkin subdiagram of $\Gamma$, and the corrective reflection $s_{\alpha + \beta}$
  for some roots $\alpha$ and $\beta$, see \S\ref{sec_longest} and Remark \ref{rem_corr_rel}.
  The map $T$ transfers sockets one to another, see \S\ref{sec_1_2_option}. The element  $T$ is said to be
  the {\it mirror map}. In Table \ref{tab_mirror_map_T} we give an explicit expression for the map $T$ for
  all mirror extensions from Tables \ref{tab_mirror_1} -- \ref{tab_mirror_2}.

\begin{table}[h]
 \Small
  \centering
  \renewcommand{\arraystretch}{1.2}
  \begin{tabular} {|c|c|c|c|c|c|}
  \hline
     Mirror Extensions & Mirror map $T$ & $\mathsf{P}v-$ & Tables & $w_0$ in  & Reference \\
                       &                & extension  &        &           &  \\
  \hline
    $\begin{array}{c}
       D_5(a_1) < E_6(a_1) \\
       D_6(a_1) < E_7(a_1) \\
       D_7(a_1) < E_8(a_1) \\
       D_4(a_1) < D_5(a_1) \text{ for } \\
       \quad D_5(a_1) \subset E_l, D_l \\
     \end{array}$
     & $\begin{array}{c} w_0 = (s_{\beta_2}s_{\alpha_1}{s}_{\alpha_2})^2 \end{array}$
     & $\mathsf{P}1$
     &  Table \ref{tab_mirror_1}
     & $W(A_3)$
     & $\begin{array}{c} \text{Proposition } \ref{prop_mirror_ext}, \\
        \text{ Remark } \ref{rem_prop_th_mirror}
       \end{array}$ \\
  \hline
     $\begin{array}{c}
       D_5 < E_6 \\
       E_6 < E_7 \\
     \end{array}$
     & $w_0$
     & $\mathsf{P}1$
     & Table \ref{tab_mirror_1}
     & $\begin{array}{c}
       W(D_5)  \\
       W(E_6)  \\
     \end{array}$
    & \S\ref{sec_max_root} \\
  \hline
   $\begin{array}{c}
       D_4(a_1) < D_5(a_1) \text{ for } \\
       \quad D_5(a_1) \subset E_l, \\
       \quad D_5(a_1) \not\subset D_l \\
     \end{array}$
     & $\begin{array}{c} T = T_1{T}_2{P}  \\
        (T_1, T_2 \text{ are corrective} \\
        \text{ reflections})
       \end{array}$
     & $\mathsf{P}1$
     & Table \ref{tab_mirror_1}
     & ---
     & Proposition \ref{prop_D4a1_toExt} \\
   \hline
     $D_{2k+2}(a_k) < D_{2k+3}(a_k)$
     &  $s_{w_0\beta_1 - \beta_1} w_0$
     & $\mathsf{P}1$
     & Table \ref{tab_mirror_1}
     & $W(A_{2k+1})$
     &  \S\ref{sect_spec_case_1} \\
  \hline
    $\begin{array}{c}
       D_5(a_1) < E_6(a_2) \\
       D_6(a_1) < E_7(a_3) \\
       D_7(a_1) < E_8(a_4) \\
     \end{array}$
     & $\begin{array}{c} w_0 = (s_{\beta_2}s_{\alpha_1}{s}_{\alpha_2})^2 \end{array}$
     & $\mathsf{P}2$
     & Table \ref{tab_mirror_2}
     & $ W(A_3)$
     & $\begin{array}{c} \text{Proposition } \ref{prop_mirror_ext}, \\
        \text{ Remark } \ref{rem_prop_th_mirror}
       \end{array}$ \\
  \hline
   $D_6(a_2) < E_7(a_4)$
   & $\begin{array}{c} T = QPw_0, \text{ where}  \\
      w_0 = (s_{\beta_2}s_{\alpha_1}{s}_{\alpha_2})^2 \\
      (Q, P \text{ are corrective} \\
        \text{ reflections}) \\
      \end{array}$
   & $\mathsf{P}3$
   & Table \ref{tab_mirror_2}
   & $ W(A_3)$
   & $\begin{array}{c}
     \text{Lemma } \ref{lem_glue_3}, \\
     \S\ref{sec_proof_E7a4} \\
     \end{array}$ \\
   \hline
\end{tabular}
  \vspace{2mm}
  \caption{The mirror map $T$ for $\mathsf{P}1$-, $\mathsf{P}2$-, $\mathsf{P}3$-extensions}
  \label{tab_mirror_map_T}
\end{table}

\clearpage
~\\
\section{\sc\bf Classification of Carter diagrams}
  \label{sec_classif}
 In this section, we add new arguments to obtain the list of Carter diagrams:
 We use the statement on intersecting cycles, Proposition \ref{prop_two_cycles};
 we exclude diagrams with cycles of length $>4$, see Theorem \ref{th_get_rid}.
 The following proposition states that
 any Carter diagram, or connection diagram, without cycles is a Dynkin diagram.
\begin{proposition}[Lemma 8, \cite{Ca72}]
  \label{prop_tree}
  Let $\Gamma$ be a Carter diagram or connection diagram.
  If $\Gamma$ is a tree, then $\Gamma$ is the Dynkin diagram.
\end{proposition}

  For the proof and examples, see \S\ref{Dynkin_diagr}. \qedsymbol
~\\

 Due to this proposition, to prove the Carter theorem (see \S\ref{sec_intro}),
 it suffices to consider only diagrams with cycles.

 \begin{remark}\rm{For $G_2$ and $A_l$, there are no Carter diagrams with cycles.
  Indeed, for $G_2$, this fact is trivial, since there at most two linearly
  independent roots; for $A_l$, see \ref{case_An}.
 }
 \end{remark}

\subsection{For the multiply-laced case, only a $4$-cycle is possible}
\begin{figure}[h]
\centering
\includegraphics[scale=0.8]{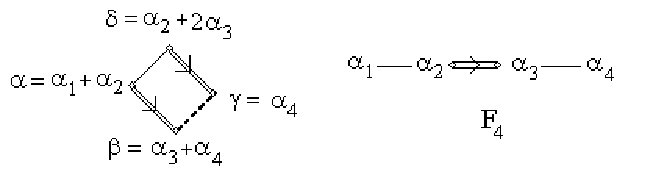}
\caption{The $4$-cycle root subset in $F_4$. The angle
$(\widehat{\beta, \gamma})$ is acute}
\label{F4_cycle}
\end{figure}

 \index{root system}
 \index{acute angle between roots}
  Consider a multiply-laced diagram containing cycles.
  If the root system $\varPhi$ contains a cycle, then  $\varPhi$  constitutes
 the $4$-cycle with one dotted edge, \cite[p. 13]{Ca72}.
 This case occurs in $F_4$, see Fig. \ref{F4_cycle}.

 If $\alpha_1, \alpha_2, \alpha_3, \alpha_4$
 are the simple roots in $F_4$, then the quadruple
 \begin{equation*}
   \alpha = \alpha_1 + \alpha_2,  \quad
   \beta = \alpha_3 + \alpha_4, \quad
   \gamma = \alpha_4, \quad
   \delta = \alpha_2 + 2\alpha_3
 \end{equation*}
  constitutes such a $4$-cycle. The values of the Tits form on the
  corresponding pairs of roots are as follows:
~\\

 \begin{equation*}
  \begin{split}
   & (\alpha, \beta) = (\alpha_1 + \alpha_2, \alpha_3 + \alpha_4) = (\alpha_2,
   \alpha_3) = -1, \\
   & (\beta, \gamma) = (\alpha_3 + \alpha_4, \alpha_4) =
    (\alpha_4, \alpha_4) - (\alpha_3, \alpha_4) = 1 - \frac{1}{2} =
    \frac{1}{2} \quad \text (dotted~edge), \\
   & (\gamma, \delta) = (\alpha_4, \alpha_2 + 2\alpha_3) =
    2(\alpha_4, \alpha_3) = -1, \\
   & (\delta, \alpha) = (\alpha_2 + 2\alpha_3, \alpha_1 + \alpha_2) =
     (\alpha_2, \alpha_2) + (\alpha_2, \alpha_1)  + 2(\alpha_2, \alpha_3) =
     2 - 1 - 2 = -1.
  \end{split}
 \end{equation*}
~\\

 In \S\ref{mult_case}, we prove that for multiply-laced cases, there are no
 other Carter diagrams with cycles.

\subsubsection{Two intersecting cycles in the simply-laced case}
  From the foregoing in this section, it suffices to consider only simply-laced diagrams.
  First of all, we discuss Carter diagrams containing
  intersecting cycles and bridges, see \S\eqref{sec_main_res}.

 \begin{proposition}[On intersecting cycles and bridges]
   \label{prop_two_cycles}

      {\rm (i)} Let $\Gamma$ be any Carter diagram, or connection diagram,
         containing two cycles with bridge $\mathcal{P}$.
         Then $\mathcal{P}$ consists of exactly one edge.

      {\rm (ii)} Let $\Gamma$ be any Carter diagram.
         Let $P_1$, $P_2 \subset \Gamma$ be two paths stemming from the opposite vertices of
         a $4$-cycle in $\Gamma$; let $\alpha_1$ (resp. $\alpha_2$) be the vertex lying in $P_1$ (resp. $P_2$). The diagram obtained from $\Gamma$ by adding the edge $\{ \alpha_1, \alpha_2 \}$   is not a Carter diagram,
         see Fig. $\ref{2paths}$.

       {\rm (iii)} Let $\Gamma$ be any Carter diagram containing two intersecting cycles.
         Then one of the cycles consists of $4$ vertices, and the other one can contain only $4$ or $6$ edges.
 \end{proposition}

 \PerfProof
  (i) Every cycle contains an odd number of dotted edges, otherwise by several reflections
  we get a cycle containing only solid edges, a case which cannot happen, see Lemma \ref{lem_must_dotted}.
  Let $n_1$ be the number of dotted edges in the top cycle:
  $\{ \alpha_1, \beta_1, \alpha_2, \beta_2, \alpha_3, \beta_n \}$,
  and $n_2$ the number of dotted edges in the bottom cycle:
  $\{ \alpha_4, \beta_1, \alpha_2, \beta_2, \alpha_5, \beta_m \}$.
  Both $n_1$ and $n_2$ are odd.
  Suppose the bridge $\mathcal{P}$ with endpoints $\beta_1$ and $\beta_2$ contains an additional vertex $\alpha_2$
  (i.e., $\mathcal{P} = \{\beta_1, \alpha_2, \beta_2\}$, see Fig. \ref{inters_cycles}$(a)$ or Fig. \ref{inters_cycles}$(b)$).
  After discarding the vertex $\alpha_2$ we get a bigger cycle
  $\mathcal{C}_3 = \{ \alpha_1, \beta_1, \alpha_4, \beta_m, \alpha_5, \beta_2, \alpha_3,  \beta_n \}$;
  in the generic case of the bridge $\mathcal{P}$, we discard from the bridge
  all vertices except $\beta_1$, $\beta_2$.
  Let $n_3$ be the number of dotted edges in the cycle  $\mathcal{C}_3$; $n_3$ is also odd.
  Therefore, $n_1 + n_2 + n_3$  is odd. On the other hand, every dotted edge enters twice,
  so $n_1 + n_2 + n_3$ is even. Thus, there is no vertex in the bridge $\{\beta_1, \beta_2 \}$
  between $\beta_1$ and $\beta_2$.
~\\

 (ii) The diagram $\Gamma \cup \{ \alpha_1, \alpha_2 \}$ contains the
 bridge $\{ \beta_1, \gamma, \beta_2 \}$ of length $2$, see
 Fig.\ref{2paths}. Thus, by (i), the diagram $\Gamma \cup \{ \alpha_1, \alpha_2 \}$ in
 Fig.\ref{2paths} is not a Carter diagram.

\begin{minipage}{7.6cm} 
  \qquad \epsfig{file=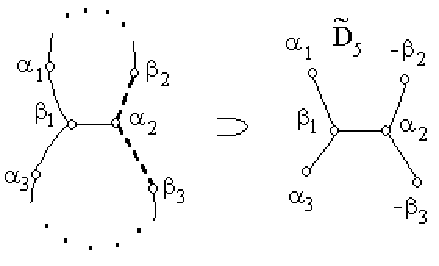, width=2.1in}
  \captionof{figure}{}
  \label{to_extD5}
\end{minipage}
\begin{minipage}{8.4cm}
\index{bridge}
  (iii) By (i) the bridge consists of one edge $\{ \beta_1, \alpha_2 \}$,
  see Fig. \ref{to_extD5}.
  Then at least one of the cycles is of length $4$.
  Otherwise, the Carter diagram contains the
  extended Dynkin diagram $\widetilde{D}_5$ contradicting
  Proposition \ref{prop_non_ext_Dynkin}.
  As above, the dotted edge may be eliminated from $\widetilde{D}_5$
  by changing the sign of one of the roots.
\end{minipage}
~\\
~\\
~\\
~\\
\begin{minipage}{5.4cm}
  \psfig{file=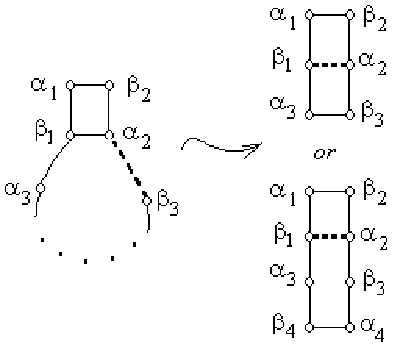, width=2.0in}
  \captionof{figure}{}
  \label{1edge_inters}
\end{minipage}
\begin{minipage}{5.2cm}
  The second cycle can be only of length $4$ or $6$ as in Fig. \ref{1edge_inters}.
  It cannot be a cycle of length $8$, otherwise the
  Carter diagram contains the extended Dynkin diagram $\widetilde{E}_7$,
  see Fig. \ref{fig_4-8-cycles}.  According to (ii),
  we cannot add edges  $\{ \alpha_1, \beta_3 \}$, $\{ \alpha_1, \beta_5 \}$,
  $\{ \beta_2, \alpha_3 \}$, or $\{ \beta_2, \alpha_5 \}$.
  \qed
\end{minipage}
\qquad
\begin{minipage}{5.4cm}
  \psfig{file=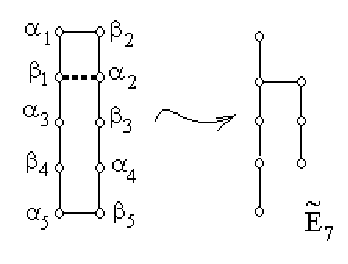, width=2.0in}
  \captionof{figure}{}
  \label{fig_4-8-cycles}
\end{minipage}

\begin{corollary}[On restrictions on quintuples of roots]
   \label{cor_numb_ep}

    {\rm (i)} Let an $\alpha$-set contain $3$ roots $\{\alpha_1, \alpha_2, \alpha_3 \}$.
    There does not exist two non-connected roots $\beta$ and $\gamma$
    connected to every $\alpha_i$ so that the vectors of the quintuple $\{\alpha_1, \alpha_2, \alpha_3, \beta, \gamma \}$
    are linearly independent.

    {\rm (ii)} Let $\{\alpha_1, \beta_1, \alpha_2, \beta_2 \}$ be
    a square in a connection diagram. There does not exist
    a root $\gamma$ connected to all vertices of the square
    so that the vectors of the quintuple $\{ \alpha_1, \beta_1, \alpha_2, \beta_2, \gamma \}$
    are linearly independent.
 \end{corollary}

\begin{figure}[h]
 \centering
 \includegraphics[scale=0.7]{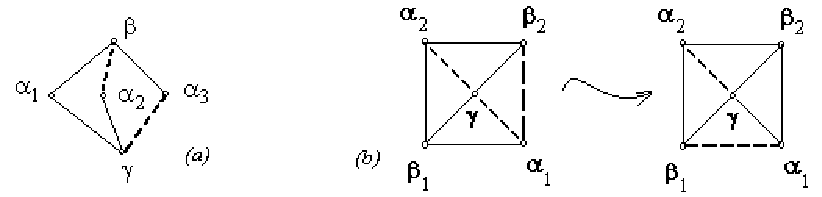}
 \caption{Every cycle should contain an odd number of dotted edges,
          a case which cannot happen}
 \label{3-cells}
 \end{figure}

 \PerfProof (i) Suppose there exist roots $\beta$ and $\gamma$ connected to every $\alpha_i$
  such that the vectors of the quintuple $\{\alpha_1, \alpha_2, \alpha_3, \beta, \gamma \}$
  are linearly independent, see Fig. \ref{3-cells}$(a)$.
  Then we have three cycles: $\{ \alpha_i, \gamma, \alpha_j, \beta \}$,
  where  $1 \leq i < j  \leq 3$.
  Every cycle should contain an odd number of dotted edges.
  Let $n_1$, $n_2$, $n_3$ be the odd numbers of dotted edges in every cycle,
  therefore $n_1 + n_2 + n_3$ is odd.
  On the other hand, every dotted edge enters twice, so $n_1 + n_2 + n_3$ is
  even, which is a contradiction.

\index{triangle of roots}

 (ii) Suppose a certain root $\gamma$ is connected
  to all vertices of the square. Then we have $5$ cycles: Four triangles
  $\{ \alpha_i, \beta_j, \gamma \}$, where $i = 1,2$ and $j = 1,2$, and the square
  $\{\alpha_1, \beta_1, \alpha_2, \beta_2 \}$, see Fig. \ref{3-cells}$(b)$.
  Every cycle should contain an odd number of dotted edges.
  Let $n_1$, $n_2$, $n_3$, $n_4$, $n_5$
  be the numbers of dotted edges in every cycle, therefore
  $n_1 + n_2 + n_3 + n_4 + n_5$ is odd.
  On the other hand, every dotted edge enters twice,
  so $n_1 + n_2 + n_3 + n_4 + n_5$ is even, which is a contradiction.
  For example, the left square in Fig. \ref{3-cells}$(b)$
  is transformed to the  right  one by the reflection $s_{\alpha_1}$,
  then the right square contains the cycle $\{\alpha_1, \beta_2, \gamma\}$
  with $3$ solid edges, i.e., the extended Dynkin diagram $\widetilde{A}_2$,
  contradicting Proposition \ref{prop_non_ext_Dynkin}.
 \qed

\subsection{Classification of simply-laced Carter diagrams with cycles}
  \label{sec_class_s-l}
  The classification of simply-laced Carter diagrams with cycles is based on the following statements:

  (i) the diagram containing any non-Dynkin diagram (in particular, any
  extended Dynkin diagram) is not a Carter diagram (Proposition \ref{prop_tree}).

  (ii) the diagram containing two cycles with a bridge of length $> 1$
  is not a Carter diagram  (Proposition \ref{prop_two_cycles}(i)).

  (iii) the diagram which can be equivalently transformed into a diagram of
  type (i) or (ii)  is not a Carter diagram. We use this fact in Lemma \ref{lem_2_impos_diagr}.

  (iv) the Carter diagrams containing cycles of length $>4$ can be excluded from Carter's list
      (Theorem \ref{th_get_rid}).
\index{bridge}

\subsubsection{The Carter diagrams with cycles on $6$ vertices}
  \label{sec_6_ver}
\begin{figure}[h]
 \centering
 \includegraphics[scale=1.6]{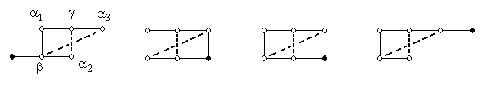}
 \caption{Not Carter diagrams on $6$-vertices}
\label{impos_6-vertices}
\end{figure}

  There are only four  $6$-vertex simply-laced Carter diagrams containing
  cycles, see Table \ref{tab_class_Carter}.
  As we show in \S\ref{sec_D6b2}, the diagram $D_6(b_2)$ is equivalent to
  $D_6(a_2)$, so  $D_6(b_2)$ can be excluded  from the list of Carter diagrams.
  The diagrams depicted in Fig.~\ref{impos_6-vertices} are not Carter diagrams.
  One should discard the bold vertex and apply Corollary \ref{cor_numb_ep}(i),
  see Fig. \ref{3-cells}.

\subsubsection{The Carter diagrams with cycles on $7$ vertices}
  \label{sec_7_ver}
\begin{figure}[h]
\centering
\includegraphics[scale=0.70]{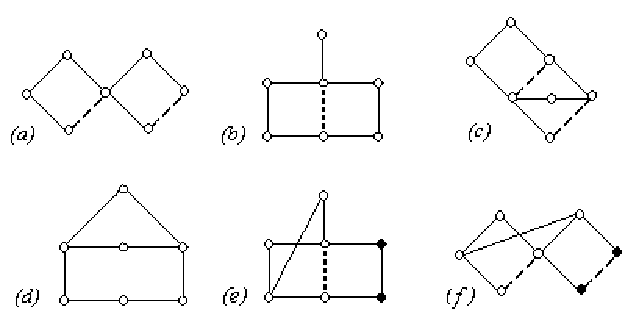}
\captionof{figure}{}
\label{fig_impos_7-vertices}
\end{figure}

 There are only six  $7$-vertex simply-laced Carter diagrams containing
 cycles, see Table \ref{tab_class_Carter}.
 According to \S\ref{sec_D6b2}, the diagram $E_7(b_2)$ is equivalent to
 $E_7(a_2)$. Thus, the diagram $E_7(b_2)$ is excluded  from the list of Carter diagrams.
 Note that the diagrams $(a)$ and $(b)$ depicted in Fig. \ref{fig_impos_7-vertices}
 are not Carter diagrams since each of
 them contains the extended Dynkin diagram $\widetilde{D}_4$. The
 diagrams $(c)$ and $(d)$ are not Carter diagrams since for each of them
 there exist two cycles with the bridge of length $> 1$,
 contradicting Proposition \ref{prop_two_cycles}. In order to see that $(e)$ and $(f)$
 are not Carter diagrams, one can discard bold vertices and apply
 Corollary \ref{cor_numb_ep}(i) as in \S\ref{sec_6_ver}.

\subsubsection{The Carter diagrams with cycles on $8$ vertices}

  There are only eleven  $8$-vertex simply-laced Carter diagrams containing
  cycles, see Table \ref{tab_class_Carter}.

\begin{figure}[h]
\centering
\includegraphics[scale=0.8]{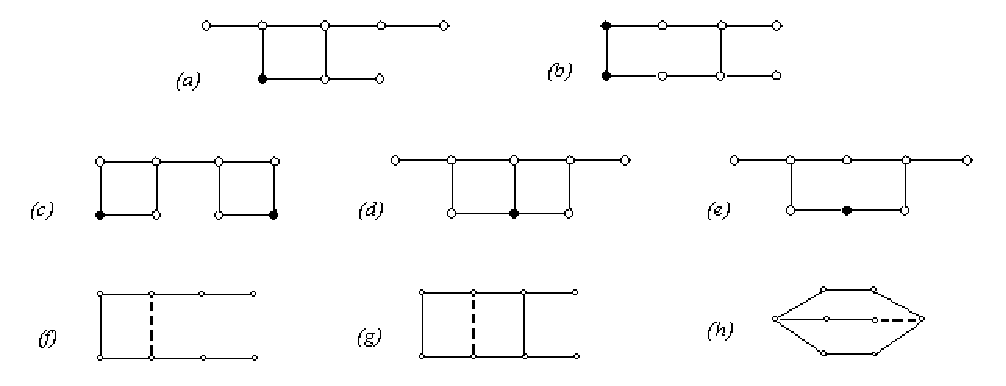}
\caption{$8$-vertex diagrams are not Carter diagrams }
\label{impos_8-vertices}
\end{figure}

 The diagrams depicted in Fig. \ref{impos_8-vertices} are not Carter diagrams.
 One can discard the bold vertices to see that each of depicted diagrams
 contains an extended Dynkin diagram. The diagram $(a)$ contains
 $\widetilde{E}_6$; $(b)$ and $(c)$ contain $\widetilde{D}_5$; $(d)$
 and $(e)$ contain $\widetilde{D}_6$. For diagrams $(f)$ and $(g)$, see
 Lemma \ref{lem_2_impos_diagr}. The diagram $(h)$ is not a Carter
 diagram since there exists the bridge of length $> 1$, see
 Proposition \ref{prop_two_cycles}\footnotemark[1]. \footnotetext[1]{
 We do not depict here the diagrams corresponding to Proposition \ref{prop_two_cycles}(ii),
 see Fig.~\ref{2paths}. For $l=6$,
 they are depicted in Fig. \ref{impos_6-vertices}; for $l=7$, see
 diagrams $(e)$, $(f)$ from \S\ref{sec_7_ver}.}

~\\

\begin{figure}[h]
\centering
\includegraphics[scale=1.5]{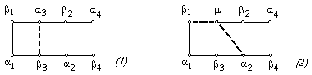}
\captionof{figure}{}
\label{eq_a_8vert}
\end{figure}

\begin{lemma}
  \label{lem_2_impos_diagr}
  Diagrams $(f)$ and $(g)$ in Fig. $\ref{impos_8-vertices}$ are not Carter diagrams.
\end{lemma}
 \PerfProof  In cases $(f)$ and $(g)$, we transform the given diagram to an equivalent one
 containing an extended Dynkin diagram.
   Let $\Gamma$ be the diagram $(f)$ in Fig. \ref{impos_8-vertices}.
   The corresponding roots are depicted in the diagram in Fig. \ref{eq_a_8vert}(1).
 Let $w$ be the $\Gamma$-associated element:
 \begin{equation*}
    w  = s_{\alpha_1}s_{\alpha_2}s_{\alpha_3}s_{\alpha_4}s_{\beta_1}s_{\beta_2}s_{\beta_3}s_{\beta_4}.
 \end{equation*}

~\\
 Since $s_{\alpha_3}s_{\beta_1}s_{\beta_3} = s_{\beta_1}s_{\beta_3}s_{\mu}$, where
 $\mu = \alpha_3 - \beta_3 + \beta_1$, we have
 \begin{equation*}
    w  = s_{\alpha_1}s_{\alpha_2}s_{\alpha_4}s_{\beta_1}s_{\beta_3}s_{\mu}s_{\beta_2}s_{\beta_4}.
 \end{equation*}
 Therefore, the element $w$ is associated with the connection diagram depicted
 in Fig. \ref{eq_a_8vert}(2).
 Discard the vertex $\beta_3$, the remaining diagram is the extended Dynkin diagram
 $\widetilde{E}_6$.
~\\

\begin{figure}[h]
\centering
\includegraphics[scale=1.8]{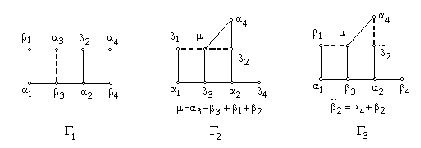}
 \caption{The equivalence transformation from $\Gamma_1$ to $\Gamma_3$}
\label{eq_b_8vert}
\end{figure}

   Let $\Gamma$ be the diagram $(g)$ in Fig. \ref{impos_8-vertices}. The
 same diagram with corresponding roots is the diagram $\Gamma_1$ depicted in Fig. \ref{eq_b_8vert}.

~\\
The $\Gamma_1$-associated element is as follows:
\begin{equation}
  \label{eq_omega_2}
   \begin{split}
    w  = & s_{\alpha_1}s_{\alpha_2}s_{\alpha_3}s_{\alpha_4}s_{\beta_1}s_{\beta_2}s_{\beta_3}s_{\beta_4} =
           s_{\alpha_1}s_{\alpha_2}s_{\alpha_4}(s_{\alpha_3}s_{\beta_1}s_{\beta_2}s_{\beta_3})s_{\beta_4} = \\
         &  s_{\alpha_1}s_{\alpha_2}s_{\alpha_4}s_{\beta_1}s_{\beta_2}s_{\beta_3}s_{\mu}s_{\beta_4},
           \text{ where } \mu = \alpha_3 - \beta_3 + \beta_1 + \beta_2.
   \end{split}
\end{equation}
 The last expression of $w$ is a $(\Gamma_2, o_2)$-associated element, where the diagram $\Gamma_2$ in Fig. \ref{eq_b_8vert}
 is the connection diagram, not a Carter diagram, and the order
 $o_2$ is given by \eqref{eq_omega_2}.
 Further,
 \begin{equation}
  \label{eq_omega_3}
    w  =  s_{\alpha_1}s_{\alpha_2}(s_{\alpha_4}s_{\beta_2})s_{\beta_1}s_{\beta_3}s_{\mu}s_{\beta_4} =
          s_{\alpha_1}s_{\alpha_2}s_{\widetilde{\beta}_2}s_{\alpha_4}s_{\beta_1}s_{\beta_3}s_{\mu}s_{\beta_4},
           \text{ where } \widetilde{\beta}_2 = \alpha_4 + \beta_2,
\end{equation}
 The obtained expression of $w$ is a $(\Gamma_3, o_3)$-associated element, where
 $\Gamma_3$ it the connection diagram in Fig. \ref{eq_b_8vert} and $o_3$ is
 the order given by \eqref{eq_omega_3}.
 The diagram $\Gamma_3$ contains the extended Dynkin diagram
 $\widetilde{D}_5 = \{\alpha_1, \alpha_2, \mu, \widetilde{\beta_2}, \beta_3, \beta_4\}$,
 but this is impossible. \qed
~\\

\subsubsection{The Carter diagrams with cycles on $l > 8$ vertices}
The Dynkin diagram $A_l$ does not contain any Carter diagrams with
cycles, see \S\ref{case_An}. For the Dynkin diagram $D_l$, we refer
to Carter's discussion in \cite[p. 13]{Ca72}.
 In this case, there are the two types of Carter diagrams (Table \ref{tab_class_Carter}, $l > 8$):

 $(1)$ pure cycles $D_l(b_{\frac{l}{2} -1})$ for $l$ is even, $l \leq n$

 $(2)$ $D_l(a_1)$,   $D_l(a_2)$, \dots, $D_l(a_{\frac{l}{2} -1})$ for $l$ is even, $l \leq n$.

 In \S\ref{sec_pure_cycle}, we will show that any pure cycle $D_l(b_{\frac{l}{2} -1})$ from $(1)$ is equivalent
 to $D_l(a_{\frac{l}{2} -1})$ from $(2)$, and hence pure cycles $D_l(b_{\frac{l}{2} -1})$ can be excluded
 from Carter's list.

\clearpage
~\\
\section{\sc\bf Exclusion of long cycles}
  \label{sec_get_rid}
In this section, we show that Carter diagrams containing cycles
of length  $n > 4$  can be discarded from the list.

\begin{theorem}[On exclusion of long cycles]
   \label{th_get_rid}
   Any Carter diagram containing $n$-cycles, where $n > 4$, is
   equivalent to another Carter diagram containing only $4$-cycles.
 \end{theorem}

\index{explicit transformation}

In all cases we construct a certain explicit transformation of the
diagram containing $n$-cycles, where $n > 4$, to a diagram containing only
$4$-cycles. The corresponding pairs of equivalent diagrams are
depicted in Table \ref{tab_eq_pairs}.

\begin{table}[t]
  \centering
  \Small
  \renewcommand{\arraystretch}{1.3}
  \begin{tabular} {|c|c|c|c|}
  \hline
           & The Carter diagram  & The equivalent                   & The characteristic \cr
           & with $n$-cycles,    & Carter diagram $\Gamma$,         & polynomial of the    \cr
           &    $n > 4$          &  only $4$-cycles                 & $\Gamma$-associated element  \\
  \hline 
       & & & \\
       1 &
       $\begin{array}{c}
       \includegraphics[scale=0.8]{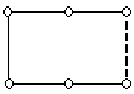}
       \end{array}$
       &
        $\begin{array}{c}
       \includegraphics[scale=0.8]{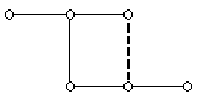}
       \end{array}$
       & $(t^3 + 1)^2$
       \cr
        & $D_6(b_2)$, $n = 6$ &  $D_6(a_2)$ &  \\
   \hline
       & & & \\
       2 &
       $\begin{array}{c}
       \includegraphics[scale=0.8]{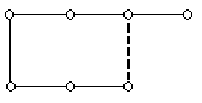}
       \end{array}$
       &
       $\begin{array}{c}
       \includegraphics[scale=0.8]{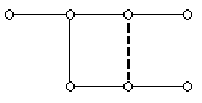}
       \end{array}$
       & $(t^4 - t^2 + 1)(t^2 - t + 1)(t+1)$
       \cr
       & $E_7(b_2)$, $n = 6$ & $E_7(a_2)$ &  \\
   \hline
       & & & \\
       3 &
       $\begin{array}{c}
       \includegraphics[scale=0.8]{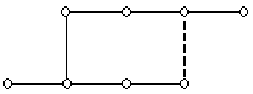}
       \end{array}$
       &
       $\begin{array}{c}
       \includegraphics[scale=0.8]{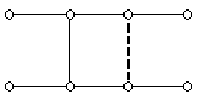}
       \end{array}$
       & $(t^4 - t^2 + 1)^2$
       \cr
       & $E_8(b_3)$, $n = 6$  & $E_8(a_3)$ & \\
   \hline
       & & & \\
       4 &
       $\begin{array}{c}
       \includegraphics[scale=0.8]{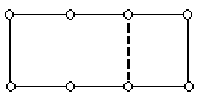}
       \end{array}$
       &
       $\begin{array}{c}
       \includegraphics[scale=0.8]{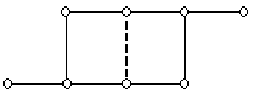}
       \end{array}$
       &  $t^8 - t^7 + t^5 - t^4 + t^3 - t^2 + 1$
       \cr
       & $E_8(b_5)$, $n = 6$  & $E_8(a_5)$  & \\
   \hline
       & & & \\
       5 &
       $\begin{array}{c}
       \includegraphics[scale=0.8]{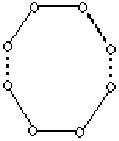}
       \end{array}$
       &
       $\begin{array}{c}
       \includegraphics[scale=0.8]{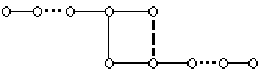}
       \end{array}$
       & $(t^{\frac{l}{2}} + 1)^2$
       \cr
       & $D_l(b_{\frac{l}{2} - 1})$, $n = l$, $l$ even &
       $D_l(a_{\frac{l}{2} - 1})$, $n = l$, $l$ even & \\
       & & & \\
  \hline
\end{tabular}
  \vspace{2mm}
  \caption{Pairs of equivalent Carter diagrams}
  \label{tab_eq_pairs}
\end{table}

 Note that the coincidence of characteristic polynomials of diagrams in
 pairs of Table \ref{tab_eq_pairs} is the necessary condition of equivalence
 of these diagrams, see \cite[Table 3]{Ca72}.
 As it is shown in Theorem \ref{th_get_rid}, this condition is also sufficient for
 the Carter diagrams.

\index{explicit transformation}

 For convenience, we consider the equivalence $D_6(b_2) \simeq D_6(a_2)$
 as a separated case, though this is a particular case of the
 pair $D_l(b_{\frac{l}{2} - 1}) \simeq D_l(a_{\frac{l}{2} - 1})$
 with $l = 6$, Table \ref{tab_eq_pairs}. The idea of explicit
 transformation connecting elements of every pair is similar for all pairs\footnotemark[1].
 \footnotetext[1]{
 Redrawing elements of pairs as the projection of
 $3$-dimensional cube in Fig. \ref{E8_a3_b3} -- Fig. \ref{E8_a5_fin}
 may give, perhaps, a hint to a geometric interpretation of these explicit transformations.}

\subsection{Equivalence $E_8(b_3) \simeq E_8(a_3)$}
  The $E_8(a_3)$-associated element $w$ is transformed as follows:
\begin{equation}
 \label{E8_mu}
  \begin{split}
   w = & s_{\alpha_1}s_{\alpha_2}s_{\alpha_3}s_{\alpha_4}s_{\beta_1}s_{\beta_2}s_{\beta_3}s_{\beta_4} = \\
       & s_{\alpha_1}s_{\alpha_4}\left (
       s_{\alpha_2}s_{\alpha_3}s_{\beta_3} \right )s_{\beta_1}s_{\beta_2}s_{\beta_4} = \\
       & s_{\alpha_1}s_{\alpha_4}s_{\mu}s_{\alpha_2}s_{\alpha_3}s_{\beta_1}
         s_{\beta_2}s_{\beta_4}, \text{ where } \mu = \beta_3 + \alpha_3 - \alpha_2.
  \end{split}
\end{equation}

\begin{figure}[h]
\centering
\includegraphics[scale=0.8]{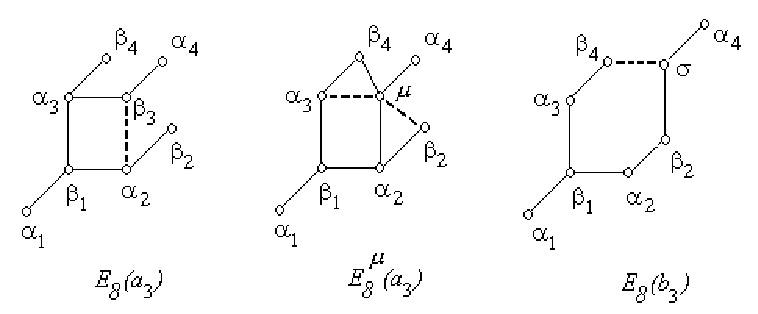}
\caption{Equivalence $E_8(b_3) \simeq E_8(a_3)$; the
 intermediate diagram $E_8^{\mu}(a_3)$ is a connection diagram, not a Carter diagram }
\label{E8_a3_b3}
\end{figure}
~\\
  The element $w$ is $(E^{\mu}_8(a_3), o)$-associated, where $E^{\mu}_8(a_3)$ is the connection diagram in
  Fig. \ref{E8_a3_b3}, the order $o$ is given by \eqref{E8_mu}.
  From \eqref{E8_mu} we have:
 \begin{equation}
  \begin{split}
   w &  \quad \stackrel{s_{\beta_2}s_{\beta_4}}{\simeq} \quad s_{\alpha_1}s_{\alpha_4}
      \left ( s_{\beta_2}s_{\beta_4}s_{\mu} \right )
      s_{\alpha_2}s_{\alpha_3}s_{\beta_1} = \\
      & \quad  s_{\alpha_1}s_{\alpha_4}s_{\sigma}s_{\beta_2}s_{\beta_4}s_{\alpha_2}s_{\alpha_3}s_{\beta_1},
         \text{ where } \sigma = \mu - \beta_2 + \beta_4.
  \end{split}
\end{equation}
~\\
So,  $\sigma = \beta_3 + \alpha_3 - \alpha_2 - \beta_2 + \beta_4$,
and it is easy to see that
 \begin{equation}
 \label{rel_E8_b3}
 \Small
  \begin{split}
    & (\sigma, \alpha_3) = (\alpha_3, \alpha_3) + (\alpha_3, \beta_4) + (\alpha_3, \beta_3)  =
       1 - \frac{1}{2} -  \frac{1}{2} = 0, \\
    & (\sigma, \alpha_2) = -(\alpha_2, \alpha_2) - (\alpha_2, \beta_2) + (\alpha_2, \beta_3)  =
       -1 + \frac{1}{2} +  \frac{1}{2} = 0, \\
    & (\sigma, \beta_1) = -(\alpha_2, \beta_1) + (\alpha_3, \beta_1) =
       \frac{1}{2} -  \frac{1}{2} = 0, \\
    & (\sigma, \alpha_1) = 0, \\
    & (\sigma, \beta_4) = (\beta_4, \beta_4) + (\alpha_3, \beta_4) = 1 - \frac{1}{2} = \frac{1}{2}, \\
    & (\sigma, \beta_2) = -(\beta_2, \beta_2) -(\alpha_2, \beta_2) = -1 + \frac{1}{2} = -\frac{1}{2}, \\
    & (\sigma, \alpha_4) = (\beta_3, \alpha_4) = -\frac{1}{2}. \\
  \end{split}
\end{equation}
Relations \eqref{rel_E8_b3} describe the Carter diagram $E_8(b_3)$,
Fig. \ref{E8_a3_b3}. We only need to check that the element $w$ is
conjugate to a product of two involutions:
 \begin{equation}
  \label{E8_a3_canon}
  \begin{split}
   w  \quad \simeq \quad & s_{\alpha_1}s_{\alpha_4}s_{\sigma}s_{\beta_2}s_{\beta_4}s_{\alpha_2}s_{\alpha_3}s_{\beta_1}
        \quad \stackrel{s_{\alpha_4}}{\simeq} \quad
        s_{\alpha_1}s_{\sigma}( s_{\beta_2}s_{\beta_4}s_{\alpha_4}) s_{\alpha_2}s_{\alpha_3}s_{\beta_1}
        \quad \stackrel{s_{\sigma}}{\simeq} \\
        & s_{\alpha_1}(s_{\beta_2}s_{\beta_4}s_{\alpha_4})(s_{\alpha_2}s_{\alpha_3}s_{\sigma})s_{\beta_1} =
          (s_{\beta_2}s_{\beta_4}s_{\alpha_4})(s_{\alpha_1}s_{\alpha_2}s_{\alpha_3}s_{\sigma})s_{\beta_1}
          \stackrel{s_{\beta_1}}{\simeq} \\
        & (s_{\beta_1}s_{\beta_2}s_{\beta_4}s_{\alpha_4})(s_{\alpha_1}s_{\alpha_2}s_{\alpha_3}s_{\sigma}).
  \end{split}
\end{equation}
  Thus, $w_1 = s_{\beta_1}s_{\beta_2}s_{\beta_4}s_{\alpha_4}$ and $w_2 = s_{\alpha_1}s_{\alpha_2}s_{\alpha_3}s_{\sigma}$
  are two involutions,
  $w = w_1 w_2$, i.e., $w$ is conjugate to the $E_8(b_3)$-associated element, which was to be proven.

  \subsection{Equivalences $E_7(b_2) \simeq E_7(a_2)$ and $D_6(b_2) \simeq D_6(a_2)$ }
    \label{sec_D6b2}
   These equivalences directly follow from the equivalence $E_8(b_3) \simeq E_8(a_3)$ that we see from Fig. \ref{E7_a2_b2} and Fig.
   \ref{D6_a2_b2}. For the equivalence $E_7(b_2) \simeq E_7(a_2)$,
   we discard $s_{\alpha_1}$ in relations \eqref{E8_mu} -- \eqref{E8_a3_canon} as follows:

\begin{figure}[h]
\centering
\includegraphics[scale=0.8]{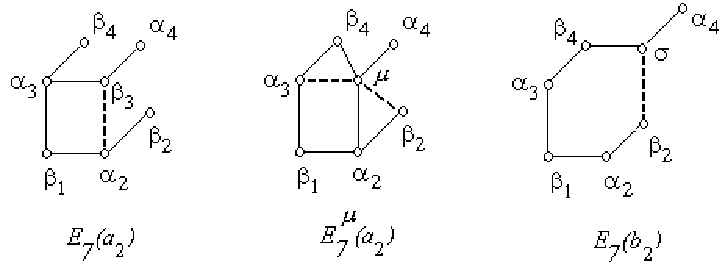}
\caption{Equivalence $E_7(b_2) \simeq E_7(a_2)$; the
intermediate diagram $E_7^{\mu}(a_2)$ is a connection diagram, not
a Carter diagram}
\label{E7_a2_b2}
\end{figure}
\begin{equation}
 \label{E7_a2_mu}
  \begin{split}
   w = & s_{\alpha_2}s_{\alpha_3}s_{\alpha_4}s_{\beta_1}s_{\beta_2}s_{\beta_3}s_{\beta_4} = \\
       & s_{\alpha_4}s_{\mu}s_{\alpha_2}s_{\alpha_3}s_{\beta_1}s_{\beta_2}s_{\beta_4}
         ~(\text{where } \mu = \beta_3 + \alpha_3 - \alpha_2) \quad  \stackrel{s_{\beta_2}s_{\beta_4}}{\simeq}\\
       &  s_{\alpha_4}s_{\sigma}s_{\beta_2}s_{\beta_4}s_{\alpha_2}s_{\alpha_3}s_{\beta_1}
         ~(\text{where } \sigma = \mu - \beta_2 + \beta_4) \quad  \stackrel{s_{\alpha_4}}{\simeq}\\
       & s_{\sigma}( s_{\beta_2}s_{\beta_4}s_{\alpha_4}) s_{\alpha_2}s_{\alpha_3}s_{\beta_1}
        \quad \stackrel{s_{\sigma}}{\simeq} \quad
       (s_{\beta_2}s_{\beta_4}s_{\alpha_4})(s_{\alpha_2}s_{\alpha_3}s_{\sigma})s_{\beta_1}
          \stackrel{s_{\beta_1}}{\simeq} \\
       & (s_{\beta_1}s_{\beta_2}s_{\beta_4}s_{\alpha_4})(s_{\alpha_2}s_{\alpha_3}s_{\sigma}).
  \end{split}
\end{equation}
  Here, $w_1 = s_{\beta_1}s_{\beta_2}s_{\beta_4}s_{\alpha_4}$ and
  $w_2 \simeq s_{\alpha_2}s_{\alpha_3}s_{\sigma}$ are two involutions, $w \simeq w_1 w_2$
  and the element $w$ is $E_7(b_2)$-associated,
  which was to be proven.  For the equivalence $D_6(b_2) \simeq D_6(a_2)$,
  we discard  $s_{\alpha_4}$ in relation \eqref{E7_a2_mu}:

\begin{figure}[h]
\centering
\includegraphics[scale=0.8]{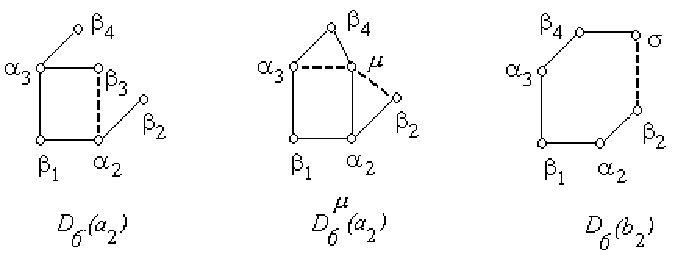}
\caption{Equivalence $D_6(b_2) \simeq D_6(a_2)$; the
 intermediate diagram $D_6^{\mu}(a_2)$ is a connection diagram, not
 a Carter diagram}
\label{D6_a2_b2}
\end{figure}

\begin{equation}
 \label{D6_a2_mu}
  \begin{split}
   w = & s_{\alpha_2}s_{\alpha_3}s_{\beta_1}s_{\beta_2}s_{\beta_3}s_{\beta_4} = \\
       & s_{\mu}s_{\alpha_2}s_{\alpha_3}s_{\beta_1}s_{\beta_2}s_{\beta_4}
         ~(\text{where } \mu = \beta_3 + \alpha_3 - \alpha_2) \quad  \stackrel{s_{\beta_2}s_{\beta_4}}{\simeq}\\
       &  s_{\sigma}s_{\beta_2}s_{\beta_4}s_{\alpha_2}s_{\alpha_3}s_{\beta_1}
         ~(\text{where } \sigma = \mu - \beta_2 + \beta_4) \quad \stackrel{s_{\sigma}}{\simeq} \quad
       (s_{\beta_2}s_{\beta_4})(s_{\alpha_2}s_{\alpha_3}s_{\sigma})s_{\beta_1}
          \stackrel{s_{\beta_1}}{\simeq} \\
       & (s_{\beta_1}s_{\beta_2}s_{\beta_4})(s_{\alpha_2}s_{\alpha_3}s_{\sigma}).
  \end{split}
\end{equation}

  Here, $w_1 = s_{\beta_1}s_{\beta_2}s_{\beta_4}$ and $w_2 = s_{\alpha_2}s_{\alpha_3}s_{\sigma}$ are two
  involutions, $w \simeq w_1 w_2$  and the element $w$ is $E_7(b_2)$-associated.

 \subsection{Equivalence $E_8(b_5) \simeq E_8(a_5)$}
  This equivalence is the most difficult.
~\\

  {\it Step 1.}  Let us transform the $E_8(b_5)$-associated element $w$ as follows:
  \begin{equation}
   \label{eq_step0}
    \begin{split}
      w  ~=~ & (s_{\beta_1}s_{\beta_2}s_{\beta_4}s_{\gamma})
             (s_{\alpha_1}s_{\alpha_2}s_{\alpha_3}s_{\alpha_4})
             \stackrel{s_{\alpha_4}}{\simeq}
             s_{\alpha_4}s_{\beta_2}s_{\beta_4}(s_{\beta_1}s_{\gamma}s_{\alpha_1}s_{\alpha_2}s_{\alpha_3}) = \\
           &
           (s_{\beta_2}s_{\beta_4}s_{\mu})(s_{\beta_1}s_{\gamma}s_{\alpha_1}s_{\alpha_2}s_{\alpha_3}),
            \text{ where } \mu = \alpha_4 - \beta_2 + \beta_4.
    \end{split}
  \end{equation}

  \begin{figure}[h]
\centering
 \includegraphics[scale=0.8]{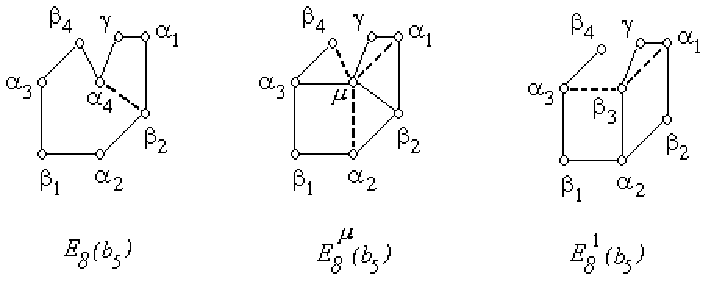}
 \caption{Transformation of $E_8(b_5)$, step $1$;
 the diagrams $E_8^{\mu}(b_5)$, $E_8^1(b_5)$ are connection ones, not
 Carter diagrams}
\label{E8_b5_step1}
\end{figure}
~\\
   We have
  \begin{equation*}
   \Small
    \begin{split}
       & (\mu, \alpha_3) = (\beta_4, \alpha_3) = -\frac{1}{2},
         \quad (\mu, \beta_4) = (\beta_4, \beta_4) + (\beta_4, \alpha_4) = 1 - \frac{1}{2} = \frac{1}{2}, \\
       & (\mu, \alpha_2) = -(\beta_2, \alpha_2) = \frac{1}{2},
         \quad (\mu, \beta_2) = -(\beta_2, \beta_2) + (\alpha_4, \beta_2) = -1 + \frac{1}{2} =
         -\frac{1}{2}, \\
       & (\mu, \alpha_1) = -(\beta_2, \alpha_1) = \frac{1}{2}.
    \end{split}
  \end{equation*}
  see $E_8^{\mu}(b_5)$ in Fig. \ref{E8_b5_step1}.  Further,
  \begin{equation}
    \label{eq_step1}
    \begin{split}
      w  ~\simeq~ & (s_{\beta_2}s_{\beta_4}s_{\beta_1})s_{\mu}s_{\alpha_2}s_{\alpha_3}
             (s_{\gamma}s_{\alpha_1}) = \\
           & (s_{\beta_2}s_{\beta_4}s_{\beta_1})(s_{\alpha_2}s_{\alpha_3}s_{\beta_3})(s_{\gamma}s_{\alpha_1}),
              \text{ where } \beta_3 = \mu - \alpha_2 + \alpha_3.
    \end{split}
  \end{equation}
  Here,
  \begin{equation*}
   \Small
    \begin{split}
       & (\beta_3, \alpha_3) = (\mu, \alpha_3) + (\alpha_3, \alpha_3) = -\frac{1}{2} + 1 = \frac{1}{2},
         \quad (\beta_3, \beta_4) = (\mu, \beta_4) + (\alpha_3, \beta_4) = \frac{1}{2} - \frac{1}{2} = 0, \\
       & (\beta_3, \alpha_2) = (\mu, \alpha_2) - (\alpha_2, \alpha_2) = \frac{1}{2} - 1 = -\frac{1}{2},
         \quad (\beta_3, \beta_2) = (\mu, \beta_2) - (\alpha_2, \beta_2) = \frac{1}{2} - \frac{1}{2} =
         0,
    \end{split}
  \end{equation*}
  see $E_8^1(b_5)$ in Fig. \ref{E8_b5_step1}.

  \begin{figure}[h]
\centering
\includegraphics[scale=0.8]{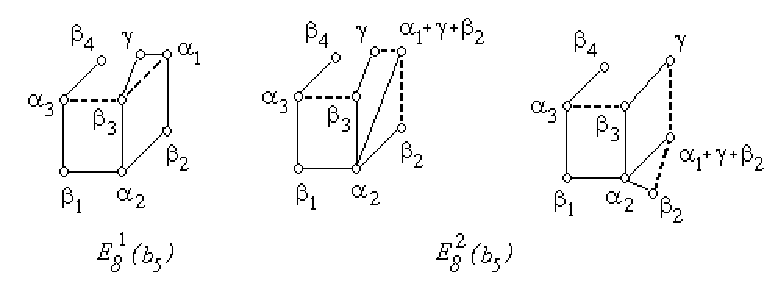}
\caption{Transformation of $E_8(b_5)$, step $2$;
diagrams $E_8^1(b_5)$, $E_8^2(b_5)$ are connection ones, not Carter
diagrams }
\label{E8_b5_step2}
\end{figure}

 {\it Step 2.} From eq. \eqref{eq_step1} we obtain
 \begin{equation}
   \label{eq_step2}
    \begin{split}
       w   \stackrel{s_{\beta_2}s_{\beta_4}s_{\beta_1}}{\simeq} & \quad
           s_{\alpha_2}s_{\alpha_3}s_{\beta_3}(s_{\alpha_1 + \gamma}s_{\gamma})
           s_{\beta_2}s_{\beta_4}s_{\beta_1}  = \\
         &  {\alpha_2}s_{\alpha_3}s_{\beta_3}s_{\alpha_1 + \gamma}
           (s_{\beta_2}s_{\beta_4}s_{\beta_1})s_{\gamma} =
            s_{\alpha_2}s_{\alpha_3}(s_{\beta_3}s_{\beta_2}s_{\beta_4}s_{\beta_1})
                s_{\alpha_1 + \gamma + \beta_2}s_{\gamma} \stackrel{s_{\alpha_2}s_{\alpha_3}}{\simeq} \\
         &  (s_{\beta_3}s_{\beta_2}s_{\beta_4}s_{\beta_1})s_{\alpha_1 + \gamma + \beta_2}
               s_{\alpha_2}s_{\alpha_3}s_{\gamma} =
            s_{\beta_2}(s_{\beta_1}s_{\beta_3}s_{\beta_4}s_{\alpha_1 + \gamma + \beta_2})
              (s_{\alpha_2}s_{\alpha_3}s_{\gamma}),
    \end{split}
  \end{equation}
  where
  \begin{equation*}
  \Small
    \begin{split}
       & (\alpha_1 + \gamma + \beta_2, \beta_3) = (\beta_3, \alpha_1) + (\beta_3, \gamma) =
       \frac{1}{2} - \frac{1}{2} = 0, \\
       & (\alpha_1 + \gamma + \beta_2, \gamma) = (\gamma, \gamma) +
       (\gamma, \alpha_1) = 1 - \frac{1}{2} = \frac{1}{2}, \\
       & (\alpha_1 + \gamma + \beta_2, \alpha_2) = (\beta_2, \alpha_2) = -\frac{1}{2}, \\
       & (\alpha_1 + \gamma + \beta_2, \beta_2) = (\beta_2, \beta_2) + (\beta_2, \alpha_1) =
          1 - \frac{1}{2} = \frac{1}{2},
    \end{split}
  \end{equation*}
   see $E_8^2(b_5)$ in Fig. \ref{E8_b5_step2}.

 \begin{figure} [h]
 \centering
 \includegraphics[scale=0.8]{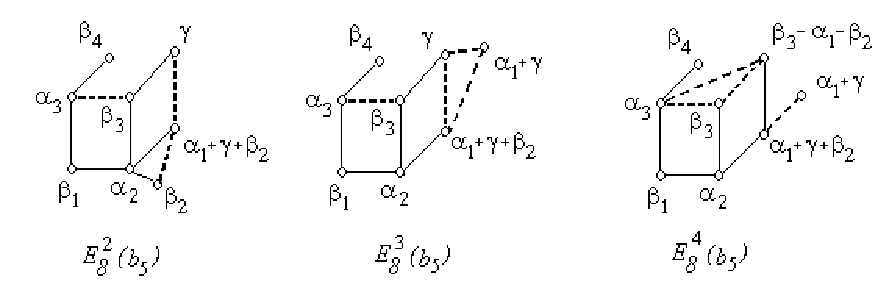}
 \caption{Transformations from $E_8^2(b_5)$ to $E_8^3(b_5)$ and from $E_8^3(b_5)$ to $E_8^4(b_5)$;
   the diagrams $E_8^2(b_5)$, $E_8^3(b_5)$, $E_8^4(b_5)$ are connection ones, not Carter diagrams}
 \label{E8_b5_step3and4}
 \end{figure}

   {\it Step 3.} Let us transform the $E_8^2(b_5)$-associated element $w$ from eq. \eqref{eq_step2}
   to a certain $E_8^3(b_5)$-associated element (where $E_8^2(b_5)$ and $E_8^3(b_5)$ are connection diagrams,
   see Fig. \ref{E8_b5_step3and4}):
  \begin{equation}
   \label{eq_step3}
    \begin{split}
    w ~\simeq~
   & s_{\beta_2}(s_{\beta_1}s_{\beta_3}s_{\beta_4}s_{\alpha_1 + \gamma + \beta_2})
    (s_{\alpha_2}s_{\alpha_3}s_{\gamma}) = \\
   & (s_{\beta_1}s_{\beta_3}s_{\beta_4})(s_{\beta_2}s_{\alpha_1 + \gamma + \beta_2})
    (s_{\alpha_2}s_{\alpha_3}s_{\gamma}) =
    (s_{\beta_1}s_{\beta_3}s_{\beta_4})(s_{\alpha_1 + \gamma + \beta_2}s_{\alpha_1 + \gamma})
    (s_{\alpha_2}s_{\alpha_3}s_{\gamma}),
    \end{split}
  \end{equation}
 where
  \begin{equation*}
   \Small
    \begin{split}
       & (\alpha_1 + \gamma,  \alpha_2) =
         (\alpha_1 + \gamma + \beta_2, \alpha_2) - (\beta_2, \alpha_2)  = -\frac{1}{2} + (\beta_2, \alpha_2) = 0, \\
       & (\alpha_1 + \gamma, \tau) = (\alpha_1 + \gamma + \beta_2, \tau) -  (\beta_2, \tau) = 0
         \text{ for } \tau = \beta_1, \beta_3, \beta_4, \alpha_3, \\
       & (\alpha_1 + \gamma, \gamma) =  -\frac{1}{2} + 1 = \frac{1}{2}.
    \end{split}
  \end{equation*}

   {\it Step 4.} Now, we transform the $E_8^3(b_5)$-associated element $w$ from eq. \eqref{eq_step3}
    into a certain $E_8^4(b_5)$-associated element ($E_8^3(b_5)$ and $E_8^3(b_5)$ are connection diagrams,
    see Fig. \ref{E8_b5_step3and4}):
  \begin{equation}
   \label{eq_step4}
    \begin{split}
    w ~\simeq~
      & (s_{\beta_1}s_{\beta_3}s_{\beta_4})(s_{\alpha_1 + \gamma + \beta_2}s_{\alpha_1 + \gamma})
        (s_{\alpha_2}s_{\alpha_3}s_{\gamma})
      \stackrel{s_{\gamma}}{\simeq} \\
      & s_{\gamma}s_{\beta_3}s_{\alpha_1 + \gamma + \beta_2}(s_{\beta_1}s_{\beta_4})s_{\alpha_1 + \gamma}
        (s_{\alpha_2}s_{\alpha_3}) =
        s_{\beta_3}s_{\alpha_1 + \gamma + \beta_2}s_{\beta_3 - \alpha_1 - \beta_2}
        s_{\alpha_1 + \gamma}s_{\beta_1}s_{\beta_4}(s_{\alpha_2}s_{\alpha_3}),  \\
    \end{split}
  \end{equation}
 since
    $\beta_3 - \alpha_1 - \beta_2 = \gamma + \beta_3 - (\alpha_1 + \gamma + \beta_2), \text{ and }
       s_{\gamma}s_{\beta_3}s_{\alpha_1 + \gamma + \beta_2} =
       s_{\beta_3}s_{\alpha_1 + \gamma + \beta_2}s_{\beta_3 - \alpha_1 - \beta_2}.$

  Here, we have
  \begin{equation*}
   \Small
    \begin{split}
       & (\beta_3 - \alpha_1 - \beta_2, \alpha_1 + \gamma) =
         (\gamma + \beta_3, \alpha_1 + \gamma)  - (\alpha_1 + \gamma + \beta_2, \alpha_1 + \gamma) =
          \frac{1}{2} - \frac{1}{2} = 0,  \\
       & (\beta_3 - \alpha_1 - \beta_2, \tau) =
         (\gamma + \beta_3, \tau) - (\alpha_1 + \gamma + \beta_2, \tau) = 0
         \text{ for } \tau  = \beta_1, \beta_4, \\
       & (\beta_3 - \alpha_1 - \beta_2, \alpha_2) = (\beta_3,  \alpha_2) - (\alpha_1 + \gamma + \beta_2, \alpha_2)  =
           \frac{1}{2} - \frac{1}{2} = 0, \\
       & (\beta_3 - \alpha_1 - \beta_2, \beta_3)  =
         (\gamma + \beta_3, \beta_3) - (\alpha_1 + \gamma + \beta_2, \beta_3) = (\gamma + \beta_3, \beta_3) =
         1 - \frac{1}{2} = \frac{1}{2},  \\
       & (\beta_3 - \alpha_1 - \beta_2, \alpha_1 + \gamma + \beta_2)  =
         (\gamma + \beta_3, \alpha_1 + \gamma + \beta_2) -
         (\alpha_1 + \gamma + \beta_2, \alpha_1 + \gamma + \beta_2) = \\
       &  \qquad \qquad (\gamma, \alpha_1 + \gamma + \beta_2) - 1 = \frac{1}{2} - 1 = -\frac{1}{2}, \\
       & (\beta_3 - \alpha_1 - \beta_2, \alpha_3)  = (\beta_3, \alpha_3) = -\frac{1}{2}.
    \end{split}
  \end{equation*}

\begin{figure} [h]
\centering
\includegraphics[scale=1.4]{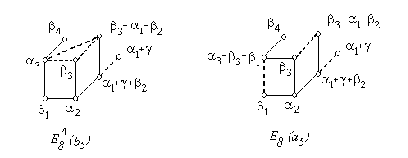}
 \caption{The last step of transformation: From $E_8^4(b_5)$ to $E_8(a_5)$}
 \label{E8_a5_fin}
 \end{figure}

   {\it Step 5.}
   The last step: From $E_8^4(b_5)$ to $E_8(a_5)$, see Fig. \ref{E8_a5_fin}.
   The $E_8^4(b_5)$-associated element $w$ from \eqref{eq_step4} is transformed as follows:

 \begin{equation}
   \label{eq_step5}
    \begin{split}
    w ~\simeq~
      & s_{\beta_3}s_{\alpha_1 + \gamma + \beta_2}s_{\beta_3 - \alpha_1 - \beta_2}
      s_{\alpha_1 + \gamma}s_{\beta_1}s_{\beta_4}
        (s_{\alpha_2}s_{\alpha_3}) =  \\
      & s_{\beta_3}s_{\beta_1}s_{\alpha_1 + \gamma + \beta_2}s_{\beta_3 - \alpha_1 - \beta_2}
        s_{\alpha_1 + \gamma}s_{\beta_4}
            (s_{\alpha_2}s_{\alpha_3})  \stackrel{s_{\alpha_3}}{\simeq} \\
      &
       (s_{\alpha_3}s_{\beta_3}s_{\beta_1})s_{\alpha_1 + \gamma + \beta_2}s_{\beta_3 - \alpha_1 - \beta_2}
        s_{\alpha_1 + \gamma}s_{\beta_4}s_{\alpha_2}  \stackrel{s_{\beta_4}}{\simeq} \\
      & (s_{\beta_4}s_{\beta_3}s_{\beta_1}s_{\alpha_1 + \gamma + \beta_2})
       (s_{\alpha_3 - \beta_3 + \beta_1}s_{\beta_3 - \alpha_1 - \beta_2}s_{\alpha_1 + \gamma}s_{\alpha_2}),
    \end{split}
  \end{equation}
 where\footnotemark[1]
 \footnotetext[1]{Recall that $\beta_3 = \mu - \alpha_2 + \alpha_3 = \alpha_4 - \beta_2 + \beta_4 - \alpha_2 + \alpha_3$,
  see \eqref{eq_step0}, \eqref{eq_step1}.}
  \begin{equation*}
    \Small
    \begin{split}
       & (\alpha_3 - \beta_3 + \beta_1,  \alpha_2) =
         -(\beta_3, \alpha_2) + (\beta_1, \alpha_2)  = -\frac{1}{2} + \frac{1}{2} = 0, \\
       & (\alpha_3 - \beta_3 + \beta_1, \tau) = (\alpha_3, \tau) -  (\beta_3, \tau)  + (\beta_1, \tau) = 0
         \text{ for } \tau = \alpha_1 + \gamma, \alpha_1 + \gamma + \beta_2, \\
       & (\alpha_3 - \beta_3 + \beta_1, \beta_3) = (\alpha_3, \beta_3) - (\beta_3, \beta_3) =
           \frac{1}{2} - 1 = - \frac{1}{2}, \\
       & (\alpha_3 - \beta_3 + \beta_1, \beta_1)  = 1 - \frac{1}{2} = \frac{1}{2}, \\
       & (\alpha_3 - \beta_3 + \beta_1, \beta_4)  = (\alpha_3, \beta_4) = - \frac{1}{2}, \\
       & (\alpha_3 - \beta_3 + \beta_1, \beta_3 - \beta_1 - \beta_2) =
          (\alpha_3, \beta_3 - \beta_1 - \beta_2) - (\beta_3, \beta_3 - \beta_1 - \beta_2)  =
           \frac{1}{2}  - \frac{1}{2}  = 0.
    \end{split}
  \end{equation*}

  Thus, \eqref{eq_step5} is a bicolored decomposition of the $E_8(a_5)$-associated element $w$,
  see Fig. \ref{E8_a5_fin}. The equivalence $E_8(b_5) \simeq E_8(a_5)$ is proven.

 \subsection{Equivalence $D_l(b_{\frac{l}{2} - 1}) \simeq D_l(a_{\frac{l}{2} - 1})$ }
  \label{sec_pure_cycle}
 \begin{figure}[h]
 \centering
 \includegraphics[scale=0.8]{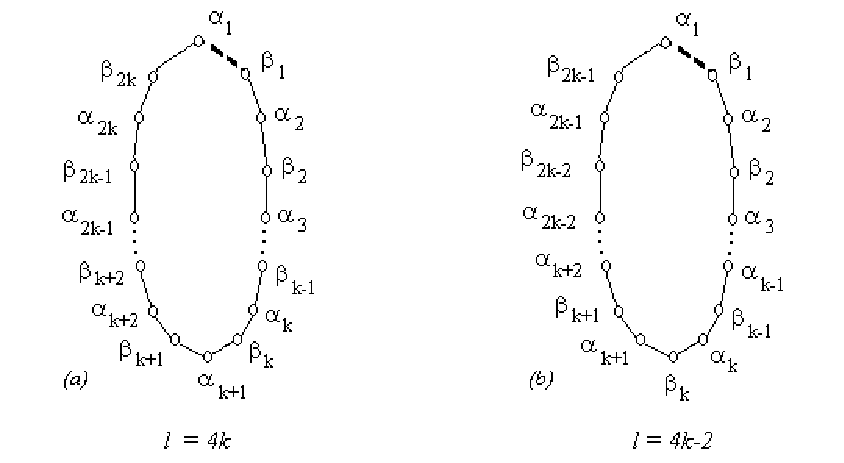}
 \caption{The two cases of even cycles $D_l(b_{\frac{l}{2} - 1})$: 1) $l = 4k$;  2)  $l = 4k-2$ }
 \label{two_circles}
 \end{figure}

 We consider the two cases of cycles $D_l(b_{\frac{l}{2} - 1})$ differing by
 length $l$, see Fig. \ref{two_circles}.

{\it  Case 1).} $l = 4k$. The opposite vertices, i.e., vertices at distance $2k$,
 are of the same type, for example, $\alpha_1$ and $\alpha_{k+1}$, see Fig. \ref{two_circles}$(a)$.

{\it  Case 2).} $l = 4k-2$. The opposite vertices, i.e., vertices at
distance $2k-1$, are of different types, for example, $\alpha_1$
and $\beta_{k}$, see Fig. \ref{two_circles}$(b)$.

\subsubsection{The case $l = 4k$}
 Consider the chains of vertices passing through the top vertex $\alpha_1$ and with endpoints
 lying on the same horizontal level, see Fig. \ref{two_circles}.
 Let $L$ (resp. $R$) be the index of the left (resp. right) end of the chain.
 Then the endpoints of these chains are as follows:
 \begin{equation}
  \label{chains}
  \begin{split}
     \{ \beta_L, \beta_R \},  & \quad L = 2k - i + 1, R = i,  \quad 1 \leq i \leq k, \text{ or }  \\
     \{ \alpha_L, \alpha_R \}, & \quad L = 2k - i + 2, R = i,  \quad 2 \leq i \leq k.
  \end{split}
 \end{equation}
 Consider the following vectors associated with chains \eqref{chains}:
 \begin{equation}
  \label{vect_chains}
  \begin{split}
     \theta(\beta_L, \beta_R) = & ~\alpha_1
              - \sum\limits_{i=1}^R{\beta_i} - \sum\limits_{i=2}^R{\alpha_i}
              + \sum\limits_{i=L}^{2k}{\beta_i} + \sum\limits_{i=L+1}^{2k}{\alpha_i},
              \quad R + L = 2k+1, \\
     \theta(\alpha_L, \alpha_R) = & ~\alpha_1
              - \sum\limits_{i=1}^{R-1}{\beta_i} - \sum\limits_{i=2}^R{\alpha_i}
              + \sum\limits_{i=L}^{2k}{\beta_i} + \sum\limits_{i=L}^{2k}{\alpha_i},
              \quad R + L = 2k+2. \\
  \end{split}
 \end{equation}
 We have the following actions on vectors \eqref{vect_chains}:
 \begin{equation}
   \label{root_actions}
  \begin{split}
       s_{\beta_1}s_{\beta_{2k}}\alpha_1 = & ~\theta(\beta_1, \beta_{2k}), \\
       s_{\alpha_2}s_{\alpha_{2k}}\theta(\beta_1, \beta_{2k}) = & ~\theta(\alpha_2, \alpha_{2k}), \\
       \dots \\
       s_{\beta_{L-1}}s_{\beta_{R}}\theta(\alpha_L, \alpha_R) = & ~\theta(\beta_{L-1}, \beta_R), \\
       s_{\alpha_L}s_{\alpha_{R+1}}\theta(\beta_L, \beta_R) = & ~\theta(\alpha_L, \alpha_{R+1}). \\
  \end{split}
 \end{equation}
 Thus, $\theta(\beta_L, \beta_R)$, $ \theta(\alpha_L, \alpha_R)$ from \eqref{vect_chains} are roots.
 The following orthogonality relations hold
 \begin{equation}
   \label{rel_orth}
  \begin{array}{lll}
     \theta(\beta_L, \beta_R) \perp \beta_i, & i \neq R, L,
        & \theta(\beta_L, \beta_R) \notperp  \beta_L, \beta_R, \\
     \theta(\beta_L, \beta_R) \perp \alpha_i, & i \neq R+1, L,
        & \theta(\beta_L, \beta_R)  \notperp \alpha_{R+1}, \alpha_L (R \neq k),
       \quad \theta(\beta_{k+1}, \beta_k) \perp \alpha_{k+1},    \\
     \theta(\alpha_L, \alpha_R) \perp \beta_i, & i  \neq L-1, R,
        & \theta(\alpha_L, \alpha_R) \notperp  \beta_{L-1}, \beta_R, \\
      \theta(\alpha_L, \alpha_R) \perp \alpha_i, &  i \neq R, L,
        & \theta(\alpha_L, \alpha_R) \notperp \alpha_L,\alpha_R.
  \end{array}
 \end{equation}
 \begin{lemma}
   \label{lem_commut}
   The following commutation relations hold:
 \begin{equation}
   \label{rel_commut}
  \begin{array}{lll}
     s_{\theta(\beta_L, \beta_R)} \prod\limits_{i=1}^{2k}{s_{\alpha_i}}
        =  \big( \prod\limits_{i=1}^{2k}{s_{\alpha_i}} \big) s_{\theta(\alpha_L, \alpha_{R+1})},
        & \text{ where } \quad L + R = 2k + 1, &  R \leq k, \\
     s_{\theta(\alpha_L, \alpha_R)}\prod\limits_{i=1}^{2k}{s_{\beta_i}}
        = \big(\prod\limits_{i=1}^{2k}{s_{\beta_i}}\big)s_{\theta(\beta_{L-1}, \beta_R)},
       & \text{ where } \quad  L + R = 2k + 2,  & R \leq k.
  \end{array}
 \end{equation}
 \end{lemma}

 \PerfProof
 According to the orthogonality relations \eqref{rel_orth}, we have
 \begin{equation*}
  \begin{split}
     s_{\theta(\beta_L, \beta_R)} \prod\limits_{i=1}^{2k}{s_{\alpha_i}} & =
       \big(\prod\limits_{\alpha_i \neq R+1, L}{s_{\alpha_i}}\big)
       s_{\theta(\beta_L, \beta_R)}s_{\alpha_{R+1}}s_{\alpha_L} = \\
     & \big(\prod\limits_{\alpha_i \neq R+1, L}{s_{\alpha_i}}\big)s_{\alpha_{R+1}}s_{\alpha_L}
        s_{\theta(\beta_L, \beta_R) - \alpha_{R+1} + \alpha_L} =
       \big(\prod\limits_{i=1}^{2k}{s_{\alpha_i}}\big)s_{\theta(\alpha_L, \alpha_{R+1})}.
  \end{split}
 \end{equation*}
  Similarly,
 \begin{equation*}
  \begin{split}
     s_{\theta(\alpha_L, \alpha_R)}\prod\limits_{i=1}^{2k}{s_{\beta_i}} & =
       \big(\prod\limits_{\beta_i \neq R, L-1}{s_{\beta_i}}\big)
       s_{\theta(\alpha_L, \alpha_R)}s_{\beta_{L-1}}s_{\beta_R} = \\
      & (\prod\limits_{\beta_i \neq R, L-1}{s_{\beta_i}}\big)s_{\beta_{L-1}}s_{\beta_R}
         s_{\theta(\alpha_L, \alpha_R) - \beta_R + \beta_{L-1}} =
        \big(\prod\limits_{i=1}^{2k}{s_{\beta_i}}\big)s_{\theta(\beta_{L-1}, \beta_R}. \qed
  \end{split}
 \end{equation*}

 \begin{proposition}
  \label{prop_main_4k}
 Let
 \begin{equation*}
    w = w_{\beta}w_{\alpha} = \prod\limits_{i=1}^{2k}{s_{\beta_i}}\prod\limits_{i=1}^{2k}{s_{\alpha_i}}
 \end{equation*}
 be the $D_l(b_{\frac{l}{2} - 1})$-associated element, where $D_l(b_{\frac{l}{2} - 1})$ is the cycle with $l = 4k$,
 see Fig. $\ref{two_circles}$.
 The element $w$ is conjugate to the element
 \begin{equation}
   \label{rel_main_4k}
    \big(\prod\limits_{i=1}^{2k}{s_{\beta_i}}\big)
    s_{\theta(\beta_{k+1}, \beta_k)}
    \big(\prod\limits_{i=2}^{2k}{s_{\alpha_i}}\big).
 \end{equation}
 \end{proposition}

 \PerfProof
  First, we have
 \begin{equation*}
  \begin{split}
   w =   \prod{s_{\beta_i}}\prod{s_{\alpha_j}}  \stackrel{s_{\alpha_1}}{\simeq}
          s_{\alpha_1}(s_{\beta_1} s_{\beta_{2k}})
      \prod\limits_{i \neq 1, 2k}{s_{\beta_i}}
      \prod\limits_{j \neq 1}{s_{\alpha_j}} =
      s_{\beta_1} s_{\beta_{2k}}s_{\theta(\beta_{2k}, \beta_{1})}
      \prod\limits_{i \neq 1, 2k}{s_{\beta_i}}
      \prod\limits_{j \neq 1}{s_{\alpha_j}}.
  \end{split}
 \end{equation*}
 By relations \eqref{rel_orth}, the elements  $s_{\theta(\beta_{1}, \beta_{2k})}$ and
  $\prod\limits_{i \neq 1, 2k}{s_{\beta_i}}$ commute, and we have:
 \begin{equation*}
   w =   s_{\beta_1} s_{\beta_{2k}} \big(\prod\limits_{i \neq 1, 2k}{s_{\beta_i}}\big)
       s_{\theta(\beta_{2k}, \beta_{1})} \prod\limits_{j \neq 1}{s_{\alpha_j}} =
         \big(\prod{s_{\beta_i}}\big)
         s_{\theta(\beta_{2k}, \beta_{1})}
         \big(\prod\limits_{j \neq 1}{s_{\alpha_j}}\big).
 \end{equation*}
 Further, we use Lemma \ref{lem_commut} to prove the equivalences:
 \begin{equation*}
  \begin{split}
   w = & \prod{s_{\beta_i}}
      \big(\prod\limits_{j \neq 1}{s_{\alpha_j}}\big)
      s_{\theta(\alpha_{2k}, \alpha_{2})}
      \quad \stackrel{s_{\theta(\alpha_{2k}, \alpha_{2})}}{\simeq} \quad
      s_{\theta(\alpha_{2k}, \alpha_{2})}\prod{s_{\beta_i}}\prod\limits_{j \neq 1}{s_{\alpha_j}} =
      \big(\prod{s_{\beta_i}}\big)s_{\theta(\beta_{2k-1}, \beta_{2})}
      \prod\limits_{j \neq 1}{s_{\alpha_j}} = \\
      & \prod{s_{\beta_i}}
      \big(\prod\limits_{j \neq 1}{s_{\alpha_j}}\big)
              s_{\theta(\alpha_{2k-1}, \alpha_{3})}
           \quad \stackrel{s_{\theta(\alpha_{2k-1}, \alpha_{3})}}{\simeq} \quad
        s_{\theta(\alpha_{2k-1}, \alpha_{3})}\prod{s_{\beta_i}}\prod\limits_{j \neq 1}{s_{\alpha_j}} =
      \big(\prod{s_{\beta_i}}\big)s_{\theta(\beta_{2k-2}, \beta_{3})}
         \prod\limits_{j \neq 1}{s_{\alpha_j}} = \\
      & \dots  \\
      & \big(\prod{s_{\beta_i}}\big)s_{\theta(\beta_{k+1}, \beta_{k})}
         \prod\limits_{j \neq 1}{s_{\alpha_j}}.
  \end{split}
 \end{equation*}
  The relation \eqref{rel_main_4k} is proved. \qed

\begin{figure}[h]
 \centering
 \includegraphics[scale=0.8]{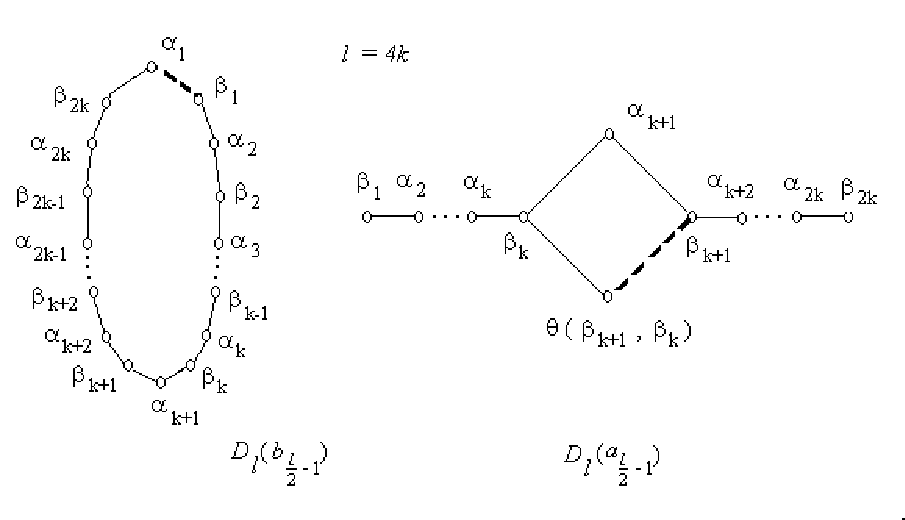}
 \caption{The equivalent diagrams $D_l(b_{\frac{l}{2} - 1})$ and
  $D_l(a_{\frac{l}{2} - 1})$, where $l = 4k$ }
 \label{circleCase_4k}
 \end{figure}

 \begin{corollary}
     The conjugacy class containing elements
 \begin{equation}
     \prod\limits_{i=1}^{2k}{s_{\beta_i}}\prod\limits_{i=1}^{2k}{s_{\alpha_i}}
     ~\simeq~
        \big(\prod\limits_{i=1}^{2k}{s_{\beta_i}}\big)
    \big(s_{\theta(\beta_{k+1}, \beta_k)}
        \prod\limits_{i=2}^{2k}{s_{\alpha_i}}\big)
 \end{equation}
     is $D_l(a_{\frac{l}{2} - 1})$-associated (as well $D_l(b_{\frac{l}{2} - 1})$-associated)
     conjugacy class for $l = 4k$, see Fig. $\ref{circleCase_4k}$.
 \end{corollary}

 \PerfProof For $i \neq k+1$, the orthogonality $\theta(\beta_{k+1}, \beta_k) \perp \alpha_i$
 follows from \eqref{rel_orth}. For $i = k+1$, it is easy to check:
 \begin{equation*}
    (\theta(\beta_{k+1}, \beta_k), \alpha_{k+1}) =
       (\beta_{k+1}, \alpha_{k+1}) - (\beta_k, \alpha_{k+1}) =  -\frac{1}{2} + \frac{1}{2}  = 0.
 \end{equation*}
 Besides, for $i \neq k, k+1$, we have $\theta(\beta_{k+1}, \beta_k) \perp \beta_i$, see \eqref{rel_orth}.
  Finally, for $i = k, k+1$, we have:
 \begin{equation*}
   \begin{split}
   &  (\theta(\beta_{k+1}, \beta_k), \beta_k) =
       (-\beta_k, \beta_k) + (-\alpha_k, \beta_k) = -1 + \frac{1}{2} = -\frac{1}{2}, \\
   &  (\theta(\beta_{k+1}, \beta_k), \beta_{k+1}) =
       (\beta_{k+1}, \beta_{k+1}) + (\alpha_{k+2}, \beta_{k+1}) = 1 -  \frac{1}{2} =  \frac{1}{2}.
       \qed
   \end{split}
 \end{equation*}

\subsubsection{The case $l = 4k-2$}
  Similarly to the case \eqref{chains}, we consider chains
 \begin{equation}
  \label{chains_2}
  \begin{split}
     \{ \beta_L, \beta_R \},  & \quad L = 2k - i, R = i,  \quad 1 \leq i \leq {k - 1}, \text{ or }  \\
     \{ \alpha_L, \alpha_R \}, & \quad L = 2k - i + 1, R = i,  \quad 2 \leq i \leq k.
  \end{split}
 \end{equation}
 Next, we consider the following vectors associated with the chains \eqref{chains_2}:
 \begin{equation}
  \label{vect_chains_2}
  \begin{split}
     \mu(\beta_L, \beta_R) = & \alpha_1
              - \sum\limits_{i=1}^R{\beta_i} - \sum\limits_{i=2}^R{\alpha_i}
              + \sum\limits_{i=L}^{2k-1}{\beta_i} + \sum\limits_{i=L+1}^{2k-1}{\alpha_i},
              \quad R + L = 2k, \\
     \mu(\alpha_L, \alpha_R) = & \alpha_1
              - \sum\limits_{i=1}^{R-1}{\beta_i} - \sum\limits_{i=2}^R{\alpha_i}
              + \sum\limits_{i=L}^{2k-1}{\beta_i} + \sum\limits_{i=L}^{2k-1}{\alpha_i},
              \quad R + L = 2k+1. \\
  \end{split}
 \end{equation}
 As above, vectors $\mu(\beta_L, \beta_R)$, $\mu(\alpha_L, \alpha_R)$ from \eqref{vect_chains_2} are roots.

 \begin{lemma}
   \label{lem_commut_2}
   The following commutation relations hold:
 \begin{equation}
   \label{rel_commut_2}
  \begin{array}{lll}
     s_{\mu(\beta_L, \beta_R)} \prod\limits_{i=1}^{2k-1}{s_{\alpha_i}}
        =  \big(\prod\limits_{i=1}^{2k-1}{s_{\alpha_i}}\big)s_{\mu(\alpha_L, \alpha_{R+1})},
         & \text{ where } \quad L + R = 2k,  & R \leq k-1, \\
     s_{\mu(\alpha_L, \alpha_R)}\prod\limits_{i=1}^{2k-1}{s_{\beta_i}}
        = \big(\prod\limits_{i=1}^{2k-1}{s_{\beta_i}}\big)s_{\mu(\beta_{L-1}, \beta_R)},
         & \text{ where } \quad L + R = 2k + 1,  & R \leq k.
  \end{array}
 \end{equation}
 \end{lemma}
  Proof is as in Lemma \ref{lem_commut}. \qed

 \begin{proposition}
   \label{prop_main_4k-2}
 Let
 \begin{equation*}
    w = w_{\beta}w_{\alpha} = \prod\limits_{i=1}^{2k-1}{s_{\beta_i}}\prod\limits_{i=1}^{2k-1}{s_{\alpha_i}}
 \end{equation*}
 be the $D_l(b_{\frac{l}{2} - 1})$-associated element, where $D_l(b_{\frac{l}{2} - 1})$ is the cycle with $l = 4k-2$,
 see Fig. $\ref{two_circles}$.
 The element $w$ is conjugate to the element
 \begin{equation}
   \label{rel_main_4k_2}
    s_{\mu(\alpha_{k+1}, \alpha_k)}(\prod\limits_{i=1}^{2k-1}{s_{\beta_i}})
    (\prod\limits_{i=2}^{2k-1}{s_{\alpha_i}}).
 \end{equation}
 \end{proposition}

 \PerfProof
  As in Proposition \ref{prop_main_4k}, we  have
 \begin{equation*}
  \begin{split}
   w =   & \prod{s_{\beta_i}}\prod{s_{\alpha_j}}  \stackrel{s_{\alpha_1}}{\simeq}
          s_{\alpha_1}(s_{\beta_1} s_{\beta_{2k-1}})
      \prod\limits_{i \neq 1, 2k-1}{s_{\beta_i}}
      \prod\limits_{j \neq 1}{s_{\alpha_j}} =
      s_{\beta_1} s_{\beta_{2k-1}}s_{\mu(\beta_{2k-1}, \beta_{1})}
      \prod\limits_{i \neq 1, 2k-1}{s_{\beta_i}}
      \prod\limits_{j \neq 1}{s_{\alpha_j}} = \\
      & s_{\beta_1} s_{\beta_{2k-1}} \big(\prod\limits_{i \neq 1, 2k-1}{s_{\beta_i}}\big)
       s_{\mu(\beta_{2k-1}, \beta_{1})} \prod\limits_{j \neq 1}{s_{\alpha_j}} =
         \big(\prod{s_{\beta_i}}\big)s_{\mu(\beta_{2k-1}, \beta_{1})}
         \big(\prod\limits_{j \neq 1}{s_{\alpha_j}}\big).
  \end{split}
 \end{equation*}

 \begin{figure}[h]
 \centering
 \includegraphics[scale=0.57]{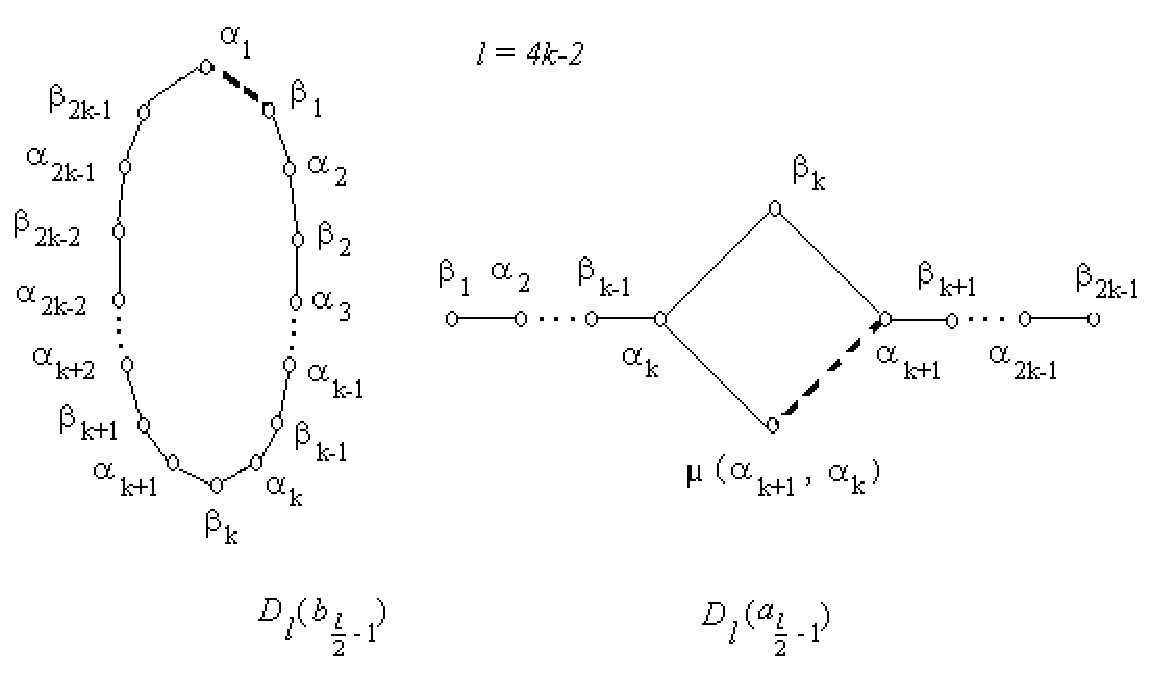}
 \caption{The equivalent diagrams $D_l(b_{\frac{l}{2} - 1})$ and
 $D_l(a_{\frac{l}{2} - 1})$, where $l = 4k-2$ }
 \label{circleCase_4k_2}
 \end{figure}

 By Lemma \ref{lem_commut_2}, we have:
 \begin{equation*}
  \begin{split}
   w = & \prod{s_{\beta_i}}
      \big(\prod\limits_{j \neq 1}{s_{\alpha_j}}\big)
      s_{\mu(\alpha_{2k-1}, \alpha_{2})}
      \quad \stackrel{s_{\mu(\alpha_{2k-1}, \alpha_{2})}}{\simeq} \quad
      s_{\mu(\alpha_{2k-1}, \alpha_{2})}\prod{s_{\beta_i}}\prod\limits_{j \neq 1}{s_{\alpha_j}} =
      \big(\prod{s_{\beta_i}}\big)s_{\mu(\beta_{2k-2}, \beta_{2})}
      \prod\limits_{j \neq 1}{s_{\alpha_j}} = \\
      & \prod{s_{\beta_i}}
        \big(\prod\limits_{j \neq 1}{s_{\alpha_j}}\big)
              s_{\mu(\alpha_{2k-2}, \alpha_{3})}
           \quad \stackrel{s_{\mu(\alpha_{2k-2}, \alpha_{3})}}{\simeq} \quad
        s_{\mu(\alpha_{2k-2}, \alpha_{3})}\prod{s_{\beta_i}}\prod\limits_{j \neq 1}{s_{\alpha_j}} =
        \big(\prod{s_{\beta_i}}\big)
        s_{\mu(\beta_{2k-3}, \beta_{3})}
         \prod\limits_{j \neq 1}{s_{\alpha_j}} = \\
      & \dots  \\
      &  s_{\mu(\alpha_{k+1}, \alpha_{k})}\prod{s_{\beta_i}}\prod\limits_{j \neq 1}{s_{\alpha_j}}.
        \qed
        \end{split}
 \end{equation*}

 \begin{corollary}
     The conjugacy class containing elements
 \begin{equation}
    \prod\limits_{i=1}^{2k-1}{s_{\beta_i}}\prod\limits_{i=1}^{2k-1}{s_{\alpha_i}}
    ~\simeq~
    \big(s_{\mu(\alpha_{k+1}, \alpha_k)}\big(\prod\limits_{i=1}^{2k-1}{s_{\beta_i}}\big)
    \prod\limits_{i=2}^{2k-1}{s_{\alpha_i}}\big).
  \end{equation}
     is $D_l(a_{\frac{l}{2} - 1})$-associated (as well $D_l(b_{\frac{l}{2} - 1})$-associated)
     conjugacy class for $l = 4k-2$, see Fig. $\ref{circleCase_4k_2}$.
 \end{corollary}

 For $i \neq k$, the orthogonality $\mu(\alpha_{k+1}, \alpha_k) \perp \beta_i$
 follows from \eqref{vect_chains_2}. For $i = k$, we have:
 \begin{equation*}
    (\mu(\alpha_{k+1}, \alpha_k), \beta_k) =
       (\alpha_{k+1}, \beta_k) - (\alpha_k, \beta_k) =  -\frac{1}{2} + \frac{1}{2}  = 0.
 \end{equation*}

 For $i \neq k, k+1$, we have $\mu(\alpha_{k+1}, \alpha_k) \perp \alpha_i$, and,
  for $i = k, k+1$, we get:
 \begin{equation*}
   \begin{split}
   &  (\mu(\alpha_{k+1}, \alpha_k), \alpha_k) =
       (-\beta_{k-1}, \alpha_k) + (-\alpha_k, \alpha_k) = \frac{1}{2} - 1 = -\frac{1}{2}, \\
   &  (\mu(\alpha_{k+1}, \alpha_k), \alpha_{k+1}) =
       (\beta_{k+1}, \alpha_{k+1}) + (\alpha_{k+1}, \alpha_{k+1}) = 1 - \frac{1}{2} =  \frac{1}{2}.
       \qed
   \end{split}
 \end{equation*}

\clearpage
~\\
\section{\sc\bf Basic patterns}
  \label{sect_basic_pat}

 \subsection{The partial Cartan matrix}
    \label{sec_partial_B}

  Similarly to the Cartan matrix associated with a given Dynkin diagram we determine
  the {\it partial Cartan matrix} corresponding to the Carter diagram $\Gamma$ as follows:

 \begin{equation}
   \label{canon_dec_2}
   B_{\Gamma} :=
      \left (
        \begin{array}{cccccc}
         (\alpha_1, \alpha_1) & \dots & (\alpha_1, \alpha_k) &
         (\alpha_1, \beta_1)  & \dots &  (\alpha_1, \beta_h) \\
         \dots                & \dots & \dots &
         \dots                & \dots & \dots \\
         (\alpha_k, \alpha_1) & \dots & (\alpha_k, \alpha_k) &
         (\alpha_k, \beta_1)  & \dots &  (\alpha_k, \beta_h) \\
         (\beta_1, \alpha_1)  & \dots & (\beta_1, \alpha_k) &
         (\beta_1, \beta_1)   & \dots & (\beta_1, \beta_h) \\
         \dots                & \dots & \dots &
         \dots                & \dots & \dots \\
         (\beta_h, \alpha_1)  & \dots & (\beta_h, \alpha_k) &
         (\beta_h, \beta_1)   & \dots & (\beta_h, \beta_h) \\
        \end{array}
      \right ),
 \end{equation}
 Here, the  $\alpha$-set $S_{\alpha} = \{\alpha_i \mid i = 1,\dots,k$\}
 and the $\beta$-set $S_{\beta} = \{\beta_j \mid j = 1,\dots,h\}$ match
 the bicolored decomposition of a certain $w \in W$
 corresponding to $\Gamma$, see \S\ref{sec_adm_diagr}.
 We will write $S = \{ \tau_1,\dots,\tau_{k+h} \}$ instead of \eqref{eq_alpha_bet} if the
 bicolored decomposition is not important.

 \index{partial Cartan matrix $B_{\Gamma}$}
 \index{Cartan matrix! - partial}
 \index{quadratic form $\mathscr{B}_{\Gamma}$}
 \index{$B^{-1}_{\Gamma}$, the inverse of the partial Cartan matrix}
 \index{$S$-associated subspace}
 \index{Weyl group}

 The symmetric bilinear form associated with
 the partial Cartan matrix $B_{\Gamma}$ is denoted by $(\cdot\hspace{0.7mm},\cdot)_{\Gamma}$ and the
 corresponding quadratic form is denoted by $\mathscr{B}_{\Gamma}$.
 The subspace $L \subseteq V$ spanned by the root subset $S$ is said to be the {\it $S$-associated subspace}.
 For $S = \{\tau_1,\dots,\tau_l\}$, we write $L =[\tau_1,\dots,\tau_l]$.

\begin{remark}{\rm
 From now on, we assume that entries $(\tau_i, \tau_j)$ of the partial Cartan matrix $B_{\Gamma}$
 for simply-laced Carter diagrams
 take values $\{0, -1, 1 \}$, not $\{0, -\frac{1}{2}, \frac{1}{2}\}$.
 In other words, we assume that the symmetric bilinear form $(\cdot\hspace{0.7mm},\cdot)_{\Gamma}$ is obtained by doubling the
 usual inner product associated with the given Weyl group $W$.
 }
\end{remark}

 \begin{proposition}
   \label{restr_forms_coincide}
 1)  The restriction of the bilinear form associated with the Cartan matrix
  ${\bf B}$ to the subspace $L$ coincides with the bilinear form associated with the
  partial Cartan matrix $B_{\Gamma}$, i.e., for any pair of vectors $v, u \in L$, we have
 \begin{equation}
   \label{restr_q}
       (v, u)_{\botG} = (v, u), \text{ and }
       \mathscr{B}_{\Gamma}(v) = \mathscr{B}(v).
 \end{equation}

 2) For every Carter diagram, the matrix $B_{\botG}$ \underline{is positive definite}.

 \end{proposition}
\PerfProof
 1) From eq. \eqref{canon_dec_2} we deduce:
 \begin{equation*}
     (v, u)_{\botG} =   (\sum\limits_i{t_i{\tau_i}}, \sum\limits_j{q_j{\tau_j}})_{\botG} =
      \sum\limits_{i,j}t_i{q}_j(\tau_i, \tau_j)_{\botG} = \sum\limits_{i,j}t_i{q}_j(\tau_i, \tau_j) =
      (v, u).
 \end{equation*}

 2) This follows from 1).
\qed

  \index{generalized Cartan matrix}
  \index{Cartan matrix! - generalized}
  \index{Coxeter element}
\begin{remark}[The classical case]
   \label{rem_classic}
 \rm{
  Recall that the $n\times{n}$ matrix $K$ such that
 \begin{equation*}
   \begin{split}
     (C1) & \quad k_{ii} = 2 \text{ for } i = 1,\dots, n, \\
     (C2) & \quad -k_{ij} \in \mathbb{Z} = \{0, 1, 2, \dots \} \text{ for } i \neq j, \\
     (C3) & \quad  k_{ij} = 0 \text{ implies } k_{ji} = 0 \text{ for } i, j = 1, \dots, n
    \end{split}
 \end{equation*}
 is called a {\it generalized Cartan matrix}, \cite{Kac80}, \cite[\S 2.1]{St08}.

  The condition (C2) \underline{does not hold} for the partial Cartan matrix:
 The values $k_{ij}$ associated with dotted edges are positive.

 If the Carter diagram does not contain any cycle,
 then the Carter diagram is the Dynkin diagram, the corresponding conjugacy class
 is the conjugacy class of the Coxeter element, and
 the partial Cartan matrix is the classical Cartan matrix, a particular case of a generalized Cartan matrix. \qed
 }
\end{remark}

\subsection{Linear dependence and maximal roots}
  \index{linearly dependent root}
  \index{root system}
 Let $S = \{\tau_1,\dots,\tau_l\}$ be a $\Gamma$-associated subset, $S'$ be another $\Gamma$-associated subset,
 and $S' = uS$ for some element $u \in W$. The matrix $B_{\Gamma}$ is the same for $S$ and $S'$,
 since $(u\tau_i, u\tau_j) = (\tau_i, \tau_j)$ for any $\tau_i, \tau_j \in S$.
 Let $\gamma$ be a root which is linearly dependent on roots of $S$ as follows:
\begin{equation}
  \label{eq_lin_depend1}
   \gamma = t_1\tau_1 + \dots + t_l\tau_l.
\end{equation}

Then, we have
 \begin{equation}
  \label{eq_lin_depend2}
  \left (
    \begin{array}{c}
      (\gamma, \tau_1) \\
      \dots \\
      (\gamma, \tau_l) \\
    \end{array}
  \right ) =  B_{\Gamma}
  \left (
    \begin{array}{c}
      t_1 \\
      \dots\\
      t_l \\
    \end{array}
   \right ) = B_{\Gamma}\gamma,  \text{ and }
     \left (
    \begin{array}{c}
      t_1 \\
      \dots \\
      t_l \\
    \end{array}
   \right ) = {B}^{-1}_{\Gamma}
    \left (
    \begin{array}{c}
      (\gamma, \tau_1) \\
      \dots \\
      (\gamma, \tau_l) \\
    \end{array}
  \right ).
 \end{equation}

 If the root $\tau_i$ is replaced by $-\tau_i$ (for any $\tau_i \in S$), then
 the coefficient $t_i$ is replaced by $-t_i$ in the decomposition \eqref{eq_lin_depend1}.

 \index{$b^{\vee}_{\eta, \eta}$ (diagonal element of $B^{-1}_{\Gamma}$)}
\begin{remark}
  \label{rem_max_root}
{\rm
   Let the vector $\gamma$ be \underline{linearly dependent on roots} of $S = \{\tau_1,\dots,\tau_l\}$,
   let $\gamma$ be connected with only one $\tau_i \in S$.
   In further considerations, we have two frequently occurring cases:
   ~\\

   (i) Suppose $\gamma$ is connected to the same point as the maximal (or minimal) root in the root system $S$.
   In other words, in eq. \eqref{eq_lin_depend2},
   the orthogonality relations $(\gamma, \tau_i)$ coincide with orthogonality relations
   for the maximal (resp. minimal) root while the edge connecting with $\gamma$
   is dotted (resp. solid). Since equation \eqref{eq_lin_depend2} has a unique solution,
   we deduce that $\gamma$ coincides with the maximal (resp. minimal) root.
   ~\\

   (ii) Consider the necessary condition that $\gamma$ is a root.
   We have
   \begin{equation}
     \label{eq_vect_inner_prod}
       \gamma^{\vee} := \left (
    \begin{array}{c}
      (\gamma, \tau_1) \\
      \dots \\
      (\gamma, \tau_{i}) \\
      \dots \\
      (\gamma, \tau_{l}) \\
    \end{array}
  \right ) =
           \left (
    \begin{array}{c}
      0 \\
      \dots \\
      \pm{1} \\
      \dots \\
      0 \\
    \end{array}
  \right ), \qquad \gamma = B_{\Gamma}^{-1}\gamma^{\vee},
   \end{equation}
  where the sign $+$ (resp. $-$) corresponds to the dotted (resp. solid)
  edge connecting $\gamma$ with $\tau_i$.
  Let $\mathscr{B}$ be the quadratic form associated with the partial Cartan matrix $B_{\Gamma}$.
  Then the value of $\mathscr{B}_{\Gamma}$ on the root $\gamma$ is as follows
\begin{equation}
  \label{eq_vect_inner_prod_2}
  \mathscr{B}_{\Gamma}(\gamma) = \langle B_{\Gamma}\gamma, \gamma \rangle =
  \langle B_{\Gamma}(B_{\Gamma}^{-1}\gamma^{\vee}, B_{\Gamma}^{-1}\gamma^{\vee}  \rangle =
  \langle \gamma^{\vee}, B_{\Gamma}^{-1}\gamma^{\vee}  \rangle = b^{\vee}_{i,i},
\end{equation}
 where $b^{\vee}_{i,i}$ is the $i$th diagonal element of $B_{\Gamma}^{-1}$. If $\gamma$ is a root, then
 $\mathscr{B}(\gamma) = 2$, and the necessary condition that $\gamma$ ia a root is the
 following simple equality:
\begin{equation}
   \label{eq_necessry_cond}
    b^{\vee}_{i,i} = 2.
\end{equation}
}
\end{remark}

 \index{attachment point}

\begin{remark}{\rm
Let us call each vertex $\tau_i$ connected with a linearly dependent
root $\gamma$ an \emph{attachment point} on $\Gamma$.
~\\

 (i) Let $\Gamma$ be a simply-laced Dynkin diagram. By Proposition
 \ref{prop_determinants} and Table \ref{tab_Cartan_E6_E7_E8}, we have
 $b^{\vee}_{i,i} = 2$ if and only if $\tau_i$ is a single attachment
 point connected with the maximal (minimal) root. Since
 $\gamma=B_{\Gamma}^{-1}\gamma^{\vee}$, then $\gamma$ is uniquely
 defined  and coincides with the maximal (minimal) root. There is
 only one attachment point. To see this, it suffices to check
 diagonal elements of $B_{\Gamma}^{-1}$.
~\\

 (ii) Let $\Gamma$ be a simply-laced Carter diagram, but not a Dynkin
 diagram. In this case, there is no such thing as a maximal (minimal)
 root. However, there exists an attachment point, not necessarily unique.
 As in (i), $\gamma=B_{\Gamma}^{-1}\gamma^{\vee}$, so it is uniquely defined.
 For example, according to Table \ref{tab_partial Cartan_1}:

 For $E_6(a_2)$, we have $b^{\vee}_{\alpha_3,\alpha_3} = 2$ and
 $b^{\vee}_{\beta_1,\beta_1} = 2$. Thus, there exists a root $\alpha$
 (resp. $\beta$) linearly dependent on vectors of $S$ and connected to $\alpha_3$ (resp. $\beta_1$).

 For $D_6(a_1)$, we have $b^{\vee}_{\alpha_1,\alpha_1} = 2$, so there is
 the root linearly dependent on vectors of $S$ and connected to $\alpha_1$.

 For $D_6(a_2)$, we have $b^{\vee}_{\beta_1,\beta_1} = 2$ and
 $b^{\vee}_{\beta_2,\beta_2} = 2$. Thus, there exists a root $\beta$ (resp. $\beta'$) linearly dependent on
 vectors of $S$ and connected to $\beta_1$ (resp. $\beta_2$). \qed
}
\end{remark}

\subsection{The orbit of the Coxeter element, and the longest element}
  \label{sec_max_root}
  \index{Coxeter element}
  \index{Coxeter groups}
  \index{orbit of the Coxeter element}

We use the following fact from the theory of Coxeter groups:
\begin{proposition}{\em(\cite[Sect. 5.7, Proposition]{Hu04})}
  \label{prop_on_leng}
 Let $w \in W$, let $\alpha$ be a positive (not necessarily simple) root.
 Let $s_{\alpha} \in W$ be the reflection corresponding to $\alpha$. Then
  \begin{equation*}
    \begin{split}
     & l(ws_{\alpha}) > l(w) \text{ if and only if } w\alpha > 0, \\
     & l(s_{\alpha}w) > l(w) \text{ if and only if } w^{-1}\alpha > 0.
    \end{split}
  \end{equation*}
\end{proposition}
  Note that the condition in the second line in Proposition \ref{prop_on_leng}
  is equivalent to the one in the first line since
  $l(ws_{\alpha}) = l(s_{\alpha}w^{-1})$ and $l(w) = l(w^{-1})$. \qed

\subsubsection{The longest elements in $W(D_l)$ and $W(E_l)$}
   \label{sec_longest}

  Let $\Gamma$ be the Dynkin diagram $D_l$ or $E_l$,
  let $w \in W$ be any $\Gamma$-associated element, and
  \begin{equation}
   \label{eq_basis_1}
    S = \{\varphi_1,\dots,\varphi_k, \delta_1,\dots,\delta_m\}
  \end{equation}
  the $\Gamma$-associated subset of linearly independent roots corresponding to the bicolored decomposition of $w$.
  Here, $w = w_{\varphi}w_{\delta}$ is the Coxeter element in the basis  \eqref{eq_basis_1}:
  \begin{equation}
    \label{decomp_w_Dl}
      w = w_{\varphi}w_{\delta}, \quad \text{ where } \quad
      w_{\varphi} = \prod\limits_{i=1}^k{s_{\varphi_i}}, \quad
      w_{\delta} = \prod\limits_{j=1}^m{s_{\delta_j}}.
  \end{equation}

 \index{Coxeter number}
 \index{longest element $w_0$}
 \index{maximal root $\mu_{max}$}
 \index{$\mu_{max}$ (maximal root)}
   Let $h$ be the Coxeter number, the order of the Coxeter element $w$.
   Since $W(A_n)$ is excluded, then $h$ is even: $h = 2g$. Let
 \begin{equation}
  \label{eq_w0_1}
     w_0 = (w_{\delta}w_{\varphi})^{\frac{h}{2}}, \quad  w_0^2 = 1.
 \end{equation}
   It is a well-known fact that $w_0$ is the longest element in $W(\Gamma)$.
   The element $w_0$ makes all positive roots negative.

\begin{proposition}[On the maximal root $\mu_{max}$ and the longest element $w_0 \in W$]
  \label{prop_conj_max_root}
  Let $\Gamma$ be the Dynkin diagram $D_l$ or $E_l$.
  Suppose $S$ is given by \eqref{eq_basis_1} and $w \in W(\Gamma)$ is the $S$-associated element.
  Let $\Gamma$ be extended to another Dynkin diagram $\widetilde{\Gamma}$
  by adding a root $\delta_{m+1}$ connected to $S$ only at $\varphi_k$
  and linearly independent of $S$;
  let $\mu_{max}$ be the maximal root in the root system $\Phi(\widetilde\Gamma)$.

  {\rm (i)} For $w_0 = w^{\frac{h}{2}} =
  (w_{\delta}w_{\varphi})^{\frac{h}{2}}$, we have
  \begin{equation}
    \label{eq_reduced_mu_max}
     w_0 \mu_{max}  ~=~  \delta_{m+1}, \quad
     \mu_{max}  ~=~ w_0\delta_{m+1}.
  \end{equation}

  {\rm (ii)} The following conjugacy relation holds:
  \begin{equation}
    \label{eq_reduced_mu_max_2}
     w{s}_{\mu_{max}}  ~\simeq~ w{s}_{\delta_{m+1}}.
  \end{equation}

\end{proposition}

   \PerfProof
   (i) Note that $w_{\varphi}$ and $w_{\delta}$ do not act on the coordinate $\delta_{m+1}$
   and this coordinate in $w_0\mu_{max}$ is positive.
   Since $w_0\mu_{max}$ is a root, we see that $w_0\mu_{max} > 0$.
   In particular, $w_0\mu_{max} \geq \delta_{m+1}$.
   Further, since $w_0$ is the longest element
   in $W(\Gamma)$, then $l(w_0{s}_{\beta}) < l(w_0)$
   for every reflection $s_{\beta} \in W(\Gamma)$.
   By Proposition \ref{prop_on_leng}, we have $w_0\beta \leq 0$ for every positive simple  root $\beta$.
   So $w_0\gamma \leq 0$  for every positive vector $\gamma \in  \mathcal{E}(\Gamma)$,
   where $\mathcal{E}$ is the linear space spanned by all roots
   $\gamma \in S$. Let $w_0\mu_{max}$ be
   decomposed into the two components:
  \begin{equation*}
      w_0\mu_{max} = \gamma + \delta_{m+1}, \text{ where } \gamma \in
      \Phi(\Gamma), \text{ and } \gamma \geq 0.
  \end{equation*}
  Since, $w_0^2 = 1$, we have
  \begin{equation*}
      \mu_{max} - w_0\gamma = w_0\delta_{m+1}.
  \end{equation*}
  Thus, $\mu_{max} - w_0\gamma$ is a root in $\Phi(\widetilde\Gamma)$. We have  $\mu_{max} - w_0\gamma \geq \mu_{max}$,
  whereas $w_0\gamma \leq 0$.  Since $\mu_{max}$ is the
  maximal root in $\Phi(\widetilde\Gamma)$, we have $w_0\gamma = 0$. Thus, $\gamma = 0$, and
  $w_0\mu_{max} = \delta_{m+1}$, as required.
~\\

   (ii)  Since $w_0 = w^{\frac{h}{2}}$, then $w_0$ and $w$
 commute. Therefore, we have
 \begin{equation}
      w s_{\mu_{max}} \quad \stackrel{w_0}{\simeq} \quad
      (w_0w{w_0})(w_0{s}_{\mu_{max}}w_0) \quad = \quad
      ws_{w_0\mu_{max}} \quad = \quad  ws_{\delta_{m+1}}.
      \qed
  \end{equation}

\begin{figure}[h]
\centering
\includegraphics[scale=1.3]{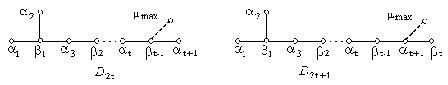}
\caption{Extending
   $D_{2t-1} < D_{2t}$ (resp. $D_{2t} < D_{2t+1}$) by adding $\alpha_{t+1}$ or $\mu_{max}$
   (resp. $\beta_t$ or  $\mu_{max}$)}
\label{Dl_extending}
\end{figure}

\begin{figure}[h]
\centering
\includegraphics[scale=0.6]{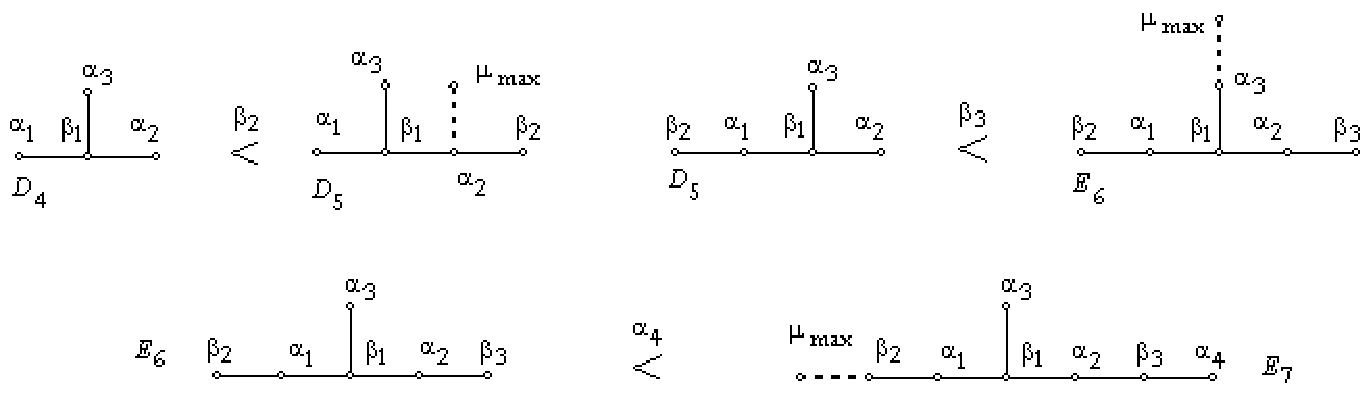}
 \caption{Extensions \quad $D_4 \stackrel{\beta_2}{<} D_5$,
                                  \quad $D_5 \stackrel{\beta_3}{<} E_6$ \quad and
                                  \quad $E_6 \stackrel{\alpha_4}{<} E_7$}
\label{D5_E6_E7_ext}
\end{figure}

\begin{remark}
   {\rm
   In Fig. \ref{Dl_extending}, we have:

    $\Gamma = D_{2t}$, $\widetilde{\Gamma} = D_{2t+1}$,  $\varphi_k = \beta_{t-1}$, $\delta_{h+1} = \alpha_{t+1}$;

    $\Gamma = D_{2t+1}$, $\widetilde{\Gamma} = D_{2t+2}$, $\varphi_k = \alpha_{t+1}$, $\delta_{h+1} =
    \beta_t$.
~\\

   In Fig. \ref{D5_E6_E7_ext}, see Table \ref{tab_orbit_max_roots}, we have:

    $\Gamma = D_4$, $\widetilde{\Gamma} = D_5$, $\varphi_k = \alpha_2$, $\delta_{h+1} = \beta_2$;

    $\Gamma = D_5$, $\widetilde{\Gamma} = E_6$, $\varphi_k = \alpha_2$, $\delta_{h+1} = \beta_3$;

    $\Gamma = E_6$, $\widetilde{\Gamma} = E_7$, $\varphi_k = \beta_3$, $\delta_{h+1} =  \alpha_4$;
    }
\end{remark}

\subsubsection{Passage from $D_{k+2}$ to $D_l(a_k)$}

\begin{figure}[h] \centering
\includegraphics[scale=1.8]{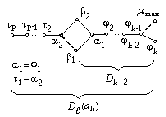}
\caption{$D_{k+2}$ and $D_l(a_k)$}
\label{fig_Dlak_and_Dk}
\end{figure}

We suppose that $D_{k+2}$ is the subdiagram of $D_l(a_k)$ as it is
depicted in Fig. \ref{fig_Dlak_and_Dk}.  For any object $O$ related
with $D_{k+2}$, we denote by $\widetilde{O}$ the corresponding
object related with $D_l(a_k)$. Let

 \begin{equation}
    S = \{\alpha_1,  \beta_1, \beta_2, \varphi_2, \dots, \varphi_{k}\} \\
 \end{equation}
 be a certain $D_{k+2}$\hspace{0.5mm}-associated root subset, let
 \begin{equation}
    \widetilde{S} = \{\alpha_1, \alpha_2, \beta_1, \beta_2, \varphi_2, \dots, \varphi_k, \tau_2, \dots, \tau_p \}
 \end{equation}
 be a $D_l(a_k)$-associated subset, $S \subset \widetilde{S}$, see Fig. \ref{fig_Dlak_and_Dk}.
 Let $V$ (resp. $\widetilde{V}$) be the space spanned by $S$ (resp. $\widetilde{S}$), $V \subset \widetilde{V}$.

\index{Coxeter number}

\begin{proposition}
  \label{prop_Dlak_and_Dk}
   Let $w$ be the $S$-associated element, $\widetilde{w}$ be the $\widetilde{S}$-associated element.
   Let $\mu_{max}$ be the maximal root in the root system $\varPhi(D_{k+2})$.
   Let $w_0$ be the longest element in $W(D_{k+2})$, i.e.,
   $w_0 = w^{\frac{h}{2}}$, where $h$ is the Coxeter number in $W(D_{k+2})$.

  {\rm (i)} For $\widetilde{w}_0 = \widetilde{w}^{\frac{h}{2}} \in W(D_l(a_k))$, we have
 \begin{equation}
   \label{eq_prop_Dlak_1}
     \widetilde{w}_0\mu_{max} = {w}_0\mu_{max} = \varphi_k.
 \end{equation}

  {\rm (ii)} The following conjugacy relation holds:
  \begin{equation}
    \label{eq_prop_Dlak_2}
     \widetilde{w}{s}_{\mu_{max}}  ~\simeq~ \widetilde{w}{s}_{\varphi_k}.
  \end{equation}
\end{proposition}

 \index{maximal root $\mu_{max}$}
 \index{$\mu_{max}$ (maximal root)}
\PerfProof (i) The maximal root $\mu_{max}$ is as follows:
 \begin{equation}
   \label{eq_max_elem_betas}
    \mu_{max} = \beta_1 + \beta_2 + 2\sum\limits_{i=1}^{k-1}\varphi_i + \varphi_k
 \end{equation}

 Note that $s_{\alpha_2}z = z$ for any $z$ having the same
 coordinates for components $\beta_1$ and  $\beta_2$.
 Indeed, if $z = b(\beta_1 + \beta_2) + \dots,$  see Fig. \ref{fig_Dlak_and_Dk}, then
 \begin{equation}
   \label{eq_alpha_touch}
       s_{\alpha_2} : z \longmapsto b(\beta_1 + \alpha_2) + b(\beta_2 - \alpha_2) +  \dots =
       b(\beta_1 + \beta_2) + \dots.
 \end{equation}
 Let $R$ (resp. $\widetilde{R}$) be the set of reflections entering the decomposition of $w$ (resp. $\widetilde{w}$),
 suppose that $\widetilde{R}\backslash{R} = \{s_{\alpha_2}, s_{\tau_2}, \dots, s_{\tau_p}\}$.
 Let $V' \subset V$ be the subspace of all vectors having the same
 coordinates for components $\beta_1$ and  $\beta_2$.
 In $\widetilde{R}\backslash{R}$, the only reflection
 that can affect $z \in V'$  is $s_{\alpha_2}$, so by \eqref{eq_alpha_touch} we have
 \begin{equation*}
  \begin{split}
    & \widetilde{w}z = wz \text{ for any } z \in V', \\
    & \widetilde{w}_0{z} = w_0{z} \text{ for any } z \in V', \\
  \end{split}
 \end{equation*}
 By \eqref{eq_max_elem_betas}, we have $\mu_{max} \in V'$. Therefore, by Proposition \ref{prop_conj_max_root}
 \begin{equation*}
   \widetilde{w}_0\mu_{max} = {w}_0\mu_{max} = \varphi_k,
 \end{equation*}
 (here, $\varphi_k$ plays the same role as $\delta_{k+1}$ in Proposition \ref{prop_conj_max_root}).
 So, \eqref{eq_prop_Dlak_1} holds.
~\\

   (ii)  Since $\widetilde{w}_0 = \widetilde{w}^{\frac{h}{2}}$, it follows that $\widetilde{w}_0$ and $\widetilde{w}$
 commute. Therefore, by \eqref{eq_prop_Dlak_1}, we have
 \begin{equation}
      \widetilde{w} s_{\mu_{max}} \quad \stackrel{\widetilde{w}_0}{\simeq} \quad
      (\widetilde{w}_0\widetilde{w}\widetilde{w}_0)(\widetilde{w}_0{s}_{\mu_{max}}\widetilde{w}_0)
        =
      \widetilde{w}s_{\widetilde{w}_0\mu_{max}}
         = \widetilde{w}s_{w_0\mu_{max}}  = \widetilde{w}s_{\varphi_k}.
  \end{equation}
 \qed

In the remaining part of \S\ref{sect_basic_pat} we consider root subsets forming the simplest diagrams:
Dipoles, triangles, squares, diamonds. For every type of these diagrams,
we describe properties helping to understand whether the given root subset is linearly independent
or not, see Lemmas \ref{lem_on_triangle_1}, \ref{lem_on_square_1}, \ref{lem_on_diamond_1}.
Further on, we consider a little more complicated diagrams obtained by gluing two simpler ones,
see \S\ref{sect_gluing_2diagr}.

\subsection{Dipoles and subsets of mutually orthogonal roots}
  \label{sec_mut_orth}

 Let $\Gamma$ be one of the diagram $D_4(a_1)$ or $D_4$.
 The pair of opposite vertices of any diagonal in $D_4(a_1)$
 or the pair of endpoints of $D_4$ is said to be a {\it dipole}.
 The pair of roots corresponding to vertices of the dipole are called a
 {\it dipole of roots}. The pair of points of the dipole are not connected by any edge.
 When no ambiguity arises we say \lq{\lq}dipole\rq\rq instead of \lq{\lq}dipole of roots\rq\rq.

\begin{figure}[h]
 \centering
\includegraphics[scale=0.8]{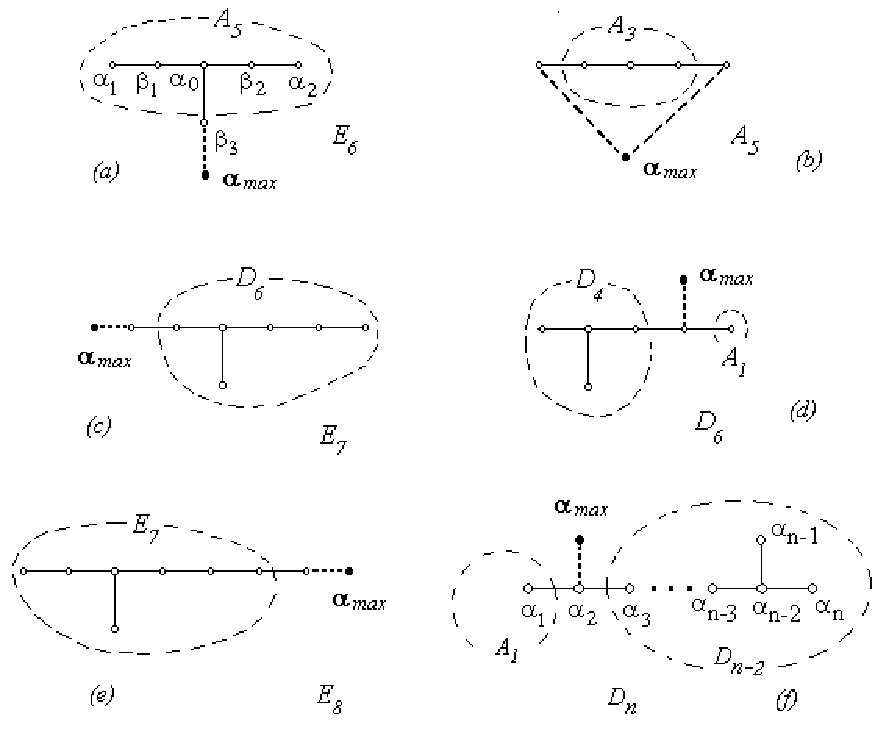}
\caption{Subsets of roots in $\varPhi(D_l)$, $\varPhi(E_l)$ orthogonal to $\alpha_{max}$}
\label{E6_E7_E8}
\end{figure}

 In the following lemma, we show that for $E_6$ (resp. $E_8$), any two
 subsets of $3$ orthogonal roots are equivalent under
 $W = W(E_6)$ (resp. $W = W(E_8))$. For $E_7$, the analogous statement does not hold
 because $E_7$ behaves differently:
 An \lq{\lq}unlucky\rq{\rq} location of the maximal root with respect to the Dynkin diagram,
 Fig. \ref{E6_E7_E8}. Recall that the location of the maximal root is the same as that of the additional vertex
 in the extended Dynkin diagram, \cite{Bo}.

\begin{lemma}
  \label{lem_mutual_ortog}
    Let $\alpha_{max}$ be the maximal root in the root system
    $\varPhi(D)$, where $D$ is  given by Table $\ref{tab_mutually_ortog}$, column $2$.
    The maximal root $\alpha_{max} \in \varPhi(D)$ is
    orthogonal to the root $\eta \in \varPhi(D)$ if and only if
    $\eta \in D'$, where $D' \subset D$ are given by Table $\ref{tab_mutually_ortog}$, column $4$.

\end{lemma}

 \index{maximal root $\mu_{max}$}
 \index{$\mu_{max}$ (maximal root)}
   \PerfProof 1) Consider line 1 of Table \ref{tab_mutually_ortog}.
     We need to prove that the maximal root $\alpha_{max}$ in $\varPhi(E_6)$ is
    orthogonal to the root $\eta \in \varPhi(E_6)$  if and only if
    $\eta \in \varPhi(A_5) = \{ \alpha_0, \alpha_1, \alpha_2, \beta_1, \beta_2\}$,  see Fig. \ref{E6_E7_E8}.
    The maximal root $\alpha_{max}$ in $\varPhi(E_6)$ is as follows:
 \begin{equation*}
 \begin{split}
  & \alpha_{max} = \alpha_1 + 2\beta_1 + 3\alpha_0 + 2\beta_2 +
  \alpha_2 + 2\beta_3, \\
  & \alpha_{max} \perp \alpha_i, \text{ for } i = 0,1,2; \quad \alpha_{max} \perp \beta_1,
  \beta_2, \text{ and } \alpha_{max} \notperp \beta_3.
 \end{split}
\end{equation*}
 Therefore, $\alpha_{max}$ is orthogonal to $\varPhi(A_5)$
 spanned by $\{ \alpha_0, \alpha_1, \alpha_2, \beta_1, \beta_2 \}$.
 Any root $\eta$ from $\varPhi(E_6)$ has the form $\eta = k_1z +
 k_2\beta_3$, where $z \in \varPhi(A_5)$. Since $\alpha_{max} \perp
 z$ and $\alpha_{max} \notperp \beta_3$, then $\alpha_{max} \perp
 \eta$ means that $k_2 = 0$, and $\eta \in \varPhi(A_5)$. If $\alpha_{max} \in \varPhi(A_5)$,
 then similar arguments show that $\eta \in \varPhi(A_5)$,
 and $\eta \perp \alpha_{max}$ if and only if $\eta \in A_3$, see Fig. \ref{E6_E7_E8}.
 The remaining cases 2) -- 6) from Table \ref{tab_mutually_ortog} are similarly considered. \qed

\begin{table}[h]
  \centering
  \renewcommand{\arraystretch}{1.8}
  \begin{tabular} {|c|c|c|c|c|}
  \hline
        & The root   & The maximal root $\alpha_{max}$
                     & The root $\eta \in \varPhi(D)$, $\eta \perp \alpha_{max}$
                     & Case in  \cr
        & system $D$ & in the root system $\varPhi(D)$
                     & if and only if $\eta \in  \varPhi(D')$, $D' \subset D$
                     & Fig. \ref{E6_E7_E8}\\
  \hline
     1 & $E_6$ & $\alpha_{max} \in \varPhi(E_6)$ &  $\eta \in \varPhi(A_5)$ & $(a)$ \\
  \hline
     2 & $A_5$ & $\alpha_{max} \in \varPhi(A_5)$ &  $\eta \in \varPhi(A_3) $ & $(b)$  \\
  \hline
     3 & $E_7$ & $\alpha_{max} \in \varPhi(E_7)$ &  $\eta \in \varPhi(D_6)$ & $(c)$ \\
  \hline
     4 & $D_6$ & $\alpha_{max} \in \varPhi(D_6)$ &  $\eta \in \varPhi(D_4) \oplus \varPhi(A_1)$ & $(d)$ \\
  \hline
     5 & $E_8$ & $\alpha_{max} \in \varPhi(E_8)$ &  $\eta \in \varPhi(E_7)$ & $(e)$ \\
  \hline
     6 & $D_l$ & $\alpha_{max} \in \varPhi(D_l)$ &  $\eta \in \varPhi(D_{l-2}) \oplus \varPhi(A_1)$ & $(f)$  \\
  \hline
\end{tabular}
  \vspace{2mm}
  \caption{Roots $\eta \in \varPhi(D)$ such that $\eta \perp \alpha_{max}$}
  \label{tab_mutually_ortog}
\end{table}

\begin{corollary}
  \label{cor_orth_roots}
  Any two dipoles in $\varPhi(D)$, where $D = E_6$, $E_7$ or $E_8$, are equivalent under $W(D)$.

\end{corollary}

 \PerfProof Let $\{\alpha_1, \alpha_2\}$ be any dipole of roots in $\varPhi(D)$, where $D = E_6, E_7, E_8$.
 We send $\alpha_1$ into $\alpha_{max} \in \varPhi(D)$ by some $u \in W(D)$.
 By Lemma \ref{lem_mutual_ortog}, for $D = E_6$ (resp. $E_7$, $E_8$), the element $u$ sends $\alpha_2$ to
 $\varPhi(D')$, where $D' = A_5$ (resp. $\varPhi(D_6)$, $\varPhi(A_5)$). By means of $W(D')$
 the root $\alpha_2$  can be sent to any other root in $\varPhi(D')$ so that the image of $\alpha_1$
 (i.e., $\alpha_{max}$) is not moved. \qed

 There are two dipoles of roots in $\varPhi(D_l)$
 which are not equivalent under $W(D_l)$, see  Fig. \ref{E6_E7_E8} and Table \ref{tab_mutually_ortog}.
 This fact is proved in Lemma \ref{lem_eq_2diag}, see also Example \ref{sec_example_4cycles}.

\begin{corollary} {\rm (\cite[Lemma 11,(i), Lemma 27]{Ca72})}
  \label{lem_mutual_ortog_2}
    {\rm (i)} Any two sets of $3$ mutually orthogonal roots in $\varPhi(E_6)$ are equivalent under $W(E_6)$.

    {\rm (ii)} Any two sets of $3$ mutually orthogonal roots in $\varPhi(E_8)$
    are equivalent under $W(E_8)$.

    {\rm (iii)} There are two sets of $3$ orthogonal roots in $\varPhi(E_7)$
     which are not equivalent under $W(E_7)$.
\end{corollary}

 \PerfProof We will show that any triple of mutually orthogonal roots in $E_6$ can
 be transformed into the triple $\{\alpha_{max}(E_6), \alpha_{max}(A_5), \alpha_{max}(A_3) \}$,
 see Fig. \ref{E6_E7_E8}, where $\alpha_{max}(D)$ means
 the maximal root in $\varPhi(D)$.
 Then any triple of mutually orthogonal roots in $E_6$
 can be transformed into each other.
 Let $\{ \varphi_1, \varphi_2, \varphi_3 \}$ be the triple of
 mutually orthogonal roots in $E_6$.  The root $\varphi_1$ can be transformed
 into any root in $W(E_6)$. We transform $\varphi_1$ into the maximal
 root $\alpha_{max}(E_6)$. By Lemma \ref{lem_mutual_ortog},  roots $\varphi_2$ and
 $\varphi_3$ are transformed into two elements in $\varPhi(A_5)$
 under $W(A_5)$. We transform $\varphi_2$ into $\alpha_{max}(A_5)$,
 see Fig. \ref{E6_E7_E8}. Again, by Lemma \ref{lem_mutual_ortog}, $\varphi_3$
 is transformed under $W(A_3)$ into the maximal root in
 $\varPhi(A_3)$. The triple of maximal roots $\{\alpha_{max}(E_6), \alpha_{max}(A_5), \alpha_{max}(A_3)\}$
 corresponds to the following chain of root subsets:
$$
   \varPhi(E_6) \supset \varPhi(A_5) \supset \varPhi(A_3).
$$
Similarly, in the case $E_8$, the triple of maximal roots
$\{\alpha_{max}(E_8), \alpha_{max}(E_7), \alpha_{max}(D_6)\}$, see Fig. \ref{E6_E7_E8}$(e)$,$(c)$ and $(d)$,
corresponds to the following chain:
$$
   \varPhi(E_8) \supset \varPhi(E_7) \supset \varPhi(D_6).
$$
In the case $E_7$, the subset associated with the third maximal root
is split up into $2$ non-connected subsets:
 \begin{equation}
  \label{split_in_E7}
   \varPhi(E_7) \supset \varPhi(D_6) \supset \varPhi(D_4) \oplus \varPhi(A_1).
 \end{equation}
The decomposition $\varPhi(D_4) \oplus \varPhi(A_1)$ in
\eqref{split_in_E7} is responsible for the presence of two
non-equivalent root subsets with $3$ mutually orthogonal roots. \qed

\subsection{Triangles, squares and diamonds}
  \label{sect_triangles}
\begin{figure}[h]
\centering
\includegraphics[scale=1.1]{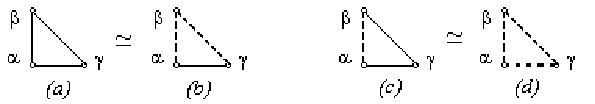}
\caption{Linearly dependent triples of roots
 $(a)$, $(b)$; linearly independent triples of roots $(c)$,$(d)$.}
\label{3cycle_1}
\end{figure}

\begin{lemma}
 \label{lem_on_triangle_1}
  Let $\Gamma$ be a $3$-cycle,  the triple of roots $S = \{\alpha, \beta, \gamma\}$ be
  a certain $\Gamma$-associated subset. The triple $S$ is linearly independent if and only if
  the number of dotted edges of $\Gamma$ is odd, see
  Fig. $\ref{3cycle_1}$$(c)$,$(d)$.

  If all edges of $\Gamma$ are solid, see Fig. $\ref{3cycle_1}$$(a)$, then
  \begin{equation}
    \label{eq_sum_0_pair}
       \alpha + \beta + \gamma = 0.
  \end{equation}

  If only one edge of $\Gamma$ is solid, for example, $\{\alpha, \gamma\}$ in Fig. $\ref{3cycle_1}$$(b)$, then
  \begin{equation}
    \label{eq_sum_0_pair_2}
       \alpha - \beta + \gamma = 0.
  \end{equation}
\end{lemma}

\begin{figure}[h]
\centering
\includegraphics[scale=1.1]{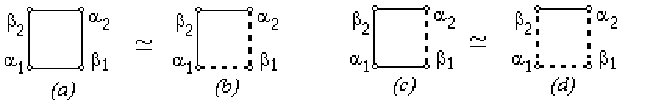}
\caption{Linearly dependent quadruples of roots
 $(a)$, $(b)$; linearly independent quadruples of roots $(c)$,$(d)$.}
\label{4cycle_1}
\end{figure}

\begin{remark}
  {\rm Note that the change of sign of any root in $S$ does not affect the linear dependence
  and does not change $\Gamma$-associated elements since $s_\alpha = s_{-\alpha}$.
  In Fig. \ref{3cycle_1}, the case $(d)$ turns into $(c)$  under the change $\gamma \longmapsto -\gamma$,
  the case $(b)$ turns into $(a)$ under the change $\beta \longmapsto -\beta$.
  In Fig. \ref{4cycle_1}, the case $(d)$ turns into $(c)$
  under the change $\alpha_1 \longmapsto -\alpha_1$, the case $(b)$ turns into $(a)$ under the change
  $\beta_1 \longmapsto -\beta_1$.}
\end{remark}

 {\it Proof of Lemma \ref{lem_on_triangle_1}.}
  If $S$ is linearly independent, then there exist $1$ or $3$
  dotted edges, see $(c)$ and $(d)$ in Fig. \ref{3cycle_1}, otherwise we have the diagram $\widetilde{A}_3$,
  contradicting Lemma \ref{lem_must_dotted}. Let $S$ be linearly
  dependent and let there be only one dotted edge as in configuration $(c)$.
  Let $\alpha = b\beta + c\gamma$. For $\Gamma_0 = \{\beta, \gamma\}$,
  by \eqref{eq_lin_depend2} we have

\begin{equation}
   \label{eq_B_lindep_triple_0}
   \Small
   \begin{split}
     B_{\Gamma_0} =
     \frac{1}{2}\left [
     \begin{array}{cc}
       2 &  -1  \\
      -1 & 2 \\
    \end{array}
    \right ],
   & \quad
     B^{-1}_{\Gamma_0} =
     \frac{2}{3}\left [
     \begin{array}{ccc}
       2 & 1  \\
       1 & 2  \\
    \end{array}
    \right ],
    \quad
   \left [
    \begin{array}{c}
     (\alpha, \beta)   \\
     (\alpha, \gamma)   \\
    \end{array}
    \right ] =
   \left [
     \begin{array}{c}
         \frac{1}{2} \\
        -\frac{1}{2} \\
         0
     \end{array}
    \right ]
    \Longrightarrow \quad
     \begin{cases}
        b = -\frac{1}{3}, \\
        c = \frac{1}{3},
     \end{cases}
   \\
   & \alpha = \frac{1}{3}(\gamma - \beta), \qquad
   \mathscr{B}(\alpha) = \frac{1}{9},
   \end{split}
  \end{equation}
  i.e., $\alpha$ is not a root. Thus, for $S$ linearly dependent,
  only cases $(a)$ and $(b)$ are possible.
  Let $\gamma = a\alpha + b\beta$.  Again, by \eqref{eq_lin_depend2} we have
\begin{equation}
   \label{eq_B_lindep_triple}
   \Small
   \begin{split}
     B_{\Gamma} =
     \frac{1}{2}\left [
     \begin{array}{cc}
       2 &  -1  \\
      -1 & 2 \\
    \end{array}
    \right ],
    & \quad
     B^{-1}_{\Gamma} =
     \frac{2}{3}\left [
     \begin{array}{ccc}
       2 & 1  \\
       1 & 2  \\
    \end{array}
    \right ],
     \quad
   \left [
    \begin{array}{c}
     (\gamma, \alpha)   \\
     (\gamma, \beta)   \\
    \end{array}
    \right ] =
   \left [
     \begin{array}{c}
        -\frac{1}{2} \\
        -\frac{1}{2} \\
        0
     \end{array}
    \right ]
    \Longrightarrow \quad
     \begin{cases}
        a = -1, \\
        b = -1,
     \end{cases} \\
     & \gamma = -\alpha - \beta.
   \end{split}
  \end{equation}
~\\
  The root $-\alpha - \beta$ is the minimal root for
  $A_2 = \{ \alpha, \beta \}$ in accordance with Remark \ref{rem_max_root}(i).
  Eq. \eqref{eq_sum_0_pair} is proved. Eq. \eqref{eq_sum_0_pair_2} is obtained from eq. \eqref{eq_sum_0_pair}
  by replacing $\beta \longmapsto -\beta$.
  \qed

\begin{lemma}
 \label{lem_on_square_1}

  Let $\Gamma$ be a $4$-cycle, the quadruple of roots
  $S = \{\alpha_1, \beta_1, \alpha_2, \beta_2 \}$ be the $\Gamma$-associated subset.
  The quadruple $S$ is linearly independent if and only if
  the number of dotted edges of $\Gamma$ is odd, see
  Fig. $\ref{4cycle_1}$$(c)$,$(d)$.

 If all edges of $\Gamma$ are solid, see Fig. $\ref{4cycle_1}$$(a)$, then:
  \begin{equation}
    \label{eq_sum_0}
       \alpha_1 + \alpha_2 + \beta_1 + \beta_2 = 0.
  \end{equation}

 If only two edges are dotted, for example,
 $\{\alpha_1, \beta_1 \}$ and $\{\alpha_2, \beta_1 \}$ in Fig. $\ref{4cycle_1}$$(b)$,
 then:
  \begin{equation}
    \label{eq_sum_0_2}
       \alpha_1 + \alpha_2 - \beta_1 + \beta_2 = 0.
  \end{equation}
\end{lemma}

 \PerfProof
  If $S$ is linearly independent, then there exist $1$ or $3$
  dotted edges, see Fig. \ref{4cycle_1}$(c)$ and $(d)$,
  otherwise we have the diagram $\widetilde{A}_4$, contradicting Lemma \ref{lem_must_dotted}.
  Let $S$ be linearly dependent and let there be only one dotted edge as in configuration $(c)$.
  Let $\beta_1 = a_1\alpha_1 + a_2\alpha_2 + b_2\beta_2$. For $\Gamma_0 = \{\alpha_1, \beta_2, \alpha_2\}$,
  by \eqref{eq_lin_depend2}, we have
\begin{equation}
   \label{eq_Gamma0_1}
   \Small
   \begin{split}
     B_{\Gamma_0} =
     \frac{1}{2}\left [
     \begin{array}{ccc}
       2 & 0 &  -1  \\
       0 & 2 &  -1  \\
      -1 & -1  & 2
    \end{array}
    \right ],
    & \quad
     B^{-1}_{\Gamma_0} =
     \frac{1}{2}\left [
     \begin{array}{cccc}
       3 & 1 & 2  \\
       1 & 3 & 2  \\
       2 & 2 &  4
    \end{array}
    \right ],
     \quad
   \left [
    \begin{array}{c}
     (\beta_1, \alpha_1)   \\
     (\beta_1, \alpha_2)   \\
     (\beta_1, \beta_2)
    \end{array}
    \right ] =
   \left [
     \begin{array}{c}
        -\frac{1}{2} \\
         \frac{1}{2} \\
        0
     \end{array}
    \right ]
    \Longrightarrow \quad
     \begin{cases}
        a_1 = -\frac{1}{2}, \\
        a_2 = \frac{1}{2}, \\
        b_1 = 0,
     \end{cases} \\
      & \beta_1 = \frac{\alpha_2 - \alpha_1}{2}, \qquad
         \mathscr{B}(\beta_1) = \frac{1}{4},
   \end{split}
  \end{equation}
  i.e., $\beta_1$ is not a root.
  Thus, for $S$ linearly dependent,
  only cases $(a)$ and $(b)$ are possible.
  Let $\beta_2 = a_1\alpha_1 + a_2\alpha_2 +  b_1\beta_1$.  Again, by \eqref{eq_lin_depend2} we have

\begin{equation}
   \label{eq_Gamma0_2}
   \Small
   \begin{split}
     B_{\Gamma} =
     \frac{1}{2}\left [
     \begin{array}{ccc}
       2 & 0 &  -1  \\
       0 & 2 &  -1  \\
      -1 & -1  & 2
    \end{array}
    \right ],
   & \quad
     B^{-1}_{\Gamma} =
     \frac{1}{2}\left [
     \begin{array}{cccc}
       3 & 1 & 2  \\
       1 & 3 & 2  \\
       2 & 2 &  4
    \end{array}
    \right ],
    \quad
   \left [
    \begin{array}{c}
     (\beta_2, \alpha_1)   \\
     (\beta_2, \alpha_2)   \\
     (\beta_2, \beta_1)
    \end{array}
    \right ] =
   \left [
     \begin{array}{c}
        -\frac{1}{2} \\
        -\frac{1}{2} \\
        0
     \end{array}
    \right ]
    \Longrightarrow \quad
     \begin{cases}
        a_1 = -1, \\
        a_2 = -1, \\
        b_1 = -1,
     \end{cases} \\
      & \beta_2 = -\alpha_1 -\alpha_2 -\beta_1.
   \end{split}
  \end{equation}
    The root $\beta_2$ is the minimal root for $A_3$ in accordance with Remark \ref{rem_max_root}(i).
    Eq. \eqref{eq_sum_0} is proved, eq. \eqref{eq_sum_0_2} is obtained from eq. \eqref{eq_sum_0}
    by replacing $\beta_1 \longmapsto -\beta_1$. \qed
~\\

\begin{figure}[h]
\centering
\includegraphics[scale=2.2]{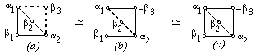}
\caption{The impossible diamonds of roots}
\label{impos_quint}
\end{figure}

 A square with a diagonal consisting of two edges,
 as in Fig. \ref{impos_quint}, is said to be a {\it diamond}.

\index{diamond}
\begin{lemma}
 \label{lem_on_diamond_1}
  Let $\Gamma = D_4$, and $S = \{ \alpha_1, \beta_1, \beta_2, \beta_3 \}$ form a linearly independent
  quadruple of roots corresponding to $\Gamma$.
  There exists no root $\alpha_2$ (even linearly dependent on $S$) completing $S$ to a diamond
  such that every $4$-cycle of the diamond contains an even number of dotted edges
  (such diamonds are depicted in Fig. $\ref{impos_quint}$).
\end{lemma}

  \PerfProof
   Suppose such a root $\alpha_2$ exists. For example, for the case $(a)$ in Fig. \ref{impos_quint},
   by \eqref{eq_sum_0_2}  we have
      \begin{equation}
        \label{eq_impos_system}
        \begin{split}
          & \alpha_1 + \beta_1 + \alpha_2 - \beta_3 = 0, \\
          & \alpha_1 + \beta_2 + \alpha_2 - \beta_3 = 0.
        \end{split}
      \end{equation}
    We derive from \eqref{eq_impos_system} that $\beta_1 = \beta_2$, which is impossible.
    The remaining cases are considered in the same way.
 \qed

\subsection{Gluing two diagrams}
 \label{sect_gluing_2diagr}
 \index{triangle of roots}

\begin{figure}[h]
\centering
\includegraphics[scale=0.64]{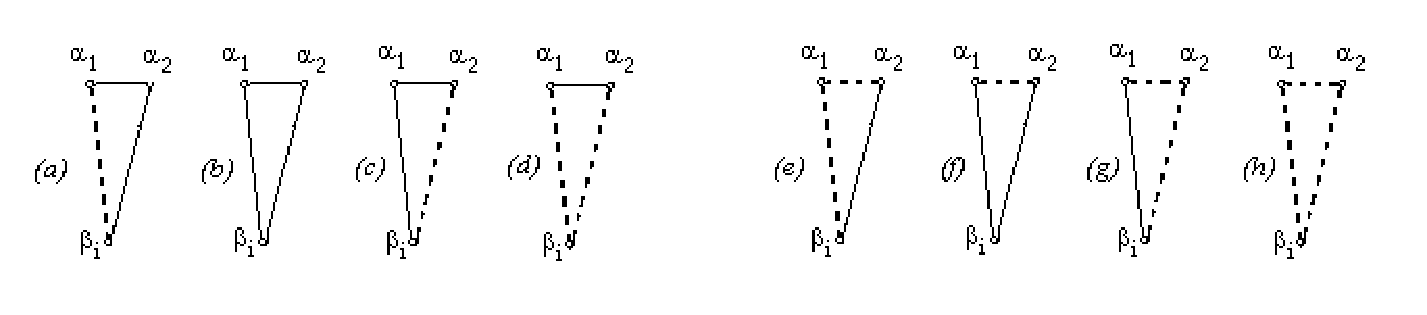}
\caption{The eight possible triangles each edge of which is either solid or dotted.}
\label{all_triangles}
\end{figure}

  Consider all possible triangles, each edge of which is
  either solid or dotted. There are $8$ such triangles, four of which, by Lemma \ref{lem_on_triangle_1},
  constitute linearly independent triples: Cases $(a)$, $(c)$, $(f)$ and $(h)$  in
  Fig. \ref{all_triangles}, and the four remaining triangles constitute linearly
  dependent triples.
 \begin{lemma}
   \label{lem_gluing_1}
    Let $\alpha_1$ and $\alpha_2$ be linearly independent roots.
    For triangles $(a)-(d)$ in Fig. $\ref{all_triangles}$,  we have $(\alpha_1, \alpha_2) = -\frac{1}{2}$,
    then $\alpha_1 + \alpha_2$ is a root and
 \begin{equation}
   \label{eq_comm_plus}
    \begin{split}
      & s_{\alpha_1 + \alpha_2}s_{\beta_i} = s_{\beta_i}s_{\alpha_1 + \alpha_2}, \\
      & s_{\alpha_1 + \alpha_2}s_{\alpha_1} = s_{\alpha_2}s_{\alpha_1 + \alpha_2},
    \end{split}
 \end{equation}
    For triangles $(e)-(h)$ in Fig. $\ref{all_triangles}$, we have $(\alpha_1, \alpha_2) = \frac{1}{2}$,
     then $\alpha_1 - \alpha_2$ is a root and
 \begin{equation}
    \label{eq_comm_minus}
    \begin{split}
      & s_{\alpha_1 - \alpha_2}s_{\beta_i} = s_{\beta_i}s_{\alpha_1 - \alpha_2}, \\
      & s_{\alpha_1 - \alpha_2}s_{\alpha_1} = s_{\alpha_2}s_{\alpha_1 - \alpha_2}.
    \end{split}
 \end{equation}
 \end{lemma}

  \PerfProof
    The second relations in  \eqref{eq_comm_plus} and
    \eqref{eq_comm_minus} follow from the relation
    $s_{\alpha + \beta}s_{\beta} = s_{\beta}s_{\alpha}$.
    Consider the first relations in  \eqref{eq_comm_plus} and
    \eqref{eq_comm_minus}.
    For cases $(a),(c)$, we have $(\alpha_1 + \alpha_2, \beta_i) = \frac{1}{2} - \frac{1}{2} = 0$,
    i.e., $s_{\alpha_1 + \alpha_2}$ and $s_{\beta_i}$ commute.
    For case $(b)$ (resp. $(d)$), by Lemma \ref{lem_on_triangle_1},
    we have $\alpha_1 + \alpha_2 + \beta_i = 0$
    (resp. $\alpha_1 + \alpha_2 - \beta_i = 0$), then $s_{\alpha_1 + \alpha_2} =
    s_{\beta_i}$. Thus, eq. \eqref{eq_comm_plus} holds.
    For cases $(f)$ and $(h)$, we have $(\alpha_1 - \alpha_2, \beta_i) = \frac{1}{2} - \frac{1}{2} = 0$,
    i.e., $s_{\alpha_1 - \alpha_2}$ and $s_{\beta_i}$ commute.
    Finally, for case $(e)$ (resp. $(g)$), we have $\alpha_2 - \alpha_1  + \beta_i = 0$
    (resp. $-\alpha_2 + \alpha_1  + \beta_i = 0$), i.e., $s_{\alpha_1 - \alpha_2} = s_{\beta_i}$.
    Eq. \eqref{eq_comm_minus} is proved.
    \qed

\begin{figure}[h]
\centering
\includegraphics[scale=0.64]{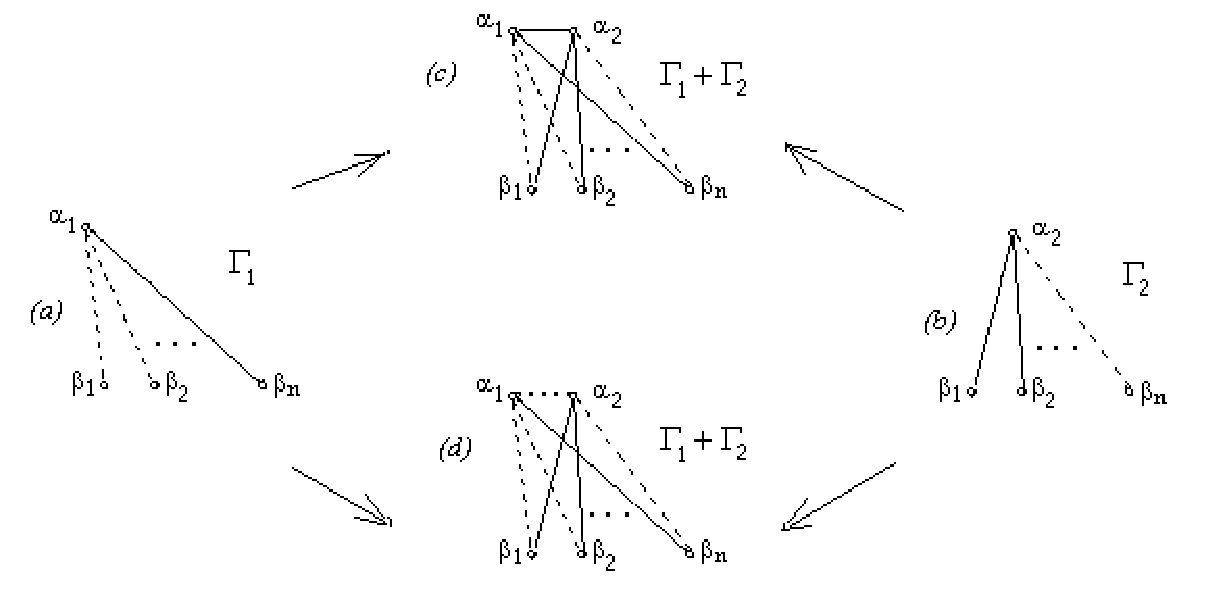}
\caption{For $(\alpha_1, \alpha_2) = -\frac{1}{2}$
(resp. $(\alpha_1, \alpha_2) = \frac{1}{2}$),
 the diagram $\Gamma_1 + \Gamma_2$ is depicted in $(c)$ (resp.
 $(d)$).}
\label{connect_alp12_3}
\end{figure}

 \begin{corollary}
  \label{cor_map_2_diagr}
   Let $\Gamma_1 = \{\alpha_1,\beta_1,\beta_2,\dots,\beta_n\}$
   (resp. $\Gamma_2 = \{\alpha_2,\beta_1,\beta_2,\dots,\beta_n\}$)
   be the star diagram with the center vertex $\alpha_1$ (resp. $\alpha_2$),
   let vertices $\{\beta_1,\beta_2,\dots,\beta_n\}$
   be common for $\Gamma_1$ and $\Gamma_2$, see Fig. $\ref{connect_alp12_3}$$(a)$,$(b)$.
   Let $\alpha_1$ and $\alpha_2$ be  connected by edge $\{\alpha_1, \alpha_2\}$,
   see Fig. $\ref{connect_alp12_3}$$(c)$,$(d)$.
   Let $w_1$ (resp. $w_2$) be $\Gamma_1$-associated (resp. $\Gamma_2$-associated) elements:
 \begin{equation}
      w_1 = s_{\alpha_1}\prod\limits_{i=1}^n{s_{\beta_i}}, \quad
      w_2 = s_{\alpha_2}\prod\limits_{i=1}^n{s_{\beta_i}}, \\
 \end{equation}
  If $(\alpha_1, \alpha_2) = -\frac{1}{2}$ (resp.  $\frac{1}{2}$),
  then $s_{\alpha_1 + \alpha_2}$ (resp. $s_{\alpha_1 - \alpha_2}$)  maps
  $\Gamma_1$ onto $\Gamma_2$. The elements $w_1$ and $w_2$ are conjugate,
 \begin{equation}
  \label{eq_corr_rel_1}
  \begin{split}
     & s_{\alpha_1 + \alpha_2} w_1 s_{\alpha_1 + \alpha_2} = w_2
     \quad \text{ for } \quad (\alpha_1, \alpha_2) = -\frac{1}{2}, \\
     & s_{\alpha_1 - \alpha_2} w_1 s_{\alpha_1 - \alpha_2} = w_2
      \quad \text{ for } \quad (\alpha_1, \alpha_2) = \frac{1}{2}, \\
  \end{split}
 \end{equation}
   The conjugacy of $w_1$ and $w_2$ is preserved also for diagrams $\Gamma_1$ and $\Gamma_2$
   containing other vertices $\tau_j$ not connected with $\alpha_1$ and $\alpha_2$.
 \end{corollary}

\begin{remark}
\label{rem_corr_rel}
{\rm
  The reflections $s_{\alpha_1 + \alpha_2}$ and  $s_{\alpha_1 - \alpha_2}$
  from \eqref{eq_corr_rel_1} behave like a map correcting the set $\{\alpha_1, \beta_1, \beta_2, \dots, \beta_n\}$ to
  the set $\{\alpha_2, \beta_1, \beta_2, \dots, \beta_n\}$.
  This is the reason to call the reflection $s_{\alpha_1 + \alpha_2}$ (resp. $s_{\alpha_1 - \alpha_2}$)
  a {\it corrective reflection}.
  In the context of the paper, the most frequently arising configuration of vertices $\Gamma_1 + \Gamma_2$
  from Corollary \ref{cor_map_2_diagr} and Fig. \ref{connect_alp12_3}
  is the configuration  with $n = 2$, see Fig. \ref{corrective_refl}$(a)$,$(c)$.
  In this case, the configuration $\Gamma_1 + \Gamma_2$ coincides with the $4$-cycle having one diagonal,
  see Fig. \ref{corrective_refl}$(b)$,$(d)$.
  \begin{figure}[h]
    \centering
    \includegraphics[scale=0.8]{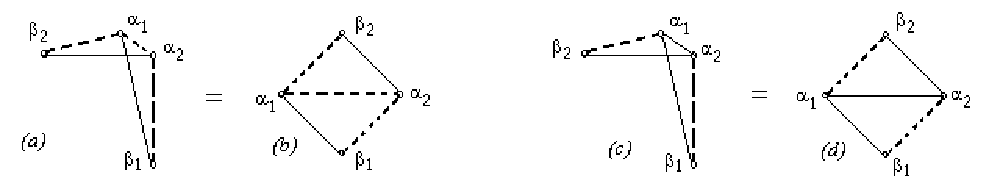}
    \caption{Corrective reflections $s_{\alpha_1 + \alpha_2}$ and $s_{\alpha_1 - \alpha_2}$}
    \label{corrective_refl}
  \end{figure}
 }
\end{remark}

\subsubsection{Gluing two $D_5(a_1)$-associated subsets}

\begin{figure}[h]
\centering
\includegraphics[scale=0.9]{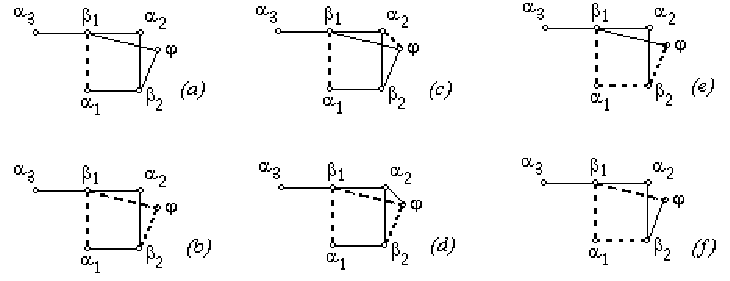}
\caption{Two $D_5(a_1)$-associated elements differing in one vertex}
\label{D5a1_with_fi_1}
\end{figure}

 \begin{lemma}[On necessarily connected roots]
   \label{lem_glue_1}
   Let $\Gamma = D_5(a_1)$; let
   $S_1 = \{\alpha_1, \alpha_2, \alpha_3, \beta_1, \beta_2\}$ and
   $S_2 = \{\alpha_1, \varphi, \alpha_3, \beta_1, \beta_2 \}$ be two
   $\Gamma$-associated subsets, the vectors of each of which being linearly
   independent, let $\{\alpha_1, \alpha_2, \beta_1, \beta_2\}$
   and $\{\alpha_1, \varphi, \beta_1, \beta_2\}$ be
   $D_4(a_1)$-associated subsets, let the root $\alpha_3$ be connected only with
   $\beta_1$, see Fig. $\ref{D5a1_with_fi_1}$.

   {\rm (i)} Configurations of Fig. $\ref{D5a1_with_fi_1}$$(a)$,$(b)$,$(e)$,$(f)$ are impossible:
   Roots $\varphi$ and $\alpha_2$ are necessarily connected, see Fig. $\ref{D5a1_with_fi_1}$$(c)$,$(d)$.

   {\rm (ii)} Let $w_1$ (resp. $w_2$) be $S_1$-associated (resp. $S_2$-associated). Then $w_1 \simeq w_2$.
 \end{lemma}

   \PerfProof (i) Suppose $\alpha_2$ and $\varphi$ are disconnected as in Fig.
   \ref{D5a1_with_fi_1}$(a),(b)$. Since $Q_0 = \{\alpha_1, \alpha_2, \beta_1, \beta_2\}$
   is $D_4(a_1)$-associated, we can assume that only the edge $\{\alpha_1, \beta_1\}$ is dotted
   in $Q_0$, otherwise we transform $Q_0$ to this case.
   Since the vectors of the quadruple $Q_1 = \{\alpha_1, \beta_1, \varphi, \beta_2 \}$ are
   linearly independent, it follows that the number of dotted edges in $Q_1$
   is either $1$ or $3$, see Fig. \ref{D5a1_with_fi_1}$(a)$,$(b)$,$(e)$,$(f)$.
   Since $Q_0$ contains an odd number of dotted edges, then the cases $(e)$ and $(f)$ are impossible.

   For cases $(a)$ and $(b)$, the quadruple $Q_2 = \{\alpha_2, \beta_1, \varphi,  \beta_2 \}$
   contains an even number of dotted edges. By Lemma \ref{lem_on_square_1}, we have
 \begin{equation}
  \begin{split}
    & \beta_1 + \beta_2 + \alpha_2 + \varphi = 0 \text{ for case } (a), \\
    & \beta_1 + \beta_2 + \alpha_2 - \varphi = 0 \text{ for case } (b),
  \end{split}
 \end{equation}
 i.e., $\varphi = \pm(\beta_1 + \beta_2 + \alpha_2)$. Thus,
 $(\varphi, \alpha_3) = (\beta_1, \alpha_3) \neq 0$, which leads to a contradiction.

  (ii) Since $\alpha_2$ and $\varphi$ are connected, see Fig. \ref{D5a1_with_fi_1}$(c)$,$(d)$,
  we can use Corollary \ref{cor_map_2_diagr}. By this corollary, the elements
  $w_1$ and $w_2$ are conjugate by means of the corrective reflection
  $s_{\varphi - \alpha_2}$ or $s_{\varphi + \alpha_2}$, see Remark \ref{rem_corr_rel}.
 \qed
~\\

\index{diamond}

\begin{figure}[h]
\centering
\includegraphics[scale=1.7]{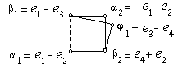}
\caption{A diamond in $D_5$. Here $\alpha_2$ and
$\varphi$ are not connected}
\label{ex_not_conn_vert}
\end{figure}

\begin{remark}{\rm
 If the edge $\{\alpha_3, \beta_1\}$ in Fig. \ref{D5a1_with_fi_1} is missing, then
 vertices $\alpha_2$ and $\varphi$ are not necessarily connected,
 see the diamond\footnotemark[1] $\{\alpha_1, \alpha_2, \beta_1, \beta_2, \varphi\}$
 in Fig. \ref{ex_not_conn_vert}.
 }
\end{remark}

 \footnotetext[1]{For definition and properties of diamonds, see \S\ref{sect_triangles}}

\subsubsection{Gluing two $E_6(a_1)$-associated subsets}

\begin{figure}[h]
\centering
\includegraphics[scale=0.9]{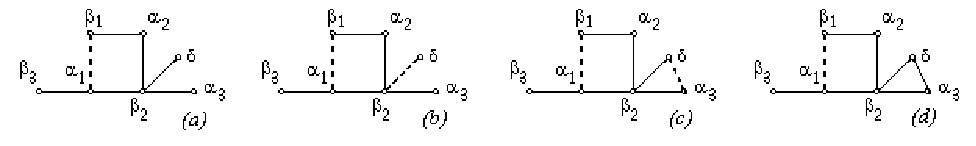}
\caption{Two $E_6(a_1)$-associated elements differing
 in one point}
\label{E6a1_with_delta}
\end{figure}

 \begin{lemma}
   \label{lem_glue_2}
   Let $\Gamma = E_6(a_1)$, let
   $S_1 = \{\alpha_1, \alpha_2, \alpha_3, \beta_1, \beta_2, \beta_3\}$ and
   $S_2 = \{\alpha_1, \alpha_2, \delta, \beta_1, \beta_2, \beta_3 \}$ be two
   $\Gamma$-associated subsets, each of which is linearly
   independent, let $\mathsf{Q}_5 = \{\alpha_1, \alpha_2, \beta_1, \beta_2, \beta_3\}$ be the $D_5(a_1)$-associated
   subset, $\alpha_3$ (resp. $\delta$) connected with $\mathsf{Q}_5$ only
   through the vertex $\beta_2$, see Fig. $\ref{E6a1_with_delta}$.

   {\rm (i)} Configurations of Fig. $\ref{E6a1_with_delta}$$(a)$,$(b)$ are  impossible --
   roots $\delta$ and $\alpha_3$ are necessarily connected, see Fig. $\ref{E6a1_with_delta}$$(c)$ or $(d)$.

   {\rm (ii)} Let $w_1$ (resp. $w_2$) be $S_1$-associated (resp. $S_2$-associated). Then $w_1 \simeq w_2$.
 \end{lemma}

\index{quintuple of roots}
\index{linearly dependent root}

  \PerfProof  (i) Consider quintuple $Q = \{\alpha_1, \alpha_2, \alpha_3, \beta_2, \delta\}$.
  If the vectors of $Q$ are linearly independent, then we get the root subsystem $\widetilde{D}_4$,
  which is impossible. If the vectors of $Q$ are linearly dependent, then $\delta$ is
  the  minimal (resp. maximal) root, i.e., $\delta = -(2\beta_2 + \alpha_1 + \alpha_2 + \alpha_3)$
  (resp. $\delta = 2\beta_2 + \alpha_1 + \alpha_2 + \alpha_3$) in the case of Fig. \ref{E6a1_with_delta}$(a)$
  (resp. Fig. \ref{E6a1_with_delta}$(b)$), see Remark \ref{rem_max_root}(i).
   Then $(\delta, \beta_3) = (\alpha_1, \beta_3) \neq 0$, but this is impossible
   since the roots $\delta$ and $\beta_3$ are disconnected.

   (ii) As in Lemma \ref{lem_glue_1}, this follows from Corollary \ref{cor_map_2_diagr}.
  \qed

\subsubsection{Corrective reflections for the diagram $E_7(a_4)$}
  \label{sec_corr_refl}
\begin{figure}[h] \centering
 \includegraphics[scale=0.64]{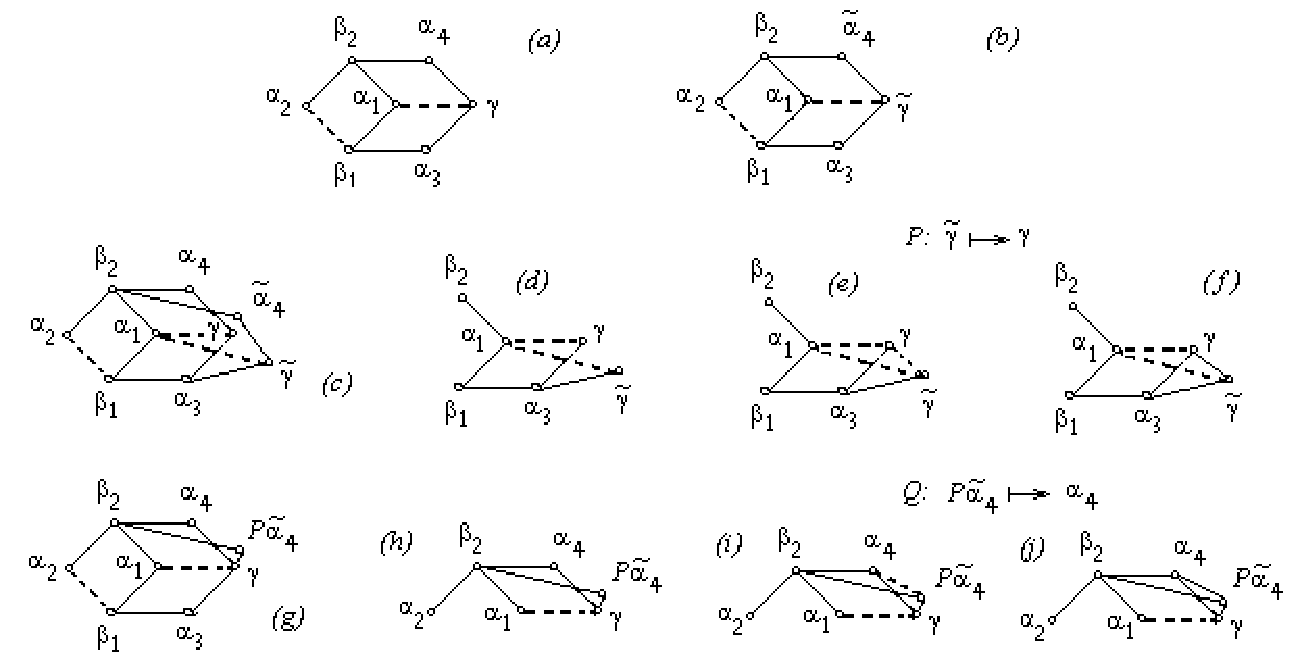}
 \caption{Two $E_7(a_4)$-associated subsets differing in two points:
 $\{\alpha_4, \gamma\}$ in the subset $(a)$ and $\{\widetilde\alpha_4, \widetilde\gamma\}$ in
 the subset $(b)$}
\label{E7a4_mir_v2_2}
\end{figure}

 \begin{lemma}
   \label{lem_glue_3}
   Let $\Gamma = E_7(a_4)$, let
   $S_1 = \{\alpha_1, \alpha_2, \alpha_3, \alpha_4, \beta_1, \beta_2, \gamma\}$ and
   $S_2 = \{\alpha_1, \alpha_2, \alpha_3, \widetilde\alpha_4, \beta_1, \beta_2, \widetilde\gamma\}$ be two
   $\Gamma$-associated subsets, see Fig. $\ref{E7a4_mir_v2_2}$$(a)$,$(b)$,
   each of which is linearly independent.
   Let
  \begin{equation*}
    \mathsf{Q}_6 = \{\alpha_1, \alpha_2, \alpha_3, \alpha_4, \beta_1, \beta_2\},
    \quad  \widetilde{\mathsf{Q}}_6 = \{\alpha_1, \alpha_2, \alpha_3, \widetilde\alpha_4, \beta_1, \beta_2\}
  \end{equation*}
   be $D_6(a_2)$-associated  subsets, see Fig. $\ref{E7a4_mir_v2_2}$$(c)$.
   Let $\alpha_4$ (resp. $\widetilde\alpha_4$) be connected with $\mathsf{Q}_6$ (resp. $\widetilde{\mathsf{Q}}_6$)
   only through $\beta_2$, let $\gamma$ (resp. $\widetilde\gamma$)
   constitute $\mathsf{P}3$-extension of $\mathsf{Q}_6$ (resp. $\widetilde{\mathsf{Q}}_6$)
   through the socket $\{\alpha_1, \alpha_3, \alpha_4\}$
   (resp. $\{\alpha_1, \alpha_3, \widetilde\alpha_4\}$), see Fig. $\ref{E7a4_mir_v2_2}$$(c)$.

   {\rm (i)} The configuration of vertices of Fig. $\ref{E7a4_mir_v2_2}$$(d)$
   (that is a fragment of Fig. $\ref{E7a4_mir_v2_2}$$(c)$)
   is  impossible: Roots $\gamma$ and $\widetilde\gamma$ are necessarily connected,
   see Fig. $\ref{E7a4_mir_v2_2}$$(e)$ or $(f)$.
   Let $P$ be the following  \underline{corrective reflection}:
   \begin{equation}
     \label{eq_corr_refl_1}
     P = s_{\gamma - \widetilde\gamma}  \text{ for Fig. }\ref{E7a4_mir_v2_2}(e), \qquad
     P = s_{\gamma +  \widetilde\gamma}  \text{ for Fig. }\ref{E7a4_mir_v2_2}(f).
   \end{equation}
   The map $P$ preserves $\{\alpha_1, \alpha_2, \alpha_3, \beta_1, \beta_2\}$, sends $S_2$ into the set
   $S^{'}_2 = \{\alpha_1, \alpha_2, \alpha_3, P\widetilde\alpha_4, \beta_1, \beta_2, \gamma\}$,
   $\widetilde\gamma$ into $\gamma$ and
   $\widetilde\alpha_4$ into $P\widetilde\alpha_4$, see Fig. $\ref{E7a4_mir_v2_2}$$(g)$.

   {\rm (ii)} The configuration of vertices Fig. $\ref{E7a4_mir_v2_2}$$(h)$
   (that is a fragment of Fig. $\ref{E7a4_mir_v2_2}$$(g)$)
   is  impossible: Roots $\alpha_4$ and $P\widetilde\alpha_4$ are necessarily connected,
   see Fig. $\ref{E7a4_mir_v2_2}$$(i)$ or $(j)$.
   Let $Q$ be the following  \underline{corrective reflection}:
    \begin{equation}
     \label{eq_corr_refl_2}
     Q = s_{\alpha_4 - P\widetilde\alpha_4}  \text{ for Fig. }\ref{E7a4_mir_v2_2}(i), \qquad
     Q = s_{\alpha_4 + P\widetilde\alpha_4}  \text{ for Fig. }\ref{E7a4_mir_v2_2}(j).
   \end{equation}
   The map $Q$ preserves $\{\alpha_1, \alpha_2, \alpha_3, \beta_1, \beta_2, \gamma\}$ and
   sends $P\widetilde\alpha_4$ into $\alpha_4$. Finally, $Q$ sends $S'_2$ into $S_1$.

   {\rm (iii)} Let $w_1$ (resp. $w_2$) be $S_1$-associated (resp. $S_2$-associated) element,
   see Fig. $\ref{E7a4_mir_v2_2}$$(a)$,$(b)$.
   Then $w_1 \simeq w_2$, the conjugacy is carried out by means of the map $QP$:
   \begin{equation}
        w_2 = (QP)^{-1}w_1{QP}.
   \end{equation}
 \end{lemma}
~\\

~\\
 \PerfProof (i) The configuration of vertices of Fig. $\ref{E7a4_mir_v2_2}$$(d)$
   is impossible due to Lemma \ref{lem_glue_1}. Therefore, $\gamma$ and $\widetilde\gamma$ are connected.
   Further, we apply the corrective reflection \eqref{eq_corr_refl_1}.

   (ii) Analogously, the configuration of vertices of Fig. $\ref{E7a4_mir_v2_2}$$(h)$
   is impossible. Therefore, $\alpha_4$ and $P\widetilde\alpha_4$ are connected.
   Finally, we apply the corrective reflection \eqref{eq_corr_refl_2}.

   (iii) This follows from (i) and (ii).  \qed
~\\

\section{\sc\bf The uniqueness of conjugacy classes for Carter diagrams from $\mathsf{C4} \coprod \mathsf{DE4}$}
     \label{sec_unique}

\subsection{The base cases of the uniqueness theorem: $D_4$ and $D_4(a_1)$}
 \label{sect_base_case}
\subsubsection{Conjugate dipoles}

 \index{non-conjugate dipoles}
 \index{Weyl group}

\begin{lemma}[On conjugate dipoles]
 \label{lem_united_2}
  For the case $\Gamma = D_4(a_1)$ (resp. $\Gamma = D_4$), let
   \begin{equation}
     \label{eq_united_2_cycles}
       \begin{cases}
           C_1 = \{\varphi_1, \varphi_2, \delta_1, \delta_2 \}, \\
           C_2 = \{\beta_1, \beta_2, \alpha_1, \alpha_2 \},
       \end{cases}
       \quad \text{ resp. } \quad
       \begin{cases}
           C_1 = \{\delta_1, \delta_2, \delta_3, \varphi_1\}, \\
           C_2 = \{\alpha_1, \alpha_2, \alpha_3, \beta_1 \}
       \end{cases}
   \end{equation}
  be two $D_4(a_1)$-associated subsets (resp. $D_4$-associated subsets) and
   \begin{equation}
     \label{eq_united_2_elem}
       \begin{cases}
          w_1 = s_{\varphi_1}s_{\varphi_2}s_{\delta_1}s_{\delta_2}, \\
          w_2 = s_{\beta_1}s_{\beta_2}s_{\alpha_1}s_{\alpha_2},
       \end{cases}
       \quad \text{ resp. } \quad
       \begin{cases}
          w_1 = s_{\delta_1}s_{\delta_2}s_{\delta_3}s_{\varphi_1}, \\
          w_2 = s_{\alpha_1}s_{\alpha_2}s_{\alpha_3}s_{\beta_1}
       \end{cases}
   \end{equation}
  be $C_i$-associated elements, where $i = 1,2$.

  \underline{Despite} the fact that there exist non-conjugate dipoles in $C_1$
  and $C_2$ (see Lemma $\ref{lem_eq_2diag}$) for any $C_1$, $C_2$ given by
  \eqref{eq_united_2_cycles}, \underline{we can choose} conjugate
  dipoles $d_1 \in C_1$ and $d_2 \in C_2$ such that $Td_2 = d_1$
  for some $T \in W$.
  For example, let  $T\{\alpha_1, \alpha_2 \} = \{\delta_1, \delta_2\}$, then
  $T^{-1}w_2{T} = s_{T\beta_1}s_{T\beta_2}s_{\delta_1}s_{\delta_2}$.
\end{lemma}

  \PerfProof
   For $E_6, E_7, E_8$, by Corollary \ref{cor_orth_roots}, there is only one equivalence class of dipoles under
  the Weyl group. Consider $\varPhi(D_l)$. By Corollary \ref{cor_conj_diag}(i)
  (resp. by Corollary \ref{cor_conj_diag}(ii))  each two $4$-cycles in  $\varPhi(D_l)$
  (resp. each two $D_4$-associated subsets in  $\varPhi(D_l)$)
  contain dipoles equivalent to each other under the action of
  $W(D_l)$, i.e., there exist dipoles $d_1 \in C_1$ and $d_2 \in C_2$,
  and the element $T\in W$ sending $d_2$ into $d_1$.
  \qed
~\\

\begin{figure}[h]
\centering
\includegraphics[scale=0.7]{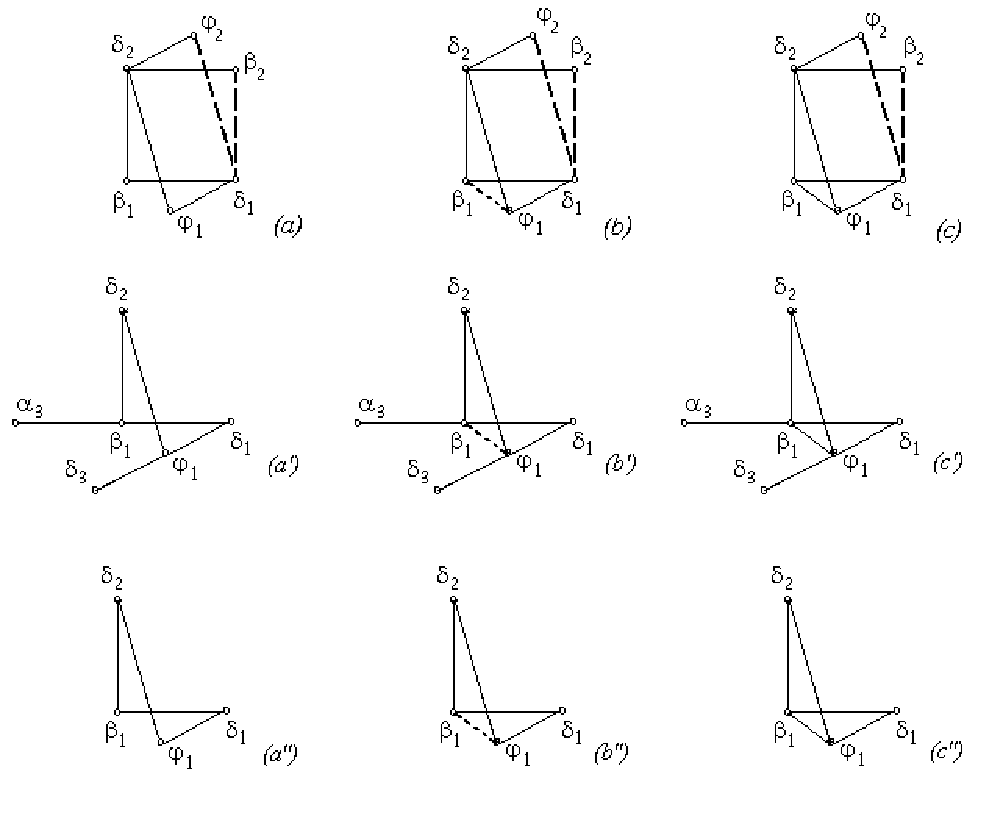}
 \caption{Two $\Gamma$-associated subsets $C_1$ and
 $C_2$ with a common dipole $\{\delta_1, \delta_2\}$, where
 $\Gamma = D_4(a_1)$ or $D_4$. There exists $T \in W$ such that
 $TC_1$ and $C_2$ coincide in $3$
 common roots, $T : \varphi_1 \longmapsto \beta_1$, see Fig. \ref{three_common_roots}}
\label{same_diagon_2}
\end{figure}

\subsubsection{From conjugate dipoles to conjugate triples}
  \label{sect_dipoles_triples}
  \index{linearly dependent root}
  In further considerations, we can assume that $T$ had already
  transformed  $d_1$ into $d_2$, and we assume that $d_1 = d_2$,
  i.e.,  $\{\alpha_1, \alpha_2 \} = \{\delta_1, \delta_2\}$,
  see Fig. \ref{same_diagon_2}.

\begin{lemma}
 \label{lem_united_3}
  For the case $\Gamma = D_4(a_1)$ (resp. $\Gamma = D_4$), let
   \begin{equation}
     \label{eq_united_3_cycles}
       \begin{cases}
           C_1 = \{\varphi_1, \varphi_2, \delta_1, \delta_2 \}, \\
           C_2 = \{\beta_1, \beta_2, \delta_1, \delta_2 \},
       \end{cases}
       \quad \text{ resp. } \quad
       \begin{cases}
           C_1 = \{\delta_1, \delta_2, \delta_3, \varphi_1\}, \\
           C_2 = \{\delta_1, \delta_2, \alpha_3, \beta_1 \}
       \end{cases}
   \end{equation}
  be two $D_4(a_1)$-associated subsets (resp. $D_4$-associated subsets) and
   \begin{equation}
     \label{eq_united_3_elem}
       \begin{cases}
          w_1 = s_{\varphi_1}s_{\varphi_2}s_{\delta_1}s_{\delta_2}, \\
          w_2 = s_{\beta_1}s_{\beta_2}s_{\delta_1}s_{\delta_2},
       \end{cases}
       \quad \text{ resp. } \quad
       \begin{cases}
          w_1 = s_{\delta_1}s_{\delta_2}s_{\delta_3}s_{\varphi_1}, \\
          w_2 = s_{\delta_1}s_{\delta_2}s_{\alpha_3}s_{\beta_1}
       \end{cases}
   \end{equation}
  be $C_i$-associated elements, where $i = 1,2$.

  There exists $T \in W$ \underline{sending $\varphi_1$ to $\beta_1$} and
  \underline{preserving $\delta_1$, $\delta_2$}.
  In other words, if two $\Gamma$-associated subsets $C_1$ and
  $C_2$ have a common dipole, then there exists $T \in W$ such that  $TC_1$ and $C_2$
  coincide in $3$ roots, see Fig. $\ref{same_diagon_2}$.
\end{lemma}

 \PerfProof
\underline{Cases $(a)$,$(a')$ in Fig. \ref{same_diagon_2}}.
  In this case, there is no connection $\{ \varphi_1, \beta_1 \}$.
  By Corollary \ref{cor_numb_ep}, $\varphi_1$
  cannot be linearly independent of $\{\beta_1, \beta_2, \delta_1, \delta_2\}$.
  By Lemma \ref{lem_on_square_1}, we have
  \begin{equation*}
     \varphi_1 + \beta_1  + \delta_1 + \delta_2 = 0,
     \quad \text{ i.e., } \quad  \varphi_1 = -(\beta_1  + \delta_1 + \delta_2),
  \end{equation*}
  For $\Gamma = D_4(a_1)$, case $(a)$ in Fig. \ref{same_diagon_2}, we have
  \begin{equation*}
    \begin{split}
      w_1 = & s_{\varphi_1}s_{\varphi_2}s_{\delta_1}s_{\delta_2} =
              s_{\beta_1  + \delta_1 + \delta_2}s_{\varphi_2}s_{\delta_1}s_{\delta_2}
              \stackrel{s_{\delta_1}s_{\delta_2}}{\simeq} \\
       & \big ( s_{\delta_1}s_{\delta_2}
        (s_{\beta_1  + \delta_1 + \delta_2}s_{\varphi_2})s_{\delta_1}s_{\delta_2} \big )
        s_{\delta_1}s_{\delta_2}  =
        s_{\beta_1}s_{\widetilde{\varphi}_2}s_{\delta_1}s_{\delta_2} =
        \widetilde{w}_1,\\
   \end{split}
 \end{equation*}
 where  $\widetilde{\varphi}_2 = T\beta_2 = \beta_2 -\delta_2 + \delta_1$.
 For  $\Gamma = D_4$, case $(a')$ in Fig. \ref{same_diagon_2}, we have
 \begin{equation*}
    \begin{split}
      w_1 = & s_{\varphi_1}s_{\delta_1}s_{\delta_2}s_{\delta_3} =
              s_{\beta_1  + \delta_1 + \delta_2}s_{\delta_1}s_{\delta_2}s_{\delta_3}
              \stackrel{s_{\delta_1}s_{\delta_2}}{\simeq} \\
       & \big ( s_{\delta_1}s_{\delta_2}
        (s_{\beta_1  + \delta_1 + \delta_2})s_{\delta_1}s_{\delta_2}s_{\delta_3} \big )
        s_{\delta_1}s_{\delta_2}  =
        s_{\beta_1}s_{\delta_1}s_{\delta_2}s_{\delta_3} = \widetilde{w}_1.
     \end{split}
  \end{equation*}
  In both cases, $\widetilde{w}_1$ is a $\widetilde{C}_1$-associated element,
  where
  \begin{equation}
   \label{eq_united_TC_1}
    \widetilde{C}_1 = TC_1 =
    \begin{cases}
         \{\beta_1, \widetilde{\varphi}_2, \delta_1, \delta_2\} \text { for } \Gamma = D_4(a_1), \\
          \{\beta_1, \delta_1, \delta_2\, \delta_3\} \text { for } \Gamma = D_4.
    \end{cases}
  \end{equation}
  Thus, $\widetilde{C}_1$ and $C_2$ coincide in $3$ roots  $\{\beta_1, \delta_1,  \delta_2\}$.
~\\

\begin{minipage}{9cm}
  \epsfig{file=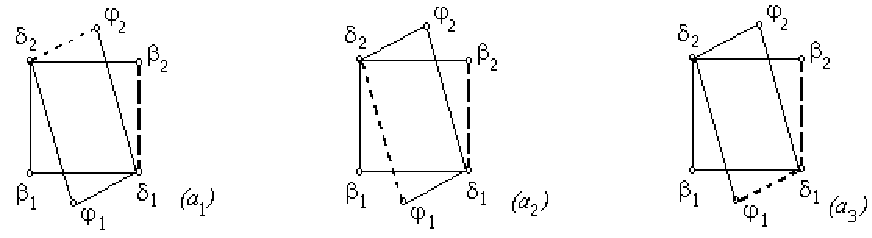, width=8.5cm}
  \captionof{figure}{}
  \label{same_diag_2_add}
\end{minipage}
\begin{minipage}{6.5cm}
\Small
 There are additional cases for gluing two diagrams $D_4(a_1)$ in the
 dipole $\{\delta_1, \delta_2\}$, see Fig. \ref{same_diag_2_add}.
 For the case $(a_1)$, we replace $\varphi_2$ by $-\varphi_2$;
 the reflection $s_{\varphi_2}$ is not changed, we get the case $(a)$ in Fig.
 \ref{same_diagon_2}. Also by replacing $\varphi_1$ by $-\varphi_1$
 we move from the case $(a_3)$ to the case $(a_2)$. The flip $\varphi_1 \longleftrightarrow
 \varphi_2$ interchanges cases $(a_2)$ and $(a_1)$.
\end{minipage}
~\\
~\\

  \underline{Cases $(b)$, $(b')$  in Fig. \ref{same_diagon_2}}. There is a dotted connection
  $(\varphi_1, \beta_1)$, i.e., $(\varphi_1, \beta_1) = \frac{1}{2}$.
  Therefore, $\beta_1 - \varphi_1$ is the root.
  By Lemma \ref{lem_gluing_1},
  for cases of Fig. \ref{same_diagon_2}$(b)$,$(b')$, the corrective reflection
  is $s_{\beta_1 - \varphi_1}$. For $\Gamma = D_4(a_1)$, we have
  \begin{equation*}
       w_1  \stackrel{s_{\beta_1 - \varphi_1}}{\simeq} s_{\beta_1 - \varphi_1}
       \big ( s_{\varphi_1}s_{\varphi_2}s_{\delta_1}s_{\delta_2} \big )s_{\beta_1 - \varphi_1}  =
        s_{\widetilde{\varphi}_2}s_{\beta_1}s_{\delta_1}s_{\delta_2} = \widetilde{w}_1, \\
  \end{equation*}
  where $\widetilde{\varphi}_2 = s_{\beta_1 - \varphi_1}(\varphi_2)$.
  For $\Gamma = D_4$, we have
  \begin{equation*}
       w_1  \stackrel{s_{\beta_1 - \varphi_1}}{\simeq} s_{\beta_1 - \varphi_1}
       \big ( s_{\varphi_1}s_{\delta_1}s_{\delta_2}s_{\delta_3} \big )s_{\beta_1 - \varphi_1}  =
        s_{\beta_1}s_{\delta_1}s_{\delta_2}s_{\delta_3} =  \widetilde{w}_1, \\
  \end{equation*}
  In both cases, $\widetilde{w}_1$ is a $\widetilde{C}_1$-associated element,
  where $\widetilde{C}_1$ is given in \eqref{eq_united_TC_1}, and
  $\widetilde{C}_1$ and $C_2$ again coincide in $3$ roots.
~\\

\begin{minipage}{9cm}
  \epsfig{file=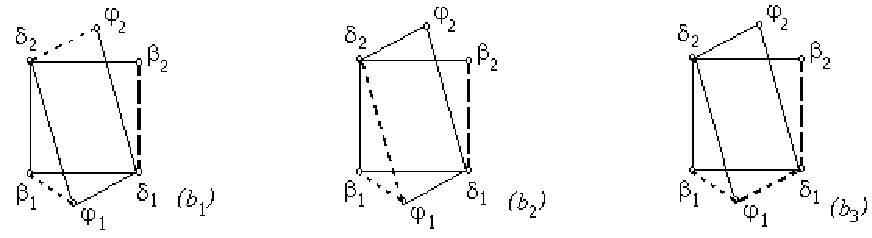, width=8.5cm}
  \captionof{figure}{}
  \label{same_diag_2_add_b}
\end{minipage}
\begin{minipage}{6.5cm}
\Small
\index{triangle of roots}
 Additional cases with the dotted edge $\{\beta_1, \varphi_1\}$ are
 depicted in Fig. \ref{same_diag_2_add_b}. For the case $(b_1)$,
 by replacing $\varphi_2$ by $-\varphi_2$ we get the case $(b)$ in Fig. \ref{same_diagon_2}.
 For the case $(b_2)$, we have the
 triangle $\{\beta_1, \varphi_1, \delta_2\}$
 and the square $\{\delta_1, \beta_2, \delta_2, \varphi_1\}$ of linearly dependent roots:
  \begin{equation}
   \label{eq_to_sum_1}
   \begin{split}
     & -\varphi_1 + \beta_1 + \delta_2 = 0, \\
     & -\delta_2  -\beta_2 + \delta_1 + \varphi_1 = 0.
   \end{split}
  \end{equation}
\end{minipage}
~\\
\begin{minipage}{16.2cm}
\Small
 see Lemmas \ref{lem_on_triangle_1} and \ref{lem_on_square_1}.
 Summing two equalities of \eqref{eq_to_sum_1} we get that
 the vectors of the triple $\{\beta_1, \beta_2, \delta_1\}$ are linearly
 dependent, which is a contradiction. For the case $(b_3)$, we have the
 triangle $\{\beta_1, \varphi_1, \delta_1\}$
 and the square $\{\delta_1, \beta_1, \delta_2, \varphi_2\}$ of linearly dependent roots:
  \begin{equation*}
   \begin{split}
     & -\varphi_1 + \beta_1 + \delta_1 = 0, \\
     & \delta_1 +  \beta_1 + \delta_2 + \varphi_2 = 0.
   \end{split}
  \end{equation*}
 Again, summing these two equalities we get that
 the vectors of the triple $\{\varphi_1, \varphi_2, \delta_2\}$ are linearly
 dependent, which is impossible.
\end{minipage}
~\\
~\\

  \underline{Cases $(c)$, $(c')$ in Fig. \ref{same_diagon_2}}.
  Suppose $\{\varphi_1, \beta_1\}$ is the solid connection, i.e., $(\varphi_1, \beta_1) = -\frac{1}{2}$.
  By Lemma \ref{lem_on_triangle_1}, each of the triples $\{\beta_1, \varphi_1, \delta_1\}$
  and $\{\beta_1, \varphi_1, \delta_2\}$ consists of linearly dependent vectors. Then by \eqref{eq_sum_0_pair} we have
  $\varphi_1 = -\delta_2 - \beta_2$ and  $\varphi_1 = -\delta_1 - \beta_2$.
  Then $\delta_1 = \delta_2$, which is impossible.
~\\

\begin{minipage}{9cm}
  \epsfig{file=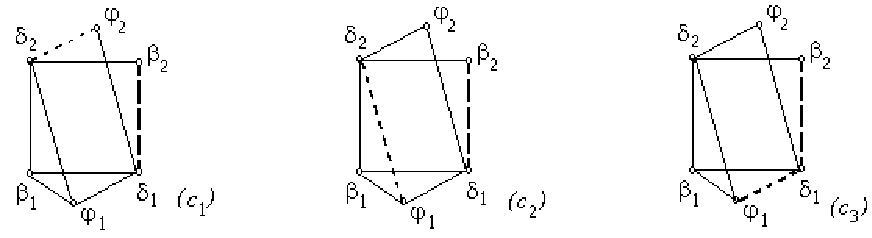, width=8.5cm}
  \captionof{figure}{}
  \label{same_diag_2_add_c}
\end{minipage}
\begin{minipage}{6.5cm}
\Small
 Additional cases with the solid edge $\{\beta_1, \varphi_1\}$ are
 depicted in Fig. \ref{same_diag_2_add_c}. For the case $(c_1)$, by replacing $\varphi_2$ by $-\varphi_2$
 we get the case $(c)$ in Fig. \ref{same_diagon_2}. For the case $(c_2)$ (resp. $(c_3)$),
 by replacing $\varphi_1$ by $-\varphi_1$ we get the case $(b_3)$ (resp. $(b_2)$)
 in Fig. \ref{same_diag_2_add_b}.
\end{minipage}
~\\
  \qed
~\\

\subsubsection{From conjugate triples to conjugate quadruples}
\begin{figure}[h]
\centering
\includegraphics[scale=0.7]{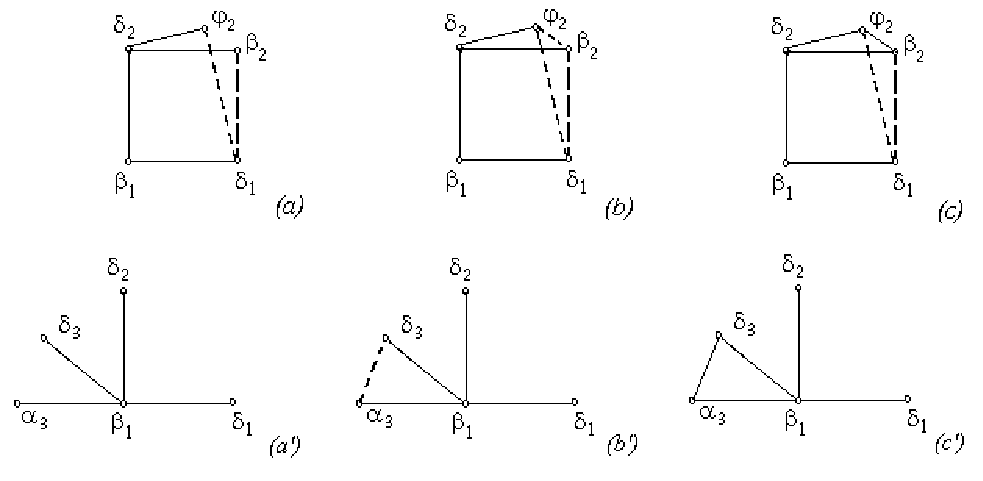}
 \caption{If two $\Gamma$-associated subsets $C_1$ and
 $C_2$ coincide in $3$ roots, then there exists $T \in W$ such that  $TC_1 = C_2$.
 For $\Gamma = D_4(a_1)$, $T : \varphi_2 \longmapsto \beta_2$;
 for $\Gamma = D_4$, $T : \delta_3 \longmapsto \alpha_3$}
\label{three_common_roots}
\end{figure}

  \index{linearly dependent root}
  Thanks to Lemma \ref{lem_united_3} we can assume that
  $\alpha_1 = \delta_1$, $\alpha_2 = \delta_2$ and $\varphi_1 =  \beta_1$,
  see Fig. \ref{three_common_roots}.

\begin{lemma}
 \label{lem_united_4}
  For the case $\Gamma = D_4(a_1)$ (resp. $\Gamma = D_4$), let
   \begin{equation}
     \label{eq_united_4_roots}
       \begin{cases}
           C_1 = \{\beta_1, \varphi_2, \delta_1, \delta_2 \}, \\
           C_2 = \{\beta_1, \beta_2, \delta_1, \delta_2 \},
       \end{cases}
       \quad \text{ resp. } \quad
       \begin{cases}
           C_1 = \{\delta_1, \delta_2, \delta_3, \beta_1\}, \\
           C_2 = \{\delta_1, \delta_2, \alpha_3, \beta_1 \}
       \end{cases}
   \end{equation}
  be two $D_4(a_1)$-associated subsets (resp. $D_4$-associated subsets) and
   \begin{equation}
     \label{eq_united_3_elem}
       \begin{cases}
          w_1 = s_{\beta_1}s_{\varphi_2}s_{\delta_1}s_{\delta_2}, \\
          w_2 = s_{\beta_1}s_{\beta_2}s_{\delta_1}s_{\delta_2},
       \end{cases}
       \quad \text{ resp. } \quad
       \begin{cases}
          w_1 = s_{\delta_1}s_{\delta_2}s_{\delta_3}s_{\beta_1}, \\
          w_2 = s_{\delta_1}s_{\delta_2}s_{\alpha_3}s_{\beta_1}
       \end{cases}
   \end{equation}
  be $C_i$-associated elements, where $i = 1,2$.

  For $\Gamma = D_4(a_1)$ (resp. $\Gamma = D_4$),
  there exists $T \in W$ \underline{sending $\varphi_2$ to $\beta_2$}
  (resp. \underline{$\alpha_3$ to $\delta_3$}) and
  \underline{preserving $\delta_1$, $\delta_2$, $\beta_1$}.
  In other words, if two $\Gamma$-associated subsets $C_1$ and
  $C_2$ coincide in $3$ roots, then there exists a certain $T \in W$ such that  $TC_1 = C_2$,
  see Fig. $\ref{three_common_roots}$.
\end{lemma}
~\\

~\\

\PerfProof

  \underline{Case $(a)$, $(a')$ in Fig. \ref{three_common_roots}.}
  Here, $(\varphi_2, \beta_2) = 0$ for case $(a)$, and $(\alpha_3, \delta_3) = 0$
  for case $(a')$. For case $(a)$,
  by Corollary \ref{cor_numb_ep} the root $\varphi_2$ and
  $\{\beta_2, \delta_1, \delta_2\}$ are linearly dependent. By Lemma \ref{lem_on_square_1}, we have
  \begin{equation*}
     \varphi_2 = -(\beta_2  - \delta_1 + \delta_2)
  \end{equation*}
  and
  \begin{equation*}
      w_2 = s_{\beta_2  - \delta_1 + \delta_2}s_{\beta_1}s_{\delta_1}s_{\delta_2} =
            s_{\beta_1}s_{\beta_2  - \delta_1 + \delta_2}s_{\delta_1}s_{\delta_2} =
            s_{\beta_1}s_{\delta_1}s_{\delta_2}s_{\beta_2}.
  \end{equation*}
  Then
  \begin{equation}
   \label{eq_2_ways_conj}
   \begin{split}
     & w_2  \stackrel{s_{\delta_1}s_{\delta_2}s_{\beta_2}}{\simeq} w_1,
              \quad \text{ and } \quad
       w_2  \stackrel{s_{\beta_1}}{\simeq} w_1.
   \end{split}
  \end{equation}
  ~\\
   For case $(a')$, the roots $\alpha_3$ and $\{\beta_1, \delta_1, \delta_2, \delta_3\}$ are linearly dependent,
   otherwise we get the extended Dynkin diagram $\widetilde{D}_4$.
   By Remark \ref{rem_max_root}(i), $\alpha_3$ is the minimal root for $D_4$, i.e.,
   $\alpha_3 = -(2\beta_1 + \alpha_1 + \alpha_2 + \alpha_3)$.
   We need to  prove that $vs_{\alpha_3} = vs_{\delta_3}$, where $v = s_{\beta_1}s_{\delta_1}s_{\delta_2}$.
   It follows from Proposition \ref{prop_conj_max_root}, eq. \eqref{eq_reduced_mu_max_2}.
~\\

  \underline{Case $(b)$, $(b')$ in Fig. \ref{three_common_roots}.}
  For $\Gamma = D_4(a_1)$, case $(b)$, we see that
  $\varphi_2$ is linearly independent of $\{\beta_2, \delta_1\}$,
  see Lemma \ref{lem_on_triangle_1}.
  Further, $(\varphi_2, \beta_2) = \frac{1}{2}$, and
  $\beta_2 - \varphi_2$ is the root. By Corollary
  \ref{cor_map_2_diagr}, we have the corrective reflection $s_{\beta_2 - \varphi_2}$:
   \begin{equation*}
      w_2  \stackrel{s_{\beta_2 - \varphi_2}}{\simeq}
      s_{\beta_2 - \varphi_2} \big ( s_{\varphi_2}s_{\beta_1}s_{\delta_1}s_{\delta_2} \big ) s_{\beta_2 - \varphi_2}  =
         s_{\beta_2}s_{\beta_1}s_{\delta_1}s_{\delta_2} = w_1.
  \end{equation*}
~\\
  For $\Gamma = D_4$, case $(b')$, we see that
  $\alpha_3$ is linearly independent of $\{\beta_1, \delta_3\}$, and
  $(\alpha_3, \delta_3) = \frac{1}{2}$. Then the corrective reflection is $s_{\alpha_3 - \delta_3}$:
   \begin{equation*}
      w_2  \stackrel{s_{\alpha_3 - \delta_3}}{\simeq}
      s_{\alpha_3 - \delta_3} \big (s_{\delta_1}s_{\delta_2}s_{\alpha_3}s_{\beta_1} \big ) s_{\alpha_3 - \delta_3}  =
         s_{\delta_1}s_{\delta_2}s_{\delta_3}s_{\beta_1} = w_1.
  \end{equation*}
~\\

\index{triangle of roots}

  \underline{Case $(c)$, $(c')$ in Fig. \ref{three_common_roots}.}
  For $\Gamma = D_4(a_1)$, case $(c)$, we have $(\varphi_2, \beta_2) = -\frac{1}{2}$,
  and $\varphi_2$ is linearly  dependent on $\{\beta_2, \delta_1\}$ and
  on $\{\beta_2, \delta_2\}$.
  By Lemma \ref{lem_on_triangle_1}, eq. \eqref{eq_sum_0_pair} we have $\varphi_2 = -\beta_2 + \delta_1$,
  and $\varphi_2 = -\beta_2 - \delta_1$, i.e., $\delta_1 = -\delta_2$, which is impossible.
  For $\Gamma = D_4$, case $(c')$, we have $(\alpha_3, \delta_3) =  -\frac{1}{2}$,
  and  $\alpha_3$ is linearly dependent on $\{\beta_1, \delta_3\}$,
  i.e., $\alpha_3 = -\delta_3 - \beta_1$. Then $(\alpha_3, \delta_1) = -(\beta_1, \delta_1) = \frac{1}{2}$,
  but this is impossible since $\alpha_3$ and $\delta_1$ are disconnected.
~\\

\begin{minipage}{9.5cm}
  \epsfig{file=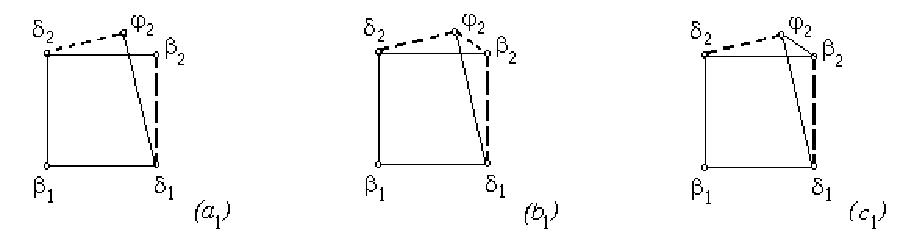, width=8.5cm}
  \captionof{figure}{}
  \label{three_com_roots_add}
\end{minipage}
\begin{minipage}{6cm}
\Small
 Additional cases with the dotted edge $\{\delta_2, \varphi_2\}$ are
 depicted in Fig. \ref{three_com_roots_add}. For the case $(a_1)$ (resp. $(b_1)$, resp. $(c_1)$),
 by replacing $\varphi_2$ by $-\varphi_2$,
 we get the case $(a)$ (resp. $(c)$, resp. $(b)$) in Fig. \ref{three_common_roots}. \qed
 \end{minipage}
~\\


\begin{corollary}[The base of induction]
  \label{cor_base_case_ind}
  For the case $\Gamma = D_4(a_1)$ (resp. $\Gamma = D_4$), let
   \begin{equation*}
       \begin{cases}
           C_1 = \{\varphi_1, \varphi_2, \delta_1, \delta_2 \}, \\
           C_2 = \{\beta_1, \beta_2, \alpha_1, \alpha_2 \},
       \end{cases}
       \quad \text{ resp. } \quad
       \begin{cases}
           C_1 = \{\delta_1, \delta_2, \delta_3, \varphi_1\}, \\
           C_2 = \{\alpha_1, \alpha_2, \alpha_3, \beta_1 \}
       \end{cases}
   \end{equation*}
  be two $D_4(a_1)$-associated subsets (resp. $D_4$-associated subsets) and let
   \begin{equation*}
       \begin{cases}
          w_1 = s_{\varphi_1}s_{\varphi_2}s_{\delta_1}s_{\delta_2}, \\
          w_2 = s_{\beta_1}s_{\beta_2}s_{\alpha_1}s_{\alpha_2},
       \end{cases}
       \quad \text{ resp. } \quad
       \begin{cases}
          w_1 = s_{\delta_1}s_{\delta_2}s_{\delta_3}s_{\varphi_1}, \\
          w_2 = s_{\alpha_1}s_{\alpha_2}s_{\alpha_3}s_{\beta_1}
       \end{cases}
   \end{equation*}
  be $C_i$-associated elements, where $i = 1,2$.

  Then there exists an element $T \in W$ such that $T^{-1}w_1{T}  = w_2$.
  In other words, any two $\Gamma$-associated elements $w_1$ and $w_2$ are conjugate.
\end{corollary}

  \PerfProof It follows from Lemmas \ref{lem_united_2}, \ref{lem_united_3}, \ref{lem_united_4}. \qed
~\\

\subsection{The uniqueness of conjugacy classes}
  \label{sec_uniq_C4DE4}

  \index{connected simply-laced Carter diagrams}

  Conjugate elements in the Weyl group $W$ are associated with the same Carter
  diagram $\Gamma$, or more exactly, with the class of diagrams equivalent to $\Gamma$.
  The Carter diagram  $\Gamma$ does not determine a single conjugacy class in $W$, \cite[Lemma 27]{Ca72}.
  Recall that $\mathsf{C4}$ (resp. $\mathsf{DE4}$) is the class of connected simply-laced
  Carter diagrams each of which \underline{contains a $4$-cycle $D_4(a_1)$}
  (resp. is a Dynkin diagram and \underline{contains $D_4$}),
  i.e., $\mathsf{C4}$ consists of Carter diagrams $E_l(a_i)$, $D_l(a_i)$ (resp. $\mathsf{DE4}$ consists of Carter diagrams
  $E_l$ for $l=6,7,8$, and $D_l$ for $l \geq 4$).
  Note that the $4$-vertex diagram $D_4(a_1)$ (resp. $D_4$) is the characterizing pattern for
  diagrams from $\mathsf{C4}$ (resp. $\mathsf{DE4}$), see Fig. \ref{tetris}.

  One of the main results of the paper is the following Theorem \ref{th_uniq_diagr} providing
  a sufficient condition of the uniqueness of conjugacy classes determining by a given Carter diagram.

 \begin{theorem}[On the conjugacy class of the diagram]
   \label{th_uniq_diagr}
   If $\Gamma$ is a connected simply-laced Carter diagram from $\mathsf{C4} \coprod \mathsf{DE4}$,
   then $\Gamma$ determines a single conjugacy class.
 \end{theorem}

  The theorem is derived from
  Corollary \ref{cor_base_case_ind} (the base of induction),
  Proposition \ref{prop_induct_step} (the induction step for homogeneous elements), and Theorem \ref{th_mirror_ext}
  (conjugacy of $\widetilde\Gamma_L$- and $\widetilde\Gamma_R$-associated elements for mirror extensions).
  \qed

\subsection{The induction step for homogeneous $\widetilde\Gamma$-associated elements}

   Two $\widetilde\Gamma$-associated elements are said to be {\it homogeneous $\widetilde\Gamma$-associated elements} if they
   are both $\widetilde\Gamma_R$-associated or both $\widetilde\Gamma_L$-associated.

\index{$\widetilde\Gamma_L$-associated element}
\index{$\widetilde\Gamma_R$-associated element}
\index{homogeneous $\widetilde\Gamma_R$-associated elements}
\index{homogeneous $\widetilde\Gamma_L$-associated elements}

  \begin{proposition} [On conjugacy of homogeneous $\Gamma$-associated elements]
    \label{prop_induct_step}
    Let $\Gamma$ be a Carter diagram such that all $\Gamma$-associated elements are
    conjugate.

    {\rm (i)} For any \underline{single-track extension} $\Gamma < \widetilde\Gamma$,
    all $\widetilde\Gamma$-associated elements are also conjugate.

    {\rm (ii)} For any \underline{mirror extensions} $\Gamma < \widetilde\Gamma_L$ and $\Gamma < \widetilde\Gamma_R$,
    all homogeneous $\widetilde\Gamma_R$-associated (resp. $\widetilde\Gamma_L$-associated)
    elements are also conjugate.

    {\rm (iii)} For any \underline{threefold extensions} $\Gamma < \widetilde\Gamma_1$,
    $\Gamma < \widetilde\Gamma_2$ and $\Gamma < \widetilde\Gamma_3$
     all homogeneous $\widetilde\Gamma_1$-associated (resp. $\widetilde\Gamma_2$-associated,
    resp. $\widetilde\Gamma_3$-associated) elements are also conjugate.

  \end{proposition}
  \index{pinholes number}
  \PerfProof
     By Proposition \ref{prop_extensions}, it suffices to
     prove that the Single-track Condition holds for the extension $\Gamma < \widetilde\Gamma$,
     i.e., $ws_\alpha \simeq ws_\gamma$ for any
     $\Gamma$-associated element $w$, regardless of the pinholes number
     of the extension $\Gamma < \widetilde\Gamma$  ($\mathsf{P}1$, $\mathsf{P}2$, $\mathsf{P}3$),
     see Fig. \ref{alp_con_ualph} and Fig. \ref{alp_perp_ualph_1}.
     The proof is based on consideration of the three principal cases:

\index{Principal Case $1$: $\alpha$ is connected to $\gamma$}
\index{linearly dependent root}

\subsubsection {Principal Case $1$: $\alpha$ is connected to $\gamma$}
  \label{sec_case_1}
\begin{figure}[h]
\centering
 $\begin{array}{c}
  \includegraphics[scale=1.4]{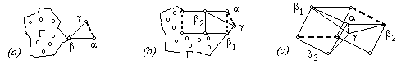} \\
  \qquad \mathsf{P}1  \qquad \qquad \qquad \qquad  \mathsf{P}2 \qquad \qquad \qquad \qquad \mathsf{P}3
 \end{array}$
\caption{Principal Case $1$: Single-track extensions}
\label{alp_con_ualph}
\end{figure}

     In this case, it does not matter whether $\gamma$ is linearly dependent on the
     subset of roots corresponding to $\widetilde\Gamma$ or not.
     First, consider single-track extensions.
     We have  $\Gamma \stackrel{\alpha}{<} \widetilde\Gamma$, ~$\Gamma \stackrel{\gamma}{<} \widetilde\Gamma$,
     see Fig. \ref{alp_con_ualph}.

\begin{figure}[h]
\centering
 \includegraphics[scale=0.62]{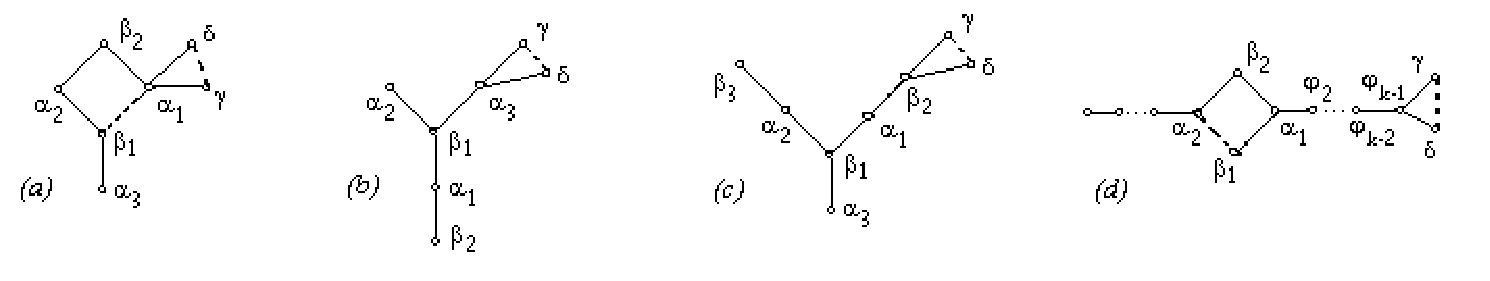}
 \caption{Principal Case $1$: Homogeneous elements for mirror $\mathsf{P}1$-extensions}
\label{mirror_P1_homogen_ext}
\end{figure}

\begin{figure}[h]
\centering
 $\begin{array}{c}
   \includegraphics[scale=0.62]{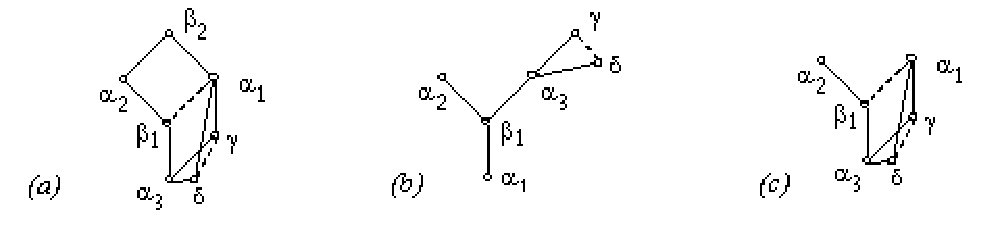} \\
    \mathsf{P}2 \qquad \qquad \qquad \qquad \mathsf{P}1  \qquad \qquad \qquad \qquad \mathsf{P}2
 \end{array}$
 \caption{Principal Case $1$: Homogeneous elements for mirror $\mathsf{P}2$-extensions
  and threefold $\mathsf{P}1$-  ad $\mathsf{P}2$-extensions}
\label{2n3_option_ext}
\end{figure}

 Corollary \ref{cor_map_2_diagr} implies that $ws_\alpha$ and $ws_\gamma$ are conjugate.
 Since $s_{\alpha - \gamma}$ and $w$ commute, then $s_{\alpha - \gamma}$ is the corrective reflection:
  \begin{equation*}
      s_{\alpha - \gamma}({w}s_{\gamma})s_{\alpha - \gamma} = w{s}_{\alpha}.
  \end{equation*}
  The same relation also holds for \underline{homogeneous associated elements} in
  the case of mirror extensions or threefold extensions (see (ii),(iii)).
  All cases in Fig. \ref{mirror_P1_homogen_ext} and case $(a)$ in Fig. \ref{2n3_option_ext} are mirror extensions.
  Cases $(b)$, $(c)$ in Fig. \ref{2n3_option_ext} are threefold extensions.

\index{Principal Case $2$: $\alpha \perp \gamma$, where $\gamma$ is linearly independent of $\widetilde\Gamma$.}

\subsubsection {Principal Case $2$: $\alpha \perp \gamma$, where $\gamma$ is linearly independent of $\widetilde\Gamma$}
  \label{sec_case_2}

 The exact wording of this case is as follows:
\begin{equation}
 \label{eq_prin_case_2}
 \begin{split}
  & \text{(i) }  \alpha  \text{ (resp. }  \gamma \text{) is the root extending } \Gamma \text{ to } \widetilde\Gamma,
      \text{ and } \alpha \perp \gamma, \\
  & \text{(ii) } \gamma \text{ is \underline{linearly independent} of vectors of } S \cup \alpha,
      \text{ where } S \text{ is a } \Gamma \text{-associated root subset}.
 \end{split}
\end{equation}
~\\

\begin{figure}[h]
\centering
 $\begin{array}{c}
 \includegraphics[scale=0.55]{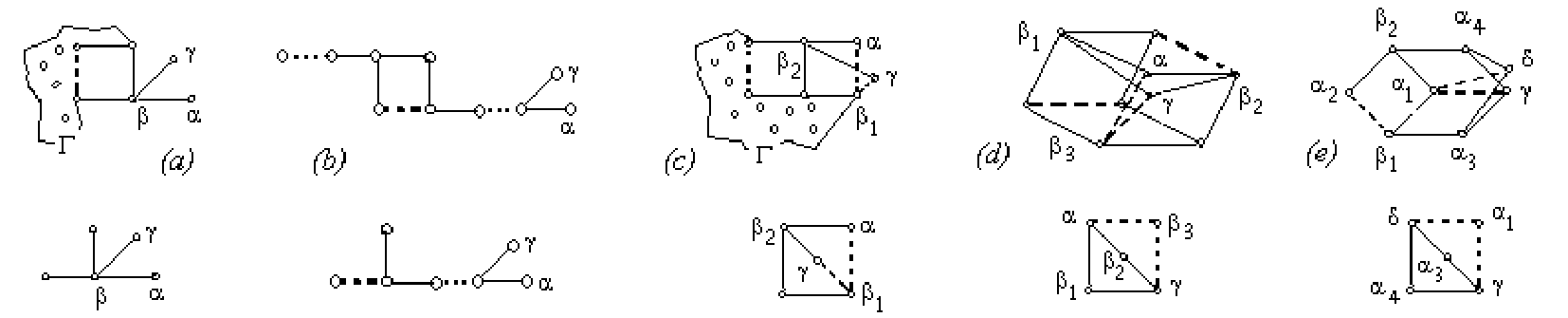} \\
  \mathsf{P}1 ~ \Small{(D\text{-joint type})}  \qquad \qquad \mathsf{P}1 ~{\Small(A\text{-joint type})}
  \qquad \qquad \qquad  \mathsf{P}2 \qquad \qquad \qquad \qquad
  \mathsf{P}3 \qquad \qquad \qquad \qquad  \mathsf{P}3 \qquad
 \end{array}$
 \caption{Principal Cases $2$ and $3$: $\alpha \perp \gamma$}
\label{alp_perp_ualph_1}
\end{figure}
~\\

\begin{figure}[h]
\centering
 \includegraphics[scale=0.62]{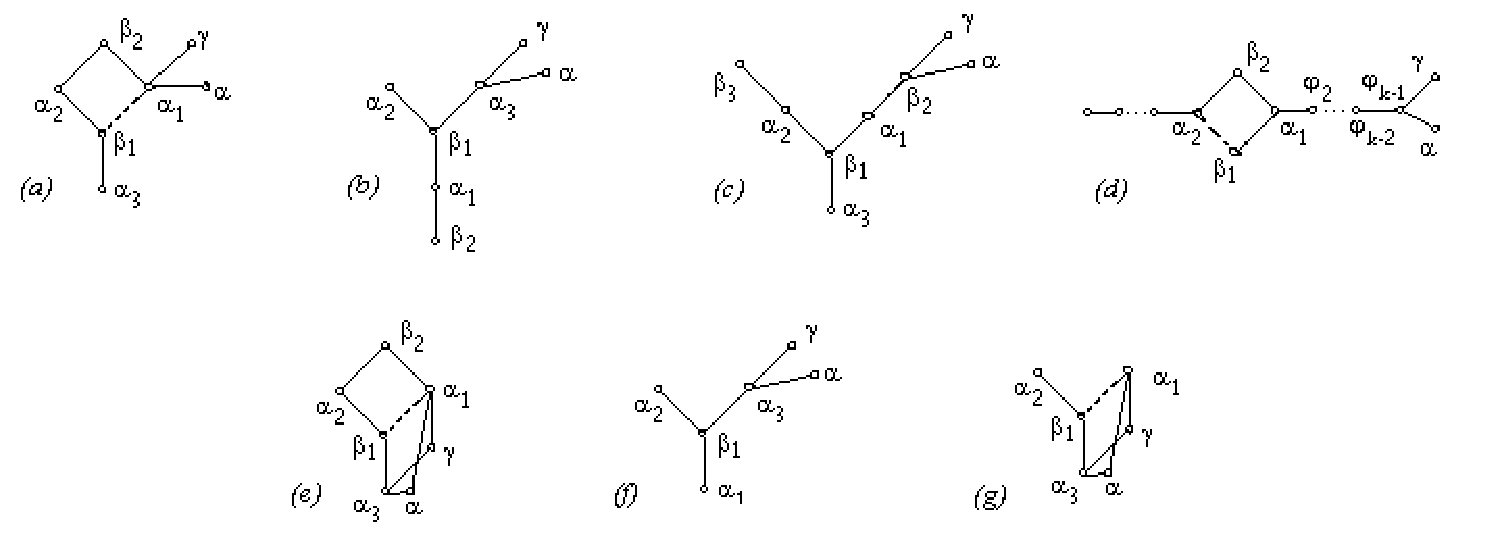}
 \caption{Principal Cases $2$ and $3$: Homogeneous elements
  for mirror and threefold $\mathsf{P}1$- and $\mathsf{P}2$-extensions}
\label{mirror_P1P2_homogen_2}
\end{figure}

 Denote by $\widetilde\Gamma^*$ the diagram obtained from $\widetilde\Gamma$ by adding the root $\gamma$,
 as depicted in Fig. \ref{alp_perp_ualph_1} or Fig. \ref{mirror_P1P2_homogen_2}.
 We will see that \underline{there does not exist a root $\gamma$ satisfying \eqref{eq_prin_case_2}}.
~\\

 \index{$\widetilde\Gamma^*$ (diagram extending $\widetilde\Gamma$)}
 For $\mathsf{P}1$-extensions of the $D$-joint type (resp. $A$-joint type), the diagram $\widetilde\Gamma^*$
 contains $\widetilde{D}_4$ (resp. $\widetilde{D}_l$, where $l > 4$) which cannot occur,
 see Fig. \ref{alp_perp_ualph_1}$(a)$,$(b)$;
 for mirror $\mathsf{P}1$-extension, see Fig. \ref{mirror_P1P2_homogen_2}$(a)$,$(b)$,$(c)$,$(d)$;
 for threefold $\mathsf{P}1$-extension, see Fig. \ref{mirror_P1P2_homogen_2}$(f)$.

\index{$\widetilde\Gamma^*$ (diagram extending $\widetilde\Gamma$)}
 For $\mathsf{P}2$- and $\mathsf{P}3$-extensions, the diagram $\widetilde\Gamma^*$ contains
 cycles with the bridge of length $> 1$ contradicting Corollary
 \ref{cor_numb_ep} for the following reasons. For $\mathsf{P}2$-extensions,
 see Fig. \ref{alp_perp_ualph_1}$(c)$; for homogeneous elements, in the case of mirror $\mathsf{P}2$-extensions,
 see Fig. \ref{mirror_P1P2_homogen_2}$(e)$; for homogeneous elements, in the case of threefold $\mathsf{P}2$-extensions,
 see Fig. \ref{mirror_P1P2_homogen_2}$(g)$.
 For $\mathsf{P}3$-extensions, see Fig. \ref{alp_perp_ualph_1}$(d)$: The path $\{\gamma, \beta_3, \alpha\}$
 forms a bridge;
 for homogeneous elements, in the case of mirror $\mathsf{P}3$-extensions, see Fig. \ref{alp_perp_ualph_1}$(e)$:
 Here, the path $\{\gamma, \alpha_3, \delta\}$ forms a bridge.

 \index{bridge}
 \index{homogeneous elements}
 \index{Principal Case $3$: $\alpha \perp \gamma$, where $\gamma$ is linearly dependent on $\widetilde\Gamma$}
 \index{quintuple of roots}

\subsubsection {Principal Case $3$: $\alpha \perp \gamma$, where $\gamma$ is
 linearly dependent on $\widetilde\Gamma$}
   \label{sec_case_lin_dep}

 The exact wording of this case is as follows:

\begin{equation}
 \label{eq_prin_case_3}
 \begin{split}
  & \text{(i) }  \alpha  \text{ (resp. }  \gamma \text{) is the root extending } \Gamma \text{ to } \widetilde\Gamma,
      \text{ and } \alpha \perp \gamma, \\
  & \text{(ii) } \gamma \text{ is \underline{linearly dependent} on vectors of } S \cup \alpha,
      \text{ where } S \text{ is the } \Gamma \text{-associated root subset}.
 \end{split}
\end{equation}
~\\

 Here, every type is considered separately.
 ~\\

 \underline{$\mathsf{P}1$-extensions, $D$-joint type}. By Lemma \ref{lem_glue_2},
 for all single-track $\mathsf{P}1$-extensions  from Table \ref{tab_single_track}
 except for the extension $D_l(a_1) < D_{l+1}(a_2)$,
 it cannot be that $\alpha \perp \gamma$ even if $\gamma$ is linearly dependent on
 the vectors of $S \cup \alpha$, see Fig. \ref{alp_perp_ualph_1}$(a)$.
 The same holds for homogeneous associated elements in
 the case of mirror $\mathsf{P}1$-extensions from Table \ref{tab_mirror_1}, see Fig. \ref{mirror_P1P2_homogen_2}$(a)$.
 For the case $D_l(a_1) < D_{l+1}(a_2)$ from Table
 \ref{tab_single_track}, the statement follows from Proposition \ref{prop_Dlak_and_Dk}(ii).
~\\

  \underline{$\mathsf{P}2$-extensions}.
  By Lemma \ref{lem_glue_1},
 for all single-track $\mathsf{P}2$-extensions from Table \ref{tab_single_track},
 it cannot be that $\alpha \perp \gamma$ even if $\gamma$ is linearly dependent on the vectors of
 $S \cup \alpha$, see Fig. \ref{alp_perp_ualph_1}$(c)$. The same holds for homogeneous associated elements in
 the case of mirror $\mathsf{P}2$-extensions from Table \ref{tab_mirror_1},
 see Fig. \ref{mirror_P1P2_homogen_2}$(e)$.
~\\

\begin{figure}[h]
$\begin{array}{c}
\includegraphics[scale=1.6]{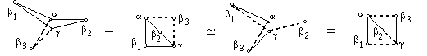} \\ \text{ Type } \mathsf{P}3
\end{array}$
 \caption{Principle Case $3$: $\mathsf{P}3$-extensions}
\label{cycles_of_E8a8}
\end{figure}

 \underline{$\mathsf{P}3$-extensions}.
   Consider the quintuple $\mathsf{Q_5} = \{\alpha, \beta_1, \beta_2, \beta_3, \gamma\}$ as a part of $E_8(a_8)$,
   see Fig. \ref{alp_perp_ualph_1}$(d)$.  By Lemma \ref{lem_on_diamond_1},
   the quintuple $\mathsf{Q_5}$ is impossible, see Fig. \ref{cycles_of_E8a8}.
   The same is true for the diagram $E_8(a_7)$ from the single-track $\mathsf{P}3$-extensions Table \ref{tab_single_track}.
   For the mirror $\mathsf{P}3$-extensions from Table \ref{tab_mirror_2}, there is only one diagram
   $E_7(a_4)$, see Fig. \ref{alp_perp_ualph_1}$(e)$. By Lemma
   \ref{lem_on_diamond_1}, the  quintuple $\{\alpha_1, \alpha_3, \alpha_4, \delta, \gamma\}$,
   a part of $E_7(a_4)$,  is also impossible.

\begin{figure}[h]
$\begin{array}{c}
\includegraphics[scale=0.62]{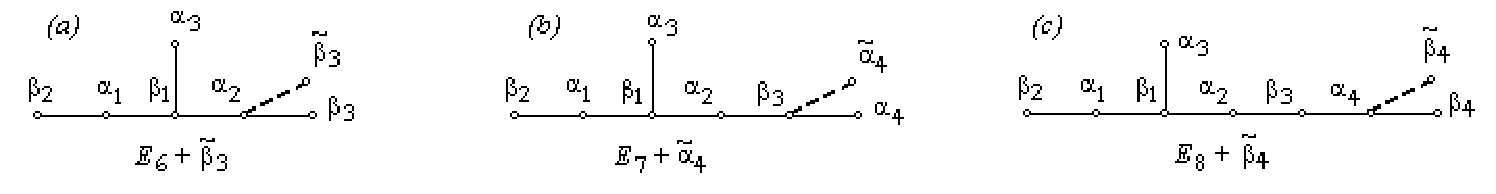}
\end{array}$
\caption{Impossible roots: $\widetilde\beta_3$ in
$E_6$; $\widetilde\alpha_4$ in $E_7$; $\widetilde\beta_4$ in $E_8$}
\label{E8_plus_Ubeta4}
\end{figure}

 \index{$b^{\vee}_{\eta, \eta}$ (diagonal element of $B^{-1}_{\Gamma}$)}
  \underline{$\mathsf{P}1$-extension, $A$-joint type}.
  1) First, we consider the extension $E_7 < E_8$ from Table \ref{tab_single_track}.
  Let there exist a root $\gamma$ connected with
  $\eta \in \widetilde\Gamma$ as in Fig. \ref{E8_plus_Ubeta4}$(c)$, where $\eta = \alpha_4$.
  We have
 \begin{equation*}
     b^{\vee}_{\alpha_4, \alpha_4} = 6 \text{ for } \Gamma = E_7 \text{ and }  \widetilde\Gamma = E_8,  \\
 \end{equation*}
  see Table \ref{tab_Cartan_E6_E7_E8} and Fig. \ref{E8_plus_Ubeta4}$(c)$.
  If there exists a root $\widetilde\beta_4$ such as in Fig. \ref{E8_plus_Ubeta4}$(c)$,
  it should be linearly dependent on the $E_8$-associated subset
   $\{\alpha_1, \alpha_2, \alpha_3, \alpha_4, \beta_1, \beta_2, \beta_3, \beta_4\}$.
  Since  $b^{\vee}_{\alpha_4, \alpha_4} > 2$, it follows by Remark \ref{rem_max_root}
  that there is no such a root $\gamma$ for $\widetilde\Gamma = E_8$.

\begin{figure}[h]
$\begin{array}{c}
\includegraphics[scale=0.62]{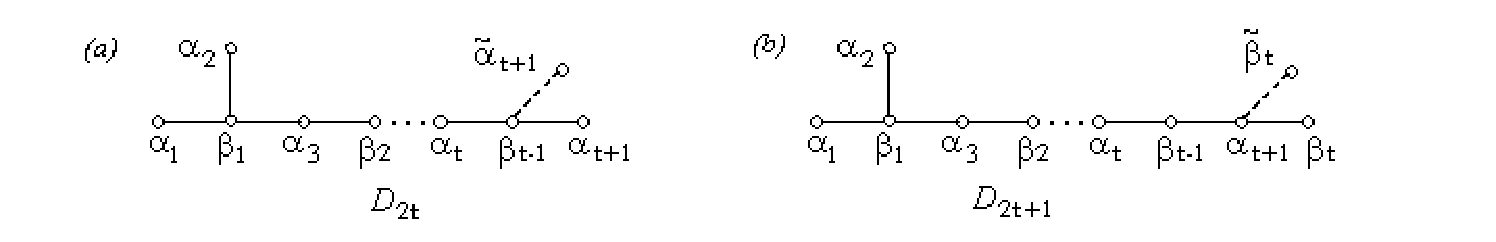}
\end{array}$
\caption{Maximal roots: $\widetilde\alpha_{t+1}$ in $D_{2t}$;
  $\widetilde\beta_t$ in $D_{2t+1}$}
\label{Dl_plus_max_root}
\end{figure}

  Similarly, there is no such a root $\gamma = \widetilde\alpha_4$  (resp. $\gamma = \widetilde\beta_3$)
  connected to $\beta_3$ (resp. $\alpha_2$) for mirror extensions $D_5 < E_6$
  (resp. $E_6 < E_7$) from Table \ref{tab_mirror_2}, see Fig. \ref{E8_plus_Ubeta4}$(a)$
  (resp. Fig. \ref{E8_plus_Ubeta4}$(b)$).
  The root $\gamma$ is linearly dependent on the $\widetilde\Gamma$-associated subset, otherwise the diagram
  $\widetilde\Gamma^*$ contains $\widetilde{D}_6$ (resp. $\widetilde{D}_7$).
  Then, for extensions $D_5 < E_6$ and $E_6 < E_7$, Proposition \ref{prop_induct_step}
  follows from Remark \ref{rem_max_root},
  we just check the value of the corresponding diagonal element of $B^{-1}_{\widetilde\Gamma}$:

 \begin{equation}
   \begin{split}
     & b^{\vee}_{\alpha_2, \alpha_2} = \frac{10}{3} > 2 \text{ for } \Gamma = D_5 \text{ and } \widetilde\Gamma = E_6, \\
     & b^{\vee}_{\beta_3, \beta_3} = 4 > 2 \text{ for } \Gamma = E_7 \text{ and } \widetilde\Gamma = E_7.
   \end{split}
 \end{equation}
~\\

 2) Consider the extension $D_l < D_{l+1}$ from Table
 \ref{tab_single_track}, where $D_{l+1}$ is obtained from $D_l$ by
  adding the endpoint $\delta = \alpha_{t+1}$ for $l = 2t$ or $\delta = \beta_t$ for $l =
  2t+1$, see Fig. \ref{Dl_plus_max_root}$(a)$,$(b)$. Let $\gamma$ be linearly dependent on $D_l + \{\delta\}$;
  let $\gamma$ also extend $D_l$ to $D_{l+1}$.
  By Remark \ref{rem_max_root}(i), such a root
  is necessarily the maximal root  $\mu_{max}$ in the root system $D_{l+1}$, see \S\ref{sec_max_root}.
  We have $\gamma = \mu_{max} = \alpha_{t+1}$ (resp. $\gamma = \mu_{max} = \beta_t$)
  for $D_{2t}$ (resp. $D_{2t+1}$).
  By Proposition \ref{prop_conj_max_root}, we have
   \begin{equation}
    \label{eq_conj_Dl_3}
    \begin{split}
      & ws_{\mu_{max}} = ws_{\gamma} \simeq ws_{\alpha_{t+1}} \text{ for the case } D_{2t}, \\
      & ws_{\mu_{max}} = ws_{\gamma} \simeq ws_{\beta_t} \text{ for the case } D_{2t+1},  \\
    \end{split}
   \end{equation}
   where $w$ is a certain $D_l$-associated element.

\begin{figure}[h]
\includegraphics[scale=0.6]{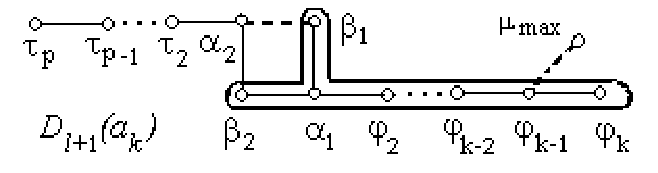}
 \caption{The Carter diagram $D_{l+1}(a_k)$ and bolded selected subdiagram $D_{k+2} \subset D_{l+1}(a_k)$}
\label{Dlak_extending_2}
\end{figure}

  3) Consider the extension $D_l(a_k) < D_{l+1}(a_k)$ from Table \ref{tab_single_track}.
     Let $D_{k+2}$ be the bolded selected subdiagram of $D_{l+1}(a_k)$, see Fig. \ref{Dlak_extending_2},
     let $S = \{\beta_1, \beta_2, \alpha_1, \varphi_1, \dots \varphi_{k-1}, \varphi_k\}$ be the $D_{k+2}$-associated subset
     of linearly independent roots.
     Let $\mu_{max}$ be the maximal root of  $D_{k+2}$.
     The root $\mu_{max}$ is linearly dependent on roots of $D_{k+2}$, and
     therefore it is linearly dependent on roots of $D_{l+1}(a_k)$. Further, $\mu_{max}$
     is orthogonal to any roots from $S$ except for $\varphi_{k-1}$, see Fig. \ref{Dlak_extending_2}
     (see also Fig. \ref{fig_Dlak_and_Dk}).
     In addition, $\mu_{max}$ is orthogonal to all remaining roots corresponding to $D_{l+1}(a_k)$ because
   \begin{equation}
     \label{eq_tau_perp}
     \begin{split}
       & \mu_{max} = \beta_1 + \beta_2 + 2\alpha_1 + 2\varphi_1 + \dots + 2\varphi_{k-1} + \varphi_k, \\
       & \mu_{max} \perp \tau_i, \text{ where } i = 2, \dots, p, \\
       & (\mu_{max}, \alpha_2) = (\beta_1 + \beta_2, \alpha_2) = \frac{1}{2} - \frac{1}{2} = 0.
      \end{split}
   \end{equation}
   see Fig. \ref{Dlak_extending_2}.  Further, if $\widetilde{w}$ is a $\Gamma$-associated
   element for $\Gamma = D_l(a_k)$, then by Proposition \ref{prop_Dlak_and_Dk}(ii), we have
   \begin{equation}
     \label{eq_conj_Dl_4}
         \widetilde{w}s_{\mu_{max}} ~\simeq~ \widetilde{w}s_{\varphi_k}. \\
   \end{equation}

   Thus, Proposition \ref{prop_induct_step} is proved also for the Principal Case $3$.
 \qed
~\\

\section{\sc\bf Mirror extensions and threefold extensions}
   \label{sec_mirror_threefold}

 For the regular extension $\Gamma \stackrel{\varphi}{<} \widetilde{\Gamma}$,
 the diagram $\widetilde{\Gamma}$ is obtained from $\Gamma$ by adding the vertex $\varphi$,
 where $\varphi$ is connected to $\Gamma$ at $n$ points, $n=1,2$ or $3$. The set of vertices of $\Gamma$
 connected to $\varphi$ is called a {\it socket}, see \S\ref{sec_1_2_option}.
 Note that each of diagrams $\Gamma$ and $\widetilde{\Gamma}$ is considered up to
 similarities \eqref{eq_1_equiv}, i.e., all edges connected to the same vertex $\alpha$ can be simultaneously
 changed in the following way: Every solid edge turns into a dotted
 one and vice versa. The Carter diagram $\widetilde{\Gamma}$ can be obtained as
 an extension of \underline{different Carter diagrams $\Gamma$}. For example,  $E_7(a_4)$
 can be obtained as the extension of $E_6(a_2)$, $D_6(a_2)$ and $D_6(b_2)$,
 see Fig. \ref{multi_sockets} and Table \ref{tab_2option_E7a4}.
 Note also that sometimes the extension
 $\Gamma < \widetilde{\Gamma}$ can be obtained \underline{from the same diagram $\Gamma$} in different
 ways ({\it multi-option extensions}, \S\ref{sec_1_2_option}) caused by symmetries of $\Gamma$.
 For example, the extension $D_4(a_1) < D_5(a_1)$, where $D_5(a_1) \subset E_l$, can be obtained in $4$ different
 ways, see Proposition \ref{prop_D4a1_toExt}.
 This section is devoted to the study of {\it mirror extensions}, see \S\ref{sec_1_2_option}.
 We will show that for a symmetric Carter diagram, all mirror extensions are conjugate. In other words, any
 $\widetilde\Gamma_L$-associated element and any $\widetilde\Gamma_R$-associated element
 are conjugate.

\index{$\widetilde\Gamma_L$-associated element $w_L$}
\index{$\widetilde\Gamma_R$-associated element $w_R$}

\begin{theorem}[On conjugacy of $\widetilde\Gamma_L$- and $\widetilde\Gamma_R$-associated elements]
 \label{th_mirror_ext}
   Let $\Gamma$ be a Carter diagram, $\Gamma < \widetilde\Gamma_L$ and
   $\Gamma < \widetilde\Gamma_R$ left and right extensions from Tables $\ref{tab_mirror_1}$ -- $\ref{tab_mirror_2}$.
   Then any $\widetilde\Gamma_L$-associated element $w_L$
   and any $\widetilde\Gamma_R$-associated element $w_R$
   are conjugate, i.e.,
   \begin{equation}
     \label{eq_mirror_conj_2}
      w_R = T^{-1}{w}_L{T} \text{ for some } T \in W.
   \end{equation}
\end{theorem}

 For any Carter diagram $\Gamma$ and
 mirror extensions $\Gamma < \widetilde\Gamma_L$ and $\Gamma < \widetilde\Gamma_R$ from
 Tables \ref{tab_mirror_1} -- \ref{tab_mirror_2},
 we explicitly construct the map $T$ for any pair of elements $w_L$ and $w_R$ from \eqref{eq_mirror_conj_2}.
 In most cases, the map $T$ is the composition of the longest element $w_0 \in W(\Gamma_0)$,
 where $\Gamma_0$ is a certain Dynkin subdiagram of $\Gamma$, and the corrective reflection $s_{\alpha + \beta}$
 for some roots $\alpha$ and $\beta$, see \S\ref{sec_longest} and Remark \ref{rem_corr_rel}.
 For the proof of Theorem \ref{th_mirror_ext}, see \S\ref{sec_th_mirror_ext}.

\subsection{Why does $A_2$ determine a single conjugacy class?}
  \label{sec_why_A2}
\begin{figure}[h]
\includegraphics[scale=1.8]{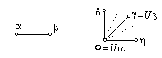}
\caption{Pairs of non-orthogonal roots $\{\varphi, \gamma\}$,
$\{\varphi, \delta\}$ and $\{\varphi, \eta\}$ are equivalent}
\label{A2_class}
\end{figure}

Let us show that for
any two pairs of non-orthogonal roots $\{ \alpha, \beta \}$ and $\{ \varphi, \delta \}$,
the elements $s_\alpha{s}_\beta$ and $s_\varphi{s}_\delta$ are conjugate:
$s_\alpha{s}_\beta \simeq s_\varphi{s}_{\delta}$.
First, there exists an element $U \in W$ sending $\alpha$ to $\varphi$. Since $\beta$ is connected to $\alpha$,
it follows that $U\beta$ is a root connected to $\varphi$, for example,  $U\beta = \gamma$, i.e.,
$U: \{ \alpha, \beta \} \longmapsto \{ \varphi, \gamma \}$, see Fig. \ref{A2_class}. This means that
\begin{equation*}
  \label{eq_conj_A2}
     U^{-1} (s_\alpha{s}_\beta) U = s_{U\alpha}{s}_{U\beta} = s_\varphi{s}_\gamma.
\end{equation*}
We construct a certain element $T \in W$ such that
\begin{equation*}
     T^{-1} (s_\varphi{s}_\gamma) T = s_\varphi{s}_\delta.
\end{equation*}

~\\
The square of the element $c = s_{\delta}s_{\gamma}s_\varphi$ sends the triple $\{\delta, \gamma, \varphi\}$
into the triple  $\{ -\gamma, -\delta, -\varphi\}$ since
\begin{equation}
 \label{eq_explicit_c}
  \begin{split}
     & c(\delta) = \varphi + \gamma,
      \quad c(\gamma) =  \varphi + \delta,
      \quad  c(\varphi) = -\varphi -\delta -\gamma, \\
      & c^2(\delta) = -\gamma,   \quad  c^2(\gamma) = -\delta, \quad c^2(\varphi) = -\varphi.
   \end{split}
\end{equation}
~\\
We put $T := c^2$. By \eqref{eq_explicit_c} we have
\begin{equation*}
     T^{-1} (s_\varphi{s}_\gamma) T = s_{c^2\varphi}{s}_{c^2\gamma} =
         s_{-\varphi}{s}_{-\delta} =  s_{\varphi}{s}_{\delta}. \qed
\end{equation*}

\index{Coxeter number}
\index{Coxeter element}

\begin{remark}
  \label{rem_T_as_w0}
{\rm
  Note that the element $c$ is the Coxeter element in the Weyl group $W(A_3)$,
  and the Coxeter number is $4$, i.e., $c^4 = 1$. The element $T = c^2$ is
  the longest element $w_0$ in $W(A_3)$.
  }
\end{remark}

\subsection{Conjugation $T$ acting as a mirror, cases $D_4(a_1)$ and $D_5(a_1)$}
  \label{sect_T_mirror}

 We will show that $T$ is the conjugation which preserves $D_4(a_1)$-associated elements and
 moves the \lq\lq{tail}\rq\rq of $D_5(a_1)$ to the opposite side of the dipole.

 \begin{figure}[h]
\includegraphics[scale=0.9]{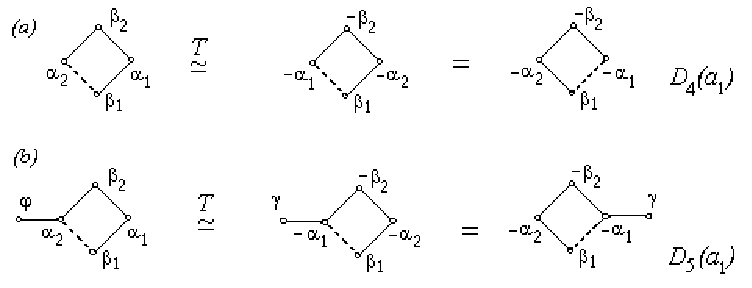}
\caption{
  Element $T = (s_{\beta_2}s_{\alpha_1}s_{\alpha_2})^2$ moves the
  \lq\lq{tail}\rq\rq $\{\varphi, \alpha_2\}$
  to $\{T\varphi, -\alpha_1\}$,
  where $\gamma := T\varphi = \varphi + \alpha_2 + \beta_2 + \alpha_1$}
\label{D5a1_map_1}
\end{figure}

 \begin{lemma}
   \label{lem_move_tail}
    {\rm (i)} Let  $\mathsf{Q}_4 = \{\alpha_1, \beta_1, \alpha_2, \beta_2\}$ be a quadruple
    corresponding to the $4$-cycle $D_4(a_1)$, see Fig. $\ref{D5a1_map_1}$$(a)$.
    Let $c: = s_{\beta_2}s_{\alpha_2}s_{\alpha_1}$, and $T: = c^2$.
    Then $T\mathsf{Q}_4 = \{-\alpha_2, \beta_1, -\alpha_1, -\beta_2 \}$,
    and the $\mathsf{Q}_4$-associated elements $v = s_{\alpha_1}s_{\alpha_2}s_{\beta_1}s_{\beta_2}$
    and $T^{-1}{v}T$ coincide:
    \begin{equation}
       \label{eq_self_eq}
         T^{-1}{v}T = v.
    \end{equation}

\index{quintuple of roots}

    {\rm (ii)}  Let  $\mathsf{Q}_5 = \{\alpha_1, \beta_1, \alpha_2, \beta_2, \varphi \}$ be a quintuple
    corresponding to $D_5(a_1)$, see Fig. $\ref{D5a1_map_1}$$(b)$. Then
    $T\mathsf{Q}_5 = \{-\alpha_2, \beta_1, -\alpha_1, -\beta_2, \gamma\}$,
    where $\gamma := T\varphi = \varphi + \alpha_2 + \beta_2 + \alpha_1$.
    Let $v{s}_{\varphi}$ be a $\mathsf{Q}_5$-associated element,
    where \underline{$\varphi$ is connected only with $\alpha_2$}.
    Then
    $$
       T^{-1}v{s}_{\varphi}T = T^{-1}vT{s}_{T\varphi} = v{s}_{T\varphi} = v{s}_{\gamma}
    $$
    is the $\mathsf{Q}_5$-associated element, where \underline{$T\varphi = \gamma$ is connected only with $\alpha_1$}.
 \end{lemma}

 \PerfProof
  (i) For the quadruple $\mathsf{Q}_4$, we have
 \begin{equation}
  \Small
    \label{eq_self_eq_3}
      \begin{split}
    & c\alpha_1 = -\alpha_1 - \beta_2, \quad c\alpha_2 = -\alpha_2 - \beta_2, \quad
      c\beta_1  = \beta_1 + \alpha_1 - \alpha_2, \quad c\beta_2 = \beta_2 + \alpha_1 + \alpha_2, \\
    & T\alpha_1 = -\alpha_2,  \quad T\alpha_2 = -\alpha_1, \quad
      T\beta_1 = \beta_1, \quad T\beta_2 = -\beta_2.
      \end{split}
  \end{equation}
  We derive from \eqref{eq_self_eq_3} that
  $$
   T^{-1}{v}T= s_{T\alpha_1}s_{T\alpha_2}s_{T\beta_1}s_{T\beta_2} =
    s_{-\alpha_2}s_{\alpha_1}s_{\beta_1}s_{-\beta_2} = v,
  $$
  i.e., eq. \eqref{eq_self_eq} holds.

  (ii) For the quintuple $\mathsf{Q}_5$, we have
 \begin{equation}
  \Small
    \label{eq_self_eq_4}
      \begin{split}
    & c\varphi = \varphi + \alpha_2 + \beta_2, \quad
      T\varphi = \varphi + \alpha_1 + \alpha_2 + \beta_2.
      \end{split}
  \end{equation}
  We derive from \eqref{eq_self_eq} and \eqref{eq_self_eq_4} that
  $$
     T^{-1}{v}s_{\varphi}T = vT^{-1}s_{\varphi}T = vs_{T\varphi}  = vs_{\gamma}. \qed
  $$

\begin{remark}{\rm
  \label{rem_tail}
 The main property that we want to get from Lemma \ref{lem_move_tail}
 is that in the case of need, the \lq\lq{tail}\rq\rq of $D_5(a_1)$
 can be transferred to the opposite side of the dipole, in other
 words, the \lq\lq{tail}\rq\rq  $s_\varphi$ in the
 $D_5(a_1)$-associated element $w$ can be transferred (by means of conjugation $T$) to $s_\gamma$:
  \begin{equation}
    \label{eq_property_T}
     \boxed{T^{-1}{v}s_{\varphi}T = vs_{\gamma}.}
  \end{equation}
  Note that, the action of $T$ on the subset $\{\alpha_1, \alpha_2, \beta_2\}$
  coincides with the action of the longest element $w_0$ in $W(A_3)$, see
  \S\ref{sec_why_A2} and Remark \ref{rem_T_as_w0}.
  }
\end{remark}

\subsection{The mirror extensions of some Carter diagrams}

\index{mirror diagrams $\widetilde\Gamma_L$, $\widetilde\Gamma_R$}
\index{mirror map $T$}

 Let $\Gamma$ be a diagram, and $\Gamma < \widetilde\Gamma_L$ be
 a mirror extension from Tables \ref{tab_mirror_1}, \ref{tab_mirror_2}.
 For different Carter diagrams, we define the element $T$ from Table
 \ref{tab_mirror_map_T}.  Let $\widetilde\Gamma_R = T\widetilde\Gamma_L$.
 We say that $\widetilde\Gamma_R$ is the {\it mirror diagram} for $\widetilde\Gamma_L$.
 The element $T$ is said to be the {\it mirror map}.
 It turns out that the property \eqref{eq_property_T} from Remark \ref{rem_tail}
 can be carried over to other mirror extensions of Carter diagrams.

\begin{proposition}
 \label{prop_mirror_ext}
   Let $\Gamma$ be one of diagrams from Tables $\ref{tab_mirror_1}$, $\ref{tab_mirror_2}$ except for the following cases:
   \begin{equation}
     \label{eq_excepr_ext}
     \begin{split}
       & D_4(a_1) < D_5(a_1), \text{ where } D_5(a_1) \subset E_l, \\
       & D_6(a_2) < E_7(a_4).
     \end{split}
   \end{equation}
  Let
 \begin{equation}
    T :=
  \begin{cases}
    w_0 = w_0(D_5) ~(\text{resp. } w_0(E_6)) & \text{ for } D_5 < E_6 ~(\text{resp. } E_6 < E_7),
               \text{ see eq. } \eqref{eq_w0_D5E6}, \\
    s_{\beta_1 - w_0\beta_1}, \text{ where } w_0 = w_0(A_{2k+1}), &
                \text{ for } D_{2k+2}(a_k) < D_{2k+3}(a_k), \text{ see } \S\ref{sect_spec_case_1},  \\
    (s_{\beta_2}s_{\alpha_1}s_{\alpha_2})^2 &
                \text{ for remaining cases of Tables } \ref{tab_mirror_1} \text{ and } \ref{tab_mirror_2} .
  \end{cases}
 \end{equation}
For any $\widetilde\Gamma_L$-associated element $w_L$, the conjugate element $w_R$ given as
   \begin{equation}
     \label{eq_mirror_conj}
      w_R := T^{-1}{w}_L{T}
   \end{equation}
   is the $\widetilde\Gamma_R$-associated element.
\end{proposition}

 \begin{remark}
   \label{rem_prop_th_mirror}
   {\rm In Proposition \ref{prop_mirror_ext}, we show that
   for any $\widetilde\Gamma_L$-associated element $w_L$, there exists
   a conjugate $\widetilde\Gamma_R$-associated element $w_R$.
   In Theorem \ref{th_mirror_ext}, we show (for the proof, see below \S\ref{sec_th_mirror_ext})
   that \underline{any} $\widetilde\Gamma_L$-associated element
   and \underline{any} $\widetilde\Gamma_R$-associated element are conjugate.
   }
 \end{remark}

 {\it Proof of Proposition \ref{prop_mirror_ext}}.
   1) If $\widetilde{\Gamma}$ is one of diagrams $E_6(a_1)$ or $ E_6(a_2)$
   (resp. $E_7(a_1)$ or $E_7(a_3)$, resp. $E_8(a_1)$ or $ E_8(a_4)$) then,
   in addition to eqs. \eqref{eq_self_eq_3} and \eqref{eq_self_eq_4},
   the conjugation $T$ also preserves the root $\alpha_3$
   (resp. $\alpha_3$, $\beta_3$, resp. $\alpha_3$, $\beta_3$, $a_4$),
   see Tables \ref{tab_mirror_1} and \ref{tab_mirror_2}.
   The conjugacy of any $\widetilde\Gamma_R$-associated element
   and any $\widetilde\Gamma_L$-associated element follows from
   ~\\

   $(a)$ the properties of the element $T$ as in Lemma \ref{lem_move_tail} and

   $(b)$ the conjugacy of homogeneous elements as in Proposition \ref{prop_induct_step}.
   ~\\

   2) Consider extensions $D_5 < E_6$ and $E_6 < E_7$, see Fig. \ref{D5_E6_E7_ext}.
      By relations \eqref{eq_w0_1} and Fig. \ref{D5_E6_E7_ext}, the longest element
      $w_0 = w_0(\Gamma)$ is as follows:
 \begin{equation}
  \Small
  \label{eq_w0_D5E6}
   \\
    w_0(D_5) =
   \left [ \begin{array}{cccccc}
      -1 & 0 & 0 & 0 & 0 & 2 \\
      0 & 0 & -1 & 0 & 0 & 2 \\
      0 & -1 & 0 & 0 & 0 & 2 \\
      0 & 0 & 0 & -1 & 0 & 3 \\
      0 & 0 & 0 & 0 & -1 & 1 \\
      0 & 0 & 0 & 0 & 0 &  1 \\
   \end{array}
     \right ]
    \begin{array}{c}
      \alpha_1 \\
      \alpha_2 \\
      \alpha_3 \\
      \beta_1 \\
      \beta_2 \\
      \beta_3 \\
   \end{array},
     \qquad
   w_0(E_6) =
   \left [
   \begin{array}{ccccccc}
      0 & -1 & 0 & 3 & 0 & 0 & 3 \\
      -1 & 0 & 0 & 3 & 0 & 0 & 3 \\
      0 & 0 & -1 & 2 & 0 & 0 & 2 \\
      0 & 0 & 0 &  1 & 0 & 0 & 0 \\
      0 & 0 & 0 &  4 & -1 & 0 & 0 \\
      0 & 0 & 0 &  2 & 0 & 0 & -1 \\
      0 & 0 & 0 &  2 & 0 & -1 & 0 \\
   \end{array}
     \right ]
   \begin{array}{c}
      \alpha_1 \\
      \alpha_2 \\
      \alpha_3 \\
      \alpha_4 \\
      \beta_1 \\
      \beta_2 \\
      \beta_3 \\
   \end{array} \\ \\
 \end{equation}
  Note that $w_0$ is the product of reflections associated with $\Gamma$, see \S\ref{sec_longest}.
  However, $w_0$ acts also on roots associated with $\widetilde{\Gamma}$. In particular, for $D_5 < E_6$,
  the element $w_0$ acts on $\beta_3$ as follows:
 \begin{equation}
    w_0\beta_3 = 2\alpha_1 + 2\alpha_2 + 2\alpha_3 + 3\beta_1 + \beta_2 + \beta_3 = \mu_{\max}(E_6),
 \end{equation}
 where $w_0 = w_0(D_5)$ in \eqref{eq_w0_D5E6} is given in the basis  $\{\alpha_1, \alpha_2, \alpha_3, \beta_1, \beta_2, \beta_3 \}$.
 For $E_6 < E_7$, the element $w_0$ acts on $\alpha_4$ as follows:
 \begin{equation}
   w_0\alpha_4 = 3\alpha_1 + 3\alpha_2 + 2\alpha_3 + \alpha_4 +  4\beta_1 + 2\beta_2 + 2\beta_3 = \mu_{\max}(E_7),
 \end{equation}
 where $w_0 = w_0(E_6)$  in \eqref{eq_w0_D5E6}
 is given in the basis $\{\alpha_1, \alpha_2, \alpha_3, \alpha_4, \beta_1, \beta_2, \beta_3 \}$.
 We put
 \begin{equation*}
   \delta =
      \begin{cases}
         \beta_3 \text{ for } D_5 < E_6, \\
         \alpha_4 \text{ for } E_6 < E_7. \\
      \end{cases}
 \end{equation*}
  By \eqref{decomp_w_Dl}, \eqref{eq_w0_1}, we have
 \begin{equation*}
  \begin{split}
   & w_L = {w}s_{\delta}, \text{ where } w = w_{\alpha}w_{\beta}, \\
   & w_0 = w^\frac{h}{2},  \text{ where } h \text{ is  the Coxeter number for } W(\Gamma).
  \end{split}
 \end{equation*}
  By Proposition \ref{prop_conj_max_root}, we have
\begin{equation*}
  w^{-1}_0w_L{w}_0  = (w^{-1}_0w{w}_0)(w^{-1}_{0}s_{\delta}{w}_0) =
  ws_{\mu_{max}} = w_R.
\end{equation*}
   Here,
\begin{equation*}
   \gamma = w_0\delta = \mu_{max} =
   \begin{cases}
      \begin{array}{ccccc}
        1 & 2 & 3 & 2 & 1 \\
        &     & 2 &     &
      \end{array} \qquad \qquad \text {for } D_5 < E_6,
      \\ \\
      \begin{array}{cccccc}
        2 & 3 & 4 & 3 & 2 & 1\\
          &   & 2 &   &   &
      \end{array} \qquad \text{ for } E_6 < E_7. \qquad\qquad \qed
   \end{cases}
\end{equation*}

 \subsubsection{Proof of Proposition \ref{prop_mirror_ext} for $D_{2k+2}(a_k) < D_{2k+3}(a_k)$}
    \label{sect_spec_case_1}
\begin{figure}[h]
\centering
\includegraphics[scale=0.7]{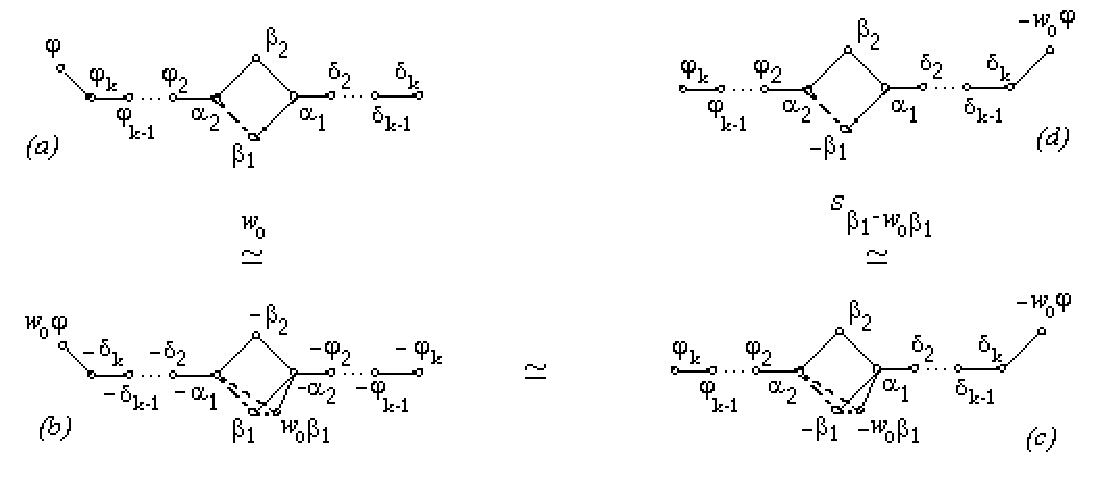}
\caption{The case $D_{2k+2}(a_k)$.
 Mirror extensions $\widetilde{\Gamma}_L$ and $\widetilde{\Gamma}_R = T\widetilde{\Gamma}_L$, where
 $T = s_{\beta_1 - w_0\beta_1} w_0$}
\label{fig_Dlak_mirror}
\end{figure}

 \index{longest element $w_0$}
 Let $w$ be a $\Gamma$-associated element, and $w_L = s_{\varphi}w$ a $\widetilde\Gamma_L$-associated element,
 where $\widetilde\Gamma_L$ is the extension of $\Gamma$ by means of the vertex $\varphi$ and the edge $\{\varphi_k, \varphi\}$,
 see Fig. \ref{fig_Dlak_mirror}$(a)$. We will show  that $w_L$ is conjugate to some $\widetilde\Gamma_R$-associated
 element $w_R$.
 Let $\Gamma'$ be the diagram $A_{2k+1}$ obtained from $\Gamma = D_{2k+2}(a_k)$ by deleting the vertex $\beta_1$
 and edges containing $\beta_1$, see Fig. \ref{fig_Dlak_mirror}.
 Let $S$ be the $\Gamma$-associated subset
\begin{equation*}
   S = \{\varphi_k,  \varphi_{k-1}, \dots, \varphi_2, \varphi_1 = \alpha_2, \beta_2, \beta_1,
   \alpha_1 = \delta_1, \delta_2, \dots, \delta_{k-1}, \delta_k \}.
\end{equation*}
 By adding the edge $\{\varphi_k, \varphi\}$, we extend $\Gamma$ to
 $\widetilde\Gamma$:
\begin{equation*}
   \Gamma' \quad \stackrel{\beta_1}{<} \quad \Gamma \quad \stackrel{\varphi}{<} \quad
   \widetilde\Gamma.
\end{equation*}
 It is a well-known fact that the longest element
 $w_0 \in W(A_{2k+1})$ determines the automorphism of the diagram
 $A_{2k+1}$, see \cite[Ch.6, \S4, $n^o$ 7, (XI)]{Bo}, and $w_0$
 acts as follows:
\begin{equation}
 \label{eq_w0_for_type_A}
 \begin{split}
    & \varphi_i ~\longmapsto~ -\delta_i, \quad \delta_i
    ~\longmapsto~ -\varphi_i \quad \text{ for } \quad i = 2,\dots,k, \\
    & \alpha_1 ~\longmapsto ~-\alpha_2, \quad \alpha_2
    ~\longmapsto~  -\alpha_1, \quad \beta_2 ~\longmapsto~ -\beta_2.
 \end{split}
\end{equation}
 For the proof of \eqref{eq_w0_for_type_A}, see \cite[\S 1.2, p. 7]{Fr01}.
 We extend eq. \eqref{eq_w0_for_type_A} from $\Gamma'$ to $\Gamma$
 and $\widetilde\Gamma$:

\begin{equation}
 \label{eq_w0_for_type_A_2}
 \begin{split}
   & w^{-1}_0({w}_L){w}_0 = w^{-1}_0(s_{\varphi}w){w}_0 = \\
   & w^{-1}_0 \big (s_{\varphi}\prod\limits_{i~even}(s_{\varphi_i}s_{\delta_i})
   s_{\alpha_1} s_{\alpha_2} \prod\limits_{i~odd}(s_{\varphi_i}s_{\delta_i}) s_{\beta_1} s_{\beta_2} \big )w_0 = \\
   & s_{w_0\varphi}\prod\limits_{i~even}({s_{-\delta_i}}{s_{-\varphi_i}})
   s_{-\alpha_2} s_{-\alpha_1} \prod\limits_{i~odd}({s_{-\delta_i}}{s_{-\varphi_i}})
   s_{w_0\beta_1} s_{-\beta_2} =  \\
   & s_{w_0\varphi}\prod\limits_{i~even}(s_{\varphi_i}s_{\delta_i})
   s_{\alpha_1} s_{\alpha_2} \prod\limits_{i~odd}(s_{\varphi_i}s_{\delta_i})
   s_{w_0\beta_1} s_{\beta_2}, \\
 \end{split}
\end{equation}
 see equivalences $(a) \stackrel{w_0}{\simeq} (b)$ and $(b) \simeq (c)$  in Fig. \ref{fig_Dlak_mirror}.
 The latter equality in \eqref{eq_w0_for_type_A_2} holds since
 every $s_{\varphi_i}$ commutes with every $s_{\delta_j}$.

 Now, we need only to fix the element $s_{w_0\beta_1}$.
 By Lemma \ref{lem_glue_1}, the vertices $\beta_1$ and
 $w_0\beta_1$ are necessarily connected. Further, we apply
 the corrective reflection $\tilde{s} = s_{\beta_1 - w_0\beta_1}$
 (resp. $ \tilde{s} = s_{\beta_1 + w_0\beta_1}$)
 if $(\beta, w_0\beta_1) = \frac{1}{2}$ (resp. $(\beta, w_0\beta_1) = -\frac{1}{2}$).
 By Corollary \ref{cor_map_2_diagr}, the element
 $\tilde{s}$ does not change any other reflection associated with roots from $S$. If, for example,
 $\tilde{s} = s_{\beta_1 - w_0\beta_1}$, we have

\begin{equation}
  \label{eq_w0_for_type_A_3}
 \begin{split}
   & s_{w_0\varphi}\prod\limits_{i~even}(s_{\varphi_i}s_{\delta_i})
   s_{\alpha_1} s_{\alpha_2} \prod\limits_{i~odd}(s_{\varphi_i}s_{\delta_i})
   s_{w_0\beta_1} s_{\beta_2}
   \quad \stackrel{s_{\beta_1 - w_0\beta_1}}{\simeq} \\
   \\
   & s_{w_0\varphi}\prod\limits_{i~even}(s_{\varphi_i}s_{\delta_i})
   s_{\alpha_1} s_{\alpha_2} \prod\limits_{i~odd}(s_{\varphi_i}s_{\delta_i})
   s_{\beta_1} s_{\beta_2} ~=~ s_{w_0\varphi}w ~=~ \quad  w_R,
 \end{split}
\end{equation}
see equivalence $(c) \simeq (d)$  in Fig. \ref{fig_Dlak_mirror}.
Then \eqref{eq_w0_for_type_A_2} and \eqref{eq_w0_for_type_A_3} yield
the equivalence $(a) \simeq (d)$  in Fig. \ref{fig_Dlak_mirror}. \qed

\subsection{Proof of Theorem \ref{th_mirror_ext} on mirror extensions}
  \label{sec_th_mirror_ext}
  For diagrams from Tables \ref{tab_mirror_1}, \ref{tab_mirror_2},
  proof of Theorem \ref{th_mirror_ext} follows from
  Propositions \ref{prop_induct_step}(ii) and \ref{prop_mirror_ext}, and Lemma \ref{lem_move_tail}.
  For proof of Theorem \ref{th_mirror_ext} for $D_6(a_2) < E_7(a_4)$,
  see \S\ref{sec_proof_E7a4}. For the case of adjacent extensions $D_4(a_1) < D_5(a_1)$,
  see \S\ref{sec_adjacent}.

 \subsubsection{Proof of Theorem \ref{th_mirror_ext} for mirror extensions $D_6(a_2) < E_7(a_4)$}
    \label{sec_proof_E7a4}
\begin{figure}[h] \centering
\includegraphics[scale=0.64]{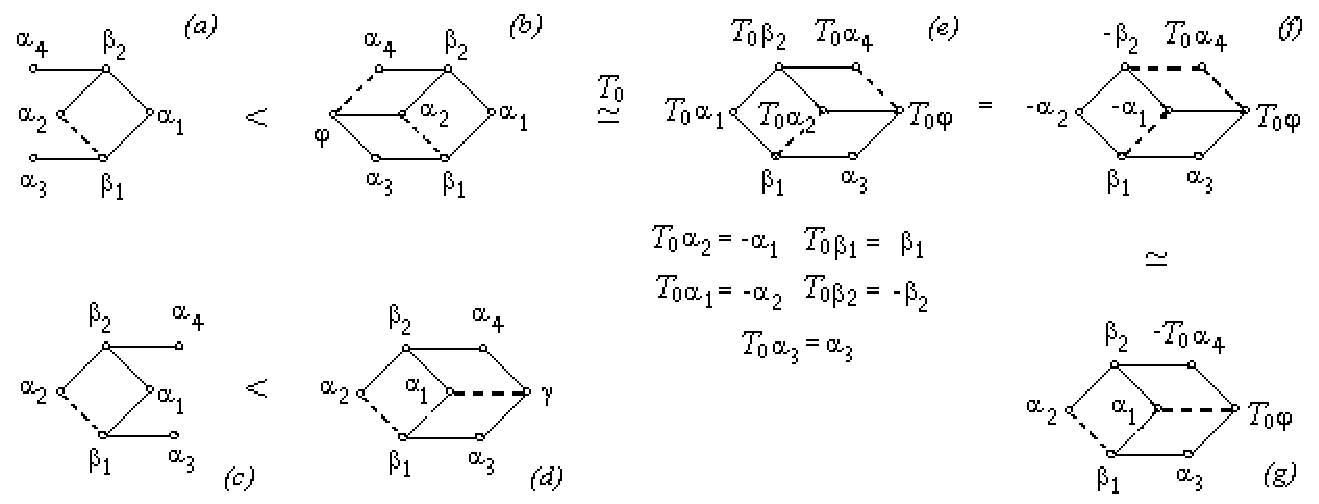}
\caption{Mirror extensions $D_6(a_2) < E_7(a_4)$}
\label{E7a4_mir_v2_1}
\end{figure}

 The Carter diagram $E_7(a_4)$ can be obtained from
 $D_6(a_2)$ by adding one vertex $\varphi$ and $3$
 edges: $\{\varphi, \alpha_2\}$, $\{\varphi, \alpha_3\}$ and
  $\{\varphi, \alpha_4\}$, see Fig. \ref{E7a4_mir_v2_1}$(a)$,$(b)$,
 or, alternatively, by adding the vertex $\gamma$
 and $3$ edges: $\{\gamma, \alpha_2\}$, $\{\gamma, \alpha_3\}$ and
  $\{\gamma, \alpha_4\}$, see Fig. \ref{E7a4_mir_v2_1}$(c)$,$(d)$.
 Let $\Gamma = D_6(a_2)$, and as
 above, in Proposition \ref{prop_mirror_ext} and
 Tables \ref{tab_mirror_1}, \ref{tab_mirror_2}, let
 $\Gamma < \widetilde{\Gamma}_L$ and $\Gamma < \widetilde{\Gamma}_R$ be two mirror
 extensions of $\Gamma$. As above, in Proposition
 \ref{prop_mirror_ext} and Tables \ref{tab_mirror_1}, \ref{tab_mirror_2}, we will show
 that the elements associated with diagrams $(b)$ and $(d)$ in Fig. \ref{E7a4_mir_v2_1} are conjugate.

 Let  $w_L$ (resp. $w_R$) be a $\widetilde{\Gamma}_L$-associated
 (resp. $\widetilde\Gamma_R$-associated) element. Let us prove that
 $w_L \simeq w_R$.  Let $T_0 := c^2$, where $c :=  s_{\beta_2}s_{\alpha_1}s_{\alpha_2}$.
  We  have
  \begin{equation*}
   \begin{split}
     & c\alpha_1 = -\alpha_1 - \beta_2, \quad c\alpha_2 = -\alpha_1 - \beta_2,
       \quad c\alpha_3 = \alpha_3, \quad c\alpha_4 = \alpha_4 + \beta_2, \\
     & c\beta_1 = \beta_1 - \alpha_1 + \alpha_2, \quad
       c\beta_2 = \beta_2 + \alpha_1 + \alpha_2,
       \quad c\varphi = \varphi + \alpha_1 + \beta_2.
   \end{split}
  \end{equation*}
 Then
  \begin{equation}
   \label{eq_T_E7a4}
   \begin{split}
     & T_0\alpha_1 = -\alpha_2, \quad T_0\alpha_2 = -\alpha_1, \quad T_0\alpha_3 = \alpha_3,
       \quad T_0\alpha_4 = 2\beta_2 + \alpha_1 + \alpha_2 + \alpha_4, \\
     & T_0\beta_1 = \beta_1, \quad T_0\beta_2 = -\beta_2,
       \quad T_0\varphi = \varphi + \alpha_1 + \alpha_2 + \beta_2.
   \end{split}
  \end{equation}

 By \eqref{eq_T_E7a4} the element $T_0$ maps the $\widetilde{\Gamma}_L$-associated
 element $w_L$, where $\widetilde{\Gamma}_L$ is the diagram on Fig. \ref{E7a4_mir_v2_1}$(b)$,
 to the $T_0\widetilde{\Gamma}_L$-associated element
 $w'_L = T_0^{-1}{w_L}T_0$, see Fig. \ref{E7a4_mir_v2_1}$(b)$,$(e),(f)$.
 By passing from Fig. \ref{E7a4_mir_v2_1}$(f)$ to Fig. \ref{E7a4_mir_v2_1}$(g)$,
 the roots $\alpha_1$, $\alpha_2$, $T_0\alpha_4$, $\beta_2$
 change the sign, the corresponding reflections are not changed, and
 $w'_L$ is preserved.  Diagrams $(d)$ and $(g)$ in Fig. \ref{E7a4_mir_v2_1}
 represent two $E_7(a_4)$-associated subsets differing only in two points: $\alpha_4$ and $\gamma$
 in the diagram $(d)$, and $T_0\alpha_4$ and $T_0\varphi$ in the diagram $(g)$.
 By Lemma \ref{lem_glue_3}, $w'_L$ and $w_R$ are conjugate, therefore $w_L$ and $w_R$
 are also conjugate. \qed

\index{adjacent extensions $D_4(a_1) < D_5(a_1)$}
\index{extension! - adjacent extensions $D_4(a_1) < D_5(a_1)$}

\subsubsection{Proof of Theorem \ref{th_mirror_ext} for adjacent extensions $D_4(a_1) < D_5(a_1)$}
   \label{sec_adjacent}

\begin{figure}[h]
\centering
\includegraphics[scale=0.7]{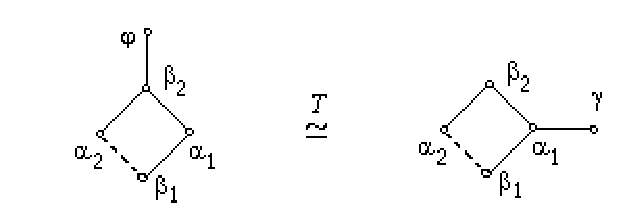}
 \caption{Adjacent extensions $D_4(a_1) < D_5(a_1)$
 obtained by adding edge $\{\beta_2, \varphi\}$ (resp. $\{\alpha_1, \gamma\}$) to
 the vertex $\beta_2$ (resp. $\alpha_1$)}
\label{adjacent_Ext}
\end{figure}

 Let $\Gamma$ be the $4$-cycle $D_4(a_1)$. The extensions $\Gamma < \Gamma_1$ and $\Gamma < \Gamma_2$,
 where $\Gamma_1$ and $\Gamma_2$ are two instances of the diagram $D_5(a_1)$,
 are said to be {\it adjacent} if they are obtained by adding edges
 in adjacent vertices of $D_4(a_1)$, see Fig. \ref{adjacent_Ext}.
 In this section, we consider extensions $D_4(a_1) < D_5(a_1)$ in
 $D_l$ and in $E_l$ and we find the conditions when these extensions may be adjacent.
~\\

\index{dipole}
\index{$2$-index dipole}
\index{$4$-index dipole}

  In what follows, we need the notion of {\it $2$-index and $4$-index dipoles in $4$-cycles}.
  The roots in the root system $D_l$ are given as
  \begin{equation}
    \label{eq_ei_ej}
      {\pm}e_i {\pm}e_j ~(1 \leq i < j \leq l), \text{ where } \{ e_i \mid i = 1, \dots l \}
      \text{ is an orthonormal basis in } \mathbb{R}^l,
  \end{equation}
  see \cite[Table IV]{Bo}.
  Any dipole of every $4$-cycle $D_4(a_1)$ in $D_l$ is determined either
  by $2$ indices or by $4$ indices. They are either
  \begin{equation}
    \label{diag_1}
      \{ e_k - e_n, e_k + e_n \}, \text{ where }  k \neq n ,
  \end{equation}
  or
  \begin{equation}
    \label{diag_2}
      \{ e_i \pm e_j, e_k \pm e_n \}, \text{ where }  i, j, k, n   \text{ are different}.
  \end{equation}
  We call the dipole of type  \eqref{diag_1} (resp. \eqref{diag_2}) the {\it $2$-index dipole}
  (resp. the {\it $4$-index dipole}).
  For properties of $2$-index and $4$-index dipoles, see Lemma \ref{lem_eq_2diag} and Corollary \ref{cor_conj_diag}
  in \S\ref{sec_2_4_index_dipole}.

 \begin{lemma}[On adjacent extensions in $D_l$]
  \label{lem_adjacent}
   Let $\Gamma = D_4(a_1)$ be a $4$-cycle in $D_l$.

   {\rm (i)} The $4$-cycle $D_4(a_1)$  can be extended to $D_5(a_1)$ only in the vertex lying in
   the $4$-index dipole of $D_4(a_1)$.

   {\rm (ii)} If both dipoles of $D_4(a_1)$
   are $4$-index dipoles, then there is no extension $D_4(a_1) < D_5(a_1)$.

   {\rm (iii)} There are no adjacent extensions $D_4(a_1) < D_5(a_1)$ for any $4$-cycle in $D_l$.
 \end{lemma}

   \PerfProof (i) Let $d = \{e_k + e_n, e_k - e_n\}$ be a $2$-index
   dipole in any $4$-cycle, see \S\ref{sec_cycles_D4} and Lemma \ref{lem_eq_2diag}.
   If $\varphi = e_i + e_j$ is connected to one of vertices of $d$, then $\{k, n\}\cap\{i, j\} \ne \emptyset$, and
   $\varphi$ is connected to another vertex of $d$. Thus, $\varphi$ does not extend
   $D_4(a_1)$ to $D_5(a_1)$.

   (ii) If $d_1$ and $d_2$ are two $4$-index dipoles of $D_4(a_1)$, then they look, for example, as follows:
   \begin{equation*}
     d_1 = \{\alpha_1 = e_k - e_n, \quad \alpha_2 = e_p - e_q\}, \quad
     d_2 = \{\beta_1 = e_n - e_q, \quad \beta_2 = -e_k - e_p\},
   \end{equation*}
   where $\{k, n, p, q\}$ are $4$ different indices.
   Let $\varphi$ be connected to one of the vertices, for example, to
   $e_k - e_n$. Then $\varphi = e_k \pm e_s$ or
   $\varphi = e_n \pm e_s$, where $s \not\in \{k, n, p, q\}$. In the first case, $\varphi$ is also connected
   to $-e_k - e_p$ , in the second case, $\varphi$ is also connected to $e_n - e_q$. In both cases,
   $\varphi$ does not extend $D_4(a_1)$ to $D_5(a_1)$.

   (iii)  By (i) and (ii) any extension of $D_4(a_1)$ exists only for the
   case where one of the dipoles is a $4$-index dipole, and another one
   is a $2$-index dipole. Since there is no extension $D_4(a_1) < D_5(a_1)$
   in any vertex lying on a $2$-index dipole, there are no adjacent extensions.
   \qed

 \begin{proposition}[On adjacent extensions in $E_l$]
   \label{prop_D4a1_toExt}

    {\rm (i)} Let $\Gamma = D_4(a_1)$ be a \underline{$4$-cycle in $E_l$}, where $l = 6,7,8$.
       For any extension $\Gamma < \Gamma_1$ $($$\Gamma_1$ is a diagram of type $D_5(a_1)$$)$,
       there exists an element $T \in W$ such that the extension $\Gamma_2 = T\Gamma_1$ constitutes the adjacent
       extension $\Gamma < \Gamma_2$.

    {\rm (ii)}Let $\Gamma = D_4(a_1)$ be a \underline{$4$-cycle in $E_l$}, where $l = 6,7,8.$
    Let $\Gamma < \Gamma_1$ and $\Gamma < \Gamma_2$ be two adjacent
    extensions, where $\Gamma_1$ and $\Gamma_2$ are diagrams of type
    $D_5(a_1)$ (diagrams $\Gamma_1$ and $\Gamma_2$ are depicted in $(a)$ and $(e)$
    in Fig. $\ref{D4a1_toExt}$, respectively).
    Then any $\Gamma_1$-associated element $w_1$ and
    any $\Gamma_2$-associated element $w_2$ are conjugate.
\begin{figure}[h]
\centering
\includegraphics[scale=0.9]{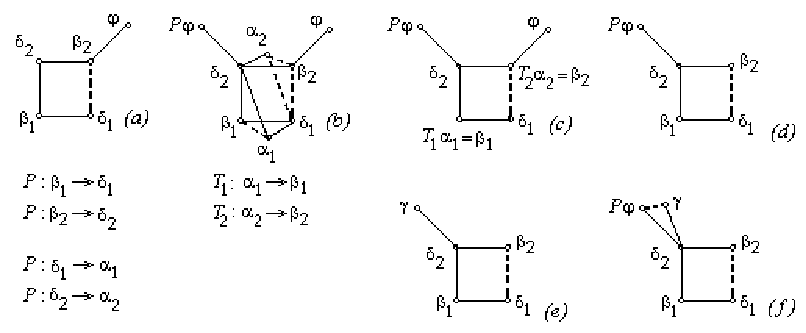}
 \caption{Adjacent extensions $(a)$ and $(d)$ of type $D_4(a_1) <
 D_5(a_1)$}
\label{D4a1_toExt}
\end{figure}

   {\rm (iii)} Let $\Gamma = D_l(a_1)$ be a  \underline{$4$-cycle in $D_l$ or $E_l$},
      let $\Gamma < \Gamma_1$ and $\Gamma < \Gamma_2$
      be two extensions in the opposite vertices of the same dipole
      (the diagrams $\Gamma_1$ and $\Gamma_2$ are depicted in $(a)$ and $(c)$ in Fig. $\ref{D5a1_1_mirror_ext}$, respectively).
      Any $\Gamma_1$-associated element $w_1$ and
      any $\Gamma_2$-associated element $w_2$ are conjugate.

\begin{figure}[h]
\centering
\includegraphics[scale=0.9]{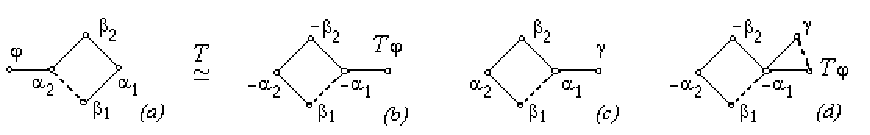}
 \caption{Mirror extensions of type $D_4(a_1) < D_5(a_1)$}
\label{D5a1_1_mirror_ext}
\end{figure}
   \end{proposition}

\index{adjacent extensions $D_4(a_1) < D_5(a_1)$}

   \PerfProof (i) For any extension $\Gamma < \Gamma_1$, we construct adjacent extension
    $\Gamma < \Gamma_2$, where $\Gamma_1$ and $\Gamma_2$ are diagrams of type $D_5(a_1)$.
    Let $\Gamma_1$ be an extension of $D_4(a_1)$ with an additional edge
    $\{\varphi, \beta_2\}$, see Fig.  \ref{D4a1_toExt}$(a)$.
    By Corollary \ref{cor_orth_roots}, for $E_6, E_7, E_8$, all dipoles are equivalent.
    Then there exists $P \in W$ mapping the dipole
    $\{\beta_1, \beta_2\}$ into the dipole $\{\delta_1, \delta_2\}$.
    Suppose $P$ maps $\{\delta_1, \delta_2\}$ to some $\{\alpha_1,\alpha_2\}$. Thus,
   \begin{equation}
   \label{eq_D4a1_toExt_1}
     \begin{split}
      & P : \{\beta_1, \beta_2, \delta_1, \delta_2\} \longrightarrow \{\delta_1, \delta_2, \alpha_1,\alpha_2\}, \\
      & P : \{\beta_2, \varphi\}  \longrightarrow  \{\delta_2,
      P\varphi\},
     \end{split}
   \end{equation}
   see Fig. \ref{D4a1_toExt}$(a)$ and $(b)$.  By Lemma \ref{lem_glue_1},
   $\alpha_1$ is connected to $\beta_1$,  and $\alpha_2$ is connected to $\beta_2$; if $\alpha_1$ is connected to
   $\beta_2$ and $\alpha_2$ is connected to $\beta_1$, we just swap $\alpha_1$ and $\alpha_2$ in Fig. \ref{D4a1_toExt}$(b)$.
   Note that if $\alpha_1$ is connected to $\beta_1$ and also to $\beta_2$, then
   $\alpha_1$ is connected to all vertices of the $4$-cycle $\{\beta_1, \beta_2, \delta_1, \delta_2\}$ contradicting
   Corollary \ref{cor_numb_ep}.  Further, by Corollary
   \ref{cor_map_2_diagr}, the corrective reflection $T_1 = s_{\alpha_1 + \beta_1}$ (resp. $T_2 = s_{\alpha_2 + \beta_2}$)
   maps $\alpha_1$ to $\beta_1$ (resp. $\alpha_2$ to $\beta_2$). Both corrective reflections preserve
   $\delta_1$, $\delta_2$ and $P\varphi$. Then
   \begin{equation}
    \begin{split}
     \label{eq_D4a1_toExt_2}
      & T_1{T}_2 : \{\delta_1, \delta_2, \alpha_1,\alpha_2\} \longrightarrow \{\delta_1, \delta_2, \beta_1, \beta_2\},\\
      & T_1{T}_2 : \{\delta_2, P\varphi\} \longrightarrow  \{\delta_2, P\varphi\},
    \end{split}
   \end{equation}
   i.e., $T_1{T}_2$ preserves the edge $\{\delta_2, P\varphi\}$. By \eqref{eq_D4a1_toExt_1} and \eqref{eq_D4a1_toExt_2}
   \begin{equation*}
    \begin{split}
      & T_1{T}_2{P} : \{\beta_1, \beta_2, \delta_1, \delta_2\} \longrightarrow \{\delta_1, \delta_2, \beta_1, \beta_2\}, \\
      & T_1{T}_2{P} : \{\beta_2, \varphi\}  \longrightarrow \{\delta_2, P\varphi\},
    \end{split}
   \end{equation*}
   i.e., $T = T_1{T}_2{P}$ preserves the $4$-cycle
   $\{\beta_1, \beta_2, \delta_1, \delta_2\}$ and maps
   extension $(a)$ to $(d)$ in Fig.  \ref{D4a1_toExt}.
~\\

   (ii) Consider the $\Gamma_1$-associated element $w_1 = s_{\beta_1}s_{\beta_2}s_{\delta_1}s_{\delta_2}s_{\varphi}$,
   and the $\Gamma_2$-associated element $w_2 = s_{\beta_1}s_{\beta_2}s_{\delta_1}s_{\delta_2}s_{P\varphi}$,
   see Fig. \ref{D4a1_toExt}$(a)$,$(e)$.
   By (i), we have
   \begin{equation*}
      ({T_1{T}_2P})^{-1}w_1({T_1{T}_2P}) = w'_1, \text{ where } w'_1 = s_{\delta_1}s_{\delta_2}s_{\beta_1}s_{\beta_2}s_{P\varphi}.
   \end{equation*}
   Since $s_{\beta_1}s_{\beta_2}$ commutes with $s_{P\varphi}$, we have
   \begin{equation*}
       w_1 \simeq w'_1 \stackrel{s_{\beta_1}s_{\beta_2}}{\simeq} w_2,
   \end{equation*}
   i.e., the Mirror Condition holds. As in Proposition \ref{prop_induct_step}, the Single-track Condition holds
   for the extension $D_4(a_1) < D_5(a_1)$. By Proposition \ref{prop_extensions}(iii)
   the $\Gamma_1$-associated elements and $\Gamma_2$-associated elements constitute the same conjugacy class.
~\\
\begin{minipage}{1cm}
 ~\\
\end{minipage}
\begin{minipage}{14.5cm}
 \Small
 ~\\
 ~\\
   The Single-track Condition for the extension $D_4(a_1) < D_5(a_1)$ is proved as follows.
   If $P\varphi$ and $\gamma$ are not connected and $P\varphi$ is
   linearly independent of $S = \{\beta_1, \beta_2, \delta_2, \gamma\}$,
   it follows that $S$ extended by $P\varphi$ constitutes $\widetilde{D}_4$, which is a contradiction.
   If $P\varphi$ and $\gamma$ are not connected, and $P\varphi$ is
   linearly dependent on $S$, then by Remark \ref{rem_max_root}(i), the root $-P\varphi$ is the maximal root
   in the root subsystem $S$:
   \begin{equation*}
      -P\varphi = 2\delta_2 + \beta_1 + \beta_2 + \gamma \quad \text{ and } \quad
      w'_1 = s_{\delta_1}s_{\delta_2}s_{\beta_1}s_{\beta_2}s_{2\delta_2 + \beta_1 + \beta_2 + \gamma}.
   \end{equation*}
   Let $w_3 :=  s_{\beta_1}s_{\beta_2}s_{\delta_1}s_{\delta_2}s_{\gamma}$.
   Since $s_{\delta_1}$ commutes with $s_{P\varphi}$, we have
   \begin{equation*}
     w'_1 ~\stackrel{s_{\delta_1}s_{\delta_2}}{\simeq}~
        s_{\beta_1}s_{\beta_2}s_{\delta_1}s_{\delta_2}s_{\delta_2 + \beta_1 + \beta_2 + \gamma}
       ~\stackrel{s_{\beta_1}s_{\beta_2}}{\simeq}~ 
       s_{\delta_1}s_{\delta_2}s_{\beta_1}s_{\beta_2}s_{\delta_2 + \gamma}
       ~\stackrel{s_{\delta_1}s_{\delta_2}}{\simeq}~
       s_{\beta_1}s_{\beta_2}s_{\delta_1}s_{\delta_2}s_{\gamma} = w_3.
   \end{equation*}
   If $P\varphi$ and $\gamma$ are connected,
   then by Corollary \ref{cor_map_2_diagr} we also have $w'_1 \simeq w_3$,
    see Fig. \ref{D4a1_toExt}$(f)$.
\end{minipage}
  ~\\
  ~\\

 (iii) We argue as in (ii). Let $w_1$ (resp. $w_2$) be the $\Gamma_1$-associated (resp. $\Gamma_2$-associated)
    element.  By \S\ref{sect_T_mirror}, $T^{-1}w_1{T} = w'_1$, where
   \begin{equation*}
    \begin{split}
      & w_1 = s_{\beta_1}s_{\beta_2}s_{\alpha_1}s_{\alpha_2}s_{\varphi},
      \quad w'_1 = s_{\beta_1}s_{\beta_2}s_{\alpha_1}s_{\alpha_2}s_{T\varphi},
      \quad  T = (s_{\beta_2}s_{\alpha_1}s_{\alpha_2})^2, \\
    \end{split}
   \end{equation*}
 i.e., the Mirror Condition holds. Further reasoning is as in (ii).
   \qed

\subsection{Threefold extensions $D_4 < D_5$}

\begin{lemma}
  There are no mirror or threefold extensions $D_4 < D_5$ in the root system $D_n$.
\end{lemma}
\begin{figure}[h]
\centering
\includegraphics[scale=0.9]{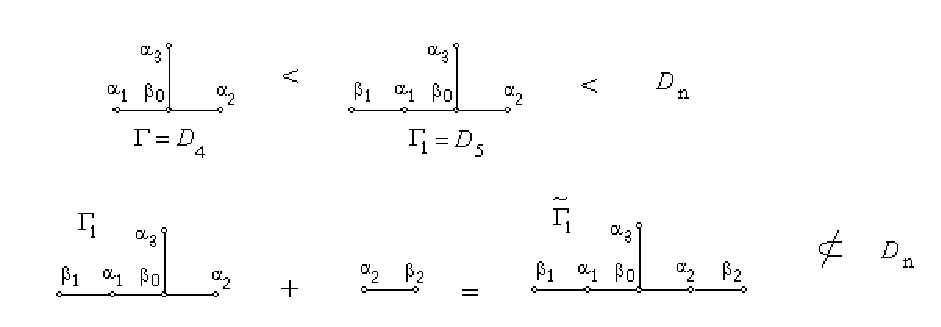}
 \caption{No extension $\Gamma_1 < \widetilde{\Gamma}_1$ is possible
  since $D_n$ does not contain $E_6$}
\label{G1G2_impos}
\end{figure}

\PerfProof  Suppose that $D_4 < D_5$ is a certain $\mathsf{P}1$-extension,
 i.e., there exists a $D_5$-associated subset in the root system $D_n$ for $n \geq 5$.
 We will show that the extension $D_5 < \widetilde{\Gamma}_1$
 obtained by adding some root $\beta_2$ and the edge $\{\alpha_2, \beta_2\}$, as in Fig. \ref{G1G2_impos},
 is impossible. This means that there are no simultaneous extensions of $\Gamma = D_4$
 by $\{\beta_1, \alpha_1\}$ and $\{\alpha_2, \beta_2\}$, i.e., there are no mirror or threefold extensions
 $D_4 < D_5$ in $D_n$.
~\\

 \underline{Case 1: $\beta_1$ and $\beta_2$ are not connected.}
 Suppose a certain root $\beta_2$ is connected with $\alpha_2$ as in Fig. \ref{G1G2_impos},
 and is linearly dependent on the $D_5$-associated subset
 $S_1 = \{\alpha_1, \alpha_2, \alpha_3, \beta_0, \beta_1\}$.
 By Remark \ref{rem_max_root}(ii) the diagonal element $b^{\vee}_{\alpha_2, \alpha_2}$ of the inverse
 matrix $B^{-1}_{\Gamma_1}$ should be equal to $2$. Actually, $b^{\vee}_{\alpha_2, \alpha_2} = \displaystyle\frac{5}{4}$.

 Suppose now that $\beta_2$ is connected with $\alpha_2$,
 and $\beta_2$ is linearly independent of $S_1$. Then the set $S' = \{S, \beta_2\}$ constitutes
 an $E_6$-associated subset. Why the root system $E_6$ is not contained in $D_n$?
 Roots of the root systems $E_6$ and $D_n$ are as follows:
 \begin{equation}
 \begin{split}
   & D_n \quad
    \begin{cases}
      \pm{e_i}  \pm{e_j} \quad (1 \le  i  < j  \le l),
    \end{cases} \\
   \\
   & E_6 \quad
    \begin{cases}
    \begin{split}
     & \pm{e_i}  \pm{e_j} \quad (1 \le  i  < j  \le 5), \\
     \\
     & \pm{\frac{1}{2}} \left ( e_8 - e_7 - e_6 + \sum\limits_{i=1}^5(-1)^{\nu(i)}e_i \right ),
     \text{ where } \sum\limits_{i=1}^5\nu(i) \text{ is even, }
    \end{split}
   \end{cases} \\
  \end{split}
 \end{equation}
 see \cite[Tables IV and V]{Bo}.
 It is clear that some roots of $E_6$ cannot be obtained as roots of $D_n$.

\begin{figure}[h]
\centering
\includegraphics[scale=0.9]{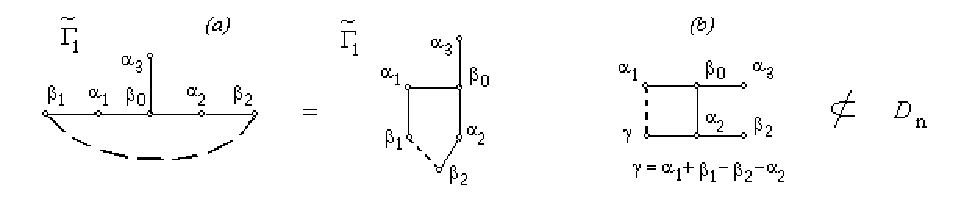}
 \caption{No extension $\Gamma_1 < \widetilde{\Gamma}_1$ is possible:
  $D_n$ does not contain $E_6(a_1)$}
\label{G1G2_impos_2}
\end{figure}

 \underline{Case 2: $\beta_1$ and $\beta_2$ are connected, see Fig. \ref{G1G2_impos_2}}. Suppose $\beta_2$ is
 linearly dependent on a certain $D_5$-associated subset $S_1$. Consider $\delta = B^{-1}_{\Gamma_1}\beta^{\vee}_2$,
 where $\beta^{\vee}_2$ is the vector obtained by a doubling the inner products
 given as in \eqref{eq_vect_inner_prod}, and $B^{-1}_{\Gamma_1}$
 is given in Table \ref{tab_Cartan_E6_E7_E8}:
\begin{equation}
 \label{eq_b2_2cases}
  \beta^{\vee}_2 =
     \left (
     \begin{array}{c}
         0 \\
         -1 \\
         0 \\
         0 \\
         1 \\
     \end{array}
   \right ), \quad
  \delta =
   \frac{1}{4}\left (
     \begin{array}{c}
         0 \\
         -3 \\
         -1 \\
         -2 \\
         2 \\
     \end{array}
   \right ) \quad
   \text{ or } \quad
  \beta^{\vee}_2 =
      \left (
     \begin{array}{c}
         0 \\
         -1 \\
         0 \\
         0 \\
         -1 \\
     \end{array}
   \right ), \quad
  \delta =
   \frac{1}{4}\left (
     \begin{array}{c}
         -8 \\
         -7 \\
         -5 \\
         -10 \\
         -6 \\
     \end{array}
   \right ).
\end{equation}
~\\
 According to \eqref{eq_vect_inner_prod_2},
 in the first case of \eqref{eq_b2_2cases}, we have
 $\mathscr{B}_{\Gamma_1}(\beta_2) = \langle \beta^{\vee}_2, \delta \rangle = \displaystyle\frac{5}{4}$
 and, in the second one, $\mathscr{B}_{\Gamma_1}(\beta_2) = \displaystyle\frac{13}{4}$, contradicting the condition
 $\mathscr{B}_{\Gamma_1}(\beta_2) = 2$, see Remark \ref{rem_max_root}(ii).

 Now, suppose that $\beta_2$ is linearly independent of $S_1$, i.e.,
 the subset $S_2 = \{\alpha_1, \alpha_2, \alpha_3, \beta_0, \beta_1, \beta_2\}$ is linearly independent.
 If we replace $\beta_1$ with $\gamma = \alpha_1 + \beta_1 - \beta_2 - \alpha_2$, we get again a set
 of linearly independent roots, which constitutes an $E_6(a_1)$-associated subset,
 see Fig. \ref{G1G2_impos_2}$(b)$, since
\begin{equation*}
\Small
 \begin{split}
   & (\gamma, \beta_0) = (\alpha_1, \beta_0) - (\alpha_2, \beta_0) = 0, \quad
    (\gamma, \alpha_1) = (\alpha_1, \alpha_1) + (\beta_1, \alpha_1) = \frac{1}{2}, \quad
    (\gamma, \alpha_3) = 0, \\
   & (\gamma, \alpha_2) = -(\alpha_2, \alpha_2) - (\alpha_2, \beta_2) = -\frac{1}{2}, \quad
     (\gamma, \beta_2) = -(\beta_2, \beta_2) + (\beta_1, \beta_2) - (\alpha_2, \beta_2) = 0. \\
 \end{split}
\end{equation*}
 However, by Lemma \ref{lem_adjacent}(iii), the root system $D_n$ does not contain $E_6(a_1)$-associated subsets.
  \qed

\begin{figure}[h]
\includegraphics[scale=0.8]{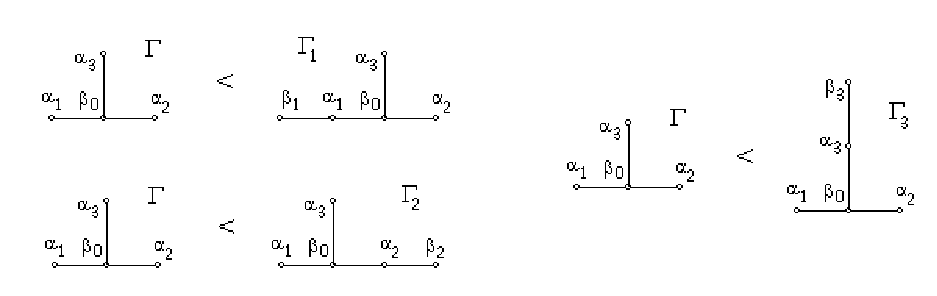}
\caption{Threefold extensions: $\Gamma < \Gamma_1$, ~$\Gamma < \Gamma_2$, ~$\Gamma < \Gamma_3$}
\label{triality_extens}
\end{figure}

  Unlike the case $D_n$,
  for the root systems $E_n$, there exist mirror and threefold extensions $D_4 < D_5$.
  In the following proposition we show that these extensions are conjugate.

 \begin{proposition}
   Let $\Gamma = D_4$ represent the conjugacy class in $W(E_n)$, where $n = 6,7,8$.
   Let $\Gamma < \Gamma_i$  be a $\mathsf{P}1$-extension with socket $\alpha_i$, where $i= 1,2,3$.
   Then any $\Gamma_i$-associated element $w_i$ and any $\Gamma_j$-associated element $w_j$
   are conjugate.
 \end{proposition}

\begin{figure}[h]
\includegraphics[scale=0.8]{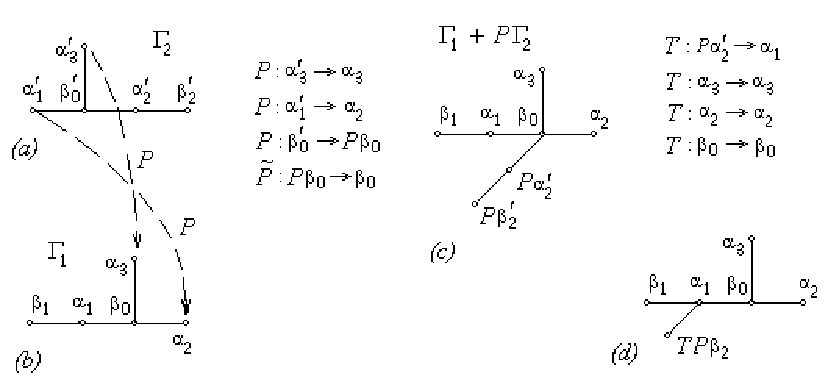}
\caption{Conjugacy of $\Gamma_1$-associated element $w_1$ and $\Gamma_2$-associated element $w_2$}
\label{triality_G1G2_conj}
\end{figure}

 \PerfProof  We will show that a given $\Gamma_1$-associated element $w_1$ and a given $\Gamma_2$-associated element $w_2$
  are conjugate, see Fig. \ref{triality_G1G2_conj}. By Corollary \ref{cor_orth_roots}, for
  $E_6$, $E_7$, $E_8$, all dipoles are equivalent. Then there exists an element $P \in W(E_n)$
  mapping the dipole $\{\alpha'_1, \alpha'_3\}$  to the dipole
  $\{\alpha_2, \alpha_3\}$, see Fig. \ref{triality_G1G2_conj}$(a)$,$(b)$.
  By Lemma \ref{lem_united_3} (see also Fig. \ref{same_diagon_2}), there exists
  a certain $\widetilde{P} \in W(E_n)$ preserving the dipole $\{\alpha_2, \alpha_3\}$ and
  sending $P\beta_0 \in \Gamma_2$ to $\beta_0 \in \Gamma_2$.  Applying $\widetilde{P}P$
  we get the diagram depicted in Fig. \ref{triality_G1G2_conj}$(c)$.
  By Lemma \ref{lem_united_4}, there exists $T \in W(E_n)$ preserving $\alpha_2$, $\alpha_3$,
  $\beta_0$ and sending $P\alpha'_2$ into $\alpha_1$, see Fig. \ref{triality_G1G2_conj}$(c)$,$(d)$.
  Thus, the Threefold Condition is proved: For a given $\Gamma_2$-associated subset
  $\{\alpha'_1, \alpha'_2, \alpha'_3, \beta'_0, \beta'_2\}$,
  the $\Gamma_1$-associated subset $\{\alpha_1, \alpha_2, \alpha_3, \beta_0, TP\beta_2\}$
  is constructed, see  Fig. \ref{triality_G1G2_conj}$(a)$,$(d)$.

\begin{minipage}{1cm}
 ~\\
\end{minipage}
\begin{minipage}{14.5cm}
 \Small
   The Single-track Condition for the extension $D_4 < D_5$ is proved as follows.
   Let $\gamma := TP\beta_2$. If $\gamma$ and $\beta_1$ are not connected, and $\gamma$ is
   linearly independent of the $D_5$-associated subset $S = \{\beta_0, \beta_1, \alpha_1, \alpha_2, \alpha_3\}$,
   then $S$ extended by $\gamma$ constitutes $\widetilde{D}_5$, but this is impossible.
   If $\gamma$ and $\beta_1$ are not connected and $\gamma$ is
   linearly dependent on $S$, then by Remark \ref{rem_max_root}(i), the root $-\gamma$
   is the maximal root in the root subsystem $D_5$.
   Let $w$ be the $D_4$-associated element, $w = s_{\alpha_1}s_{\alpha_2}s_{\alpha_3}s_{\beta_0}$.
   We need to prove that $ws_\gamma \simeq w\beta_1$.
   This follows from Proposition \ref{prop_conj_max_root}, eq. \eqref{eq_reduced_mu_max_2}.
   If $\gamma$ and $\beta_1$ are connected, then by Corollary \ref{cor_map_2_diagr}
   we also have $ws_\gamma \simeq w\beta_1$.
\end{minipage}
~\\
~\\
   Since the Threefold Condition and Single-track Condition hold
   for the extension $D_4 < D_5$, then by Proposition \ref{prop_extensions}(iv)
   the set of $\Gamma_1$-associated elements and the set of $\Gamma_2$-associated elements
   constitute the same conjugacy class.
  \qed

\appendix
\section{\sc\bf Cycles}
 \subsection{The ratio of lengths of roots}
   \label{sec_trees}
   \index{obtuse angle between roots}
  Let $\Gamma$ be a Dynkin diagram, and
  $\sqrt{t}$ be the ratio of the length of any long root to the length of any short root.
  The inner product between two long roots is
$$
   (\alpha, \beta) = \sqrt{t}\cdot\sqrt{t}\cdot\cos(\widehat{\alpha, \beta}) =
   \sqrt{t}\cdot\sqrt{t}\cdot \big ( \pm\frac{1}{2} \big ) =
   \pm\displaystyle\frac{t}{2}.
$$
 By Remark \ref{rem_tree}, we may put $(\alpha, \beta) = -\displaystyle\frac{t}{2}$.
The inner product between two short roots is
$$
   (\alpha, \beta) = \cos(\widehat{\alpha, \beta}) = \pm\displaystyle\frac{1}{2}.
$$
 Again, by Remark \ref{rem_tree}, we may put $(\alpha, \beta) =  -\displaystyle\frac{1}{2}$.
 The inner product $(\alpha, \beta)$ between roots of different lengths is
$$
   (\alpha, \beta) = 1\cdot\sqrt{t}\cdot\cos(\widehat{\alpha, \beta}) =
   1\cdot\sqrt{t}\cdot \big ( \pm\frac{\sqrt{t}}{2} \big ) =
   \pm\displaystyle\frac{t}{2}.
$$
As above, we choose the obtuse angle and put
  $(\alpha, \beta) = -\displaystyle\frac{t}{2}$.
~\\
We can summarize:
\begin{equation}
 \label{eq_all_inner_prod}
  (\alpha, \beta) =
  \begin{cases}
    -\frac{1}{2} \quad  \text{for} & \norm{\alpha} = \norm{\beta} = 1, \\
    \\
    -1 \quad  \text{for} & \norm{\alpha}  = \norm{\beta} =  2, \quad \text{or} \quad \norm{\alpha}  = 1, ~\norm{\beta} = 2, \\
    \\
    -\frac{3}{2} \quad  \text{for} & \norm{\alpha} = \norm{\beta} = 3, \quad \text{or} \quad \norm{\alpha}  = 1, ~\norm{\beta} = 3, \\
  \end{cases}
\end{equation}
where all angles $\widehat{\alpha, \beta}$ are obtuse.

\subsection{Cycles in the simply-laced case}

\subsubsection{The Carter and connection diagrams for trees}
  \label{Dynkin_diagr}

  \begin{lemma}
    \label{lem_must_dotted}
  There is no root subset (in the root system associated with a Dynkin diagram)
  forming a simply-laced cycle containing only solid edges.
  Every cycle in the Carter diagram or in the connection diagram contains at least one solid edge
  and at least one dotted edge.
  \end{lemma}
  \PerfProof Suppose a subset $S = \{ \alpha_1, \dots,  \alpha_n\} \subset \varPhi $
 forms a cycle containing only solid edges. Consider the vector
 \begin{equation*}
     v = \sum\limits_{i=1}^n \alpha_i.
 \end{equation*}
 The value of the quadratic Tits form $\mathscr{B}$ (see \cite{St08}) on $v$ is equal to
 \begin{equation*}
    \mathscr{B}(v) =  \sum\limits_{i \in \Gamma_0} 1 - \sum\limits_{i \in \Gamma_1} 1 =  n - n = 0,
 \end{equation*}
 where $\Gamma_0$ (resp. $\Gamma_1$) is the set of all vertices (resp. edges)
 of the diagram $\Gamma$ associated with $S$.
 Therefore, $v = 0$ and elements of the root subset $S$ are linearly dependent. \qed
 \\

   The following proposition is true only for trees.
\begin{proposition}[Lemma 8, \cite{Ca72}]
  \label{prop_non_ext_Dynkin}
  Let $S = \{ \alpha_1, \dots,  \alpha_n\} $ be a subset of linearly independent (not necessarily simple)
  roots of the root system $\varPhi$ associated with a certain Dynkin diagram $\Gamma$, and
  $\Gamma_S$  the Carter diagram or the connection diagram  associated with $S$.
  If $\Gamma_S$ is a tree, then $\Gamma_S$ is a Dynkin diagram.

\end{proposition}
\PerfProof
 If $\Gamma_S$ is not a Dynkin diagram, then $\Gamma_S$ contains an
 extended Dynkin diagram $\widetilde{\Gamma}$ as a subdiagram. Since $\Gamma_S$
 is a tree, we can turn all dotted edges to solid ones\footnotemark[1],
 see Remark \ref{rem_tree}.
\footnotetext[1]{This fact is not true for cycles, since by Lemma \ref{lem_must_dotted} we cannot
 eliminate all dotted edges.}
 Further, we consider the vector
 \begin{equation}
    \label{eq_lin_dep}
        v = \sum\limits_{i \in \widetilde{\Gamma}_0}{t_i\alpha_i},
 \end{equation}
 where $\widetilde{\Gamma}_0$ is the set of all vertices of $\widetilde{\Gamma}$, and $t_i$ (where $i \in \widetilde{\Gamma}_0$)
 are the coefficients of the nil-root, see \cite{Kac80}. Let the remaining coefficients corresponding
 to  $\Gamma_S \backslash  \widetilde{\Gamma}$ be equal to $0$.
 Let $\mathscr{B}$ be the positive definite quadratic Tits form (see \cite{St08}) associated with the diagram $\Gamma$,
 and $(\cdot\hspace{0.7mm},\cdot)$ the symmetric bilinear form associated with $\mathscr{B}$.
 Let  $\{\delta_i \mid i \in \widetilde{\Gamma}_0\}$
 be the set of simple roots associated with vertices $\widetilde{\Gamma}_0$.
 For all $i,j \in \widetilde{\Gamma}_0$, we have $(\alpha_i, \alpha_j) = (\delta_i, \delta_j)$,
 since this value is described by edges of $\widetilde{\Gamma}$. Therefore,
 \begin{equation*}
    \mathscr{B}(v) =
    \sum\limits_{i,j \in \widetilde{\Gamma}_0}t_i{t}_j(\alpha_i, \alpha_j) =
    \sum\limits_{i,j \in \widetilde{\Gamma}_0}t_i{t}_j(\delta_i, \delta_j) =
    \mathscr{B}(\sum\limits_{i \in \widetilde{\Gamma}_0}{t_i\delta_i}) = 0.
 \end{equation*}
 Since $\mathscr{B}$ is a positive definite form, we have $v = 0$, i.e.,
 vectors $\alpha_i$ are linearly dependent. This contradicts the definition of the set $S$.
 \qed

\begin{example}[multiply-laced cases]
{\rm
 On Figs. \ref{fig_F41_F42b}, \ref{fig_B2_C2b}, \ref{fig_B3_C3b} and \ref{fig_G21_G22b},
 let the coefficients of linear dependence be as in the proof of Proposition \ref{prop_non_ext_Dynkin}.
 The labels at vertices are coordinates of the nil-root of the corresponding
 extended Dynkin diagrams, \cite{Kac80}.
 In all cases below, inner products are calculated in accordance with \S\ref{sec_trees} and eq. \eqref{eq_all_inner_prod}.
 In all these calculations except for $\widetilde{G}_{21}$ and $\widetilde{G}_{22}$, $t = 2$.
 ~\\
 ~\\
\begin{minipage}{7cm}
 \qquad
  \epsfig{file=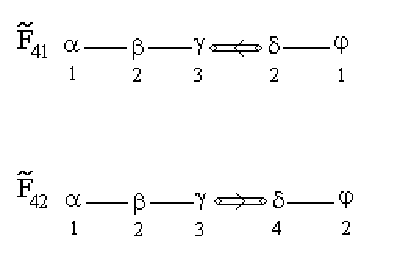, width=2.1in}
  \captionof{figure}{}
  \label{fig_F41_F42b}
\end{minipage}
 \quad
\begin{minipage}{8cm}
Case $\widetilde{F}_{41}$. We have $v = \alpha + 2\beta + 3\gamma +
2\delta + \varphi$. Then
 \begin{equation*}
   \begin{split}
   \norm{v}
    & = 1 + 4 + 9 + 4t + t - 1\cdot{2} - 2\cdot{3} - 3\cdot{2t} - 2\cdot{t} \\
    & = 6 - 3t = 0.
   \end{split}
 \end{equation*}
\vspace{1mm} \\
 Case $\widetilde{F}_{42}$. Here, $v = \alpha + 2\beta + 3\gamma +
4\delta + 2\varphi$, and
 \begin{equation*}
 \begin{split}
   \norm{v}
  & = (1 + 4 + 9)t + 16 + 4 - 1\cdot{2t} - 2\cdot{3t}  \\
  & - 3\cdot{4t} - 4\cdot{2} = 12 - 6t = 0.
 \end{split}
 \end{equation*}
 \vspace{2mm}
\end{minipage}
 ~\\
 ~\\
 ~\\
\begin{minipage}{10cm}
Case $\widetilde{C}_{2}$. We have $v = \alpha + t\beta + \gamma$,
where $t = 2$. Then
$$
  \norm{v}
   = t + t^2 + t - t\cdot{t} - t\cdot{t} = 2t - t^2
   = 0.
$$
\vspace{1mm} \\
Case $\widetilde{B}_{2}$.
 In this case, $v = \alpha + \beta + \gamma$. Then
 $$
  \norm{v}
 = 1 + 1 +  t - t - t =  2 - t = 0.
 $$
\end{minipage}
\begin{minipage}{6cm}
  \qquad \qquad
  \psfig{file=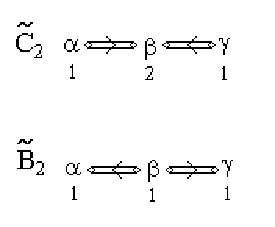, width=1.3in}
  \captionof{figure}{}
  \label{fig_B2_C2b}
\end{minipage}
 ~\\
 ~\\
 ~\\
\begin{minipage}{6cm}
  \qquad
  \psfig{file=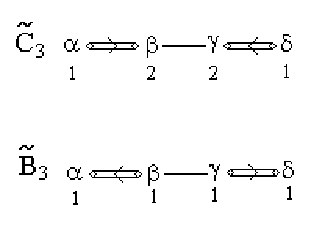, width=1.6in}
  \captionof{figure}{}
    \label{fig_B3_C3b}
\end{minipage}
 \qquad
\begin{minipage}{10cm}
Case $\widetilde{C}_{3}$. Here, $v = \alpha + t\beta + t\gamma +
\delta$, and
$$
  \norm{v}
   = t + t +  t^2 + t^2 - t^2 - t^2 - t^2 = 2t - t^2 = 0.
$$
\vspace{1mm} \\
 Case $\widetilde{B}_{3}$.
 We have
 $v = \alpha + \beta + \gamma + \delta$. Then
 $$
  \norm{v}
 = 1 + 1 +  t + t  - t - t - t =  2 - t = 0.
 $$
\end{minipage}
~\\
~\\
~\\
 For $\widetilde{C}_n$, where $n \geq 4$, we have
$v = \alpha_1 + t\alpha_2 + \dots + t\alpha_{n} + \alpha_{n+1}.$ Any
new short edge adds $t^2 - t^2$, i.e., $\norm{v} = 0$.
For $\widetilde{B}_n$, where $n \geq 4$, we have $v = \alpha_1 + \alpha_2 + \dots +
\alpha_{n} + \alpha_{n+1}.$ Any new long edge adds $t - t$, i.e., $\norm{v}$ = 0.
 ~\\
\begin{minipage}{10cm}
Case $\widetilde{G}_{21}$. Here, $v = \alpha + 2\beta + \gamma$. Here, $t = 3$.
 Then
 $$
  \norm{v}
 = 1 + 4 + t - 1\cdot{2} - 2\cdot{t} = 3 - t = 0.
 $$
\vspace{1mm} \\
Case $\widetilde{G}_{22}$. In this case, $v = \alpha + 2\beta + 3\gamma$.
Again, $t = 3$. Then
 $$
  \norm{v}
 = t + 4\cdot{t} + 9 - 2\cdot{t} - 2\cdot{3}\cdot{t} = 9 - 3t = 0. \qquad \qed
 $$
\end{minipage}
\begin{minipage}{6cm}
  \qquad  \qquad
  \psfig{file=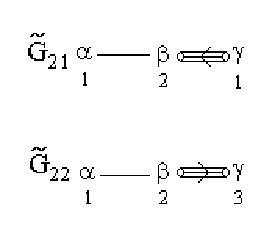, width=1.3in}
  \captionof{figure}{}
  \label{fig_G21_G22b}
\end{minipage}
 }
\end{example}
~\\

\subsubsection{There are no cycles for the root system $A_n$}
  \label{case_An}
  \index{similarity transformation for Carter diagrams}
   Recall that any root in $A_n$ is of the form $\pm(e_i - e_j)$, where $1 \leq i < j \leq n+1$.
 Then, up to the similarity $\alpha \longmapsto -\alpha$,
 a cycle of roots is of one of the followings forms:
 \begin{equation*}
  \begin{split}
   & \{ e_{i_1} - e_{i_2},  e_{i_2} - e_{i_3}, \dots,
    e_{i_{k-1}} - e_{i_k}, e_{i_{k}} - e_{i_1} \}, \\
   & \{ e_{i_1} - e_{i_2},  e_{i_2} - e_{i_3}, \dots,
    e_{i_{k-1}} - e_{i_k}, -(e_{i_{k}} - e_{i_1}) \}.
   \end{split}
 \end{equation*}
 In the first case, the sum of all these roots is equal to $0$, and roots are linearly dependent.
 In the second case, the sum of the $k-1$ first roots is equal to the last one, and
 roots are also linearly dependent. Thus, for $A_n$, there are no cycles of linearly independent roots.
~\\

\subsection{Cycles in the multiply-laced case}
  \label{mult_case}
  \index{obtuse angle between roots}

\subsubsection{There are no $4$-cycles with all angles obtuse}
 No root system $R$ containing a $4$-cycle with all angles obtuse
 can occur. Suppose this is possible, so a quadruple of roots
 $\{\alpha, \beta, \gamma, \delta\}$  yields pairs with the
 following values of the Tits form:

\begin{equation*}
   (\alpha, \beta) = -1, \quad (\beta, \gamma) =  -\frac{1}{2} ,
   \quad (\gamma, \delta) = -1, \quad (\delta, \alpha) = -1,
 \end{equation*}
 see Fig. \ref{fig_multiply_laced_sq}.  Then we see that

\begin{equation*}
 \begin{split}
  & s_\alpha = \left (
    \begin{array}{cccc}
      -1 &  1 &  0 &  1 \\
      0  &  1 & 0  & 0 \\
      0  &  0 & 1  & 0 \\
      0  &  0 & 0  & 1 \\
    \end{array} \right), \qquad
   s_\beta = \left (
    \begin{array}{cccc}
      1 &  0 &  0 &  0 \\
      2  &  -1 & 1  & 0 \\
      0  &  0 & 1  & 0 \\
      0  &  0 & 0  & 1 \\
    \end{array} \right), \\
      \\
  & s_\gamma = \left (
    \begin{array}{cccc}
      1 &  0 &  0 &  0 \\
      0  &  1 & 0  & 0 \\
      0  &  1 & -1  & 2 \\
      0  &  0 & 0  & 1
    \end{array} \right), \qquad
   s_\delta = \left (
    \begin{array}{cccc}
      1 &  0 &  0 &  0 \\
      0  &  1 & 0  & 0 \\
      0  &  0 & 1  & 0 \\
      1  &  0 & 1  & -1
    \end{array} \right).
   \end{split}
\end{equation*}

 \index{semi-Coxeter element}

   Then the semi-Coxeter element ${\bf C} = s_\alpha s_\beta s_\gamma s_\delta$
   in the Weyl group generated by the quadruple
   $\{s_\alpha, s_\beta, s_\gamma, s_\delta\}$, and its characteristic polynomial
   is as follows:
\begin{equation*}
 {\bf C}  = s_\alpha s_\beta s_\gamma s_\delta = \left (
    \begin{array}{cccc}
      4  &  0 &  2 &  -3 \\
      4  &  0 &  1 &  -2 \\
      2  &  1 &  1  & -2 \\
      1  &  0 &  1  & -1 \\
    \end{array} \right), \qquad
    \chi({\bf C})  = x^4 - 4x^3 - x^2 -4x + 1.
\end{equation*}
~\\
\begin{minipage}{8cm}
 Since the maximal root of $\chi({\bf C})$ is $\lambda \approx 4.419 > 1$,
 then the semi-Coxeter element ${\bf C}$ is of infinite order,
 but this is impossible.
\end{minipage}
\begin{minipage}{7cm}
  \qquad \qquad \psfig{file=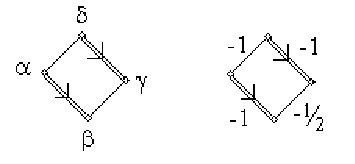, width=2.3in}
  \captionof{figure}{}
   \label{fig_multiply_laced_sq}
\end{minipage}

\subsubsection{More of impossible cases of multiply-laced cycles}

 \index{acute angle between roots}
 We consider several patterns (of multiply-laced diagrams) that
 are not a part of any Carter diagram. First of all, the
 arrows  on the double
~\\
\begin{minipage}{8cm} 
  edges connecting roots of
 different lengths should be directed face to face, otherwise we have
 $3$ different lengths of roots, as depicted in Fig. \ref{fig_3_diff_lengths}:
  $$
    \norm{\delta}  > \norm{\gamma}  =  \norm{\beta}  > \norm{\alpha}.
  $$
  ~\\
\end{minipage}
\begin{minipage}{8cm} 
  \qquad \qquad \qquad \qquad
  \psfig{file=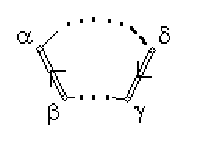, width=1.0in}
  \captionof{figure}{}
  \label{fig_3_diff_lengths}
\end{minipage}
~\\
~\\
\begin{minipage}{8cm} 
  \qquad \qquad
  \psfig{file=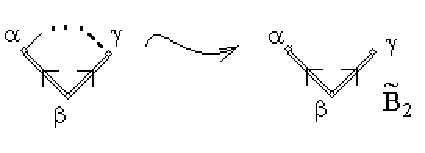, width=2in}
  \captionof{figure}{}
  \label{fig_1long_2short}
\end{minipage}
\begin{minipage}{8.3cm} 
 Further, two double edges connecting roots of different lengths
 cannot be adjacent, as depicted in Fig. \ref{fig_1long_2short}. Otherwise,
 the root subset contains the extended Dynkin diagram
 $\widetilde{B}_2$ or $\widetilde{C}_2$ which cannot occur.
\end{minipage}
~\\
\begin{minipage}{8cm}
 For cycles of length $5$ or more, the diagram contains
 the extended Dynkin diagram of type $\widetilde{B}_n$ or $\widetilde{C}_n$
 which cannot happen. If the acute angle (resp. the dotted edge)
 lies on the part corresponding to $\widetilde{B}_n$ (or $\widetilde{C}_n$), this
 obstacle can be easily eliminated by changing certain roots with
 their opposites; the procedure of eliminating the acute angle
 may be applied to any tree regardless of whether it
 contains roots of different lengths or not.
\end{minipage}
\begin{minipage}{7cm}
  \qquad \qquad  \psfig{file=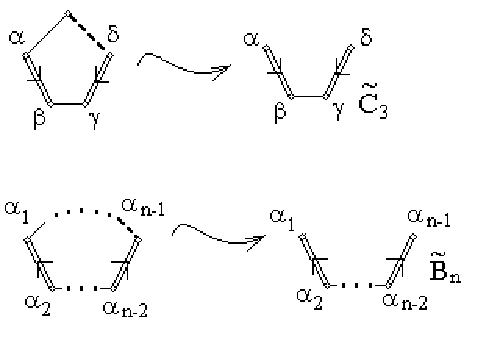, width=2.3in}
  \captionof{figure}{}
\end{minipage}
~\\
~\\
\begin{minipage}{8cm}
  \psfig{file=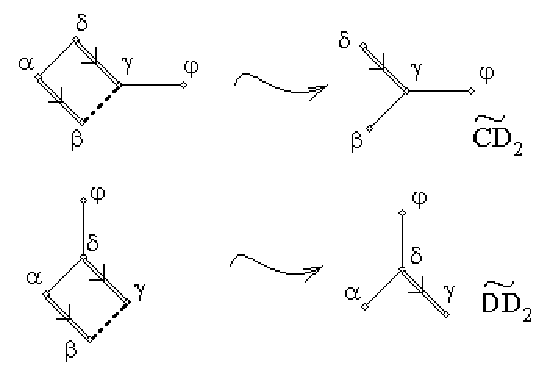, width=3in}
  \captionof{figure}{}
  \label{fig_5cycles_1}
\end{minipage}
\begin{minipage}{8.3cm}
 There are no \lq\lq{kites}\rq\rq, i.e., cycles of length $4$ with an additional fifth edge, since
 any such subset contains the extended Dynkin diagram $\widetilde{CD}_2$ or $\widetilde{DD}_2$
 which cannot be, see Fig. \ref{fig_5cycles_1}. One should note that every cycle in
 the Carter diagram contains, by definition, an even number of vertices,
 so the connection like $\{\varphi, \alpha \}$
 or $\{\varphi, \beta \}$  forming a triangle cannot occur.
 \end{minipage}
~\\

~\\
~\\
\section{\sc\bf The $4$-cycles in the root system $D_l$}

\subsection{The $D_4(a_1)$-associated and $D_4$-associated subsets}
  \label{sec_cycles_D4}

 In this section, we consider $D_4(a_1)$-associated subsets (i.e., $4$-cycles) and
 $D_4$-associated subsets for $W(D_l)$.

\subsubsection{The $2$-index and $4$-index dipoles in $4$-cycles}
  \label{sec_2_4_index_dipole}
  For definition of $2$-index and $4$-index dipoles, see \S\ref{sec_adjacent}.

\begin{lemma}
  \label{lem_eq_2diag}
    {\rm (i)} Let $\mathcal{C}$ be a $4$-cycle with a certain $2$-index (resp. $4$-index) dipole $d$.
    Let $w$ be any element from $W$.  Then $wd$ is also a $2$-index (resp. $4$-index) dipole.

    {\rm (ii)} If $d_1$ (resp. $d_2$) is a $2$-index dipole of some $4$-cycle $\mathcal{C}_1$
    (resp. $\mathcal{C}_2$), then there exists an element $w \in W$ such that $wd_1 = d_2$.

    {\rm (iii)} If $d_1$ (resp. $d_2$) is a $4$-index dipole from some $4$-cycle $\mathcal{C}_1$
    (resp. $\mathcal{C}_2$), then there exists an element $w \in W$ such that $wd_1 = d_2$.

    {\rm (iv)} For any $4$-cycle $\mathcal{C}$, one of dipoles is a $4$-index dipole.
       It can happen that both dipoles are $4$-index dipoles.

\end{lemma}

\begin{figure}[h]
 \centering
\includegraphics[scale=0.9]{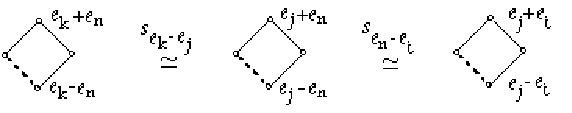}
 \caption{Conjugation of $2$-index dipoles (of $4$-cycles).}
\label{2index_diagonals}
\end{figure}

  \PerfProof
   (i) It suffices to prove this fact for reflections.
   Let $d = \{e_n + e_k, e_n - e_k \}$ be a $2$-index dipole.
   By \eqref{eq_ei_ej},
   any reflection looks as $s_\alpha$, where $\alpha$ is the root  $\pm e_i \pm e_j$ for some $e_i$, $e_j$.
   If $e_i, e_j \not\in \{e_n, e_k \}$, then $\alpha \perp \{e_n + e_k, e_n - e_k \}$ and
   $s_\alpha$ acts trivially on $d$.
   If $s_\alpha$ acts nontrivially on $e_n - e_k$, then $e_i$ or $e_j$ belongs to $\{ e_n, e_k \}$.
   Let $\alpha = e_k - e_j$. We have
   \begin{equation*}
     \begin{split}
       & (\alpha, e_n - e_k) = (e_k - e_j, e_n - e_k ) = -\frac{1}{2}, \\
       & (\alpha, e_n + e_k) = (e_k - e_j, e_n + e_k ) = \frac{1}{2}.
     \end{split}
   \end{equation*}
   Then
   \begin{equation*}
     \begin{split}
       & s_\alpha(e_n - e_k) = (e_n - e_k) + (e_k - e_j) = e_n - e_j, \\
       & s_\alpha(e_n + e_k) = (e_n + e_k) - (e_k - e_j) = e_n + e_j.
     \end{split}
   \end{equation*}
   Thus, we get the $2$-index dipole $\{ e_n + e_j, e_n - e_j \}$, see Fig. \ref{2index_diagonals}.

   We get the same result for $-\alpha = -e_k + e_j$.
   For $\alpha = e_k + e_j$, we get the $2$-index dipole $\{ e_n + e_j, e_n - e_j\}$.
   Let both $e_i$ and $e_j$ belong to $\{e_n, e_k \}$, for example, $\alpha = e_n - e_k$.
   Then $e_n - e_k \perp e_n + e_k$ and
   $s_\alpha(e_n - e_k) = e_k - e_n$, and $s_\alpha(e_n + e_k) = e_n + e_k$.
   We get the $2$-index dipole $\{ e_k - e_n, e_k + e_n\}$.
   For $\alpha = e_n + e_k$, we get the $2$-index dipole $\{e_n - e_k, -e_k - e_n\}$.

   It is clear that a $4$-index dipole cannot be sent to a $2$-index dipole by any $w \in W$,
   otherwise the inverse element $w^{-1}$ sends a $2$-index dipole to a $4$-index dipole,
   contradicting the proved above.

  (ii) This fact follows from Fig. \ref{2index_diagonals}. The reflection $s_{e_k - e_j}$ changes
  one index: $k \longmapsto j$, and the reflection $s_{e_n - e_i}$ changes another index:
  $n \longmapsto i$.

\begin{figure}[h]
 \centering
\includegraphics[scale=0.9]{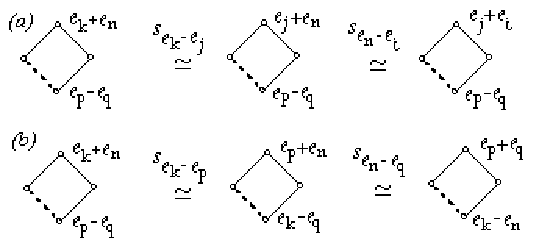}
 \caption{Conjugation of $4$-index dipoles (of $4$-cycles).}
\label{4index_diagonals}
\end{figure}

  (iii) Here, every reflection changes only one index and does not touch
  three other indices as in Fig. \ref{4index_diagonals}$(a)$, or permutes two indices
  and does not touch two other indices as in Fig. \ref{4index_diagonals}$(b)$.
  Consistently using these reflections we send $d_1$ into $d_2$.

  (iv) Suppose both dipoles $d_1 = \{ e_i + e_j, e_i - e_j \}$ and $d_2 = \{ e_k + e_n, e_k - e_n \}$ are  $2$-index dipoles
  in $\mathcal{C}$. Let $g$ be a number of coinciding indices from
  $\{k, n\}$ and  $\{i, j \}$, i.e.,
  \begin{equation*}
    g = | \{ k, n \} \cap \{ i, j \} |.
  \end{equation*}
  If $g = 0$, then $e_k + e_n \perp e_i + e_j$, contradicting the fact that
  the vertices of $d_1$ and $d_2$ are pairwise connected.
  If $g = 2$, then $d_1$ = $d_2$, which is impossible. So $g = 1$. For example, $e_n = e_j$, i.e.,
   \begin{equation*}
    d_1 = \{ \alpha_1 = e_i + e_j,  \alpha_2 = e_i - e_j \}, \quad
    d_2 = \{ \beta_1 = e_k + e_j, \beta_2 = e_k - e_j \}.
  \end{equation*}
  We have $\alpha_1 - \alpha_2 = \beta_1 - \beta_2$, contradicting the fact of linear independence of
  $\{ \alpha_1,  \alpha_2, \beta_1, \beta_2 \}$.  Therefore, one of dipoles is a $4$-index dipole.

  The following $4$-cycle gives an example of two $4$-index dipoles:
   \begin{equation*}
    d_1 = \{ \alpha_1 = e_1 + e_2,  \alpha_2 = e_3 + e_4 \}, \quad
    d_2 = \{ \beta_1 = e_3 - e_1, \beta_2 = -e_2 - e_4 \}. \qed
  \end{equation*}

\begin{corollary}
   \label{cor_conj_diag}
   {\rm (i)} Let $\mathcal{C}_1$ and $\mathcal{C}_2$ be $4$-cycles.
   There exist dipoles $d_1 \in \mathcal{C}_1$ and $d_2 \in \mathcal{C}_2$, and
   an element $w \in W$ such that $wd_1 = d_2$.

   {\rm (ii)}
   Let $\mathcal{C}_1$ and $\mathcal{C}_2$ be $D_4$-associated subsets.
   There exist dipoles $d_1 \in \mathcal{C}_1$ and $d_2 \in \mathcal{C}_2$, and
   an element $w \in W$ such that $wd_1 = d_2$.
\end{corollary}

\PerfProof
  (i) By Lemma \ref{lem_eq_2diag}(iv) there exist $4$-index dipoles $d_i \in \mathcal{C}_i$ for $i = 1,2$;
  by (iii) they are conjugate, i.e., there exists an element $w \in W$ such that $wd_1 = d_2$.

  (ii) It cannot be that all $3$ dipoles of $\mathcal{C}_1$ (resp. $\mathcal{C}_2$) are $2$-index
  dipoles. Thus, there exist $4$-index dipoles $d_i \in \mathcal{C}_i$ for $i = 1,2$. Further, as in (i).
\qed

\subsubsection{Example of equivalent $4$-cycles}
  \label{sec_example_4cycles}
 Let us take the diagram $D_6$ with simple roots:
 \begin{equation}
   \begin{split}
     & \alpha_1 = e_1 - e_2, \quad \alpha_2 = e_2 - e_3, \quad \alpha_3 = e_3 - e_4, \\
     & \alpha_4 = e_4 - e_5, \quad \alpha_5 = e_5 - e_6, \quad \alpha_2 = e_5 + e_6,
   \end{split}
 \end{equation}
 and the maximal root $\alpha = e_1 + e_2$, see Fig. \ref{example_D6_2sq}$(a)$.
 Consider the following two $4$-cycles:
 \begin{equation}
   \begin{split}
     & \mathcal{C}_1 = \{ e_1 + e_2,  e_4 - e_1, e_1 - e_2, e_2 - e_3  \}, \\
     & \mathcal{C}_2 = \{ e_1 + e_2,  e_4 - e_1, e_3 - e_4, e_2 - e_3  \}.
   \end{split}
 \end{equation}

 Since the $4$-index dipole can be mapped only onto
 a $4$-index dipole, see Lemma \ref{lem_eq_2diag}, and $\mathcal{C}_2$ consists of two $4$-index dipoles,
 it follows that $\mathcal{C}_1$ and $\mathcal{C}_2$ can not be conjugate.
 However, $\mathcal{C}_1$ and $\mathcal{C}_2$
 have the \underline{common second dipole}  $\{ e_4 - e_1, e_2 - e_3\}$.
 Moreover, let $w_1$ (resp. $w_2$) be $\mathcal{C}_1$-associated (resp. $\mathcal{C}_2$-associated):
 \begin{equation}
   \begin{split}
     & w_1 = s_{e_1 + e_2}s_{e_1 - e_2}s_{e_4 - e_1}s_{e_2 - e_3}, \\
     & w_2 = s_{e_1 + e_2}s_{e_3 - e_4}s_{e_4 - e_1}s_{e_2 - e_3},
   \end{split}
 \end{equation}
then elements $w_1$ and $w_2$ are conjugate:
 \begin{equation}
   \begin{split}
     w_1 = & s_{e_1 + e_2}s_{e_1 - e_2}s_{e_4 - e_1}s_{e_2 - e_3} =
             s_{e_1 + e_2}(s_{e_1 - e_2}s_{e_2 - e_3})s_{e_4 - e_1} = \\
           & s_{e_1 + e_2}s_{e_2 - e_3}s_{(e_2 - e_3) + (e_1 - e_2)}s_{e_4 - e_1} =
             s_{e_1 + e_2}s_{e_2 - e_3}(s_{e_1 - e_3}s_{e_4 - e_1}) = \\
           & s_{e_1 + e_2}s_{e_2 - e_3}s_{e_4 - e_1}s_{(e_4 - e_1) + (e_1 - e_3)} =
             s_{e_1 + e_2}s_{e_2 - e_3}s_{e_4 - e_1}s_{e_4 - e_3} \stackrel{s_{e_4 - e_3}}{\simeq}   \\
           & s_{e_1 + e_2}s_{e_3 - e_4}s_{e_4 - e_1}s_{e_2 - e_3} = w_2.  \qed
   \end{split}
 \end{equation}

\subsection{The Carter diagrams corresponding to different conjugacy classes}

\subsubsection{Intersection of $E$-supports}
 \label{sec_inters_supports}

  The simple roots in $D_l$ are as follows:
 \begin{equation*}
     \alpha_i = e_i - e_{i+1}, \text{ for } 1 \leq i \leq {l-1},  \quad \alpha_l = e_{l-1} + e_l,
 \end{equation*}
 any other roots are of the form
  \begin{equation*}
     \beta = e_k \pm e_n, \text{ for  different } i \neq j,
 \end{equation*}
 see \cite[Table IV]{Bo}.
 The set of the vectors $e_i$ entering any root $\beta$ is said to be the {\it $E$-support}
 of $\beta$. We denote this set by $E(\beta)$. The $E$-supports for the simple and arbitrary roots are as follows:
\begin{equation*}
  \begin{split}
    &  E(\alpha_i) = \{e_i, e_{i+1} \} \text{ for } 1 \leq i \leq {l-2},  \\
    &  E(\alpha_l) = E(\alpha_{l-1}) = \{e_{l-1}, e_l \}, \\
    &  E(\beta) =  \{e_k, e_n\}, \text{ where } \beta = e_k - e_n.
  \end{split}
\end{equation*}
For each pair of roots $\alpha$ and $\beta$, we set
\begin{equation}
  \label{eq_def_inters_supp}
    E(\alpha, \beta) :=  E(\alpha) \cap E(\beta).
\end{equation}

 Note that $s_{e_p - e_q}$  acts non-trivially on the pair of indices
 $\{k, n \}$ of the root $e_k - e_n$ if and only if
 $|\{p, q \} \cap \{k, n \}| = 1$. If $|\{p, q \} \cap \{k, n \}| = 1$, let,
 for example, $p = k$. Then $s_{e_p - e_q}$ acts on $\{k, n \}$
 by replacing $k$ with $q$:
\begin{equation}
 \label{act_e_supp_1}
    s_{e_p - e_q}(\{k, n \}) =
        \begin{cases}
            \{k, n \}  \text{ if } |\{p, q \} \cap \{k, n \}| = 0 \text{ or } 2, \\
            \{q, n \} \text{ if }  |\{p, q \} \cap \{k, n \}| = 1 ~(\text{for } k = p), \\
        \end{cases}
\end{equation}
  For the action \eqref{act_e_supp_1} on the $E$-support $\{k, n \}$, there is
  the corresponding action by conjugation on the reflections $s_{e_k \pm e_n}$:
\begin{equation}
 \label{act_e_supp_1A}
      s_{e_p - e_q}s_{e_k \pm e_n}s^{-1}_{e_p - e_q} =
        \begin{cases}
            s_{e_k \pm e_n}  \text{ if } |\{p, q \} \cap \{k, n \}| = 0 \text{ or } 2, \\
            s_{e_q \pm e_n}  \text{ if } |\{p, q \} \cap \{k, n \}| = 1 ~(\text{for } k = p).
        \end{cases}
\end{equation}

\begin{lemma}
 \label{lem_inters_E}
 Let $\gamma_k = e_{i_k} - e_{j_k}$ (or $e_{i_k} + e_{j_k}$), where $k = 1,2,3$,
 be $3$ linearly independent roots in the root system $D_l$, and
    \begin{equation}
        \label{inters_E_supp}
       E(\gamma_1, \gamma_2, \gamma_3) =  E({\gamma_1}) \cap E({\gamma_2}) \cap E({\gamma_3}).
    \end{equation}
 The cardinality $|E(\gamma_1, \gamma_2, \gamma_3)|$ is invariant under the action \eqref{act_e_supp_1},
 and therefore, under the conjugation \eqref{act_e_supp_1A}.
\end{lemma}
 \PerfProof Since roots $\gamma_1$, $\gamma_2$, $\gamma_3$ are
 linearly independent, we have
$$
  |E(\gamma_1, \gamma_2, \gamma_3)| < 2.
$$
Let $|E(\gamma_1, \gamma_2, \gamma_3)| = 1$; then the $E$-supports are
\begin{equation*}
  E(\gamma_1) = \{i, j_1 \}, \quad
  E(\gamma_2) = \{i, j_2 \}, \quad
  E(\gamma_3) = \{i, j_3 \}.
\end{equation*}
 By \eqref{act_e_supp_1}, the action $s_\beta = s_{e_p - e_i}$ moves
 the $E$-supports of $\gamma_1$, $\gamma_2$, $\gamma_3$ as follows:
\begin{equation}
 \label{act_e_supp_2}
   \{ \{i, j_1 \},  \{i, j_2 \}, \{i, j_3 \} \} \longmapsto
        \begin{cases}
            \{\{i, j_1\}, \{j_1, j_2\}, \{j_1, j_3\}\} & \text{ if } s_\beta = s_{e_{j_1} - e_i}, \\
            \{\{p, j_1\}, \{p, j_2\}, \{p, j_3\}\} & \text{ if } s_\beta = s_{e_p - e_i}, p \not\in \{j_1, j_2, j_3\}\\
            \{\{i, p\},   \{i, j_2\}, \{i, j_3\}\} & \text{ if } s_\beta = s_{e_p - e_{j_1}}, p \neq i \\
            \{\{i, j_1\}, \{i, j_2\}, \{i, j_3\}\} & \text{ if } s_\beta = s_{e_p - e_q},
                p, q \not\in \{i, j_1, j_2, j_3 \}. \\
        \end{cases}
\end{equation}
 In all cases of \eqref{act_e_supp_2}, $|E(s_\beta\gamma_1, s_\beta\gamma_2, s_\beta\gamma_3)| = 1$.
 Let $|E(\gamma_1, \gamma_2, \gamma_3)| = 0$, then the $E$-supports are
\begin{equation}
 \label{act_e_supp_3}
 \begin{split}
   &  \{ \{i, j \},  \{j, k \}, \{k, i \} \}, i.e.,
     E(\gamma_f) \cap E(\gamma_g) \neq 0, \text{ where } 1 \leq f < g \leq 3,
   \text{ or } \\
   & \{ \{i, j \},  \{k, n \}, \{k, r\} \}, i.e.,
     E(\gamma_1) \cap E(\gamma_g) = 0,  E(\gamma_2) \cap E(\gamma_3) \neq 0,
      \text{ where } g = 2,3,
     \text{ or } \\
   &  \{ \{i, j \},  \{k, n \}, \{t, r\} \}, i.e.,
     E(\gamma_f) \cap E(\gamma_g) = 0, \text{ where } 1 \leq f < g \leq 3. \\
 \end{split}
\end{equation}
 Consider the action by $s_\beta = s_{e_p - e_i}$. For $p \not\in \{j,k,r,t \}$,
 the reflection $s_\beta$ only replaces $i$ with $p$. For $p \in \{j,k,r,t \}$,
 for example, $p = k$, the reflection $s_\beta$ swaps $i$ and $k$. In all cases of
 \eqref{act_e_supp_3}, $|E(s_\beta\gamma_1, s_\beta\gamma_2, s_\beta\gamma_3)|$ is not changed.
\qed

\subsubsection{The Carter diagram $A_3$ determines different conjugacy classes in $W(D_l)$}
  \label{sec_A3}

\begin{proposition}
  \label{prop_non_conj_1}
    Let $D$ be a Dynkin diagram  from Table $\ref{tab_non_conj_elem}$,
    let $\Gamma$ be the Carter diagram associated with a root subset in $\varPhi(D)$,
    see columns $2$ and $3$. The two elements $w_1$, $w_2$  corresponding to $\Gamma$
    and given by column $4$ are not conjugate.
\end{proposition}

\begin{table} 
  \centering
  \renewcommand{\arraystretch}{1.9}
  \begin{tabular} {|c|c|c|c|}
  \hline
        & Dynkin      & Carter diagram $\Gamma$         & Non-conjugate \cr
        & diagram $D$ & in the root system $\varPhi(D)$ & elements $w_1$ and $w_2$ \\
  \hline
     1 & $\begin{array}{c} \includegraphics[scale=1.4]{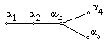}  \\ D_5  \end{array}$
       & $\begin{array}{c} \includegraphics[scale=1.4]{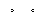}  \\ A_1 + A_1  \end{array}$
       & $\begin{array}{c} w_1 = s_{\alpha_1}{s}_{\alpha_3}   \\ w_2 = s_{\alpha_4}{s}_{\alpha_5} \end{array}$
       \\
  \hline
     2 & $\begin{array}{c} \includegraphics[scale=1.4]{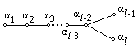}  \\ D_l  \end{array}$
       & $\begin{array}{c} \includegraphics[scale=1.4]{A1_p_A1.eps}  \\ A_1 + A_1  \end{array}$
       & $\begin{array}{c} w_1 = s_{\alpha_1}{s}_{\alpha_3}   \\ w_2 = s_{\alpha_{l-1}}s_{\alpha_l} \end{array}$ \\
  \hline
     3 & $\begin{array}{c} \includegraphics[scale=1.4]{Dl_3_2.eps}  \\ D_l  \end{array}$
       & $\begin{array}{c} \includegraphics[scale=1.4]{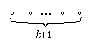}  \\ (k+1)A_1  \end{array}$
       & $\begin{array}{c} w_1 = s_{\alpha_1}s_{\alpha_3}\dots{s}_{\alpha_{2k+1}}   \\
                           w_2 = s_{\alpha_{l-1}}s_{\alpha_l}s_{\alpha_1}s_{\alpha_3}\dots{s}_{\alpha_{2k-3}} \end{array}$ \\
  \hline
     4 & $\begin{array}{c} \includegraphics[scale=1.4]{Dl_3_2.eps}  \\ D_l  \end{array}$
       & $\begin{array}{c} \includegraphics[scale=1.4]{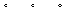}  \\ A_3  \end{array}$
       & $\begin{array}{c} w_1 = s_{\alpha_1}s_{\alpha_3}s_{\alpha_2}   \\
                           w_2 = s_{\alpha_{l-1}}s_{\alpha_l}s_{\alpha_{l-2}} \end{array}$ \\
  \hline
\end{tabular}
  \vspace{2mm}
  \caption{Non-conjugate elements $w_1$ and $w_2$ corresponding to the Carter diagram $\Gamma$}
  \label{tab_non_conj_elem}
\end{table}

\PerfProof  \underline{Line 1) in Table \ref{tab_non_conj_elem}}.
  Here,  $D = D_5$, and  $W = W(D_5)$;
  $\Gamma = A_1 + A_1$ is the Carter diagram consisting of two not connected points;
  $w_1 = s_{\alpha_1}s_{\alpha_3}$ and  $w_2 = s_{\alpha_4}s_{\alpha_5}$
  the two elements corresponding to $\Gamma$. By \eqref{eq_def_inters_supp}, we have
\begin{equation}
  \begin{split}
    \label{inters_are_invar}
    & E(\alpha_1, \alpha_3) = \{e_1, e_2 \} \cap \{e_3, e_4\}, \qquad \abs{E(\alpha_1, \alpha_3)} = 0, \\
    & E(\alpha_4, \alpha_5) = \{e_4, e_5 \} \cap \{e_4, e_5\},  \qquad \abs{E(\alpha_4, \alpha_5)} = 2.
  \end{split}
\end{equation}
~\\
 To prove the statement, it suffices to prove that
 cardinalities $|E_(\alpha_1, \alpha_3)|$ and $|E_{45}|$ from
 eq. \eqref{inters_are_invar} do not change under conjugation by any
 reflection. Let us conjugate $w_2$ by means of $s_{e_i - e_j}$, where $i,j \not\in  \{4, 5\}$, then
\begin{equation*}
  s_{e_i - e_j}w_2{s}_{e_i - e_j} = s_{e_i - e_j}s_{e_4 - e_5}s_{e_4 + e_5}{s}_{e_i - e_j} = s_{e_4 - e_5}s_{e_4 + e_5} = w_2
\end{equation*}
 If, for example, $i = 4$ , we have
\begin{equation*}
  s_{e_4 - e_j}w_2{s}_{e_4 - e_j} = s_{e_4 - e_j}s_{e_4 - e_5}s_{e_4 + e_5}{s}_{e_4 - e_j} = s_{e_j - e_5}s_{e_j + e_5}
\end{equation*}
 Thus, the intersection $E(\alpha_4, \alpha_5)$ turns into $\{e_j, e_5 \} \cap \{e_j, e_5\}$ with the same cardinality $2$.

  \underline{Line 2) in  Table \ref{tab_non_conj_elem}}.
  Here, $D = D_l$, and $W = W(D_l)$;  $\Gamma = A_1 + A_1$,
  $w_1 = s_{\alpha_1}s_{\alpha_3}$, $w_2 = s_{\alpha_{l-1}}s_{\alpha_l}$.

In this case,
\begin{equation*}
    w_1 = s_{e_1 - e_2}s_{e_3 - e_4}, \qquad w_2 = s_{e_{l-1} - e_l}s_{e_{l-1} + e_l}.
\end{equation*}
The proof is the same as for \underline{Line 1)}.

 \underline{Line 3) in  Table \ref{tab_non_conj_elem}}.
 Here, $D = D_l$, and $W = W(D_l)$;  $\Gamma = (k+1)A_1$
 is the Carter diagram consisting of $(k+1)$ not connected points, where $2k+1 \leq l-2$.
 Let $w_1 = s_{\alpha_1}s_{\alpha_3}\dots{s}_{\alpha_{2k+1}}$ and
      $w_2 = s_{\alpha_{l-1}}s_{\alpha_l}s_{\alpha_1}s_{\alpha_3}\dots{s}_{\alpha_{2k-3}}$
 be two elements corresponding to $\Gamma$.
 We have
\begin{equation*}
  \begin{split}
  \Small
   & w_1 = s_{e_1 - e_2}s_{e_3 - e_4}\dots{s}_{e_{2k-1} - e_{2k}}{s}_{e_{2k+1} - e_{2k+2}}, \\
   & w_2 = s_{e_{l-1} - e_l}s_{e_{l-1} + e_l}s_{e_1 - e_2}s_{e_3 - e_4}\dots{s}_{e_{2k-3} - e_{2k-2}}.
  \end{split}
\end{equation*}
As above, we put
\begin{equation*}
  \begin{split}
  \Small
 E(\alpha_i, \alpha_{i+2}) = E(\alpha_i) \cap E(\alpha_{i+2}) &= \{ e_i, e_{i-1} \} \cap \{ e_{i+2}, e_{i+1}\}
 \text{ for } i+2 \leq 2k+1, \\
 E(\alpha_{l-1},\alpha_l) = E(\alpha_{l-1}) \cap E(\alpha_{l}) &=  \{ e_i, e_{i-1} \}. \\
  \end{split}
\end{equation*}
Therefore,
\begin{equation*}
 \abs{E(\alpha_i, \alpha_{i+2})} = 0  \text{ and } \abs{E(\alpha_{l-1},\alpha_l)} = 2.
\end{equation*}

 As for \underline{Line 1)},  any conjugation of $w_2$ preserves the
$\abs{E(\alpha_{l-1},\alpha_l)}$ and it cannot be changed to $0$.
~\\

 \index{ambiguity of conjugacy classes}
 \index{conjugacy class! - ambiguity}

 \underline{Line 4) in  Table \ref{tab_non_conj_elem}}.
 Consider the following two elements in the Weyl group $D_l$:
\begin{equation}
  \label{eq_wi_A3}
  \begin{split}
     w_1 = s_{\alpha_1}s_{\alpha_3} & s_{\alpha_2}, \quad
     w_2 = s_{\alpha_{l-1}}s_{\alpha_l}s_{\alpha_{l-2}}, \text{ where } \\
    & \alpha_i = e_i - e_{i+1} \text{ for } 1 \leq i \leq l-1, \text{ and } \alpha_l = e_{l-1} + e_l,
  \end{split}
\end{equation}
see Fig. \ref{Dl_w1_w2}.
The elements $w_1$ and $w_2$ are $A_3$-associated. From \eqref{eq_wi_A3} we get
\begin{equation*}
   |E(\alpha_1, \alpha_2, \alpha_3)| = 0, \qquad
   |E(\alpha_{l-1}, \alpha_l, \alpha_{l-2})| = 1.
\end{equation*}
By Lemma \ref{lem_inters_E}, any conjugation of $w_1$ and $w_2$ preserves
intersections of $E$-supports $|E(\alpha_1, \alpha_2, \alpha_3)|$ and $|E(\alpha_{l-1}, \alpha_l, \alpha_{l-2})|$.
Thus $w_1$ and $w_2$ are not conjugate.
\qed

\begin{remark}
  {\rm
  Note that the $A_3$-associated subset of roots $\{ \alpha_{l-2}, \alpha_{l-1}, \alpha_l \}$
  cannot be extended to any $A_4$-associated subset.
  Suppose $\{\alpha_{l-2}, \alpha_{l-1}, \alpha_l, \beta\}$ constitutes an $A_4$-associated subset.
  Then the support $E(\beta)$ does not contain $e_l$ and $e_{l-1}$, otherwise
  we get a cycle,  see \eqref{eq_wi_A3}. Therefore, $\beta$ is not connected to
  $\alpha_{l-1}$ and $\alpha_l$.
  }
\end{remark}

\subsubsection{The Carter diagram $A_1 + A_1$ determines several conjugacy classes in $W(D_l)$}
   \label{sec_A1A1}

\index{connected simply-laced Carter diagrams}

 In this paper we consider only \underline{connected} simply-laced Carter diagrams except for
 the diagram $A_1 + A_1$. This exception is very important to us for the following reason.
 The diagram $A_1 + A_1$ is the essential subdiagram (dipole) of the two diagrams, from which
 we start our considerations. These diagrams are $D_4(a_1)$ and $D_4$, see Fig. \ref{tetris}.
 The diagrams $D_4(a_1)$ and $D_4$ lie in the base of the
 two non-intersecting chains of Carter diagrams: $\mathsf{C4}$ and $\mathsf{DE4}$, see Fig. \ref{diagram_tree};
 $\mathsf{C4}$ consists of the Carter diagrams $E_l(a_i)$ and $D_l(a_i)$;
 $\mathsf{DE4}$ consists of the Carter diagrams $E_l$ for $n=6,7,8$, and $D_l$ for $l \geq 4$.

 In the case $D_4(a_1)$, both dipoles are diagrams $A_1 + A_1$; in the case $D_4$,
 the three dipoles are diagrams of type $A_1 + A_1$. {\it Despite the fact that
 $A_1 + A_1$ determines several conjugacy classes in $W(D_l)$, the diagrams
 $D_4(a_1)$ and $D_4$  (which contain $A_1 + A_1$) determine only one conjugacy classes in $W(D_l)$}. For $W(E_l)$,
 the situation is simpler: $A_1 + A_1$ determines a single conjugacy class, and also
 $D_4(a_1)$ and $D_4$ determine only one conjugacy classes in $W(E_l)$.
 These facts constitute the base of induction in the proof of the main theorem,
 see Theorem \ref{th_uniq_diagr}, stating that every connected Carter diagram containing
 $D_4(a_1)$ and $D_4$ determines a single conjugacy class in a given Weyl group $W$.
 The base case is proved is \S\ref{sect_base_case}.

Consider the following two elements in $W(D_l)$, see Fig. \ref{Dl_w1_w2}:
\begin{equation}
  \begin{split}
  \label{eq_wi_A1A1}
     & w_1 = s_{\alpha_1}s_{\alpha_3}, \text{ where }
               \alpha_1 = e_1 - e_2, \quad \alpha_3 = e_3 - e_4, \\
     & w_2 = s_{\alpha_{l-1}}s_{\alpha_l},  \text{ where }
               \alpha_{l-1} = e_{l-1} - e_l,  \quad \alpha_l = e_{l-1} + e_l. \\
  \end{split}
\end{equation}
  The elements $w_1$ and $w_2$ are associated
 with the Carter diagram $\Gamma = A_1 + A_1$.
 From \eqref{eq_wi_A1A1} we deduce that
\begin{equation*}
   |E(\alpha_1, \alpha_3)| = 0, \quad |E(\alpha_{l-1}, \alpha_l)| = 2.
\end{equation*}
 As above, in the proof of Proposition \ref{prop_non_conj_1} (see \underline{Line 3)})
 any conjugation of $w_1$ preserves the cardinality of $E(\alpha_1, \alpha_3)$ and it cannot be changed to $2$.
 Therefore, $w_1$ and $w_2$ are not conjugate.

~\\
\section{\sc\bf Diagonal elements}

\subsection{The diagonal elements of  $B^{-1}_{\Gamma}$ for $A_l$, $D_l$, $D_l(a_k)$}
  \index{$b^{\vee}_{\eta, \eta}$ (diagonal element of $B^{-1}_{\Gamma}$)}

\begin{figure}[h]
\centering
\includegraphics[scale=0.5]{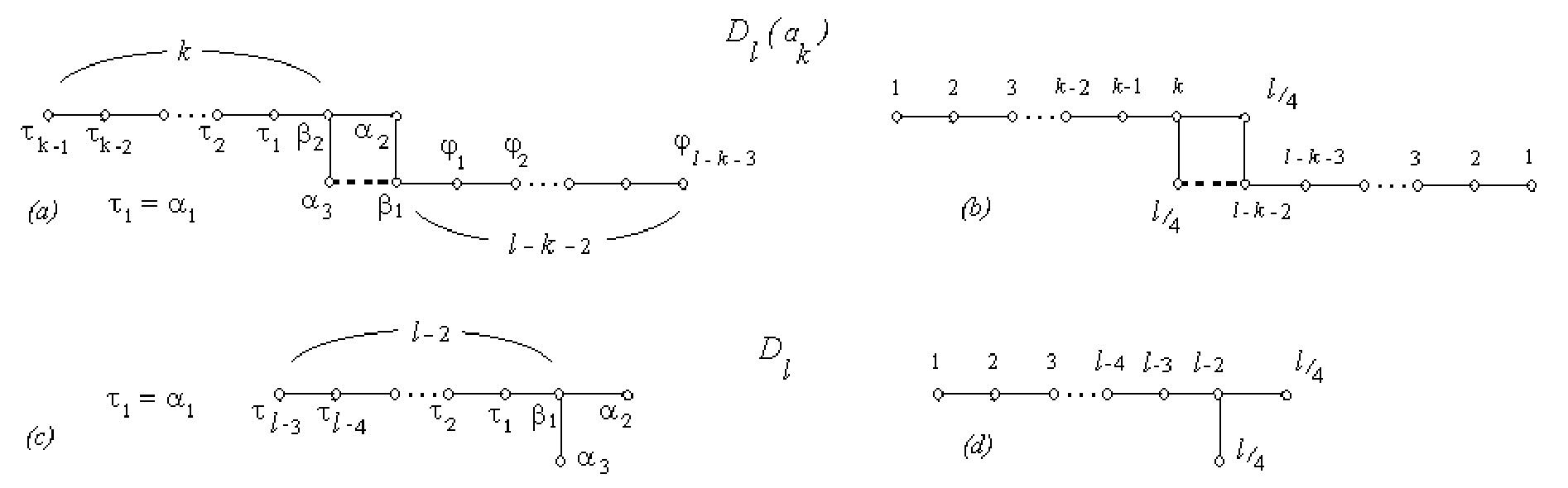}
\caption{The numerical labels in the diagrams in the right hand side
 are the diagonal elements of $B^{-1}_{\Gamma}$}
\label{Dl_ak_Carter_diagr_in_prop}
\end{figure}

 Let $\eta$ be one of vertices of $\Gamma = D_l(a_k)$ (resp. $\Gamma = D_l$) such that $\eta \neq \alpha_2, \alpha_3$,
 see Fig. \ref{Dl_ak_Carter_diagr_in_prop},
 \begin{equation}
        \label{diag_elem_1}
        \eta \in
        \begin{cases}
           \{ \tau_1, \dots, \tau_{k-1}, \varphi_1, \dots, \varphi_{l-k-3},
          \alpha_2, \alpha_3, \beta_1, \beta_2\}  & \text{ for } \Gamma = D_l(a_k),  \\
           \{ \tau_1, \dots, \tau_{l-3},
          \alpha_2, \alpha_3, \beta_1, \beta_2\}  & \text{ for } \Gamma = D_l,
        \end{cases}
 \end{equation}
 Removing the vertex $\eta$ with its bonds from $\Gamma$ we get the diagram $\Gamma'$ which is decomposed,
 except for $\eta = \tau_{k-1}$ and $\eta = \varphi_{l-k-3}$,
 into the  union of two connected subdiagrams:
  \begin{equation}
    \label{eq_one_of_diagr}
      \Gamma' =
        \begin{cases}
           A_{d-1} \oplus D_{l-p-1},  & \text{ where } d = l-k-2 \text{~(resp. } d = k), \\
                                  & \text{ for } \eta = \beta_1 ~(\text{resp. } \eta = \beta_2), \\
           A_{d-1} \oplus D_{l-p-1(a_i)},  & \text{ where } d = k-i, \text{~(resp. } d = l-k-2-i), \\
                & \text{ for } \eta = \tau_i ~(1 \leq i \leq k-1) ~(\text{resp. }
                               \eta = \varphi_i ~(1 \leq i \leq l-k-3)).  \\
        \end{cases}
  \end{equation}
  For $\eta = \tau_{k-1}$ (resp. $\eta = \varphi_{l-k-3}$),
  we have $d = 1$ and $\Gamma' = D_{l-1}(a_{k-1})$ (resp. $\Gamma' = D_{l-1}(a_{l-k-3})$).
  In Fig. \ref{Dl_ak_Carter_diagr_in_prop}$(b)$ and $(d)$, we associate
  the numerical label $d$ with the corresponding vertex $\eta$.

 \begin{proposition}
   \label{prop_determinants}
    {\rm (i)} The determinant of the partial Cartan matrix $B_{\Gamma}$ is as follows:
  \begin{equation}
    \label{eq_det_AD}
      \det{B_{\Gamma}} =
        \begin{cases}
           l + 1 & \text{for } A_l, \text{ where } l \geq 2; \text{ here } B_{\Gamma} = {\bf B}, \\
           4     & \text{for } D_l \text{ and }  D_l(a_k), \text{ where } l \geq 4. \\
        \end{cases}
  \end{equation}

    {\rm  (ii)} Let $b^{\vee}_{\eta,\eta}$ be diagonal elements of $B_{\Gamma}^{-1}$,
     where $\eta$ is given by \eqref{diag_elem_1}.
     For $\Gamma = D_l(a_k)$ or $\Gamma = D_l$,
      \begin{equation}
        \label{bi_eq_1}
     b^{\vee}_{\eta,\eta} =
     \begin{cases}
         \displaystyle\frac{l}{4} \text{ for } \eta = \alpha_2, \text{ or } \eta = \alpha_3, \\
         d \text{ for } \eta \text{ given in }  \eqref{diag_elem_1},
     \end{cases}
     \end{equation}
 where $d$ is given in the vertex $\eta$
 in Fig. \ref{Dl_ak_Carter_diagr_in_prop}$(b)$ and $(d)$, and by the relation \eqref{eq_one_of_diagr}.

   {\rm  (iii)} For $\Gamma = A_l$, let $b^{\vee}_{\eta,\eta}$ be diagonal elements of $B_{\Gamma}^{-1}$, where
 \begin{equation}
   \label{seq_vetices}
   \eta \in
  \begin{cases}
     \alpha_1, \beta_1, \alpha_2, \beta_2, \dots, \alpha_n, \beta_n & \text{ for } l = 2n,  \\
     \alpha_1, \beta_1, \alpha_2, \beta_2, \dots, \alpha_n, \beta_n, \alpha_{n+1} & \text{ for } l = 2n + 1,
  \end{cases}
 \end{equation}
 see Fig. $\ref{fig_A_l}$. Then
   \begin{equation}
     \label{eq_diag_elem_0}
      b^{\vee}_{\eta,\eta} =  d - \frac{d^2}{l+1},
   \end{equation}
 where $d$ is the sequential number of the vertex in eq. \eqref{seq_vetices} or in Fig. $\ref{fig_A_l}$.

\begin{figure}[h]
\centering
\includegraphics[scale=0.8]{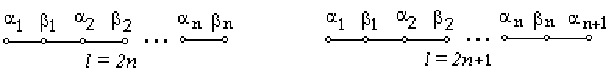}
\caption{The Dynkin diagrams $A_l$}
\label{fig_A_l}
\end{figure}

\end{proposition}

\PerfProof  (i) This statement is easily verified for $A_1$, $A_2$,
$D_4$, $D_5$ and $D_4(a_1)$.
   By induction
 \begin{equation*}
   \begin{split}
    & \det{{\bf B}(A_{l+1})} = 2\det{{\bf B}(A_{l})} -  \det{{\bf B}(A_{l-1})} =
        2(l+1) - l = l + 2, \\
    & \det{{\bf B}(D_{l+1})} = 2\det{{\bf B}(D_{l})} -  \det{{\bf B}(D_{l-1})} = 4,
   \end{split}
 \end{equation*}
 where  ${\bf B}(A_{l})$ (resp. ${\bf B}(D_{l})$)  is the Cartan matrix for $A_{l}$ (resp. $D_{l}$).

 To get the determinant of the partial Cartan matrix $B_{\Gamma}$ for the diagram $\Gamma = D_{l+1}(a_k)$,
 we expand $\det B_{\Gamma}(D_{l+1}(a_k))$ with respect
 to the minors corresponding to the $i$-th line (associated with the vertex $i$ of $\Gamma$).
 By induction, we have
 \begin{equation*}
   \begin{array}{lll}
    & \det B_{\Gamma}(D_{l+1}(a_k)) = 2\det B_{\Gamma}(D_{l}(a_{k-1})) -  \det B_{\Gamma}(D_{l-1}(a_{k-2})) = 4
       & \text{ for } k > 2, \\
    & \det B_{\Gamma}(D_{l+1}(a_{l - k - 1})) = 2\det B_{\Gamma}(D_{l}(a_{l - k - 2})) -
                                                   \det B_{\Gamma}(D_{l-1}(a_{l - k - 3})) = 4
       & \text{ for } l - k  > 3, \\
    & \det B_{\Gamma}(D_{l+1}(a_2)) = 2\det B_{\Gamma}(D_{l}(a_1)) -  \det {\bf B}(D_{l-1}) = 4
       & \text{ for } k = 2 \text{ or } l - k = 3.
    \end{array}
 \end{equation*}
 In the case $k < 2$ and $l - k  < 3$, we have $l < 5$, i.e., $\Gamma = D_4(a_1)$.
~\\

  (ii) For $\Gamma = D_l(a_k)$, we have
   \begin{equation}
    b^{\vee}_{\eta,\eta} =
    \begin{cases}
    ~ \displaystyle\frac{\det {B_{\Gamma}(D_{l}(a_k))}}{\det {B_{\Gamma}(D_{l+1}(a_k))}} = \frac{4}{4} = 1
        & \text{ for } $d = 1$, \\
    & \\
    ~ \displaystyle\frac{\det{\bf B}(A_{d-1}) \det B_{\Gamma}(D_{l-d}(a_k))}
                {\det B_{\Gamma}(D_{l+1}(a_k))} = \frac{d \cdot 4}{4} = d
        & \text{ for } d > 1, \\
    & \\
    ~  \displaystyle\frac{\det {{\bf B}(A_{l-1})}}{\det {B_{\Gamma}(D_{l+1}(a_k))}} = \frac{l}{4}
       & \text{ for } \eta = \alpha_2, \alpha_3.
     \end{cases}
   \end{equation}
   For $\Gamma = D_l$, we have
   \begin{equation}
    b^{\vee}_{\eta,\eta} =
    \begin{cases}
    ~ \displaystyle\frac{\det {{\bf B}(D_{l})}}{\det {{\bf B}(D_{l+1})}} = \frac{4}{4} = 1
        & \text{ for } $d = 1$, \\
    & \\
    ~ \displaystyle\frac{\det{\bf B}(A_{d-1}) \det {\bf B}(D_{l-d})}
                {\det {\bf B}(D_{l+1})} = \frac{d \cdot 4}{4} = d
        & \text{ for } d > 1, \\
    & \\
    ~  \displaystyle\frac{\det {{\bf B}(A_{l-1})}}{\det {{\bf B}(D_{l+1})}} = \frac{l}{4}
        & \text{ for } \eta = \alpha_2, \alpha_3.
     \end{cases}
   \end{equation}
~\\

  (iii) For $\Gamma = A_l$, we have:
  \begin{equation}
   \Small
        b^{\vee}_{\eta,\eta} = \frac{\det \overline{B}_{\eta,\eta}}{\det {\bf B}(A_l)},
  \end{equation}
  where  $\overline{B}_{\eta,\eta}$ is the matrix obtained from $B_{\Gamma}$ by deleting the $\eta$th column and $\eta$th row.
  The matrix $\overline{B}_{\eta,\eta}$ splits into the direct sum:
  \begin{equation*}
        \overline{B}_{\eta,\eta} = {\bf B}(A_{d-1}) \oplus {\bf B}(A_{l - d}).
  \end{equation*}
   Hence,
  \begin{equation}
    \label{eq_diag_elem_1}
        b^{\vee}_{\eta,\eta} = \frac{\det {\bf B}(A_{d-1}) \times \det {\bf B}(A_{l-d})}{\det {\bf B}(A_l)}.
  \end{equation}
  By \eqref{eq_det_AD} $\det {\bf B}(A_{l-1}) = l$ and by \eqref{eq_diag_elem_1} we have
  \begin{equation*}
        b^{\vee}_{\eta,\eta} = \frac{d \times (l + 1 -d)}{l + 1} =  d - \frac{d^2}{l+1}. \qed
  \end{equation*}
  ~\\

   For expressions ${\bf B^{-1}}$ for the Dynkin diagram $A_l$, see \cite[p. 295]{OV90}.
   For completeness sake, we give the partial  Cartan matrices and its inverse for some Cartan diagrams
   in Tables \ref{tab_Cartan_E6_E7_E8} and \ref{tab_partial Cartan_1}.
   For the Dynkin diagrams, the Cartan matrix ${\bf B}$ coincides with the partial Cartan matrix
   $B_{\Gamma}$.

 \begin{table}[h]
  \tiny
  \centering
  \renewcommand{\arraystretch}{1.5}
  \begin{tabular} {|c|c|c|}
  \hline
      The diagram $\Gamma$  & The Cartan matrix $B_{\Gamma}$   &   The inverse matrix $B^{-1}_{\Gamma}$    \cr
  \hline  
    & & \\
     $\begin{array}{c} \includegraphics[scale=0.4]{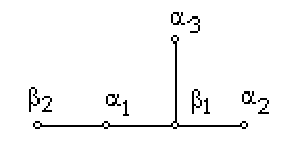} \\
        {\bf D_5} \end{array}$ &
     $\begin{array}{c}
       \left [
   \begin{array}{cccccc}
     2 & 0 & 0  &   -1 & -1  \\
     0 & 2 & 0  &   -1 & 0 \\
     0 & 0 & 2  &   -1 & 0  \\
    -1 & -1 & -1 &   2 & 0  \\
    -1 & 0 & 0 &    0 & 2 \\
  \end{array}
  \right ]
     \end{array}$   &
      $\begin{array}{c}
    \frac{1}{4} \left [ \begin{array}{cccccc}
     8 & 4 & 4 & 8 & 4 \\
     4 & 5 & 3 & 6 & 2 \\
     4 & 3 & 5 & 6 & 2 \\
     8 & 6 & 6 & 12 & 4 \\
     4 & 2 & 2 & 4 & 4 \\
  \end{array} \right ]
     \end{array}$    \\
    & & \\
  \hline  
    & & \\
     $\begin{array}{c} \includegraphics[scale=0.6]{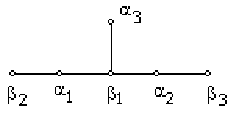} \\
        {\bf E_6} \end{array}$ &
     $\begin{array}{c}
       \left [
   \begin{array}{cccccc}
     2 & 0 & 0  &   -1 & -1 & 0  \\
     0 & 2 & 0  &   -1 &  0 & -1 \\
     0 & 0 & 2  &   -1 &  0 & 0  \\
    -1 & -1 & -1 &   2 & 0 & 0  \\
    -1 &  0 & 0  &   0 & 2 & 0 \\
      0 & -1 & 0 &   0 & 0 & 2 \\
  \end{array}
  \right ]
     \end{array}$   &
      $\begin{array}{c}
    \frac{1}{3} \left [ \begin{array}{cccccc}
     10 & 8 & 6 & 12 & 5 & 4  \\
     8 & 10 & 6 & 12 & 4 & 5 \\
     6 & 6 &  6  & 9 & 3 & 3  \\
    12 & 12 &  9 & 18 & 6 & 6 \\
     5 & 4 &  3  & 6 & 4 & 2 \\
     4 & 5 & 3 & 6 & 2 & 4 \\
  \end{array} \right ]
     \end{array}$    \\
    & & \\
     \hline
    & & \\
     $\begin{array}{c} \includegraphics[scale=0.6]{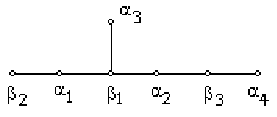} \\
        {\bf E_7} \end{array}$ &
     $\begin{array}{c}
       \left [
   \begin{array}{ccccccc}
     2 & 0 & 0  &  0 & -1 & -1 & 0  \\
     0 & 2 & 0  &  0 & -1 &  0 & -1 \\
     0 & 0 & 2  &  0 & -1 &  0 & 0  \\
     0 & 0 & 0  &  2 &  0 &  0 & -1 \\
    -1 & -1 & -1 & 0 &  2 &  0 & 0  \\
    -1 &  0 &  0 & 0 &  0 &  2 & 0 \\
     0 & -1 &  0 & -1 &  0 &  0 & 2 \\
  \end{array}
  \right ]
     \end{array}$   &
      $\begin{array}{c}
  \frac{1}{2} \left [ \begin{array}{ccccccc}
     12 & 12 & 8 & 4  & 16 & 6 & 8  \\
     12 & 15 & 9 & 5  & 18 & 6 & 10 \\
     8  & 9  & 7 & 3  & 12 & 4 & 6 \\
     4  & 5  & 3 & 3  & 6  & 2 & 4 \\
     16 & 18 & 12 & 6 & 24 & 8 & 12 \\
     6  & 6  & 4 & 2  & 8  & 4 & 4 \\
     8  & 10 & 6 & 4 & 12  & 4 & 8 \\
  \end{array} \right ]
     \end{array}$    \\
    & & \\
     \hline
    & & \\
     $\begin{array}{c} \includegraphics[scale=0.5]{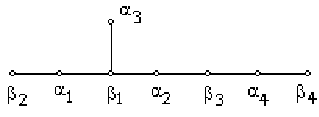} \\
        {\bf E_8} \end{array}$ &
     $\begin{array}{c}
       \left [
   \begin{array}{cccccccc}
     2 & 0  &  0  &  0  & -1 & -1 & 0  &  0 \\
     0 & 2  &  0  &  0  & -1 &  0 & -1 &  0\\
     0 & 0  &  2  &  0  & -1 &  0 & 0  &  0 \\
     0 & 0  &  0  &  2  &  0 &  0 & -1 &  -1\\
    -1 & -1 & -1  &  0  &  2 &  0 & 0  &  0\\
    -1 &  0 &  0  &  0  &  0 &  2 & 0  &  0 \\
     0 & -1 &  0  & -1  &  0 &  0 & 2  &  0 \\
     0 &  0 &  0  &  -1  &  0 &  0 & 0  &  2 \\
  \end{array}
  \right ]
     \end{array}$   &
      $\begin{array}{c}
     \left [ \begin{array}{cccccccc}
     14 & 16 & 10 & 8  & 20 & 7 & 12 & 4  \\
     16 & 20 & 12 & 10 & 24 & 8 & 15 & 5  \\
     10 & 12  & 8 & 6  & 15 & 5 & 9  & 3 \\
     8  & 10  & 6 & 6  & 12 & 4 & 8  & 3 \\
     20 & 24 & 15 & 12 & 30 & 10 & 18 & 6\\
     7  & 8  & 5 & 4  & 10  & 4 & 6  & 2 \\
     12 & 15 & 9 & 8  & 18  & 6 & 12 & 4 \\
     4  & 5 & 3  & 3  & 6  &  2 & 4  & 2 \\
  \end{array} \right ]
     \end{array}$    \\
    & & \\
     \hline
\end{tabular}
  \vspace{2mm}
  \caption{\Small The partial Cartan matrix $B_{\Gamma}$ and its inverse matrix $B_{\Gamma}^{-1}$}
  \label{tab_Cartan_E6_E7_E8}
\end{table}

 \begin{table}[h]
\tiny
  \centering
  \renewcommand{\arraystretch}{1.5}
  \begin{tabular} {|c|c|c|}
  \hline
      The Carter    & The partial Cartan         &   The inverse          \cr
   diagram $\Gamma$   &  matrix $B_{\Gamma}$       &   matrix $B^{-1}_{\Gamma}$    \\
  \hline  
    & & \\
     $\begin{array}{c} \includegraphics[scale=0.6]{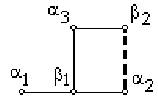} \\
        {\bf D_5(a_1) = D_5(a_2)} \end{array}$ &
     $\begin{array}{c}
       \left [
   \begin{array}{ccccc}
     2 & 0 &   0 & -1 & 0  \\
     0 & 2 &   0 & -1 & 1\\
     0 & 0 &   2 & -1 & -1  \\
     -1 & -1 &  -1 & 2 & 0  \\
      0 &  1 &   -1  & 0 & 2 \\
  \end{array}
  \right ]
     \end{array}$   &
      $\begin{array}{c}
    \frac{1}{4} \left [ \begin{array}{ccccc}
     4 & 2 &   2 & 4 & 0  \\
     2 & 5 &   1 & 4 & -2\\
     2 & 1 &   5  & 4 & 2  \\
     4 & 4 &   4 & 8 & 0  \\
     0 & -2 &  2  & 0 & 4 \\
  \end{array} \right ]
     \end{array}$    \\
    & & \\
    \hline
    & & \\
     $\begin{array}{c} \includegraphics[scale=0.6]{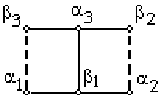} \\
        {\bf E_6(a_2)} \end{array}$ &
      $\begin{array}{c}
      \left [
   \begin{array}{cccccc}
     2 & 0 &   0 & -1 & 0  & 1 \\
     0 & 2 &   0 & -1 & 1  & 0\\
     0 & 0 &   2 & -1 & -1 & -1 \\
     -1 & -1 &  -1  & 2 & 0 & 0 \\
      0 &  1 &  -1  & 0 & 2 & 0 \\
      1 &  0 &  -1   & 0 & 0 & 2 \\
  \end{array}
  \right ]
     \end{array}$   &
           $\begin{array}{c}
      \frac{1}{3}\left [
   \begin{array}{cccccc}
     4 & 2 & 0 &  3 & -1 & -2 \\
     2 & 4 & 0 &  3 & -2 & -1 \\
     0 & 0 & 6 &  3 & 3 &  3 \\
     3 &  3  &  3 &  6 & 0 & 0 \\
     -1 & -2 &  3 &  0 & 4 & 2 \\
     -2 & -1 &  3 &  0 & 2 & 4 \\
  \end{array}
  \right ]
     \end{array}$  \\
    & & \\
    \hline  
     & & \\
     $\begin{array}{c} \includegraphics[scale=0.6]{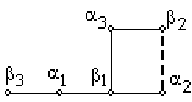} \\
        {\bf D_6(a_1) = D_6(a_3)} \end{array}$  &
     $\begin{array}{c}
       \left [
   \begin{array}{cccccc}
     2 & 0 &   0 & -1 & 0  & -1 \\
     0 & 2 &   0 & -1 & 1  & 0\\
     0 & 0 &   2 & -1 & -1 & 0 \\
     -1 & -1 &  -1  & 2 & 0 & 0 \\
      0 &  1 &  -1  & 0 & 2 & 0 \\
     -1 &  0 &  0   & 0 & 0 & 2 \\
  \end{array}
  \right ]
     \end{array}$   &
    $\begin{array}{c}
    \frac{1}{2} \left [
       \begin{array}{cccccc}
     4 & 2 &  2 &   4  & 0 & 2 \\
     2 & 3 &  1 &   3 & -1 & 1  \\
     2 & 1 &  3 &   3 &  1 & 1 \\
     4 & 3 &  3 &   6 &  0 & 2 \\
     0 &  -1 & 1 &  0 &  2 & 0 \\
     2 &  1 &  1 &  2 &  0 & 2 \\
  \end{array}
 \right ]
     \end{array}$    \\
    & & \\
     \hline
     & & \\
     $\begin{array}{c} \includegraphics[scale=0.6]{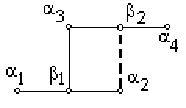} \\
        {\bf D_6(a_2)} \end{array}$  &
     $\begin{array}{c}
       \left [
   \begin{array}{cccccc}
      2 & 0  &   0  &  0  & -1  & 0 \\
      0 & 2  &   0  &  0  & -1  & 1 \\
      0 & 0  &   2  &  0  & -1  & -1 \\
      0 & 0  &   0  &  2  &  0  & -1 \\
     -1 & -1 &  -1  &  0  &  2  & 0 \\
      0 &  1 &  -1  & -1  &  0  & 2 \\
  \end{array}
  \right ]
     \end{array}$   &
    $\begin{array}{c}
    \frac{1}{2} \left [
       \begin{array}{cccccc}
     2 & 1 &  1  &  0  & 2 & 0 \\
     1 & 3 &  0  &  -1  & 2 & -2  \\
     1 & 0 &  3 &   1 &  2  &  2 \\
     0 &  -1 &  1 &  2 & 0  & 2 \\
     2 &  2 &  2 &  0 & 4 & 0 \\
     0 & -2 &  2 &  2 & 0 & 4 \\
  \end{array}
 \right ]
     \end{array}$    \\
    & & \\
     \hline
\end{tabular}
  \vspace{2mm}
  \caption{\small The partial Cartan matrix $B_{\Gamma}$ and the inverse matrix $B^{-1}_{\Gamma}$ for Carter
   diagrams with the number of vertices $l < 7$}
  \label{tab_partial Cartan_1}
\end{table}

\renewcommand{\appendixname}{}

\newpage
~\\~\\
\clearpage
~\\

\printindex

\end{document}